\documentclass[11pt,reqno,a4paper,dvipsnames]{amsart}

\usepackage[utf8]{inputenc}
\usepackage[dvipsnames]{xcolor}
\usepackage[textwidth=14cm]{geometry}
\usepackage{xfrac}

\usepackage{amsmath,amssymb}
\usepackage{enumitem}

\usepackage{tikz-cd}

\usepackage{mathtools}

\usepackage{amsthm}
\usepackage{etexcmds}
\usepackage{thmtools}
\usepackage{bbm}

\usepackage[stretch=10]{microtype}

\usepackage{mathtools}
\usepackage{amssymb}
\usepackage{amsthm}
\usepackage{amsmath}
\usepackage{graphicx}
\usepackage{xcolor}
\usepackage{enumerate}
\usepackage[normalem]{ulem}

\theoremstyle{plain}
\newtheorem{theorem}{Theorem}[subsection]
\newtheorem{lemma}[theorem]{Lemma}
\newtheorem{proposition}[theorem]{Proposition}

\theoremstyle{definition}

\newtheorem{remark}[theorem]{Remark}

\numberwithin{equation}{subsection}

\usepackage[hidelinks, colorlinks=false, hypertexnames=false]{hyperref}

\newcommand{\N}{{\mathbb{N}}}
\newcommand{\Z}{{\mathbb{Z}}}
\newcommand{\Q}{{\mathbb{Q}}}
\newcommand{\R}{{\mathbb{R}}}

\DeclareMathOperator{\dens}{\operatorname{dens}}
\DeclareMathOperator{\supp}{\operatorname{supp}}

\DeclareMathOperator{\Lip}{\operatorname{Lip}}

\newcommand{\fp}{{\mathfrak p}}
\newcommand{\fP}{{\mathfrak P}}
\newcommand{\fu}{{\mathfrak u}}
\newcommand{\fU}{{\mathfrak U}}

\newcommand{\fq}{{\mathfrak q}}

\newcommand{\fT}{{\mathfrak T}}
\newcommand{\fL}{{\mathfrak L}}
\newcommand{\fC}{{\mathfrak C}}
\newcommand{\pc}{{\mathrm{c}}}
\newcommand{\ps}{{\mathrm{s}}}

\newcommand{\fc}{{\Omega}}

\newcommand{\scI}{{\mathcal{I}}} 
\newcommand{\tQ}{{Q}}
\newcommand{\mfa}{{\vartheta}}
\newcommand{\mfb}{{\theta}}
\newcommand{\Mf}{{\Theta}}
\newcommand{\fcc}{{\mathcal{Q}}}

\newcommand{\proves}[1]{}
\newcommand{\lean}[1]{}
\newcommand{\leanok}{}
\ExplSyntaxOn \NewDocumentCommand{\uses}{m}
 {\clist_map_inline:nn{#1}{\vphantom{\ref{##1}}}%
  \ignorespaces}
\ExplSyntaxOff

\def\C{\mathbb{C}}

\usepackage[nameinlink, capitalize]{cleveref}

\usepackage[style=trad-alpha]{biblatex}
\bibliography{bibliography.bib}

\author[Becker et al.]{Lars Becker}
\address{Lars Becker,
	University of Bonn}
\email{becker@math.uni-bonn.de}

\author[]{Mar\'ia In\'es de Frutos-Fern\'andez}
\address{Mar\'ia In\'es de Frutos-Fern\'andez,
	University of Bonn}
\email{midff@math.uni-bonn.de}

\author[]{Leo Diedering}
\address{Leo Diedering,
	University of Bonn}
\email{leo.diedering@uni-bonn.de}

\author[]{Floris van Doorn}
\address{Floris van Doorn,
	University of Bonn}
\email{vdoorn@math.uni-bonn.de}

\author[]{S\'ebastien Gou\"ezel}
\address{S\'ebastien Gou\"ezel,
	University of Rennes}
\email{sebastien.gouezel@univ-rennes1.fr}

\author[]{Asgar Jamneshan}
\address{Asgar Jamneshan,
	University of Bonn}
\email{ajamnesh@math.uni-bonn.de}

\author[]{Evgenia Karunus}
\address{Evgenia Karunus,
  University of Bonn}
\email{lakesare@gmail.com}

\author[]{Edward van de Meent}
\address{Edward van de Meent,
  University of Utrecht}
\email{edwardvdmeent@gmail.com}

\author[]{Pietro Monticone}
\address{Pietro Monticone,
University of Trento}
\email{pit.monticone@gmail.com}

\author[]{Jasper Mulder-Sohn}
\address{Jasper Mulder-Sohn,
	The Hague}
\email{jasper.mulder@planet.nl}

\author[]{Jim Portegies}
\address{Jim Portegies,
  Eindhoven University of Technology}
\email{j.w.portegies@tue.nl}

\author[]{Joris Roos}
\address{Joris Roos,
  University of Massachusetts Lowell}
\email{joris\_roos@uml.edu}

\author[]{Michael Rothgang}
\address{Michael Rothgang,
	University of Bonn}
\email{rothgang@math.uni-bonn.de}

\author[]{Rajula Srivastava}
\address{Rajula Srivastava,
	University of Bonn}
\email{rajulas@math.uni-bonn.de}

\author[]{James Sundstrom}
\address{James Sundstrom,
	Baruch College}
\email{james.sundstrom@baruch.cuny.edu}

\author[]{Jeremy Tan}
\address{Jeremy Tan,
	National University of Singapore}
\email{jtjierui@comp.nus.edu.sg}

\author[]{Christoph Thiele}
 \address{Christoph Thiele,
	University of Bonn}
\email{thiele@math.uni-bonn.de}

\begin{document}

%
%

\title[Formalization of Carleson's theorem]{A blueprint for the formalization of Carleson's theorem on convergence of Fourier series}

\date{\today}

\begin{abstract}
This paper is the blueprint underlying the Lean formalization of the proof of Carleson's classical result \cite{carleson} asserting almost everywhere convergence of Fourier series of continuous functions. We break up the proof into two steps, a reduction of the classical result to a new theorem that appears in a sibling communication (to be submitted) and a proof of this new theorem, which is also detailed as blueprint in this paper. An early version of this blueprint was used to initiate the Lean formalization. During the formalization, many contributors elaborated the blueprint with minor corrections, modifications and extensions.
The final version is presented here as a guide through the accompanying Lean code.
\end{abstract}

\maketitle

\tableofcontents


\section{Introduction}

Trigonometric series represent functions as possibly infinite linear combinations of pure frequencies.
They gained particular prominence through the work of J. Fourier, who used them in his analytical theory of heat \cite{Fourier}, see also \cite{MR3470070}, thereby establishing them as a tool for solving partial differential equations.
Fourier also made the groundbreaking claim that a wide range of functions could be represented using trigonometric series. This sparked the interest of many mathematicians,
including Dirichlet, who gave some rigorous conditions for convergence of Fourier series, as trigonometric series are now called. Dirichlet also opened a branch of analytic number theory partially inspired by the ideas of
Fourier.
Nowadays, Fourier analysis plays an important role in many areas of mathematics.

With Euler's formula to represent pure frequencies in mind, a trigonometric polynomial can be expressed as
\begin{equation}\label{eq:trig-series}
    s_N(x):= \sum_{n=-N}^N c_n e^{inx} \ .
\end{equation}
The Fourier series is then defined as the limit $f$
of such a sequence $s_N$ as $N$ tends to $\infty$.
Fourier's vision to represent rather general functions raises two fundamental questions.
The first question is to identify the appropriate choice of coefficients
$c_n$ to use to represent a given $f$. The second question addresses the convergence of $s_N$ to $f$.

The first question has a fairly canonical and standard answer, provided by the Fourier integral formula:
\begin{equation}\label{eq:fourier-coefficients}
    c_n:=\widehat{f}_n:=\frac 1{2\pi}\int_0^{2\pi}f(x) e^{- i nx}\, dx,
\end{equation}
 where the precise interpretation of the integral depends on the chosen theory of integration. For a continuous function $f$,
Riemann's notion of the integral suffices. If $f$ is integrable in the Lebesgue sense,
$f\in L^1(0,2\pi)$, the Lebesgue integral is appropriate. More generally, if
$f$ is a distribution in the sense of Schwartz, supported in $[0,2\pi]$
the integral can be understood as an evaluation against the periodic test function
$e^{-i nx}$. In each case, the more general definition reduces to the simpler one within the respective more restrictive domain.
Hence, the Fourier coefficients given by \eqref{eq:fourier-coefficients} serve as a universal choice.
This choice is unique in several respects, in particular if one is to exactly reproduce
a trigonometric polynomial $f$ in the form
\eqref{eq:trig-series}.

The second question of convergence bifurcates into the question of pointwise convergence
of the series \eqref{eq:trig-series} (with coefficients given by \eqref{eq:fourier-coefficients}) for a given $x$ on the one hand and convergence of the functions $s_N$ to the function $f$ in a suitable
function space with corresponding topology on the other hand. There are at
least as many function spaces for
the question of convergence as there are different definitions of the integral elaborated earlier.
There are some very
canonical answers to the convergence question in function spaces, albeit not known at the time of Fourier and Dirichlet.
One example is convergence in the Hilbert space sense for $f$ in $L^2(0,2\pi)$,
as discovered in the first decade of the twentieth century
as a consequence of the rapid development of Lebesgue integration theory.
Another canonical example is convergence in the sense of distributions for general distributional $f$, as discovered a few decades after Lebesgue integration. For some other natural spaces, such as $L^1(0,2\pi)$, there is no guarantee of convergence in the
norm of that space even if $f$ is in the space.

In contrast to these examples of function spaces with a very natural theory of
convergence of Fourier series in the topology of the function space, there are no similarly elegant solutions to the
characterization of pointwise convergence.
In particular, the space of functions $f$ such that the sequence $s_N(x)$ converges for every $x$ does not have a good characterization in terms of $f$ itself.
Similarly, the space of all functions $f$ such that the sequence of coefficients $\widehat{f}_n$ is absolutely summable has also no good characterization.

When the Fourier integral is defined in the Lebesgue sense and
$f\in L^1(0,2\pi)$, then the function
$f$ itself is meaningful not everywhere but only pointwise almost everywhere
in the Lebesgue sense. The question of pointwise convergence to $f$ for all $x$ then becomes
meaningless, and instead one asks for almost everywhere convergence.
Such convergence
was conjectured by N. Luzin \cite{Luzin13}
for the space
$L^2(0,2\pi)$, see also the collected works \cite{Luzin53}. Kolmogorov's example \cite{Kolmogorov} of an $L^1$ function whose Fourier series diverges at almost every point seemed to provide some evidence against Luzin's conjecture.
Indeed, in the 1960s, L. Carleson tried to construct a counterexample to Luzin's conjecture.
In his own recollection, his efforts led him to better and better understand how such a potential counterexample should look like.
In the end, the counterexample had to satisfy so many properties that the properties started to contradict each other,
and Carleson realized that a counterexample could not exist.
Thus, he had proved Luzin's conjecture \cite{carleson}.
In particular, he had proven the more elementary statement

\begin{theorem}[classical Carleson]
    \label{classical-carleson}
    \leanok
    \lean{classical_carleson}
    \uses{exceptional-set-carleson}
    Let $f$ be a $2\pi$-periodic complex-valued continuous function on $\mathbb{R}$.
    Then for almost all $x \in \mathbb{R}$ we have
    \begin{equation}\label{eq:fourier-limit}
      \lim_{N\to\infty}s_N f(x) = f(x),
    \end{equation}
    where $S_N f$ is the $N$-th partial Fourier sum of $f$ defined in \eqref{eq:trig-series}
    with coefficients \eqref{eq:fourier-coefficients}.
\end{theorem}

Here, almost every $x$ means in the Lebesgue sense, i.e., for every $\epsilon>0$ the set of $x$ where convergence fails can be covered by a sequence of intervals
such that the sum of the lengths of these
intervals is less than $\epsilon$. While Carleson had proven the more general Luzin conjecture for functions in $L^2[0,2\pi]$, even the more elementary statement for continuous functions was not known before Carleson's work.
Moreover, until now, the elementary statement has not seen any substantially easier proof than those generalizing to $L^2$,
partially because there is no readily usable criterion on the level of
Fourier coefficients to distinguish between continuous functions and $L^2$ functions.

Shortly after Carleson's breakthrough, Hunt \cite{MR238019} proved the analogous result for $L^p$ functions
with $p>1$ and for functions in $L^1\log(L)^2$.
Billard \cite{MR217510} adapted Carleson's arguments to prove almost everywhere convergence of Walsh-Fourier series
for functions in $L^2$.

In \cite{fefferman}, C. Fefferman gave an alternative proof of Carleson's theorem via an a priori
bound for Carleson's operator, the maximally modulated singular integral
\begin{equation}\label{eq:feffermancarleson}
    Tf(x):=\sup_{N} \int_{-\pi}^\pi e^{iNy}f(x-y)\frac 1y {dy}\, .
\end{equation}
In \cite{lacey-thiele}, a duality
between the approaches by Carleson and Fefferman was pointed out and a more symmetric and self-dual proof was presented.

Various strengthenings of Fefferman's estimates for Carleson's operator have appeared since, such as bounds for a higher dimensional Carleson operator in the isotropic
\cite{MR336222}, \cite{MR2007237} and anisotropic \cite{MR4014801} setting.
The supremum norm in the variable $N$ in \eqref{eq:feffermancarleson} was strengthened
to maximal multiplier norms in \cite{MR2420509} with applications to return times theorems in ergodic theory.
It was strengthened to variation norms $V^r$ with $r>2$ in \cite{MR2881301} with applications to nonlinear analysis.

Stein and Wainger \cite{stein-wainger} proposed to study variants of \eqref{eq:feffermancarleson} replacing the
linear modulation phase $Ny$ by polynomial phases in $y$ and proved a result for polynomials without a linear term.
This was generalized by Lie to general polynomial Carleson operators, first for a supremum over all quadratic
polynomials \cite{lie-quadratic} and then for a supremum over all polynomials in \cite{lie-polynomial}.
The polynomial Carleson operator was further generalized to higher dimensions and H\"older regular kernels in \cite{zk-polynomial}.
The development of polynomial Carleson operators has led to further generalisations, for example to almost everywhere convergence of Malmquist-Takenaka series
\cite{mnatsakanyan} and maximally modulated singular Radon transforms \cite{ramos}, \cite{becker2024maximal}.

Surveys to Carleson's theorem can be found in \cite{MR2091007}, \cite{MR3334208}.

In this blueprint, we prove Theorem \ref{classical-carleson} as a corollary to a further generalization of the polynomial Carleson operator towards doubling metric measure spaces,
Theorem \ref{metric-space-Carleson} below,
which is new and appears in a sibling communication (to be submitted).
Theorem \ref{metric-space-Carleson} is an axiomatic approach to
Carleson type theorems on doubling metric measure spaces.
This axiomatic approach is suitable for formalization and a good route towards the classical theorem.

Early drafts of the sibling communication (to be submitted) existed in summer 2023.
Based on this, the first draft of the present blueprint was written in the first half of 2024,
containing a much more detailed proof, which involved increasing the size by a factor of four,
and adding the derivation of Carleson's classical result.
In June 2024, the fourth author launched a \href{https://florisvandoorn.com/carleson/}{public website}
to post the blueprint, using the open-source software \texttt{leanblueprint} developed by Patrick Massot,
calling for contributions to formalize the proof.
The goal was to formalize the blueprint in the Lean proof assistant~\cite{moura2021lean},
building on top of its mathematical library \texttt{Mathlib}~\cite{mathlib}.
The work was split up into about 180 tasks, to be claimed by individual contributors.
Most tasks were to formalize the proof of a single lemma from the blueprint, and some were to develop basic theory or refactor existing code.
The contributors adapted the blueprint to fix some gaps found during the formalization and gave feedback that led to discussions about the proof.
This even resulted in a few changes to the general setup and the main theorems.
All of the gaps found required only fairly localized changes to the blueprint, indicating that the initial blueprint was already of high quality.
The formalization was completed in July 2025. It is attached to this arXiv posting, and the latest version can be found \href{https://github.com/fpvandoorn/carleson}{on Github}.

Everyone that completed a substantial amount of tasks is included as a coauthor of the blueprint.
The authors acknowledge contributions in the form of small formalization additions,
pointing out corrections to the blueprint,
or supplying ideas to the Lean efforts by the following people:
Michel Alexis,
Bolton Bailey,
Julian Berman,
Joachim Breitner,
Martin Dvořák,
Georges Gonthier,
Aaron Hill,
Austin Letson,
Bhavik Mehta,
Eric Paul,
Clara Torres,
Dennis Tsar,
Andrew Yang,
Ruben van de Velde.

\noindent \textit{Acknowledgement.}
L.B., M.I.d.F.F., L.D., F.v.D., M.R., R.S., and C.T. were funded by the Deutsche For\-schungs\-gemein\-schaft (DFG, German Research Foundation) under Germany's Excellence Strategy -- EXC-2047/1 -- 390685813.
L.B. , R.S., and C.T. were also supported by SFB 1060.
A.J. is funded by the T\"UBITAK (Scientific and Technological Research Council of T\"urkiye) under Grant Number 123F122.
J.R. was supported in part by NSF grant DMS-2154835 and a HIM fellowship for the Fall 2024 trimester program in Bonn.

\subsection{Statement of the metric space Carleson theorems}

Let
\begin{equation}
    a\ge 4
\end{equation}
be a natural number that as it gets larger will allow for more general applications of Theorem \ref{metric-space-Carleson}
but will worsen the constants in the estimates.

A doubling metric measure space $(X,\rho,\mu, a)$ is a complete
and locally compact metric space $(X,\rho)$
equipped with a non-zero locally finite Borel measure $\mu$ that satisfies the doubling condition that for all $x\in X$ and all $R>0$ we have
\begin{equation}\label{doublingx}
    \mu(B(x,2R))\le 2^a\mu(B(x,R))\,,
\end{equation}
where we have denoted by $B(x,R)$ the open ball of radius $R$ centered at $x$:
\begin{equation}\label{eq-define-ball}
 B(x,R):=\{y\in X: \rho(x,y)<R\}. \end{equation}
In a doubling metric measure space the closed balls are compact and $\mu$ is positive on all non-empty open sets.

A collection $\Mf$ of real valued continuous functions on the doubling metric measure space $(X,\rho,\mu,a)$ is called compatible,
if there is a point $o\in X$ where all the functions are equal to $0$,
and if there exists for each ball $B \subset X$ a metric $d_B$ on $\Mf$,
such that the following five properties \eqref{osccontrol}, \eqref{firstdb}, \eqref{monotonedb}, \eqref{seconddb}, and \eqref{thirddb} are satisfied.
For every ball $B \subset X$
\begin{equation}\label{osccontrol}
    \sup_{x,y\in B}|\mfa(x)-{\mfa(y)}- \mfb(x)+{\mfb(y)}| \le d_{B}(\mfa,\mfb)\,.
\end{equation}
For any two balls $B_1=B(x_1,R)$, $B_2= B(x_2,2R)$ in $X$ with $x_1\in B_2$ and any $\mfa,\mfb\in \Mf$,
\begin{equation}\label{firstdb}
    d_{B_2}(\mfa,\mfb)\le 2^a d_{B_1}(\mfa,\mfb) .
\end{equation}
For any two balls $B_1, B_2$ in $X$ with $B_1 \subset B_2$ and any $\mfa, \mfb \in \Mf$
\begin{equation}\label{monotonedb}
    d_{B_1}(\mfa,\mfb) \le d_{B_2}(\mfa, \mfb)
\end{equation}
and for any two balls
$B_1=B(x_1,R)$, $B_2= B(x_2,2^aR)$
with $B_1\subset B_2$, and $\mfa,\mfb\in \Mf$,
\begin{equation}\label{seconddb}
    2d_{B_1}(\mfa,\mfb)
\le d_{B_2}(\mfa,\mfb) .
\end{equation}
For every ball $B$ in $X$ and every $d_B$-ball $\tilde B$ of radius $2R$ in $\Mf$, there is a collection $\mathcal{B}$ of
at most $2^a$ many $d_B$-balls of radius $R$ covering $\tilde B$, that is,
\begin{equation}\label{thirddb}
    \tilde B\subset \bigcup \mathcal{B}.
\end{equation}

Further, a compatible collection $\Mf$ is called cancellative, if
for any ball $B$ in $X$ of radius $R$, any Lipschitz function $\varphi: X\to \C$
supported on $B$, and any $\mfa,\mfb\in \Mf$ we have
\begin{equation}
    \label{eq-vdc-cond}
    |\int_B e(\mfa(x)-{\mfb(x)}) \varphi(x) d\mu(x)|\le 2^a \mu(B)\|\varphi\|_{\Lip(B)}
(1+d_B(\mfa,\mfb))^{-\frac{1}{a}},
\end{equation}
where $\|\cdot\|_{\Lip(B)}$ denotes the inhomogeneous Lipschitz norm on $B$:
$$
    \|\varphi\|_{\Lip(B)} = \sup_{x \in B} |\varphi(x)| + R \sup_{x,y \in B, x \neq y} \frac{|\varphi(x) - \varphi(y)|}{\rho(x,y)}\,.
$$

A one-sided Calder\'on--Zygmund kernel $K$ on the doubling metric measure space $(X, \rho, \mu, a)$ is a measurable function
\begin{equation}\label{eqkernel0}
  K:X\times X\to \mathbb{C}
\end{equation}
such that for all $x,y',y\in X$ with $x\neq y$, we have
\begin{equation}\label{eqkernel-size}
    |K(x,y)| \leq \frac{2^{a^3}}{V(x,y)}
\end{equation}
and if $2\rho(y,y') \leq \rho(x,y)$, then
\begin{equation}
  \label{eqkernel-y-smooth}
  |K(x,y) - K(x,y')| \leq \left(\frac{\rho(y,y')}{\rho(x,y)}\right)^{\frac{1}{a}}\frac{2^{a^3}}{V(x,y)},
\end{equation}
where \[V(x,y):=\mu(B(x,\rho(x,y))).\]
Define the maximally truncated non-tangential singular integral $T_{*}$ associated with $K$ by
\begin{equation}
    \label{def-tang-unm-op}
    T_{*}f(x):=\sup_{R_1 < R_2} \sup_{\rho(x,x')<R_1} \left|\int_{R_1< \rho(x',y) < R_2} K(x',y) f(y) \, \mathrm{d}\mu(y) \right|\,.
\end{equation}
We define the generalized Carleson operator $T$ by
\begin{equation}
    \label{def-main-op}
    Tf(x):=\sup_{\mfa\in\Mf} \sup_{0 < R_1 < R_2}\left| \int_{R_1 < \rho(x,y) < R_2} K(x,y) f(y) e(\mfa(y)) \, \mathrm{d}\mu(y) \right|\, ,
\end{equation}
where $e(r)=e^{ir}$.

Our first main result is the following restricted weak type estimate for $T$ in the range $1<q\le 2$,
which by interpolation techniques recovers $L^q$ estimates for the open range $1<q<2$.
\begin{theorem}[metric space Carleson]
\label{metric-space-Carleson}
\uses{linearised-metric-Carleson,int-continuous}
\leanok
\lean{metric_carleson}
    For all integers $a \ge 4$ and real numbers $1<q\le 2$ the following holds.
    Let $(X,\rho,\mu,a)$ be a doubling metric measure space.
    Let $\Mf$ be a cancellative compatible collection of functions and let $K$ be a one-sided Calder\'on--Zygmund kernel on $(X,\rho,\mu,a)$.
    Assume that for every bounded measurable function $g$ on $X$ supported on a set of finite measure we have
    \begin{equation}\label{nontanbound}
        \|T_{*}g\|_{2} \leq 2^{a^3} \|g\|_2\,,
    \end{equation}
    where $T_{*}$ is defined in \eqref{def-tang-unm-op}.
    Then for all Borel sets $F$ and $G$ in $X$ and all Borel functions $f:X\to \C$ with
    $|f|\le \mathbf{1}_F$, we have, with $T$ defined in \eqref{def-main-op},
    \begin{equation}
        \label{resweak}
        \left|\int_{G} T f \, \mathrm{d}\mu\right| \leq \frac{2^{443a^3}}{(q-1)^6} \mu(G)^{1-\frac{1}{q}} \mu(F)^{\frac{1}{q}}\, .
    \end{equation}
\end{theorem}

In some applications, such as the Walsh-case of Carleson's theorem, the kernel $K$ naturally depends also on the modulation functions $\mfa$.
The fact that we don't assume H\"older-continuity of the kernel $K$ in the first argument allows us to also capture this situation.

We now state another Theorem with a slightly weaker replacement of the assumption \eqref{nontanbound}.
It implies the Walsh-case of Carleson's theorem.
For a Borel function $\tQ:X\to \Mf$, and $\mfa \in \Mf$ and $x\in X$ define
\begin{equation}
    R_{\tQ}(\mfa,x)=\sup\{r>0:d_{B(x,r)}(\mfa, \tQ(x))<1\}
\end{equation}
and define further
\begin{equation}
    \label{def-lin-star-op}
    T_{\tQ}^\mfa f(x):=\sup_{R_1<R_2} \ \sup_{\rho(x,x')<R_1}
    \left|\int_{R_1< \rho(x',y) < \min\{R_2, R_{\tQ}(\mfa,x')\}} K(x',y) f(y) \, \mathrm{d}\mu(y) \right|
\end{equation}
Define further the linearized generalized Carleson operator $T_\tQ$ by
\begin{equation}
    \label{def-lin-main-op}
    T_\tQ f(x):= \sup_{0 < R_1 < R_2}\left| \int_{R_1 < \rho(x,y) < R_2} K(x,y) f(y) e(\tQ(x)(y)) \, \mathrm{d}\mu(y) \right|\, ,
\end{equation}
where again $e(r)=e^{ir}$.

\begin{theorem}[linearised metric Carleson]
\label{linearised-metric-Carleson}
\uses{R-truncation,int-continuous}
\leanok
\lean{linearized_metric_carleson}
    For all integers $a \ge 4$ and real numbers $1<q\le 2$ the following holds.
    Let $(X,\rho,\mu,a)$ be a doubling metric measure space. Let $\Mf$ be a
    cancellative compatible collection of functions.
    Let $\tQ:X\to \Mf$ be a Borel function with finite range.
    Let $K$ be a one-sided Calder\'on--Zygmund kernel on $(X,\rho,\mu,a)$. Assume that for every $\mfa\in \Mf$ and every bounded measurable function $g$ on $X$ supported on a set of finite measure we have
    \begin{equation}\label{linnontanbound}
        \|T_{\tQ}^\mfa g\|_{2} \leq 2^{a^3} \|g\|_2\,,
    \end{equation}
    where $T_{\tQ}^\mfa$ is defined in \eqref{def-lin-star-op}.
    Then for all Borel sets $F$ and $G$ in $X$ and all Borel functions $f:X\to \C$ with
    $|f|\le \mathbf{1}_F$, we have, with $T_\tQ$ defined in \eqref{def-lin-main-op},
    \begin{equation}
        \label{linresweak}
        \left|\int_{G} T_\tQ f \, \mathrm{d}\mu\right| \le \frac{2^{443a^3}}{(q-1)^6} \mu(G)^{1-\frac{1}{q}} \mu(F)^{\frac{1}{q}}\, .
    \end{equation}
\end{theorem}

\begin{remark}
The value of the constant factor $2^{443a^3}$ in Theorems~\ref{metric-space-Carleson}
and~\ref{linearised-metric-Carleson} is by no means sharp.
One source of non-sharpness is our choice to write for readability most constants in the form $2^{na^3}$ for some explicit constant $n$.

An a posteriori byproduct of the Lean formalization of this document is that Theorems~\ref{metric-space-Carleson}
and~\ref{linearised-metric-Carleson} remain true with $2^{44a^3}$ instead of $2^{443a^3}$. This is obtained by
changing the parameter $D$ introduced in~\eqref{defineD} from $2^{100a^3}$ to $2^{7a^3}$, checking that exactly the same
proof goes through, tracking the constants, and getting $2^{44a^3}$ in the end. This value is again by no means sharp.
In the formalization we define the parameter $D$ as $2^{\mathbbm{c}a^3}$
and the constants in this blueprint are obtained by setting $\mathbbm{c} = 100$.
\end{remark}

\section{Proof of Metric Space Carleson, overview}
\label{overviewsection}

This section organizes the proof of \Cref{metric-space-Carleson} into sections \ref{thmfromproplinear}, \ref{christsection}, \ref{proptopropprop}, \ref{antichainboundary}, \ref{treesection}, \ref{liphoel}, and \ref{sec-hlm}. These sections are mutually independent except for referring to the statements formulated in the present section. \Cref{thmfromproplinear} proves the main \Cref{metric-space-Carleson}, while sections \ref{christsection}, \ref{proptopropprop}, \ref{antichainboundary}, \ref{treesection}, \ref{liphoel}, and \ref{sec-hlm} each prove one proposition that is stated in the present section. The present section also introduces all definitions used across these sections.

\Cref{global-auxiliary-lemmas} proves some auxiliary lemmas that are used in more than one of the sections 3-9.

Let $a, q$ be given as in \Cref{metric-space-Carleson}.

Define
\begin{equation}\label{defineD}
    D:= 2^{100 a^2}\, ,
\end{equation}
\begin{equation}\label{definekappa}
    \kappa:= 2^{-10a}\,,
\end{equation}
and
\begin{equation}
    \label{defineZ}
    Z := 2^{12a}\,.
\end{equation}
Let
 $\psi:\R \to \R$ be the unique compactly supported, piece-wise linear, continuous function with corners precisely at $\frac 1{4D}$, $\frac 1{2D}$, $\frac 14$ and $\frac 12$ which satisfies
 \begin{equation}
    \label{eq-psisum}
    \sum_{s\in \mathbb{Z}} \psi(D^{-s}x)=1
\end{equation}
for all $x>0$. This function vanishes outside $[\frac1{4D},\frac 12]$, is constant one on
$[\frac1{2D},\frac 14]$, and is Lipschitz
with constant $4D$.

Let a doubling metric measure space $(X,\rho,\mu, a)$ be given.
Let a cancellative compatible collection $\Mf$ of functions on $X$ be given.
Let $o\in X$ be a point such that $\mfa(o)=0$
for all $\mfa\in \Mf$.

Let a Borel measurable function $\tQ:X\to \Mf$ with finite range be given.
Let a one-sided Calder\'on--Zygmund kernel $K$ on $X$ be given so that for every $\mfa\in \Mf$ the operator $T_{\tQ}^{\mfa}$ defined in \eqref{def-lin-star-op}
satisfies \eqref{linnontanbound}.

For $s\in\mathbb{Z}$, we define
\begin{equation}\label{defks}
    K_s(x,y):=K(x,y)\psi(D^{-s}\rho(x,y))\,,
\end{equation}
so that for each $x, y \in X$ with $x\neq y$ we have
$$K(x,y)=\sum_{s\in\mathbb{Z}}K_s(x,y).$$

In \Cref{thmfromproplinear}, we prove \Cref{metric-space-Carleson} and \Cref{linearised-metric-Carleson}
from the finitary version, \Cref{finitary-Carleson} below. Recall
that a function from a measure space to a finite set is measurable if the pre-image of each of the elements in the range is measurable.

\begin{proposition}[finitary Carleson]
\label{finitary-Carleson}
\leanok
\lean{finitary_carleson}
\uses{discrete-Carleson, grid-existence, tile-structure, tile-sum-operator}
Let ${\sigma_1},\sigma_2\colon X\to \mathbb{Z}$ be measurable functions with finite range and ${\sigma_1}\leq \sigma_2$.  Let $F,G$ be bounded Borel sets in $X$. Then there is a Borel set $G'$ in $X$ with $2\mu(G')\leq \mu(G)$ such that
for all Borel functions $f:X\to \C$ with $|f|\le \mathbf{1}_F$.
\begin{equation*}
    \int_{G \setminus G'} \left|\sum_{s={\sigma_1}(x)}^{{\sigma_2}(x)} \int K_s(x,y) f(y) e(\tQ(x)(y)) \, \mathrm{d}\mu(y) \right| \mathrm{d}\mu(x)
\end{equation*}
\begin{equation}
    \label{eq-linearized}
    \le \frac{2^{442a^3}}{(q-1)^5} \mu(G)^{1-\frac{1}{q}} \mu(F)^{\frac 1 q}\,.
\end{equation}
\end{proposition}
Let measurable functions ${\sigma_1}\leq \sigma_2\colon X\to \mathbb{Z}$ with finite range be given.
Let bounded Borel sets $F,G$ in $X$ be given.
Let $S$ be the smallest integer such that the ranges of
$\sigma_1$ and $\sigma_2$ are contained in $[-S,S]$ and $F$ and $G$ are contained
in the ball $B(o, \frac{1}{4}D^S)$.

In \Cref{christsection}, we prove \Cref{finitary-Carleson} using a
bound for a dyadic model formulated in \Cref{discrete-Carleson} below.

A grid structure $(\mathcal{D}, c, s)$ on $X$ consists of a finite collection $\mathcal{D}$ of pairs $(I, k)$ of Borel
sets in $X$ and integers $k \in [-S, S]$, the projection $s\colon \mathcal{D}\to [-S, S], (I, k) \mapsto k$ to the second component which is assumed to be surjective and
called scale function, and a function $c:\mathcal{D}\to X$
called center function such that the five properties
\eqref{coverdyadic}, \eqref{dyadicproperty}, \eqref{subsetmaxcube},
\eqref{eq-vol-sp-cube}, and \eqref{eq-small-boundary} hold. We call the elements of $\mathcal{D}$ dyadic cubes. By abuse of notation, we will usually write just $I$ for the cube $(I,k)$, and we will write $I \subset J$ to mean that for two cubes $(I,k), (J, l) \in \mathcal{D}$ we have $I \subset J$ and $k \le l$.

For each dyadic cube $I$ and each $-S\le k<s(I)$ we have
\begin{equation}\label{coverdyadic}
    I\subset \bigcup_{J\in \mathcal {D}: s(J)=k}J\, .
\end{equation}
Any two non-disjoint dyadic cubes $I,J$ with $s(I)\le s(J)$ satisfy
\begin{equation}\label{dyadicproperty}
    I\subset J.
\end{equation}
There exists a $I_0 \in \mathcal{D}$ with $s(I_0) = S$ and $c(I_0) = o$
and for all cubes $J \in \mathcal{D}$, we have
\begin{equation}\label{subsetmaxcube}
    J \subset I_0\,.
\end{equation}
For any dyadic cube $I$,
\begin{equation}
    \label{eq-vol-sp-cube}
    c(I)\in B(c(I), \frac{1}{4} D^{s(I)}) \subset I \subset B(c(I), 4 D^{s(I)})\,.
\end{equation}
For any dyadic cube $I$ and any $t$ with $tD^{s(I)} \ge D^{-S}$,
\begin{equation}
    \label{eq-small-boundary}
    \mu(\{x \in I \ : \ \rho(x, X \setminus I) \leq t D^{s(I)}\}) \le 2 t^\kappa \mu(I)\,.
\end{equation}

A tile structure $(\fP,\scI,\fc,\fcc,\pc,\ps)$
for a given grid structure $(\mathcal{D}, c, s)$
is a finite set $\fP$ of elements called tiles with five maps
\begin{align*}
    \scI&\colon \fP\to {\mathcal{D}}\\
    \fc&\colon \fP\to \mathcal{P}(\Mf) \\
    \fcc &\colon \fP\to \tQ(X)\\
    \pc &\colon \fP\to X\\
    \ps &\colon \fP\to \mathbb{Z}
\end{align*}
with $\scI$ surjective and $\mathcal{P}(\Mf)$ denoting the power set of $\Mf$ such that the five properties \eqref{eq-dis-freq-cover}, \eqref{eq-freq-dyadic},
\eqref{eq-freq-comp-ball}, \eqref{tilecenter}, and
\eqref{tilescale} hold.
For each dyadic cube $I$, the restriction of the map $\Omega$ to the set
\begin{equation}\label{injective}
    \fP(I)=\{\fp: \scI(\fp) =I\}
\end{equation}
is injective
and we have the disjoint covering property (we use the union symbol with dot on top to denote a disjoint union)
\begin{equation}\label{eq-dis-freq-cover}
    \tQ(X)\subset \dot{\bigcup}_{\fp\in \fP(I)}\fc(\fp).
\end{equation}
For any tiles $\fp,\fq$ with $\scI(\fp)\subset \scI(\fq)$ and $\fc(\fp) \cap \fc(\fq) \neq \emptyset$ we have
\begin{equation} \label{eq-freq-dyadic}
    \fc(\fq)\subset \fc(\fp) .
\end{equation}
For each tile $\fp$,
\begin{equation}\label{eq-freq-comp-ball}
    \fcc(\fp)\in B_{\fp}(\fcc(\fp), 0.2) \subset \fc(\fp) \subset B_{\fp}(\fcc(\fp),1)\,,
\end{equation}
where
\begin{equation}
    B_{\fp} (\mfa, R) := \{\mfb \in \Mf \, : \, d_{\fp}(\mfa, \mfb) < R\,\} ,
\end{equation}
and
\begin{equation}\label{defdp}
    d_{\fp} := d_{B(\pc(\fp),\frac 14 D^{\ps(\fp)})}\, .
\end{equation}
We have for each tile $\fp$
\begin{equation}\label{tilecenter}
    \pc(\fp)=c(\scI(\fp)),
\end{equation}
\begin{equation}\label{tilescale}
    \ps(\fp)=s(\scI(\fp)).
\end{equation}

\begin{proposition}[discrete Carleson]
\label{discrete-Carleson}
\leanok
\lean{discrete_carleson}
\uses{exceptional-set, forest-union, forest-complement}
Let $(\mathcal{D}, c, s)$ be a grid structure and
\begin{equation*}
    (\fP,\scI,\fc,\fcc,\pc,\ps)
\end{equation*}
a tile structure for this grid structure.
Define for $\fp\in \fP$
\begin{equation}\label{defineep}
    E(\fp)=\{x\in \scI(\fp): \tQ(x)\in \fc(\fp) , {\sigma_1}(x)\le \ps(\fp)\le {\sigma_2}(x)\}
\end{equation}
and
\begin{equation}\label{definetp}
    T_{\fp} f(x)= \mathbf{1}_{E(\fp)}(x) \int K_{\ps(\fp)}(x,y) f(y) e(\tQ(x)(y)-\tQ(x)(x))\, d\mu(y).
\end{equation}
Then there exists a Borel set $G'$ with $2\mu(G') \leq \mu(G)$ such that for all Borel functions $f:X\to \C$ with $|f|\le \mathbf{1}_F$
we have
\begin{equation}
    \label{disclesssim}
   \int_{G \setminus G'} \left| \sum_{\fp \in \fP} T_{\fp} f (x) \right| \, \mathrm{d}\mu(x) \le \frac{2^{442a^3}}{(q-1)^5} \mu(G)^{1-\frac{1}{q}} \mu(F)^{\frac{1}{q}}\,.
\end{equation}
\end{proposition}

The proof of \Cref{discrete-Carleson} is done in \Cref{proptopropprop}
by a reduction to two further propositions that we state below.

Fix a grid structure $(\mathcal{D}, c, s)$ and a tile structure $(\fP,\scI,\fc,\fcc,\pc,\ps)$
for this grid structure.

We define the relation
\begin{equation}\label{straightorder}
    \fp\le \fp'
\end{equation}
 on $\fP\times \fP$ meaning
$\scI(\fp)\subset \scI(\fp')$ and
$\Omega(\fp')\subset \Omega(\fp)$.
We further define for $\lambda,\lambda' >0$
the relation
\begin{equation}\label{wiggleorder}
    \lambda \fp \lesssim \lambda' \fp'
\end{equation}
on $\fP\times \fP$ meaning
$\scI(\fp)\subset \scI(\fp')$ and
\begin{equation}
    B_{\fp'}(\fcc(\fp'),\lambda') \subset B_{\fp}(\fcc(\fp),\lambda)\, .
\end{equation}

Define for a tile $\fp$ and $\lambda>0$
\begin{equation}\label{definee1}
    E_1(\fp):=\{x\in \scI(\fp)\cap G: \tQ(x)\in \fc(\fp)\}\, ,
\end{equation}
\begin{equation}\label{definee2}
    E_2(\lambda, \fp):=\{x\in \scI(\fp)\cap G: \tQ(x)\in B_{\fp}(\fcc(\fp), \lambda)\}\, .
\end{equation}

Given a subset $\fP'$ of $\fP$, we define
$\fP(\fP')$ to be the set of
all $\fp \in \fP$ such that there exist $\fp' \in \fP'$ with $\scI(\fp)\subset \scI(\fp')$. Define the densities
\begin{equation}\label{definedens1}
    {\dens}_1(\fP') := \sup_{\fp'\in \fP'}\sup_{\lambda \geq 2} \lambda^{-a} \sup_{\fp \in \fP(\fP'), \lambda \fp' \lesssim \lambda \fp}
    \frac{\mu({E}_2(\lambda, \fp))}{\mu(\scI(\fp))}\, ,
\end{equation}
\begin{equation}\label{definedens2}
    {\dens}_2(\fP') := \sup_{\fp'\in \fP'}
    \sup_{r\ge 4D^{\ps(\fp)}}
    \frac{\mu(F\cap B(\pc(\fp),r))}{\mu(B(\pc(\fp),r))}\, .
\end{equation}

An antichain is a subset $\mathfrak{A}$
of $\fP$ such that for any distinct $\fp,\fq\in \mathfrak{A}$ we do not have $\fp\le \fq$.

The following proposition is proved in \Cref{antichainboundary}.

\begin{proposition}[antichain operator]
\label{antichain-operator}
\leanok
\lean{antichain_operator}

\uses{dens2-antichain,dens1-antichain}
For any antichain $\mathfrak{A} $ and for all $f:X\to \C$ with $|f|\le \mathbf{1}_F$ and all $g:X\to\C$ with $|g| \le \mathbf{1}_G$
\begin{equation} \label{eq-antiprop}
    |\int \overline{g(x)} \sum_{\fp \in \mathfrak{A}} T_{\fp} f(x)\, d\mu(x)|
\end{equation}
\begin{equation}
    \le \frac{2^{117a^3}}{q-1} \dens_1(\mathfrak{A})^{\frac {q-1}{8a^4}}\dens_2(\mathfrak{A})^{\frac 1{q}-\frac 12} \|f\|_2 \|g\|_2\, .
\end{equation}
\end{proposition}

Let $n\ge 0$.
An $n$-forest is a pair $(\fU, \mathfrak{T})$
where $\fU$ is a subset of $\fP$
and $\mathfrak{T}$ is a map assigning to
each $\fu\in \fU$ a nonempty set $\fT (\fu)\subset \fP$ called tree
such that the following properties
\eqref{forest1}, \eqref{forest2},
\eqref{forest3},
\eqref{forest4},
\eqref{forest5}, and
\eqref{forest6}
hold.

For each $\fu\in \fU$ and each $\fp\in \fT(\fu)$
we have $\scI(\fp) \ne \scI(\fu)$ and
\begin{equation}\label{forest1}
    4\fp\lesssim \fu.
\end{equation}
For each $\fu \in \fU$ and each $\fp,\fp''\in \fT(\fu)$ and $\fp'\in \fP$
we have
\begin{equation}\label{forest2}
    \fp, \fp'' \in \mathfrak{T}(\fu), \fp \leq \fp' \leq \fp'' \implies \fp' \in \mathfrak{T}(\fu).
\end{equation}
We have
\begin{equation}\label{forest3}
   \|\sum_{\fu\in \fU} \mathbf{1}_{\scI(\fu)}\|_\infty \leq 2^n\,.
\end{equation}
We have for every $\fu\in \fU$
\begin{equation}\label{forest4}
    \dens_1(\fT(\fu))\le 2^{4a + 1-n}\, .
\end{equation}
We have for $\fu, \fu'\in \fU$ with $\fu\neq \fu'$ and $\fp\in \fT(\fu')$ with $\scI(\fp)\subset \scI(\fu)$ that
\begin{equation}\label{forest5}
    d_{\fp}(\fcc(\fp), \fcc(\fu))>2^{Z(n+1)}\, .
\end{equation}
We have for every $\fu\in \fU$ and $\fp\in \fT(\fu)$ that
\begin{equation}\label{forest6}
    B(\pc(\fp), 8D^{\ps(\fp)})\subset \scI(\fu).
\end{equation}

The following proposition is proved in \Cref{treesection}.
\begin{proposition}[forest operator]
\label{forest-operator}
\leanok
\lean{forest_operator}
\uses{forest-row-decomposition,row-bound,row-correlation,disjoint-row-support}
For any $n\ge 0$ and any $n$-forest $(\fU,\fT)$ we have for all $f,g: X \to \mathbb{C}$ with $|f| \le \mathbf{1}_F$ and $|g| \le \mathbf{1}_G$
$$
    | \int \overline{g(x)} \sum_{\fu\in \fU} \sum_{\fp\in \fT(\fu)} T_{\fp} f(x) \, \mathrm{d}\mu(x)|
$$
$$
    \le
    2^{440a^3}2^{-\frac{q-1}{q} n} \dens_2\left(\bigcup_{\fu\in \fU}\fT(\fu)\right)^{\frac{1}{q}-\frac{1}{2}} \|f\|_2 \|g\|_2 \,.
$$
\end{proposition}

\Cref{metric-space-Carleson} is formulated at the level of generality
for general kernels satisfying the mere H\"older regularity condition \eqref{eqkernel-y-smooth}. On the other hand, the cancellative condition \eqref{eq-vdc-cond} is a testing condition against more regular,
namely Lipschitz functions. To bridge the gap, we follow \cite{zk-polynomial} to observe a variant of \eqref{eq-vdc-cond} that we formulate
in the following proposition proved in \Cref{liphoel}.

Define
\begin{equation}
    \tau:=\frac 1a\, .
\end{equation}
Define for any open ball $B$ of radius $R$ in $X$ the $L^\infty$-normalized $\tau$-H\"older norm by
\begin{equation}
    \label{eq-Holder-norm}
    \|\varphi\|_{C^\tau(B)} = \sup_{x \in B} |\varphi(x)| + R^\tau \sup_{x,y \in B, x \neq y} \frac{|\varphi(x) - \varphi(y)|}{\rho(x,y)^\tau}\,.
\end{equation}

\begin{proposition}[Holder van der Corput]
    \label{Holder-van-der-Corput}
    \leanok
    \lean{holder_van_der_corput}
    \uses{Lipschitz-Holder-approximation}
     Let $z\in X$ and $R>0$ and set $B=B(z,R)$.
     Let $\varphi: X \to \mathbb{C}$ be
     supported on $B$ and satisfy $\|{\varphi}\|_{C^\tau(B(z, 2R))}<\infty$.
     Let $\mfa, \mfb \in \Mf$. Then
    \begin{equation}
        \label{eq-vdc-cond-tau-2}
        |\int e(\mfa(x)-{\mfb(x)})\varphi(x) dx|\le
         2^{7a} \mu(B) \|{\varphi}\|_{C^\tau(B(z, 2R))}
       (1 + d_{B}(\mfa,\mfb))^{-\frac{1}{2a^2+a^3}}
    \,.
    \end{equation}
\end{proposition}

We further formulate a classical Vitali covering result
and maximal function estimate that we need throughout several sections.
This following proposition will typically be applied to the absolute value of a complex valued function and be proved in \Cref{sec-hlm}. By a ball $B$ we mean a set $B(x,r)$ with $x\in X$
and $r>0$ as defined in \eqref{eq-define-ball}.
For a finite collection $\mathcal{B}$ of balls in $X$
and $1\le p< \infty$ define the measurable function $M_{\mathcal{B},p}u$ on $X$ by
\begin{equation}\label{def-hlm}
    M_{\mathcal{B},p}u(x):=\left(\sup_{B\in \mathcal{B}} \frac{\mathbf{1}_{B}(x)}{\mu(B)}\int _{B} |u(y)|^p\, d\mu(y)\right)^\frac 1p\, .
\end{equation}
Define further $M_{\mathcal{B}}:=M_{\mathcal{B},1}$.

\begin{proposition}[Hardy--Littlewood]
\label{Hardy-Littlewood}
\leanok
\lean{Finset.measure_biUnion_le_lintegral, hasStrongType_maximalFunction,
MeasureTheory.AEStronglyMeasurable.globalMaximalFunction, laverage_le_globalMaximalFunction, hasStrongType_globalMaximalFunction}
\uses{layer-cake-representation,covering-separable-space}
   Let $\mathcal{B}$ be a finite collection of balls in $X$.
If for some $\lambda>0$ and some measurable function $u:X\to [0,\infty)$ we have
\begin{equation}\label{eq-ball-assumption}
\int_{B} u(x)\, d\mu(x)\ge \lambda \mu(B)
\end{equation}
for each $B\in \mathcal{B}$, then
\begin{equation}\label{eq-besico}
    \lambda \mu(\bigcup \mathcal{B}) \le 2^{2a}\int_X u(x)\, d\mu(x)\, .
\end{equation}
For every measurable function $v$
and $1\le p_1<p_2$ we have
\begin{equation}\label{eq-hlm}
    \|M_{\mathcal{B},p_1} v\|_{p_2}\le 2^{2a}\frac{p_2}{p_2-p_1} \|v\|_{p_2}\, .
\end{equation}
Moreover, given any measurable bounded function $w: X \to \C$ there exists a measurable function $Mw: X \to [0, \infty)$ such that the following \eqref{eq-ball-av} and \eqref{eq-hlm-2} hold. For each ball $B \subset X$ and each $x \in B$
\begin{equation}
    \label{eq-ball-av}
    \frac{1}{\mu(B)} \int_{B} |w(y)| \, \mathrm{d}\mu(y) \le Mw(x)
\end{equation}
and for all $1 \le p_1 < p_2 \le \infty$
\begin{equation}
    \label{eq-hlm-2}
    \|M(w^{p_1})^{\frac{1}{p_1}}\|_{p_2} \le 2^{4a} \frac{p_2}{p_2-p_1}\|w\|_{p_2}\,.
\end{equation}

\end{proposition}

This completes the overview of the proof of \Cref{metric-space-Carleson}.

\subsection{Auxiliary lemmas}
\label{global-auxiliary-lemmas}
We close this section by recording some auxiliary lemmas about the objects defined in \Cref{overviewsection}, which will be used in multiple sections to follow.

First, we record an estimate for the metrical entropy numbers of balls in the space $\Mf$ equipped with any of the metrics $d_B$, following from the doubling property \eqref{thirddb}.

\begin{lemma}[ball metric entropy]
    \label{ball-metric-entropy}
    \leanok
    \lean{Θ.finite_and_mk_le_of_le_dist}
    Let $B' \subset X$ be a ball. Let $r > 0$, $\mfa \in \Mf$ and $k \in \mathbb{N}$. Suppose that $\mathcal{Z} \subset B_{B'}(\mfa, r2^k)$ satisfies that $\{B_{B'}(z,r)\mid z \in \mathcal{Z}\}$ is a collection of pairwise disjoint sets. Then
    $$ |\mathcal{Z}| \le 2^{ka}\,. $$
\end{lemma}

\begin{proof}
    \leanok
    By applying property \eqref{thirddb} $k$ times, we obtain a collection $\mathcal{Z}' \subset \Mf$ with $|\mathcal{Z}'| = 2^{ka}$ and
    $$
        B_{B'}(\mfa,r2^k) \subset \bigcup_{z' \in \mathcal{Z}'} B_{B'}(z', \frac{r}{2})\,.
    $$
    Then each $z \in \mathcal{Z}$ is contained in one of the balls $B(z', \frac{r}{2})$, but by the separation assumption no such ball contains more than one element of $\mathcal{Z}$. Thus $|\mathcal{Z}| \le |\mathcal{Z}'| = 2^{ka}$.
\end{proof}

The next lemma concerns monotonicity of the metrics $d_{B(c(I), \frac 14 D^{s(I)})}$ with respect to inclusion of cubes $I$ in a grid.

\begin{lemma}[monotone cube metrics]
    \label{monotone-cube-metrics}
    \leanok
    \lean{Grid.dist_mono, Grid.dist_strictMono}
    Let $(\mathcal{D}, c, s)$ be a grid structure. Denote for cubes $I \in \mathcal{D}$
    $$
        I^\circ := B(c(I), \frac{1}{4} D^{s(I)})\,.
    $$
    Let $I, J \in \mathcal{D}$ with $I \subset J$.
    Then for all $\mfa, \mfb \in\Mf$ we have
    $$
        d_{I^\circ}(\mfa, \mfb) \le d_{J^\circ}(\mfa, \mfb)\,,
    $$
    and if $I \ne J$ then we have
    $$
        d_{I^\circ}(\mfa, \mfb) \le 2^{-95a} d_{J^\circ}(\mfa, \mfb)\,.
    $$
\end{lemma}

\begin{proof}
    \leanok
    If $s(I) \ge s(J)$ then \eqref{dyadicproperty} and the assumption $I\subset J$ imply $I = J$. Then the lemma holds by reflexivity.

    If $s(J) \ge s(I)+1$, then using the monotonicity property \eqref{monotonedb}, \eqref{defineD} and \eqref{seconddb}, we get
    \begin{equation}
    \label{eq-dIJ-est}
        d_{I^\circ}(\mfa, \mfb) \le d_{B(c(I), 4 D^{s(I)})}(\mfa, \mfb) \le 2^{-100a} d_{B(c(I), 4D^{s(J)})}(\mfa, \mfb)\,.
    \end{equation}
    Using \eqref{eq-vol-sp-cube}, together with the inclusion $I \subset J$, we obtain
    $$
        c(I) \in I \subset J \subset B(c(J), 4 D^{s(J)})
    $$
    and consequently by the triangle inequality
    $$
        B(c(I), 4 D^{s(J)}) \subset B(c(J), 8 D^{s(J)})\,.
    $$
    Using this together with the monotonicity property \eqref{monotonedb} and \eqref{firstdb} in \eqref{eq-dIJ-est}, we obtain
    \begin{align*}
        d_{I^\circ}(\mfa, \mfb) &\le 2^{-100a} d_{B(c(J), 8D^{s(J)})}(\mfa, \mfb)\\
        &\le 2^{-100a + 5a} d_{B(c(J), \frac{1}{4}D^{s(J)})}(\mfa, \mfb)\\
        &= 2^{-95a}d_{J^\circ}(\mfa, \mfb)\,.
    \end{align*}
    This proves the second inequality claimed in the Lemma, from which the first follows since $a \ge 4$ and hence $2^{-95a} \le 1$.
\end{proof}

We also record the following basic estimates for the kernels $K_s$.

\begin{lemma}[kernel summand]
\label{kernel-summand}
\leanok
\lean{dist_mem_Icc_of_Ks_ne_zero, enorm_Ks_le, enorm_Ks_sub_Ks_le}
    Let $-S\le s\le S$ and $x,y,y'\in X$.
    If $K_s(x,y)\neq 0$, then we have
    \begin{equation}\label{supp-Ks}
      \frac{1}{4} D^{s-1} \leq \rho(x,y) \leq \frac{1}{2} D^s\, .
    \end{equation}
    We have
    \begin{equation}
       \label{eq-Ks-size}
        |K_s(x,y)|\le \frac{2^{102 a^3}}{\mu(B(x, D^{s}))}\,
    \end{equation}
    and
    \begin{equation}
        \label{eq-Ks-smooth}
        |K_s(x,y)-K_s(x, y')|\le \frac{2^{127a^3}}{\mu(B(x, D^{s}))}
        \left(\frac{ \rho(y,y')}{D^s}\right)^{\frac 1a}\,.
    \end{equation}
\end{lemma}

\begin{proof}
    \leanok
    By Definition \eqref{defks}, the function $K_s$ is the product of
    $K$ with a function which is supported in the set of all
    $x,y$ satisfying \eqref{supp-Ks}. This proves
    \eqref{supp-Ks}.

    Using \eqref{eqkernel-size} and the lower bound in \eqref{supp-Ks}
    we obtain
    \begin{equation}
        \label{eqkernel-size-Ks}
        |K_s(x,y)|\le |K(x,y)|\le \frac{2^{a^3}}{\mu(B(x,\frac 14 D^{s-1}))}
    \end{equation}
    Using $D=2^{100a^2}$
    and the doubling property \eqref{doublingx} $2 +100a^2$ times estimates
    the last display by
    \begin{equation}
        \label{eq-Ks-aux}
        \le \frac{2^{2a+101a^3}}{\mu(B(x, D^{s}))}\, .
    \end{equation}
    Using $a\ge 4$ proves \eqref{eq-Ks-size}.

    To prove \eqref{eq-Ks-smooth} when $2\rho(y,y') > \rho(x,y)$, use the lower bound in
    \eqref{supp-Ks} and the inequality $2\rho(y,y') > \frac{1}{4}D^{s-1}$. Then \eqref{eq-Ks-smooth} follows from
    the triangle inequality, \eqref{eq-Ks-size} and $a \ge 4$.

    If $2\rho(y,y') \le \rho(x,y)$, we rewrite $|K_s(x,y)-K_s(x, y')|$ as
    \begin{equation}
        |(K(x,y)-K(x,y')) \psi(D^{-s}\rho(x,y)) +
        K(x,y)(\psi(D^{-s}\rho(x,y))-\psi(D^{-s}\rho(x,y')))|\,.
    \end{equation}
    An upper bound for $|K(x,y)-K(x, y')|$ is obtained similarly to the proof of
    \eqref{eq-Ks-size}, using \eqref{eqkernel-y-smooth} and the lower bound in \eqref{supp-Ks}
    \begin{equation}
        |K(x,y)-K(x, y')|\le \frac{2^{a^3}}{\mu(B(x, \frac 14 D^{s-1}))}
        \left(\frac{ \rho(y,y')}{\frac 14 D^{s-1}}\right)^{\frac 1a}\,.
    \end{equation}
    As above, this is estimated by
    \begin{equation}
       \le \frac{4D 2^{2a+101a^3}}{\mu(B(x, D^{s}))}
        \left(\frac{ \rho(y,y')}{D^{s}}\right)^{\frac 1a}
         = \frac{2^{2+2a+100a^2+101a^3}}{\mu(B(x, D^{s}))}
        \left(\frac{ \rho(y,y')}{D^{s}}\right)^{\frac 1a}\,.
    \end{equation}
    We have the trivial bound $|\psi(D^{-s}\rho(x,y))| \leq 1$, and \eqref{eq-Ks-aux}
    provides a bound for $|K(x,y)|$. Finally, we show that
    \begin{equation}
        |\psi(D^{-s}\rho(x,y))-\psi(D^{-s}\rho(x,y'))|\le
        4D \left(\frac{\rho(y, y')}{D ^ s}\right)^{\frac 1a}
    \end{equation}
    by considering separately the cases $\rho(y,y')/D^s \ge 1$ and $\rho(y,y')/D^s < 1$. In the
    former case, the inequality is trivial; in the latter case, it follows from the
    fact that $\psi$ is Lipschitz with constant $4D$.

    Combining the above bounds and using $a\ge 4$ proves \eqref{eq-Ks-smooth} when
    $2\rho(y,y') \le \rho(x,y)$.
\end{proof}

\section{Proof of Metric Space Carleson}
\label{thmfromproplinear}

In this section we prove \Cref{metric-space-Carleson} and \Cref{linearised-metric-Carleson}.

Let $(X, \rho, \mu, a)$ be a doubling metric measure space and $\Theta$ a cancellative, compatible collection of functions on $X$. Let $K$ be a one-sided Calderon-Zygmund kernel.

We begin by proving some continuity properties of the integrand in \eqref{def-main-op}.

\begin{lemma}[int continuous]
    \label{int-continuous}
    \leanok
    \lean{continuous_carlesonOperatorIntegrand, rightContinuous_carlesonOperatorIntegrand, leftContinuous_carlesonOperatorIntegrand, measurable_carlesonOperatorIntegrand, enorm_carlesonOperatorIntegrand_le}
    Let $f$ be a measurable function with $|f| \le 1$. Then the function
    \[
        G: X \times \Theta \times (0,\infty) \times (0, \infty) \to \mathbb{C}
    \]
    \[
        G(x, \mfa, R_1, R_2) := \int_{R_1 < \rho(x,y) < R_2} K(x,y) f(y) e(\mfa(y)) \, \mathrm{d}\mu(y)
    \]
    is continuous in $\mfa$ for fixed $x, R_1, R_2$. It is right-continuous in $R_1$ for fixed $x, \vartheta, R_2$ and left-continuous in $R_2$ for fixed $x, \vartheta, R_1$. Finally, it is measurable in $x$ and bounded for fixed $\vartheta, R_1, R_2$.
\end{lemma}

\begin{proof}
    \leanok
    Measurability in $x$ follows from joint measurability of
    \[
        K(x,y) \mathbf{1}_{B(x,R_2) \setminus B(x,R_1)}(y)
    \]
    in $x$ and $y$ and (part of the proof of) Fubini's theorem.

    (Partial) continuity in $R_1$ and $R_2$ is also a standard lemma in measure theory. It follows for example by splitting the integrand as $F_1 - F_{-1} + iF_i - iF_{-i}$ for positive, disjointly supported functions $F_{-}$ and applying the monotone convergence theorem to each summand.

    For continuity in $\mfa$ note that
    \[
        |G(x, \mfa, R_1, R_2) - G(x, \mfa', R_1, R_2)|
    \]
    \[
        = \left| \int_{R_1 < \rho(x,y) < R_2} K(x,y) f(y) (e(\mfa(y)) - e(\mfa'(y))) \, \mathrm{d}\mu(y) \right|
    \]
    \[
        \le \int_{R_1 < \rho(x,y) < R_2} |K(x,y)| |f(y)| |e(\mfa(y) - \mfa'(y) - \mfa(o) + \mfa'(o)) - 1| \, \mathrm{d}\mu(y).
    \]
    By $1$-Lipschitz continuity of $e^{ix}$, the property \eqref{osccontrol} of the metrics $d$, the kernel upper bound \eqref{eqkernel-size} and $|f| \le 1$ this is
    \[
        \le \mu(B(x, R_2)) \sup_{R_1 < \rho(x,y) < R_2} \frac{2^{a^3}}{V(x,y)} d_{B(x,\rho(o,x)+R_2)}(\vartheta, \vartheta').
    \]
    If $R_1 < \rho(x,y) < R_2$ then there exists $n$ with $B(x,R_2) \subset B(x, 2^n \rho(x,y))$ and $2^n \le 2 R_2/R_1$. Hence, by the doubling property \eqref{doublingx},
    \[
        V(x,y) = \mu(B(x, \rho(x,y))) \ge 2^{-an} \mu(B(x, R_2)) \ge (2R_2 / R_1)^{-a} \mu(B(x, R_2)).
    \]
    Hence
    \[
        |G(x, \mfa, R_1, R_2) - G(x, \mfa', R_1, R_2)| \le 2^{a^3} \Big(\frac{2R_2}{R_1}\Big)^a d_{B(x, \rho(o,x)+R_2)}(\mfa, \mfa').
    \]
    Since the topology on $\Mf$ is the one induced by any of the metrics $d_B$, continuity follows.

    Finally, for boundedness as a function of $x$ note that we also have by a similar computation using $|e(\mfa)|=1$
    \[
        |G(x, \mfa, R_1, R_2)| \le 2^{a^3} \Big(\frac{2R_2}{R_1}\Big)^a.
    \]
\end{proof}

We now prove Theorem \ref{metric-space-Carleson} using Theorem \ref{linearised-metric-Carleson}.

\begin{proof}[Proof of \Cref{metric-space-Carleson}]
    \leanok
    \proves{metric-space-Carleson}
    Let Borel sets $F$, $G$ in $X$ be given. Let a Borel function $f: X \to \mathbb{C}$ with $f \le \mathbf{1}_F$ be given.

    Let $\Mf' \subset \Mf$ be a countable dense set. By Lemma \ref{int-continuous} we have
    \[
        \sup_{\mfa\in\Mf} \sup_{0 < R_1 < R_2}\left| \int_{R_1 < \rho(x,y) < R_2} K(x,y) f(y) e(\mfa(y)) \, \mathrm{d}\mu(y) \right|
    \]
    \[
        = \sup_{\mfa\in\Mf'} \sup_{0 < R_1 < R_2, R_i \in \mathbb{Q}}\left| \int_{R_1 < \rho(x,y) < R_2} K(x,y) f(y) e(\mfa(y)) \, \mathrm{d}\mu(y) \right|.
    \]
    Consider an enumeration of $\Mf'$ and let $\Mf_n$ be the finite set consisting of the first $n+1$ functions in the enumeration. Then by the monotone convergence theorem
    \[
        \Big| \int_G Tf \, \mathrm{d}\mu\Big|
    \]
    \[
        = \int_G \sup_{\mfa\in\Mf'} \sup_{0 < R_1 < R_2, R_i \in \mathbb{Q}}\left| \int_{R_1 < \rho(x,y) < R_2} K(x,y) f(y) e(\mfa(y)) \, \mathrm{d}\mu(y) \right| \, \mathrm{d}\mu(x)
    \]
    \[
        = \lim_{n \to \infty} \int_G \sup_{\mfa\in\Mf_n} \sup_{0 < R_1 < R_2, R_i \in \mathbb{Q}}\left| \int_{R_1 < \rho(x,y) < R_2} K(x,y) f(y) e(\mfa(y)) \, \mathrm{d}\mu(y) \right| \, \mathrm{d}\mu(x).
    \]
    For each $n$, let $Q_n(x)$ be the measurable function specifying the maximizer in the supremum in $\mfa$ in the previous line. It is possible to construct such a measurable function of $x$, because the function inside the supremum is measurable in $x$, being a countable supremum of measurable functions, and because $\Mf_n$ is finite. (For example one may fix some order on the finite set $\Mf_n$ and pick the smallest maximizer with respect to that order. This is a measurable choice function.) Then the previous display becomes
    \[
        \le \lim_{n \to \infty} \int_G \sup_{0 < R_1 < R_2} \left| \int_{R_1 < \rho(x,y) < R_2} K(x,y) f(y) e(Q_n(x)(y)) \, \mathrm{d}\mu(y) \right| \, \mathrm{d}\mu(x).
    \]
    \[
        = \lim_{n \to \infty} \int_G T_{Q_n} f \, \mathrm{d}\mu.
    \]
    It remains to verify the assumptions of Theorem \ref{linearised-metric-Carleson}, which when applied here completes the proof.

    The assumptions of Theorem \ref{metric-space-Carleson} and Theorem \ref{linearised-metric-Carleson} are the same,
    with the exception of the assumption \eqref{linnontanbound}.
    The assumption \eqref{linnontanbound} however is weaker than the assumption \eqref{nontanbound}
    of Theorem \ref{metric-space-Carleson}.
    Indeed, setting for fixed $x, \mfa$ the outer radius $R_2$ in \eqref{def-tang-unm-op} to
    $\min\{R_2', R_Q(\mfa, x')\}$, where $R_2'$ is the outer radius in \eqref{def-lin-star-op},
    shows that \eqref{def-lin-star-op} is smaller than or equal to \eqref{def-tang-unm-op}.
    Thus we can apply Theorem \ref{linearised-metric-Carleson}, which completes the proof.
\end{proof}

We continue with the proof of Theorem \ref{linearised-metric-Carleson}, via a series of reductions to simpler lemmas.

Let a measurable function $Q$ with finite range be given.

\begin{proof}[Proof of \Cref{linearised-metric-Carleson}]
\leanok
\proves{linearised-metric-Carleson}
Let Borel sets $F$, $G$ in $X$ with finite measure be given. Let a Borel function $f: X \to \mathbb{C}$ with $f \le \mathbf{1}_F$ be given.

For each $0 < R_1,R_2,R$, we define $T_{R_1,R_2,R}f$ as in \eqref{TRR}.
By Lemma \ref{int-continuous},
$T_{R_1,R_2,R}f$
is measurable and bounded, and we clearly have for each $x \in X$
\[
    T_Q f(x) = \lim_{n \to \infty} \sup_{2^{-n} < R_1 < R_2 < 2^n} T_{R_1, R_2, 2^n}f(x).
\]
For each $x$ and all $f$, the functions $\sup_{2^{-n} < R_1 < R_2 < 2^n} T_{R_1, R_2, 2^n} f(x)$ are measurable by Lemma \ref{int-continuous} and form an increasing sequence in $n$. By the monotone convergence theorem, the claimed estimate \eqref{linresweak} then follows from Lemma \ref{R-truncation}.
\end{proof}

\begin{lemma}[R truncation]
    \label{R-truncation}
    \leanok
    \lean{R_truncation}
    \uses{S-truncation}

    Let $F$, $G$ be Borel sets in $X$. Let $f:X\to \C$ be a Borel function with $|f|\le 1_F$. Then for all $R\in 2^\N$ we have
    \begin{equation} \label{Rcut}
        \int \mathbf{1}_G \sup_{1/R<R_1<R_2<R} |T_{R_1,R_2,R} f(x)|\, d\mu(x) \le \frac{2^{443a^3}}{(q-1)^6} \mu(G)^{1-\frac{1}{q}} \mu(F)^{\frac{1}{q}},
    \end{equation}
    where
    \begin{equation}\label{TRR}
        T_{R_1,R_2,R} f(x)= \mathbf{1}_{B(o,R)}(x)
        \int_{R_1 < \rho(x,y) < R_2} K(x,y) f(y) e(Q(x)(y)) \, \mathrm{d}\mu(y) .
    \end{equation}
\end{lemma}

\begin{proof}
\leanok
Let $F,G,f$ as in the lemma be given. Let $R\in 2^\N$ be given.
By replacing $G$ with $G\cap B(o,R)$ if necessary, a replacement that does not change the conclusion of the Lemma \ref{R-truncation}, it suffices to show \eqref{Rcut} under the assumption that $G$ is contained in $B(o,R)$ and thus bounded. We make this assumption. For every $x\in G$ and $R_2 < R$, the domain of integration in \eqref{TRR} is contained in $B(o,2R)$. By replacing $F$ with $F\cap B(o,2R)$ if necessary, and correspondingly restricting $f$ to $B(o, 2R)$, it suffices to show \eqref{Rcut} under the assumption that $F$ is contained in $B(o,2R)$ and thus bounded. We make this assumption.

Using the definition \eqref{defks} of $K_s$ and the partition of unity \eqref{eq-psisum},
we observe that
for fixed $R_1<R_2$ we have
\begin{equation}
    \label{KKs}
    K(x,y)\mathbf{1}_{R_1<\rho(x,y)<R_2}=\sum_{s=s_1-2}^{s_2+2} K_s(x,y) \mathbf{1}_{R_1<\rho(x,y)<R_2},
\end{equation}
where $s_1=\lfloor\log_D 2R_1\rfloor+3$ and $s_2=\lceil\log_D 4R_2\rceil-2$.
To obtain the identity \eqref{KKs}, we have used that on the set where $R_1<\rho(x,y)<R_2$ the kernels $K_s$ vanish unless $s$ is inside the interval of summation in \eqref{KKs}. Similarly, we observe
\begin{equation}
    \label{KsrhoKs}
    \sum_{s=s_1}^{s_2} K_s(x,y) \mathbf{1}_{R_1<\rho(x,y)<R_2} = \sum_{s=s_1}^{s_2} K_s(x,y),
\end{equation}
because on the support of $K_s$ with $s_1\le s\le s_2$ we have necessarily $R_1<\rho(x,y)<R_2$. We thus express \eqref{TRR} as the sum of
\begin{equation}\label{middles}
T_{s_1,s_2}f(x):=\sum_{s_1 \le s\le s_2}
\int K_s(x,y) f(y) e(Q(x)(y)) \, \mathrm{d}\mu(y)
\end{equation}
and
\begin{equation}\label{boundarys}
\sum_{s=s_1-2,s_1-1, s_2+1, s_2+2}
\int_{R_1 < \rho(x,y) < R_2} K_s(x,y) f(y) e(Q(x)(y)) \,
 \mathrm{d}\mu(y),
\end{equation}

We apply the triangle inequality and estimate \eqref{Rcut} separately with $T_{R_1,R_2,R}$ replaced by \eqref{middles} and by \eqref{boundarys}.
To handle the case \eqref{middles}, we employ \Cref{S-truncation}. Here, we utilize the fact that if $\frac 1R\le R_1\le R_2\le R$, then $s_1$ and $s_2$ as in \eqref{middles} are in an interval $[-S,S]$ for some sufficiently large $S$ depending on $R$.
To handle the case \eqref{boundarys}, we use the triangle inequality and the properties \eqref{supp-Ks}, \eqref{eq-Ks-size} of $K_s$  and $|f| \le \mathbf{1}_F$ to obtain for arbitrary $s$ the inequality
\begin{multline}
\left|\int_{R_1 < \rho(x,y) < R_2} K_s(x,y) f(y) e(Q(x)(y)) \, \mathrm{d}\mu(y)\right|\\
\leq \frac{2^{102 a^3}}{\mu(B(x, D^{s}))}
 \int_{B(x, D^s)} \mathbf{1}_F(y) \, \mathrm{d}\mu(y)
\leq 2^{102 a^3} M\mathbf{1}_F(x),
\end{multline}
where $M\mathbf{1}_F$ is as defined in \Cref{Hardy-Littlewood}.
Now the left-hand side of \eqref{Rcut}, with $T_{R_1,R_2,R}$ replaced by \eqref{boundarys}, can be estimated using H\"older's inequality and \Cref{Hardy-Littlewood} by
$$
    2^{102a^3+2}\int \mathbf{1}_{G}(x) M\mathbf{1}_F(x)\, d\mu(x)
    \le \frac{2^{102a^3+4a+3}q}{q-1}\mu(G)^{1-\frac{1}{q}} \mu(F)^{\frac{1}{q}}\,.
$$
Summing the contributions from \eqref{middles} and \eqref{boundarys} completes the proof.
\end{proof}

\begin{lemma}[S truncation]
    \label{S-truncation}
    \leanok
    \lean{S_truncation}
    \uses{Hardy-Littlewood, linearized-truncation}
    Let $F$, $G$ be bounded Borel sets in $X$.
    Let $f:X\to \C$ be a Borel function with $|f|\le 1_F$. Then for all $S\in\Z$ we have
    \begin{equation} \label{Scut}
        \int \mathbf{1}_G(x) \sup_{-S\le s_1\le s_2\le S} |T_{s_1,s_2} f(x)|\, d\mu(x)
        \le \frac{2^{442a^3+2}}{(q-1)^6} \mu(G)^{1-\frac{1}{q}} \mu(F)^{\frac{1}{q}},
    \end{equation}
    where
    \begin{equation}\label{Tss}
        T_{s_1,s_2} f(x) = \sum_{s_1\le s \le s_2} \int_X K_s(x,y) f(y) e(Q(x)(y)) \, \mathrm{d}\mu(y).
    \end{equation}
\end{lemma}

\begin{proof}[Proof of \Cref{S-truncation}]
    \leanok
    We reduce \Cref{S-truncation} to \Cref{linearized-truncation}. For each $x$, let $\sigma_1(x)$ be the
    minimal element $s'\in [-S,S]$ such that
    \[
    \max_{s'\le s_2\le S} |T_{s',s_2} f(x)| = \max_{-S\le s_1\le s_2\le S} |T_{s_1,s_2} f(x)| =: T_{1,x}.
    \]
    Similarly, let ${\sigma}_2(x)$ be the minimal element $s''\in [-S,S]$ such that
    \[
    |T_{\sigma_2(x), s''} f(x)| = T_{1,x}\,.
    \]
    With these choices, and noting that with the definition of $T_{2, \sigma_1, \sigma_2}$ from \eqref{middles1}
    \begin{equation*}
    T_{\sigma_1(x),\sigma_2(x)} f(x)=T_{2,\sigma_1,\sigma_2} f(x),
    \end{equation*}
    we conclude that the left-hand side of \eqref{Scut} and \eqref{Sqlin} are equal.
    Therefore, \Cref{S-truncation} follows from \Cref{linearized-truncation}.
\end{proof}

\begin{lemma}[linearized truncation]
    \label{linearized-truncation}
    \leanok
    \lean{linearized_truncation}
    \uses{finitary-Carleson}
    Let $\sigma_1,\sigma_2\colon X\to \mathbb{Z}$ be measurable functions with finite range and $\sigma_1\leq \sigma_2$.
    Then we have
    \begin{equation} \label{Sqlin}
        \int \mathbf{1}_{G}(x)
        \left| {T}_{2,\sigma_1,\sigma_2}f(x)\right|\, d\mu(x)
        \le \frac{2^{442a^3+2}}{(q-1)^6} \mu(G)^{1-\frac{1}{q}} \mu(F)^{\frac{1}{q}},
    \end{equation}
    with
    \begin{equation}\label{middles1}
        {T}_{2,\sigma_1,\sigma_2}f(x)=\sum_{\sigma_1(x) \le s\le \sigma_2(x)}
        \int K_s(x,y) f(y) e(\tQ(x)(y)) \, \mathrm{d}\mu(y)\,.
    \end{equation}
\end{lemma}

\begin{proof}[Proof of \Cref{linearized-truncation}]
\leanok
Fix $\sigma_1$, $\sigma_2$ and $\tQ$ as in the lemma. Applying \Cref{finitary-Carleson} recursively, we obtain a sequence of sets $G_n$ with $G_0=G$ and, for each $n\ge 0$, $G_{n+1} \subset G_n$, $\mu(G_{n})\le 2^{-n} \mu(G)$ and
\begin{equation*}
    \int \mathbf{1}_{G_{n}\setminus G_{n+1}}(x)
    \left| {T}_{2,\sigma_1,\sigma_2} f(x) \right|\, d\mu(x)
\end{equation*}
\begin{equation}
    \le \frac{2^{442a^3}}{(q-1)^5} \mu(G_n)^{1 - \frac{1}{q}} \mu(F)^{\frac{1}{q}}.
\end{equation}
Adding the first $n$ of these inequalities, we obtain by bounding a geometric series
    \begin{equation} \label{Sqcut2}
    \int \mathbf{1}_{G\setminus G_{n}}(x)
\left| {T}_{2,\sigma_1,\sigma_2}f(x) \right|\, d\mu(x)
\le \frac{2^{442a^3+2}}{(q-1)^6} \mu(G)^{1-\frac{1}{q}} \mu(F)^{\frac{1}{q}}.
\end{equation}
As the integrand is non-negative and non-decreasing in $n$, we obtain by the monotone convergence theorem
 \begin{equation} \label{Sqcut3}
    \int \mathbf{1}_{G}(x)
\left| {T}_{2,\sigma_1,\sigma_2}f(x) \right|\, d\mu(x)
\le \frac{2^{442a^3+2}}{(q-1)^6} \mu(G)^{1-\frac{1}{q}} \mu(F)^{\frac{1}{q}}.
\end{equation}
This completes the proof of \Cref{linearized-truncation}
and thus \Cref{metric-space-Carleson}.
\end{proof}

\section{Proof of Finitary Carleson}
\label{christsection}

To prove Proposition
\ref{finitary-Carleson}, we already fixed in \Cref{overviewsection}
measurable functions ${\sigma_1},\sigma_2, \tQ$ and Borel sets $F,G$. We have also
defined $S$ to be the smallest
integer such that the ranges of
$\sigma_1$ and $\sigma_2$ are contained in $[-S,S]$ and $F$ and $G$ are contained
in the ball $B(o, \frac 14 D^S)$.

The proof of the next lemma is done in \Cref{subsecdyadic},
following the construction of dyadic cubes in \cite[\S 3]{christ1990b}.

\begin{lemma}[grid existence]
    \label{grid-existence}
    \leanok
    \lean{grid_existence}
    \uses{counting-balls,boundary-measure} There exists a grid structure $(\mathcal{D}, c,s)$.
\end{lemma}

The next lemma, which we prove in \Cref{subsectiles}, should be compared
with the construction in \cite[Lemma 2.12]{zk-polynomial}.

\begin{lemma}[tile structure]
    \label{tile-structure}
    \leanok
    \lean{tile_existence}
    \uses{ball-metric-entropy,frequency-ball-cover,disjoint-frequency-cubes,frequency-cube-cover}
        For a given grid structure $(\mathcal{D}, c,s)$, there exists a tile structure
        $(\fP,\scI,\fc,\fcc,\pc,\ps)$.
\end{lemma}

Choose a grid structure $(\mathcal{D}, c,s)$ with \Cref{grid-existence} and a tile structure for this
grid structure $(\fP,\scI,\fc,\fcc,\pc,\ps)$ with \Cref{tile-structure}.
Applying \Cref{discrete-Carleson}, we obtain a Borel set $G'$ in $X$ with $2\mu(G')\leq \mu(G)$ such that for all Borel functions $f:X\to \C$ with $|f|\le \mathbf{1}_F$
we have \eqref{disclesssim}.

\begin{lemma}[tile sum operator]
    \label{tile-sum-operator}
    \leanok
    \lean{tile_sum_operator, integrable_tile_sum_operator}
    We have for all $x\in G\setminus G'$
    \begin{equation}\label{eq-sump}
        \sum_{\fp\in \fP}T_{\fp} f(x)= \sum_{s=\sigma_1(x)}^{\sigma_2(x)}
        \int K_{s}(x,y) f(y) e(\tQ(x)(y)-\tQ(x)(x))\, d\mu(y).
    \end{equation}
\end{lemma}
\begin{proof}
    \leanok
    Fix $x\in G\setminus G'$.
    Sorting the tiles $\fp$ on the left-hand-side of \eqref{eq-sump} by the value $\ps(\fp)\in [-S,S]$,
    it suffices to prove for every $-S\le s\le S$ that
    \begin{equation}\label{outsump}
        \sum_{\fp\in \fP: \ps(\fp)=s}T_{\fp} f(x)=0
    \end{equation}
    if $s\not\in [\sigma_1(x), \sigma_2(x)]$ and
    \begin{equation}\label{insump}
        \sum_{\fp\in \fP: \ps(\fp)=s}T_{\fp} f(x)=
        \int K_{s}(x,y) f(y) e(\tQ(x)(y) - \tQ(x)(x))\, d\mu(y).
    \end{equation}
    if $s\in [\sigma_1(x),\sigma_2(x)]$.
    If $s\not\in [\sigma_1(x), \sigma_2(x)]$, then by definition of $E(\fp)$ we have
    $x\not\in E(\fp)$ for any $\fp$ with $\ps(\fp)=s$ and thus $T_{\fp} f(x)=0$. This proves
    \eqref{outsump}.

    Now assume $s\in [\sigma_1(x),\sigma_2(x)]$.
    By \eqref{coverdyadic}, \eqref{subsetmaxcube}, \eqref{eq-vol-sp-cube}, the fact that $c(I_0) = o$ and $G\subset B(o,\frac 14 D^S)$, there is at least
    one $I\in \mathcal{D}$ with $s(I)=s$ and $x\in I$.
    By \eqref{dyadicproperty}, this $I$ is unique. By \eqref{eq-dis-freq-cover}, there is precisely one $\fp\in \fP(I)$ such that
    $\tQ(x)\in \fc(\fp)$. Hence there is precisely one $\fp\in \fP$ with $\ps(\fp)=s$ such that
    $x\in E(\fp)$. For this $\fp$, the value $T_{\fp}(x)$ by its definition in \eqref{definetp}
    equals the right-hand side of \eqref{insump}. This proves the lemma.
\end{proof}

We use this to prove \Cref{finitary-Carleson}.
\begin{proof}[Proof of \Cref{finitary-Carleson}]
\proves{finitary-Carleson}
\leanok
We now estimate with \Cref{tile-sum-operator} and \Cref{discrete-Carleson}
\begin{equation}
 \int_{G \setminus G'} \left|\sum_{s={\sigma_1}(x)}^{{\sigma_2}(x)} \int K_s(x,y) f(y) e(\tQ(x)(y)) \, \mathrm{d}\mu(y)\right| \mathrm{d}\mu(x)
\end{equation}
\begin{equation}
 =\int_{G \setminus G'} \left|\sum_{s={\sigma_1}(x)}^{{\sigma_2}(x)} \int K_s(x,y) f(y) e(\tQ(x)(y) - \tQ(x)(x))\mathrm{d}\mu(y)\right| \mathrm{d}\mu(x)
\end{equation}
\begin{equation}
 =\int_{G \setminus G'} \left|\sum_{\fp\in \fP}T_{\fp} f(x)\right| \mathrm{d}\mu(x)
 \le \frac{2^{442a^3}}{(q-1)^5} \mu(G)^{1 - \frac{1}{q}} \mu(F)^{\frac{1}{q}} \,.
\end{equation}
This proves \eqref{eq-linearized} for the chosen set $G'$ and arbitrary $f$ and thus completes the proof of Proposition
\ref{finitary-Carleson}.
\end{proof}

\subsection{Proof of Grid Existence Lemma}
\label{subsecdyadic}

We begin with the construction of the centers of the dyadic cubes.
\begin{lemma}[counting balls]
    \label{counting-balls}
    \leanok
    \lean{counting_balls}
    Let $-S\le k\le S$. Consider $Y\subset X$ such that for any $y\in Y$, we have
    \begin{equation}\label{ybinb}
    y\in B(o,4D^S-D^k),
    \end{equation}
    furthermore, for any $y'\in Y$ with $y\neq y'$, we have
    \begin{equation} \label{eq-disj-yballs}
        B(y,D^k)\cap B(y',D^k)=\emptyset.
    \end{equation}
    Then the cardinality of $Y$ is bounded by
    \begin{equation}\label{boundY}
        |Y|\le 2^{3a + 200Sa^3}\, .
    \end{equation}
\end{lemma}

\begin{proof}
    \leanok
Let $k$ and $Y$ be given. By applying the doubling property \eqref{doublingx} inductively, we have for each integer $j\ge 0$
\begin{equation}\label{jballs}
    \mu(B(y,2^{j}D^k))\le 2^{aj} \mu(B(y,D^k))\, .
\end{equation}
Since $X$ is the union of the balls $B(y,2^{j}D^k)$ and $\mu$ is not zero, at least one of the balls $B(y,2^{j}D^k)$ has positive measure, thus $B(y,D^k)$ has positive measure.

Applying \eqref{jballs} for $j' = \ln_2(8D^{2S}) = 3 + 2S \cdot 100a^2$ by \eqref{defineD}, using $-S\le k\le S$, $y\in B(o,4D^S)$, and the triangle inequality, we have
\begin{equation}
    B(o, 4D^S) \subset B(y, 8D^S) \subset B(y,2^{j'}D^k) \, .
\end{equation}
Using the disjointedness of the balls in \eqref{eq-disj-yballs}, \eqref{ybinb}, and the triangle inequality for $\rho$, we obtain
\begin{equation}
|Y|\mu(B(o,4D^S))\le 2^{j'a}\sum_{y\in Y}\mu(B(y,D^k))
\end{equation}
\begin{equation}
\le
2^{j'a}\mu(\bigcup_{y\in Y}B(y,D^k))
\le 2^{j'a}\mu(o,4D^S)\, .
\end{equation}
As $\mu(o,4D^S)$ is not zero, the lemma follows.
\end{proof}

For each $-S\le k\le S$, let $Y_k$ be a set of
maximal cardinality in $X$ such that $Y=Y_k$ satisfies
the properties \eqref{ybinb} and \eqref{eq-disj-yballs} and such that $o \in Y_k$.
By the upper bound of \Cref{counting-balls}, such a set exists.

For each $-S\le k\le S$, choose an enumeration of the points in the finite set $Y_k$ and thus a total
order $<$ on $Y_{k}$.

\begin{lemma}[cover big ball]
    \label{cover-big-ball}
    \leanok
    \lean{cover_big_ball}
    For each $-S\le k\le S$, the ball
    $B(o, 4D^S-D^k)$ is contained
    in the union of the balls $B(y,2D^k)$ with $y\in Y_k$.
\end{lemma}

\begin{proof}
\leanok
Let $x$ be any point of $B(o, 4D^S-D^k)$. By maximality of $|Y_k|$, the ball
$B(x, D^k)$ intersects one of the balls
$B(y, D^k)$ with $y\in Y_k$. By the triangle
inequality, $x\in B(y,2D^k)$.
\end{proof}

Define the set
\begin{equation}
    \mathcal{C}:= \{(y,k): -S\le k\le S, y\in Y_k\}\,
\end{equation}
We totally order the set $\mathcal{C}$ lexicographically by setting
$(y,k)<(y',k')$ if $k< k'$ or both $k=k'$ and $y<y'$.
In what follows, we define recursively in the sense of this order a function
\begin{equation}
    (I_1,I_2,I_3): \mathcal{C}\to \mathcal{P}(X)\times \mathcal{P}(X)\times \mathcal{P}(X)\, .
\end{equation}

Assume the sets ${I}_j(y',k')$ have already been defined for $j=1,2,3$ if $k'<k$ and if $k=k'$ and $y'<y$.

If $k=-S$, define for $j\in \{1,2\}$ the set
${I}_j(y,k)$ to be $B(y,jD^{-S})$.
If $-S<k$, define for $j\in \{1,2\}$
and $y\in Y_k$ the set ${I}_j(y,k)$ to be
\begin{equation}\label{defineij}
\bigcup\{I_3(y',k-1):
y'\in Y_{k-1}\cap B(y,jD^k)\}.
\end{equation}
Define for {$-S\leq k\leq S$} and $y\in Y_k$
\begin{equation}\label{definei3}
I_3(y,k):={I_1}(y,k)\cup \left[{I_2}(y,k)\setminus \left[X_k\cup \bigcup\{I_3(y',k):y'\in Y_{k}, y'<y\})\right]\right]
\end{equation}
with
\begin{equation}
      X_{k}:=\bigcup\{I_1(y', k):y'\in Y_{k}\}.
\end{equation}

\begin{lemma}[basic grid structure]
    \label{basic-grid-structure}
    \uses{cover-big-ball}
    \leanok
    \lean{I1_prop_1,I3_prop_1,I2_prop_2,I3_prop_2,I3_prop_3_1,I3_prop_3_2}
    For each $-S\le k\le S$ and $1\le j\le 3$
    the following holds.

    If $j\neq 2$ and for some $x\in X$ and $y_1,y_2\in Y_k$ we have
    \begin{equation}\label{disji}
        x\in I_j(y_1,k)\cap I_j(y_2,k),
    \end{equation}
    then $y_1=y_2$.

    If $j\neq 1$, then
    \begin{equation}\label{unioni}
    B(o, 4D^S-2D^k)\subset \bigcup_{y\in Y_k} I_j(y,k)\, .
    \end{equation}
    We have for each $y\in Y_k$,
    \begin{equation}\label{squeezedyadic}
        B(y,\frac 12 D^k) \subset I_3(y,k)\subset
        B(y,4D^k).
    \end{equation}
\end{lemma}

\begin{proof}
\leanok
We prove these statements simultaneously by induction on the ordered set of pairs $(y,k)$. Let $-S\le k\le S$.

We first consider \eqref{disji} for $j=1$. If $k=-S$, disjointedness of the sets $I_1(y,-S)$ follows by definition of $I_1$ and $Y_k$. If $k>-S$, assume $x$ is in $I_1(y_m,k)$ for $m=1,2$. Then, for $m=1,2$, there is $z_m\in Y_{k-1}\cap B(y_m,D^k)$ with $x\in I_3(z_m,k-1)$. Using \eqref{disji} inductively for $j=3$, we conclude $z_1=z_2$. This implies that the balls $B(y_1, D^k)$ and $B(y_2, D^k)$ intersect. By construction of $Y_k$, this implies $y_1=y_2$. This proves \eqref{disji} for $j=1$.

We next consider \eqref{disji} for $j=3$. Assume $x$ is in $I_3(y_m,k)$ for $m=1,2$ and $y_m\in Y_k$. If $x$ is in $X_k$, then by definition \eqref{definei3}, $x\in I_1(y_m,k)$ for $m=1,2$. As we have already shown \eqref{disji} for $j=1$, we conclude $y_1=y_2$. This completes the proof in case $x\in X_k$, and we may assume $x$ is not in $X_k$. By definition \eqref{definei3}, $x$ is not in $I_3(z,k)$ for any $z$ with $z<y_1$ or $z<y_2$. Hence, neither $y_1<y_2$ nor $y_2<y_1$, and by totality of the order of $Y_k$, we have $y_1=y_2$. This completes the proof of \eqref{disji} for $j=3$.

We show \eqref{unioni} for $j=2$. In case $k=-S$, this follows from \Cref{cover-big-ball}. Assume $k>-S$. Let $x$ be a point of $B(o, 4D^S-2D^k)$. By induction, there is $y'\in Y_{k-1}$ such that $x\in I_3(y',k-1)$. Using the inductive statement \eqref{squeezedyadic}, we obtain $x\in B(y',4D^{k-1})$. As $D>4$, by applying the triangle inequality with the points, $o$, $x$, and $y'$, we obtain that $y'\in B(o, 4D^S-D^k)$. By \Cref{cover-big-ball}, $y'$ is in $B(y,2D^k)$ for some $y\in Y_k$. It follows that $x\in I_2(y,k)$. This proves \eqref{unioni} for $j=2$.

We show \eqref{unioni} for $j=3$. Let $x\in B(o, 4D^S-2D^k)$. In case $x\in X_k$, then by definition of $X_k$ we have $x\in I_1(y,k)$ for some $y\in Y_k$ and thus $x\in I_3(y,k)$. We may thus assume $x\not\in X_k$. As we have already seen \eqref{unioni} for $j=2$, there is $y\in Y_k$ such that $x\in I_2(y,k)$. We may assume this $y$ is minimal with respect to the order in $Y_k$. Then $x\in I_3(y,k)$. This proves \eqref{unioni} for $j=3$.

Next, we show the first inclusion in \eqref{squeezedyadic}. Let $x\in B(y,\frac 12 D^k)$. As $I_1(y,k)\subset I_3(y,k)$, it suffices to show $x\in I_1(y,k)$. If $k=-S$, this follows immediately from the assumption on $x$ and the definition of $I_1$. Assume $k>-S$. By the inductive statement \eqref{unioni} and $D>4$, there is a $y'\in Y_{k-1}$ such that $x\in I_3(y',k-1)$. By the inductive statement \eqref{squeezedyadic}, we conclude $x\in B(y',4D^{k-1})$. By the triangle inequality with points $x$, $y$, $y'$, and $D\geq8$, we have $y'\in B(y,D^k)$. It follows by definition \eqref{defineij} that $I_3(y',k-1)\subset I_1(y,k)$, and thus $x\in I_3(y,k)$. This proves the first inclusion in \eqref{squeezedyadic}.

We show the second inclusion in \eqref{squeezedyadic}. Let $x\in I_3(y,k)$. As $I_1(y,k)\subset I_2(y,k)$ directly from the definition \eqref{defineij}, it follows by definition \eqref{definei3} that $x\in I_2(y,k)$. By definition \eqref{defineij}, there is $y'\in Y_{k-1}\cap B(y,2D^k)$ with $x\in I_3(y',k-1)$. By induction, $x\in B(y', 4D^{k-1})$. By the triangle inequality applied to the points $x,y',y$ and $D>4$, we conclude $x\in B(y,4D^k)$. This shows the second inclusion in \eqref{squeezedyadic} and completes the proof of the lemma.
\end{proof}

\begin{lemma}[cover by cubes]
    \label{cover-by-cubes}
    \leanok
    \lean{cover_by_cubes}
    Let $-S\le l\le k\le S$ and
    $y\in Y_k$.
    We have
    \begin{equation}\label{3coverdyadic}
    I_3(y,k)\subset \bigcup_{y'\in Y_l} I_3(y',l)\, .
    \end{equation}
\end{lemma}
\begin{proof}
    \leanok
    Let $-S\le l\le k\le S$ and $y\in Y_k$.
    If $l=k$, the inclusion \eqref{3coverdyadic}
    is true from the definition of set union.
    We may then assume inductively that $k>l$ and the statement of the lemma is true if $k$ is replaced by $k-1$.
    Let $x\in I_3(y,k)$.
    By definition \eqref{definei3}, $x\in I_j(y,k)$
    for some $j\in \{1,2\}$. By \eqref{defineij},
    $x\in I_3(w,k-1)$ for some $w\in Y_{k-1}$.
    We conclude \eqref{3coverdyadic} by induction.
\end{proof}

\begin{lemma}[dyadic property]
    \label{dyadic-property}
    \uses{basic-grid-structure, cover-by-cubes}
    \leanok
    \lean{dyadic_property}
    Let $-S\le l\le k\le S$ and $y\in Y_k$ and $y'\in Y_l$ with
    $I_3(y',l)\cap I_3(y,k)\neq \emptyset$. Then
    \begin{equation}
        \label{3dyadicproperty}
    I_3(y',l)\subset I_3(y,k).
    \end{equation}
\end{lemma}

\begin{proof}
\leanok
Let $l,k,y,y'$ be as in the lemma. Pick $x\in I_3(y',l)\cap I_3(y,k)$. Assume first $l=k$. By \eqref{disji} of \Cref{basic-grid-structure}, we conclude $y'=y$, and thus \eqref{3dyadicproperty}. Now assume $l<k$. By induction, we may assume that the statement of the lemma is proven for $k-1$ in place of $k$.

By \Cref{cover-by-cubes}, there is a $y''\in Y_{k-1}$ such that $x\in I_3(y'',k-1)$. By induction, we have $I_3(y',l)\subset I_3(y'',k-1)$. It remains to prove
\begin{equation}\label{wyclaim}
I_3(y'',k-1)\subset I_3(y,k).
\end{equation}
We make a case distinction and assume first $x\in X_k$. By Definition \eqref{definei3}, we have $x\in I_1(y,k)$. By Definition \eqref{defineij}, there is a $v\in Y_{k-1}\cap B(y,D^k)$ with $x\in I_3(v,k-1)$. By \eqref{disji} of \Cref{basic-grid-structure}, we have $v=y''$. By Definition \eqref{defineij}, we then have $I_3(y'',k-1)\subset I_1(y,k)$. Then \eqref{wyclaim} follows by Definition \eqref{definei3} in the given case.

Assume now the case $x\notin X_k$. By \eqref{definei3}, we have $x\in I_2(y,k)$. Moreover, for any $u<y$ in $Y_k$, we have $x\not\in I_3(u,k)$. Let $u<y$. By transitivity of the order in $Y_k$, we conclude $x\not \in I_2(u,k)$. By \eqref{defineij} and the disjointedness property of \Cref{basic-grid-structure}, we have $I_3(y'',k-1)\cap I_2(u,k)= \emptyset$. Similarly, $I_3(y'',k-1)\cap I_1(u,k)= \emptyset$. Hence $I_3(y'',k-1)\cap I_3(u,k)=\emptyset$. As $u<y$ was arbitrary, we conclude with \eqref{definei3} the claim in the given case. This completes the proof of \eqref{wyclaim}, and thus also \eqref{3dyadicproperty}.
\end{proof}

For $-S\le k'\le k\le S$ and $y'\in Y_{k'}$, $y\in Y_k$
write $(y',k'|y,k)$ if $I_3(y',k')\subset I_3(y,k)$ and
\begin{equation}\label{bdcond}
    \inf_{x\in X\setminus I_3(y,k)}\rho(y',x)<6D^{k'}\, .
\end{equation}

\begin{lemma}[transitive boundary]
    \label{transitive-boundary}
    \uses{dyadic-property}
    \leanok
    \lean{transitive_boundary}
    Assume $-S\le k''< k'< k\le S$ and
    $y''\in Y_{k''}$, $y'\in Y_{k'}$, $y\in Y_k$.
    Assume there is $x\in X$ such that
    \begin{equation}
    x\in I_3(y'',k'')\cap I_3(y',k')\cap I_3(y,k)\, .
    \end{equation}
    If $(y'',k''|y,k)$, the also
    $(y'',k''|y',k')$ and $(y',k'|y,k)$
\end{lemma}

\begin{proof}
    \leanok
    As $x\in I_3(y'',k'')\cap I_3(y',k')$ and $k''< k'$, we have by \Cref{dyadic-property} that
    $I_3(y'',k'')\subset I_3(y',k')$. Similarly,
    $I_3(y',k')\subset I_3(y,k)$.
    Pick $x'\in X\setminus I_3(y,k)$ such that
    \begin{equation}\label{yppxp}
        \rho(y'',x')< 6D^{k''}\, ,
    \end{equation}
    which exists as $(y'',k''|y,k)$. As
    $x'\in X\setminus I_3(y',k')$ as well, we conclude
    $(y'',k''| y',k')$.
    By the triangle inequality, we have
    \begin{equation}
    \rho(y',x')\le \rho(y',x)+\rho(x,y'')+\rho(y'',x')
    \end{equation}
    Using the choice of $x$ and \eqref{squeezedyadic}
    as well as \eqref{yppxp}, we estimate this by
    \begin{equation}
    < 4D^{k'}+4D^{k''}+6D^{k''}\le 6D^{k'}\, ,
    \end{equation}
    where we have used $D>5$ and $k''<k'$.
    We conclude $(y',k'|y,k)$.
\end{proof}

\begin{lemma}[small boundary]
    \label{small-boundary}
    \uses{transitive-boundary}
    \leanok
    \lean{small_boundary}
    Let $K = 2^{4a+1}$. For each $-S+K\le k\le S$ and $y\in Y_k$ we have
      \begin{equation}
            \label{new-small-boundary}
            \sum_{z\in Y_{k-K}: (z,k-K|y,k)}\mu(I_3(z,k-K)) \le \frac 12 \mu(I_3(y,k))\,.
        \end{equation}
\end{lemma}

\begin{proof}
\leanok
Let $K$ be as in the lemma. Let $-S+K\le k\le S$ and $y\in Y_k$.

Pick $k'$ so that $k-K\le k'\le k$.
For each $y''\in Y_{k-K}$ with $(y'',k-K| y,k)$,
by \Cref{cover-by-cubes} and \Cref{dyadic-property}, there is a unique $y'\in Y_{k'}$ such that
\begin{equation}\label{4.31}
    I_3(y'',k-K)\subset I_3(y',k')\subset I_3(y,k)\, .
\end{equation}
Using \Cref{transitive-boundary}, $(y',k'|y,k)$.

We conclude using the disjointedness property of
\Cref{basic-grid-structure} that
\begin{equation}\label{scalecompare}
    \sum_{y'': (y'',k-K|y,k)}\mu(I_3(y'',k-K))
   \le
\sum_{y': (y',k'|y,k)}
\mu(I_3(y',k')) \, .
   \end{equation}
Adding over $k-K<k'\le k$, and using
\[\mu(I_3(y',k'))\le 2^{4a} \mu(B(y', \frac 14 D^{k'}))\]
from the doubling property \eqref{doublingx} and
\eqref{squeezedyadic} gives
\begin{equation}\label{sumcompare}
    K\sum_{y'': (y'',k-K|y,k)}
    \mu(I_3(y'',k-K))
\end{equation}
\begin{equation}\label{sumcompare1}
    \le 2^{4a} \sum_{k-K<k'\le k}
   \left[ \sum_{y': (y',k'|y,k)}
\mu(B(y', \frac 14 D^{k'}))\right]
\end{equation}
Each ball $B(y', \frac 14 D^{k'})$ occurring in
\eqref{sumcompare1} is contained in $I_3(y',k')$
by \eqref{squeezedyadic} and in turn contained in
$I_3(y,k)$ by \eqref{4.31}. Assume for the moment all these balls are pairwise disjoint. Then
by additivity of the measure,
\begin{equation}
    K\sum_{y'': (y'',k-K|y,k)}
    \mu(I_3(y'',k-K))
    \le 2^{4a}
\mu(I_3(y,k))\,
\end{equation}
which by $K=2^{4a+1}$ implies \eqref{new-small-boundary}.

It thus remains to prove that the balls
occurring in
\eqref{sumcompare1} are pairwise disjoint.
Let $(u,l)$ and $(u',l')$ be two parameter pairs occurring
in the sum of \eqref{sumcompare1} and let
$ B(u, \frac 14 D^l)$ and $B({u'}, \frac 14 D^{l'})$
be the corresponding balls. If
$l=l'$, then the balls are equal or disjoint by
\eqref{squeezedyadic} and \eqref{disji} of \Cref{basic-grid-structure}. Assume then without loss of generality that $l'<l$. Towards a contradiction, assume that
\begin{equation}\label{bulbul}
    B(u, \frac 14 D^l)\cap B({u'}, \frac 14 D^{l'})\neq \emptyset
\end{equation}
As $(u',l'|y,k)$, there is a point $x$ in
$X\setminus I_3(y,k)$ with $\rho(x,u')<6D^{l'}$.
Using $D>25$, we conclude from the triangle inequality and \eqref{bulbul} that
$x\in B(u,\frac 12D^l)$. However, $B(u,\frac 12 D^l)\subset I_3(u,l)$,
and $I_3(u,l)\subset I_3(y,k)$, a contradiction to
$x\not\in I_3(y,k)$.
This proves the lemma.
\end{proof}

\begin{lemma}[smaller boundary]
    \label{smaller-boundary}
    \uses{small-boundary}
    \leanok
    \lean{smaller_boundary}
    Let $K = 2^{4a+1}$
    and let $n\ge 0$ be an integer. Then
    for each $-S+nK\le k\le S$ we have
    \begin{equation}
            \label{very-new-small}
            \sum_{y'\in Y_{k-nK}: (y',k-nK|y,k)}\mu(I_3(y',k-nK)) \le 2^{-n} \mu(I_3(y,k))\,.
        \end{equation}
\end{lemma}
\begin{proof}
    \leanok
    We prove this by induction on $n$. If $n=0$,
    both sides of \eqref{very-new-small} are equal to
    $\mu(I_3(y,k))$ by \eqref{disji}. If $n=1$, this follows from \Cref{small-boundary}.

    Assume $n>1$ and \eqref{very-new-small} has been proven for $n-1$.
    We write \eqref{very-new-small}
    \begin{equation}
        \sum_{y''\in Y_{k-nK}: (y'',k-nK|y,k)}\mu(I_3(y'',k-nK))
    \end{equation}
    \begin{equation}
        = \sum_{y'\in Y_{k-K}:(y',k-K|y,k)} \left[ \sum_{y''\in Y_{k-nK}: (y'',k-nK|y',k-K)}\mu(I_3(y'',k-nK)) \right]
    \end{equation}
    Applying the induction hypothesis, this is bounded by
    \begin{equation}
        = \sum_{y'\in Y_{k-K}:(y',k-K|y,k)} 2^{1-n}\mu(I_3(y',k-K))
    \end{equation}
    Applying \eqref{new-small-boundary} gives
    \eqref{very-new-small}, and proves the lemma.
\end{proof}

\begin{lemma}[boundary measure]
    \label{boundary-measure}
    \uses{smaller-boundary}
    \leanok
    \lean{boundary_measure}
    For each $-S\le k\le S$ and $y\in Y_k$ and $0<t<1$
    with $tD^k\ge D^{-S}$ we have
    \begin{equation}
        \label{old-small-boundary}
        \mu(\{x \in I_3(y,k) \ : \ \rho(x, X \setminus I_3(y,k)) \leq t D^{k}\}) \le 2 t^\kappa \mu(I_3(y,k))\,.
    \end{equation}
\end{lemma}
\begin{proof}
\leanok
Let $x\in I_3(y,k)$ with $\rho(x, X \setminus I_3(y,k)) \leq t D^{k}$. Let $K = 2^{4a+1}$ as in \Cref{smaller-boundary}.
Let $n$ be the largest integer such that
$D^{nK} \le \frac{1}{t}$, so that $tD^k \le D^{k-nK}$ and
\begin{equation}
\label{eq-n-size}
    D^{nK} > \frac{1}{tD^K}\,.
\end{equation}
Let $k' = k - nK$, by the assumption $tD^k \ge D^{-S}$, we have $k' \ge -S$. By \eqref{3coverdyadic}, there exists $y' \in Y_{k'}$ with $x \in I_3(y',k')$. By the squeezing property \eqref{squeezedyadic} and the assumption on $x$, we have
$$
    \rho(y', X \setminus I_3(y,k)) \le \rho(x,y') + \rho(x, X \setminus I_3(y,k)) \le 4 D^{k'} + t D^{k}\,.
$$
By the assumption on $n$ and the definition of $k'$, this is
$$
    \le 4D^{k'} + D^{k - nK} < 6D^{k'}\,.
$$
Together with \eqref{3dyadicproperty} thus $(y',k'|y,k)$. We have shown that
$$
    \{x \in I_3(y,k) \ : \ \rho(x, X \setminus I_3(y,k)) \leq t D^{k}\}
$$
$$
    \subset \bigcup_{y'\in Y_{k-nK}: (y',k-nK|y,k)}I_3(y',k-nK)\,.
$$
Using monotonicity and additivity of the measure and \Cref{smaller-boundary}, we obtain
$$
    \mu(\{x \in I_3(y,k) \ : \ \rho(x, X \setminus I_3(y,k)) \leq t D^{k}\}) \le 2^{-n} \mu(I_3(y,k))\,.
$$
By \eqref{eq-n-size} and the definition \eqref{defineD} of $D$, this is bounded by
$$
    2 t^{1/(100a^2 K)}  \mu(I_3(y,k))\,,
$$
which completes the proof by the definition \eqref{definekappa} of $\kappa$.
\end{proof}

Let $\tilde{\mathcal{D}}$ be the set of all $I_3(y,k)$ with $k\in [-S,S]$ and
$y\in Y_k$. Define
\begin{equation}
s(I_3(y,k)):=k
\end{equation}
\begin{equation}
c(I_3(y,k)):=y\,.
\end{equation}
We define $\mathcal{D}$ to be the set of all $I \in \tilde{\mathcal{D}}$ such that $I \subset I_3(o, S)$.

\begin{proof}[Proof of \Cref{grid-existence}]
\proves{grid-existence}
\leanok
We first show that $(\tilde{\mathcal{D}},c,s)$ satisfies properties \eqref{coverdyadic}, \eqref{dyadicproperty}, \eqref{eq-vol-sp-cube} and \eqref{eq-small-boundary}. Property \eqref{eq-vol-sp-cube}
follows from \eqref{squeezedyadic}, while \eqref{eq-small-boundary} follows from \Cref{boundary-measure}.

Let $x\in B(o, D^S)$.
We show properties
\eqref{coverdyadic} and
\eqref{dyadicproperty}
for $(\tilde{\mathcal{D}},c,s)$ and $x$.

We first show \eqref{coverdyadic} for $(\tilde{\mathcal{D}},c,s)$ by contradiction. Then there is an $I$ violating the conclusion of
\eqref{coverdyadic}. Pick such $I=I_3(y,l)$ such that $l$ is minimal.
{By assumption, we have $-S\le k<l$; in particular $-S<l$}.
By definition, $I_3(y, l)$ is contained in $I_1(y, l)\cup I_2(y, l)$, which is contained in the union of $I_3(y',l-1)$ with $y'\in Y_{l-1}$.
By minimality of $l$, each such $I_3(y',l-1)$ is contained in the union of
all $I_3(z,k)$ with $z\in Y_k$. This proves \eqref{coverdyadic}.

We now show \eqref{dyadicproperty} for $(\tilde{\mathcal{D}},c,s)$. Assume to get a contradiction that
there are non-disjoint $I, J\in \tilde{\mathcal{D}}$ with $s(I)\le s(J)$
and $I \not \subset J$. We may assume the existence of such $I$ and $J$ with minimal
$s(J)-s(I)$. Let $k=s(I)$. Assume first $s(J)=k$. Let $I=I_3(y_1,k)$ and $J=I_3(y_2,k)$ with $y_1,y_2\in Y_k$.
If $y_1=y_2$, then $I=J$, a contradiction to $I\not \subset J$.
If $y_1\neq y_2$, then $I\cap J=\emptyset$ by \eqref{disji}, a contradiction to the non-disjointedness of $I,J$.
Assume now $s(J)>k$. Choose $y\in I\cap J$. By property \eqref{coverdyadic},
there is $K\in \tilde{\mathcal{D}}$ with $s(K)=s(J)-1$ and $y\in K$. By construction
of $J$, and pairwise disjointedness of all $I_3(w,s(J)-1)$ that we have already seen,
we have $K\subset J$. By minimality of $s(J)$, we have $I\subset K$.
This proves $I\subset J$ and thus \eqref{dyadicproperty}.

Now note that properties \eqref{dyadicproperty}, \eqref{eq-vol-sp-cube} and \eqref{eq-small-boundary} immediately carry over to $(\mathcal{D},c, s)$ by restriction. \eqref{subsetmaxcube} is true for $(\mathcal{D}, c, s)$ by definition, and \eqref{coverdyadic} follows from \eqref{coverdyadic} and \eqref{dyadicproperty} for $(\tilde{\mathcal{D}}, c, s)$.
\end{proof}

\subsection{Proof of Tile Structure Lemma}
\label{subsectiles}
Choose a grid structure $(\mathcal{D}, c, s)$ with \Cref{grid-existence}
Let $I \in \mathcal{D}$. Suppose that
\begin{equation}
    \label{eq-tile-Z}
    \mathcal{Z} \subset \tQ(X)
\end{equation}
is such that for any $\mfa, \mfb \in \mathcal{Z}$ with $\mfa\ne \mfb$ we have
\begin{equation}
    \label{eq-tile-disjoint-Z}
    B_{I^\circ}(\mfa, 0.3) \cap B_{I^\circ}(\mfb, 0.3) \cap \tQ(X) = \emptyset\,.
\end{equation}
Since $\tQ(X)$ is finite, there exists a set $\mathcal{Z}$ satisfying both \eqref{eq-tile-Z} and \eqref{eq-tile-disjoint-Z} of maximal cardinality among all such sets. We pick for each $I \in \mathcal{D}$ such a set $\mathcal{Z}(I)$.

\begin{lemma}[frequency ball cover]
    \label{frequency-ball-cover}
    \leanok
    \lean{frequency_ball_cover}
    For each $I \in \mathcal{D}$, we have
    \begin{equation}
        \label{eq-tile-cover}
        \tQ(X) \subset \bigcup_{z \in \mathcal{Z}(I)} B_{I^\circ}(z, 0.7)\,.
    \end{equation}
\end{lemma}

\begin{proof}
    \leanok
    Let $\mfb \in \bigcup_{\mfa \in \tQ(X)} B_{I^\circ}(\mfa, 1)$. By maximality of $\mathcal{Z}(I)$, there must be a point $z \in \mathcal{Z}(I)$ such that $B_{I^\circ}(z, 0.3) \cap B_{I^\circ}(\mfb, 0.3) \ne \emptyset$. Else, $\mathcal{Z}(I) \cup \{\mfb\}$ would be a set of larger cardinality than $\mathcal{Z}(I)$ satisfying \eqref{eq-tile-Z} and \eqref{eq-tile-disjoint-Z}. Fix such $z$, and fix a point $z_1 \in B_{I^\circ}(z, 0.3) \cap B_{I^\circ}(\mfb, 0.3)$. By the triangle inequality, we deduce that
    $$
        d_{I^\circ}(z,\mfb) \le d_{I^\circ}(z,z_1) + d_{I^\circ}(\mfb, z_1) < 0.3 + 0.3 = 0.6\,,
    $$
    and hence $\mfb \in B_{I^\circ}(z, 0.7)$.
\end{proof}

We define
$$
    \fP = \{(I, z) \ : \ I \in \mathcal{D}, z \in \mathcal{Z}(I)\}\,,
$$
$$\scI((I, z)) = I\qquad \text{and} \qquad \fcc((I, z)) = z.$$ We further set $$\ps(\fp) = s(\scI(\fp)),\qquad \qquad \pc(\fp) = c(\scI(\fp)).$$ Then \eqref{tilecenter}, \eqref{tilescale} hold by definition.

It remains to construct the map $\Omega$, and verify properties \eqref{eq-dis-freq-cover}, \eqref{eq-freq-dyadic} and
\eqref{eq-freq-comp-ball}. We first construct an auxiliary map $\Omega_1$. For each $I \in \mathcal{D}$, we pick an enumeration of the finite set $\mathcal{Z}(I)$
$$
    \mathcal{Z}(I) = \{z_1, \dotsc, z_M\}\,.
$$
We define {$\Omega_1:\fP \mapsto \mathcal{P}(\Mf) $ as below}. Set
$$
    \Omega_1((I, z_1)) = B_{I^\circ}(z_1, 0.7) \setminus \bigcup_{z \in \mathcal{Z}(I)\setminus \{z_1\}} B_{I^\circ}(z, 0.3)
$$
and then define iteratively
\begin{equation}
    \label{eq-def-omega1}
    \Omega_1((I, z_k)) = B_{I^\circ}(z_k, 0.7) \setminus \bigcup_{z \in \mathcal{Z}(I) \setminus \{z_k\}} B_{I^\circ}(z, 0.3) \setminus \bigcup_{i=1}^{k-1} \Omega_1((I, z_i))\,.
\end{equation}
\begin{lemma}[disjoint frequency cubes]
    \label{disjoint-frequency-cubes}
    \leanok
    \lean{Construction.disjoint_frequency_cubes}
    For each $I \in \mathcal{D}$, and $\fp_1, \fp_2\in \fP(I)$,
    if $$\Omega_1(\fp_1)\cap \Omega_1(\fp_2)\neq \emptyset,$$ then $\fp_1=\fp_2$.
\end{lemma}

\begin{proof}
    \leanok
    By the definition of the map $\scI$, we have
    $$
        \fP(I) = \{(I, z) \, : \, z \in \mathcal{Z}(I)\}\,.
    $$
    By \eqref{eq-def-omega1}, the set $\Omega_1((I, z_k))$ is disjoint from each $\Omega_1((I, z_i))$ with $i < k$. Thus the sets $\Omega_1(\fp)$, $\fp \in \fP(I)$ are pairwise disjoint.
\end{proof}

\begin{lemma}[frequency cube cover]
    \label{frequency-cube-cover}
    \leanok
    \lean{Construction.iUnion_ball_subset_iUnion_Ω₁, Construction.ball_subset_Ω₁, Construction.Ω₁_subset_ball}
    For each $I \in \mathcal{D}$, it holds that
    \begin{equation}
    \label{eq-omega1-cover}
            \bigcup_{z \in \mathcal{Z}(I)} B_{I^\circ}(z, 0.7)\subset \bigcup_{\fp \in \fP(I)} \Omega_1(\fp)\,.
    \end{equation}
    For every $\fp \in \fP$, it holds that
    \begin{equation}
        \label{eq-omega1-incl}
        B_{\fp}(\fcc(\fp), 0.3) \subset \Omega_1(\fp) \subset B_{\fp}(\fcc(\fp), 0.7)\,.
    \end{equation}
\end{lemma}

\begin{proof}
    \leanok
    For \eqref{eq-omega1-incl} let $\fp = (I, z)$.
    The second inclusion in \eqref{eq-omega1-incl} then follows from \eqref{eq-def-omega1} and the equality $B_{\fp}(\fcc(\fp), 0.7) = B_{I^\circ}(z, 0.7)$, which is true by definition.
    For the first inclusion in \eqref{eq-omega1-incl} let $\mfa \in B_{\fp}(\fcc(\fp),0.3)$. Let $k$ be such that $z = z_k$ in the enumeration we chose above. It follows immediately from \eqref{eq-def-omega1} and \eqref{eq-tile-disjoint-Z} that
    $\mfa \notin \Omega_1((I, z_i))$ for all $i < k$. Thus, again from \eqref{eq-def-omega1}, we have
    $\mfa \in \Omega_1((I,z_k))$.

    To show \eqref{eq-omega1-cover} let $\mfa \in \bigcup_{z \in \mathcal{Z}(I)} B_{I^\circ}(z,0.7)$.
    If there exists $z \in \mathcal{Z}(I)$ with $\mfa \in B_{I^\circ}(z,0.3)$, then
    $$
        z \in \Omega_1((I, z)) \subset \bigcup_{\fp \in \fP(I)} \Omega_1(\fp)
    $$
    by the first inclusion in \eqref{eq-omega1-incl}.

    Now suppose that there exists no $z \in \mathcal{Z}(I)$ with $\mfa \in B_{I^\circ}(z, 0.3)$. Let $k$ be minimal such that $\mfa \in B_{I^\circ}(z_k, 0.7)$. Since $\Omega_1((I, z_i)) \subset B_{I^\circ}(z_i, 0.7)$ for each $i$ by \eqref{eq-def-omega1}, we have that $\mfa \notin \Omega_1((I, z_i))$ for all $i < k$. Hence $\mfa \in \Omega_1((I, z_k))$, again by \eqref{eq-def-omega1}.
\end{proof}

Now we are ready to define the function $\Omega$. We define for all $\fp \in \fP(I_0)$
\begin{equation}
    \label{eq-max-omega}
    \fc(\fp) = \Omega_1(\fp)\,.
\end{equation}
For all other cubes $I \in \mathcal{D}$, $I \ne I_0$, there exists, by \eqref{dyadicproperty} and \eqref{subsetmaxcube}, $J \in \mathcal{D}$ with $I \subset J$ and $I \ne J$. On such $I$ we define $\Omega$ by recursion. We can pick an inclusion minimal $J \in \mathcal{D}$ among the finitely many cubes such that $I \subset J$ and $I \ne J$. This $J$ is unique: Suppose that $J'$ is another inclusion minimal cube with $I \subset J'$ and $I \ne J'$. Without loss of generality, we have that $s(J) \le s(J')$. By \eqref{dyadicproperty}, it follows that $J \subset J'$. Since $J'$ is minimal with respect to inclusion, it follows that $J = J'$. Then we define
\begin{equation}
    \label{eq-it-omega}
    \fc(\fp) = \bigcup_{z \in \mathcal{Z}(J) \cap \Omega_1(\fp)} \Omega((J, z)) \cup B_{\fp}(\fcc(\fp),0.2)\,.
\end{equation}

We now verify that $(\fP,\scI,\fc,\fcc,\pc,\ps)$ forms a tile structure.
\begin{proof}[Proof of \Cref{tile-structure}]
    \leanok
    \proves{tile-structure}
    First, we prove \eqref{eq-freq-comp-ball}. If $I =I_0$, then \eqref{eq-freq-comp-ball} holds for all $\fp \in \fP(I)$ by \eqref{eq-max-omega} and \eqref{eq-omega1-incl}. Now suppose that $I$ is not maximal in $\mathcal{D}$ with respect to set inclusion. Then we may assume by induction that for all $J \in \mathcal{D}$ with $I \subset J$ and all $\fp' \in \fP(J)$, \eqref{eq-freq-comp-ball} holds. Let $J$ be the unique minimal cube in $\mathcal{D}$ with $I \subsetneq J$.

    Suppose that $\mfa \in \Omega(\fp)$. If $\mfa\in B_{\fp}(\mathcal{Q}(\fp), 0.2)$, then since
    \begin{equation*}
        B_{\fp}(\mathcal{Q}(\fp), 0.2)\subset B_{\fp}(\mathcal{Q}(\fp), 1)\, ,
    \end{equation*}
    we conclude that $\mfa\in B_{\fp}(\mathcal{Q}(\fp), 0.7)$. If not, by \eqref{eq-it-omega}, there exists $z \in \mathcal{Z}(J) \cap \Omega_1(\fp)$ with $\mfa \in \Omega(J,z)$. Using the triangle inequality and \eqref{eq-omega1-incl}, we obtain
    $$
        d_{I^\circ}(\fcc(\fp),\mfa) \le d_{I^\circ}(\fcc(\fp), z) + d_{I^\circ}(z, \mfa) \le 0.7 + d_{I^\circ}(z, \mfa)\,.
    $$
    By \Cref{monotone-cube-metrics} and the induction hypothesis, this is estimated by
    $$
        \le 0.7 + 2^{-95a} d_{J^\circ}(z,\mfa) \le 0.7 + 2^{-95a}\cdot 1 < 1\,.
    $$
    This shows the second inclusion in \eqref{eq-freq-comp-ball}. The first inclusion is immediate from \eqref{eq-it-omega}.

    Next, we show \eqref{eq-dis-freq-cover}. Let $I \in \mathcal{D}$.

    If $I = I_0$, then disjointedness of the sets $\fc(\fp)$ for $\fp \in \fP(I)$ follows from the definition \eqref{eq-max-omega} and \Cref{disjoint-frequency-cubes}. To obtain the inclusion in \eqref{eq-dis-freq-cover} one combines the inclusions \eqref{eq-tile-cover} and \eqref{eq-omega1-cover}
    of \Cref{frequency-cube-cover} with \eqref{eq-max-omega}.

    Now we turn to the case where there exists $J \in \mathcal{D}$ with $I \subset J$ and $I\ne J$. In this case we use induction: It suffices to show \eqref{eq-dis-freq-cover} under the assumption that it holds for all cubes $J \in \mathcal{D}$ with $I \subset J$. As shown before definition \eqref{eq-it-omega}, we may choose the unique inclusion minimal such $J$. To show disjointedness of the sets $\fc(\fp), \fp \in \fP(I)$ we pick two tiles $\fp, \fp' \in \fP(I)$ and $\mfa \in \fc(\fp) \cap \fc(\fp')$.
    Then we are by \eqref{eq-it-omega} in one of the following four cases.

    1. There exist $z \in \mathcal{Z}(J) \cap \Omega_1(\fp)$ such that $\mfa \in \Omega(J, z)$, and there exists $z' \in \mathcal{Z}(J) \cap \Omega_1(\fp')$ such that $\mfa \in \Omega(J, z')$. By the induction hypothesis, that \eqref{eq-dis-freq-cover} holds for $J$, we must have $z = z'$. By \Cref{disjoint-frequency-cubes}, we must then have $\fp = \fp'$.

    2. There exists $z \in \mathcal{Z}(J) \cap \Omega_1(\fp)$ such that $\mfa \in \Omega(J,z)$, and $\mfa \in B_{\fp'}(\fcc(\fp'), 0.2)$. Using the triangle inequality, \Cref{monotone-cube-metrics} and \eqref{eq-freq-comp-ball}, we obtain
    $$
        d_{\fp'}(\fcc(\fp'),z) \le d_{\fp'}(\fcc(\fp'), \mfa) + d_{\fp'}(z, \mfa) \le 0.2 + 2^{-95a} \cdot 1 < 0.3\,.
    $$
    Thus $z \in \Omega_1(\fp')$ by \eqref{eq-omega1-incl}. By \Cref{disjoint-frequency-cubes}, it follows that $\fp = \fp'$.

    3. There exists $z' \in \mathcal{Z}(J) \cap \Omega_1(\fp')$ such that $\mfa \in \Omega(J,z')$, and $\mfa \in B_{\fp}(\fcc(\fp), 0.2)$. This case is the same as case 2., after swapping $\fp$ and $\fp'$.

    4. We have $\mfa \in B_{\fp}(\fcc(\fp), 0.2) \cap B_{\fp'}(\fcc(\fp'), 0.2)$. In this case it follows that $\fp = \fp'$ since the sets $B_{\fp}(\fcc(\fp), 0.2)$ are pairwise disjoint by the inclusion \eqref{eq-omega1-incl} and \Cref{disjoint-frequency-cubes}.

    To show the inclusion in \eqref{eq-dis-freq-cover}, let $\mfa \in \tQ(X)$. By the induction hypothesis, there exists $\fp \in \fP(J)$ such that $\mfa \in \Omega(\fp)$. By definition of the set $\fP$, we have $\fp = (J, z)$ for some $z \in \mathcal{Z}(J)$.
    Thus, by \eqref{eq-tile-cover}, there exists $z' \in \mathcal{Z}(I)$ with $z \in B_{I^\circ}(z', 0.7)$. Then by Lemma \eqref{frequency-cube-cover} there exists $\fp' \in \fP(I)$ with $z \in \mathcal{Z}(J) \cap \Omega_1(\fp')$. Consequently, by \eqref{eq-it-omega}, $\mfa \in \fc(\fp')$. This completes the proof of \eqref{eq-dis-freq-cover}.

    Finally, we show \eqref{eq-freq-dyadic}. Let $\fp, \fq \in \fP$ with $\scI(\fp) \subset \scI(\fq)$ and $\fc(\fp) \cap \fc(\fq) \ne \emptyset$. If we have $\ps(\fp) \ge \ps(\fq)$, then it follows from \eqref{dyadicproperty} that $I = J$, thus $\fp, \fq \in \fP(I)$. By \eqref{eq-dis-freq-cover} we have then either $\fc(\fp) \cap \fc(\fq) = \emptyset$ or $\fc(\fp) = \fc(\fq)$. By the assumption in \eqref{eq-freq-dyadic} we have $\fc(\fp) \cap \fc(\fq) \ne \emptyset$, so we must have $\fc(\fp) = \fc(\fq)$ and in particular $\fc(\fq) \subset \fc(\fp)$.

    So it remains to show \eqref{eq-freq-dyadic} under the additional assumption that $\ps(\fq) > \ps(\fp)$. In this case, we argue by induction on $\ps(\fq)-\ps(\fp)$. By \eqref{coverdyadic}, there exists a cube $J \in \mathcal{D}$ with $s(J) = \ps(\fq) - 1$ and $J \cap\scI(\fp) \ne \emptyset$. We pick one such $J$. By \eqref{dyadicproperty}, we have $\scI(\fp) \subset J \subset \scI(\fq)$.

    Thus, by \eqref{eq-tile-cover}, there exists $z' \in \mathcal{Z}(J)$ with $\fcc(\fq) \in B_{J^\circ}(z', 0.7)$. Then by \Cref{frequency-cube-cover} there exists $\fq' \in \fP(J)$ with $\fcc(\fq) \in\Omega_1(\fq')$.
    By \eqref{eq-it-omega}, it follows that $\Omega(\fq) \subset \Omega(\fq')$. Note that then $\scI(\fp) \subset \scI(\fq')$ and $\fc(\fp) \cap \fc(\fq') \ne \emptyset$ and $\ps(\fq') - \ps(\fp) = \ps(\fq) - \ps(\fp) - 1$. Thus, we have by the induction hypothesis that $\Omega(\fq') \subset \Omega(\fp)$. This completes the proof.
\end{proof}

\section{Proof of discrete Carleson}
\label{proptopropprop}

Let a grid structure $(\mathcal{D}, c, s)$ and a tile structure $(\fP, \scI, \fc, \fcc)$ for this grid structure be given. In \Cref{subsectilesorg}, we decompose the set $\fP$ of tiles into subsets. Each subset will be controlled by one of three methods. The guiding principle of the decomposition is to be able to apply the forest estimate of \Cref{forest-operator} to the final subsets defined in \eqref{defc5}. This application is done in \Cref{subsecforest}. The miscellaneous subsets along the construction of the forests will either be thrown into exceptional sets, which are defined and controlled in \Cref{subsetexcset}, or will be controlled by the antichain estimate of \Cref{antichain-operator}, which is done in \Cref{subsecantichain}. \Cref{subsec-lessim-aux} contains some auxiliary lemmas needed for the proofs in Subsections \ref{subsecforest}-\ref{subsecantichain}.

\subsection{Organisation of the tiles}\label{subsectilesorg}

In the following definitions, $k, n$, and
$j$ will be nonnegative integers.
Define
$\mathcal{C}(G,k)$ to be the set of $I\in \mathcal{D}$
such that there exists a $J\in \mathcal{D}$ with $I\subset J$
and
\begin{equation}\label{muhj1}
   {\mu(G \cap J)} > 2^{-k-1}{\mu(J)}\, ,
\end{equation}
but there does not exist a $J\in \mathcal{D}$ with $I\subset J$ and
\begin{equation}\label{muhj2}
   {\mu(G \cap J)} > 2^{-k}{\mu(J)}\,.
\end{equation}
Let
\begin{equation}
    \label{eq-defPk}
    \fP(k)=\{\fp\in \fP \ : \ \scI(\fp)\in \mathcal{C}(G,k)\}
\end{equation}
Define $ {\mathfrak{M}}(k,n)$ to be the set of $\fp \in \fP(k)$ such that
 \begin{equation}\label{ebardense}
    \mu({E_1}(\fp)) > 2^{-n} \mu(\scI(\fp))
 \end{equation}
and there does not exist $\fp'\in \fP(k)$ with
$\fp'\neq \fp$ and $\fp\le \fp'$ such that
 \begin{equation}\label{mnkmax}
    \mu({E_1}(\fp')) > 2^{-n} \mu(\scI(\fp')).
 \end{equation}
Define for a collection $\fP'\subset \fP(k)$
\begin{equation}
    \label{eq-densdef}
   \dens_k' (\fP'):= \sup_{\fp'\in \fP'}\sup_{\lambda \geq 2} \lambda^{-a} \sup_{\fp \in \fP(k): \lambda \fp' \lesssim \lambda \fp}
    \frac{\mu({E}_2(\lambda, \fp))}{\mu(\scI(\fp))}\,.
\end{equation}
Sorting by density, we define
\begin{equation}
    \label{def-cnk}
    \fC(k,n):=\{\fp\in \fP(k) \ : \
    2^{4a}2^{-n}< \dens_k'(\{\fp\}) \le
    2^{4a}2^{-n+1}\}\,.
\end{equation}
Following Fefferman \cite{fefferman}, we
define for $\fp \in \fC(k,n)$
\begin{equation}\label{defbfp}
     \mathfrak{B}(\fp) := \{ \mathfrak{m} \in \mathfrak{M}(k,n) \ : \ 100 \fp \lesssim \mathfrak{m}\}
\end{equation}
and
\begin{equation}\label{defcnkj}
       \fC_1(k,n,j) := \{\fp \in \fC(k,n) \ : \ 2^{j} \leq |\mathfrak{B}(\fp)| < 2^{j+1}\}\,.
\end{equation}
and
\begin{equation}\label{defl0nk}
       \fL_0(k,n) := \{\fp \in \fC(k,n) \ : \ |\mathfrak{B}(\fp)| <1\}\,.
\end{equation}
Together with the following removal of minimal layers, the splitting into $\fC_1(k,n,j)$ will lead to a separation of trees.
Define recursively for $0\le l\le Z(n+1)$
\begin{equation}
    \label{eq-L1-def}
    \fL_1(k,n,j,l)
\end{equation}
to be the set of minimal elements with respect to $\le$ in
\begin{equation}
    \fC_1(k,n,j)\setminus \bigcup_{0\le l'<l}
\fL_1(k,n,j,l')\, .
\end{equation}
Define
\begin{equation}
    \label{eq-C2-def}
    \fC_2(k,n,j):= \fC_1(k,n,j)\setminus \bigcup_{0\le l'\le Z(n+1)}
\fL_1(k,n,j,l')\, .
\end{equation}

The remaining tile organization will be relative to
prospective tree tops, which we define now.
Define
\begin{equation}\label{defunkj}
     \fU_1(k,n,j)
\end{equation}
to be the set of all
$\fu \in \fC_1(k,n,j)$ such that
for all $\fp \in \fC_1(k,n,j)$
with $\scI(\fu)$ strictly contained in
$\scI(\fp)$ we have $B_{\fu}(\fcc(\fu), 100) \cap B_{\fp}(\fcc(\fp),100) = \emptyset$.

We first remove the pairs that are outside the immediate reach of any of the prospective tree tops.
Define
\begin{equation}
\label{eq-L2-def}
\fL_2(k,n,j)
\end{equation}
to be the set of all $\fp\in \fC_2(k,n,j)$ such that there
does not exist
$\fu\in \fU_1(k,n,j)$
with $\scI(\fp)\neq \scI(\fu)$ and $2\fp\lesssim \fu$.
Define
\begin{equation}
\label{eq-C3-def}
\fC_3(k,n,j):=\fC_2(k,n,j)
  \setminus \fL_2(k,n,j)\, .
\end{equation}

We next remove the maximal layers.
Define recursively for $0 \le l \le Z(n+1)$
\begin{equation}
    \label{eq-L3-def}
    \fL_3(k,n,j,l)
\end{equation}
to be the set of all maximal elements with respect to $\le$ in
\begin{equation}
    \fC_3(k,n,j) \setminus \bigcup_{0 \le l' < l} \fL_3(k,n,j,l')\,.
\end{equation}
Define
\begin{equation}
\label{eq-C4-def}
\fC_4(k,n,j):=\fC_3(k,n,j)
  \setminus \bigcup_{0 \le l \le Z(n+1)} \fL_3(k,n,j,l)\,.
\end{equation}

Finally, we remove the boundary pairs relative to the prospective tree tops. Define
\begin{equation}
    \label{eq-L-def}
    \mathcal{L}(\fu)
\end{equation}
to be the set of all $I \in \mathcal{D}$ with $I \subset \scI(\fu)$ and $s(I) = \ps(\fu) - Z(n+1) - 1$ and
\begin{equation}
    B(c(I), 8 D^{s(I)})\not \subset \scI(\fu)\, .
\end{equation}
Define
\begin{equation}
    \label{eq-L4-def}
    \fL_4(k,n,j)
\end{equation}
to be the set of all $\fp\in \fC_4(k,n,j)$ such that there exists
$\fu\in \fU_1(k,n,j)$
with $\scI(\fp) \subset \bigcup \mathcal{L}(\fu)$, and define
\begin{equation}\label{defc5}
\fC_5(k,n,j):=\fC_4(k,n,j)
  \setminus \fL_4(k,n,j)\, .
\end{equation}

We define three exceptional sets.
The first exceptional set $G_1$ takes into account the ratio of the measures of $F$ and $G$.
Define $\fP_{F,G}$ to be the set of all $\fp\in \fP$
with
\begin{equation}
    \dens_2(\{\fp\})> 2^{2a+5}\frac{\mu(F)}{\mu(G)}\,.
\end{equation}
Define
\begin{equation}\label{definegone}
    G_1:=\bigcup_{\fp\in \fP_{F,G} }\scI(\fp)\, .
\end{equation}
For an integer $\lambda\ge 0$, define $A(\lambda,k,n)$ to be the set
of all $x\in X$ such that
\begin{equation}
    \label{eq-Aoverlap-def}
    \sum_{\fp \in \mathfrak{M}(k,n)}\mathbf{1}_{\scI(\fp)}(x)>\lambda 2^{n+1}
\end{equation}
and define
\begin{equation}\label{definegone2}
    G_2:=
\bigcup_{k\ge 0}\bigcup_{k\le n}
A(2n+6,k,n)\, .
\end{equation}
Define
    \begin{equation}\label{defineg3}
        G_3 :=
        \bigcup_{k\ge 0}\, \bigcup_{n \geq k}\,
        \bigcup_{0\le j\le 2n+3}
        \bigcup_{\fp \in \fL_4 (k,n,j)}
        \scI(\fp)\, .
     \end{equation}
Define $G'=G_1\cup G_2 \cup G_3$. The following bound of the measure of $G'$ will be proven in \Cref{subsetexcset}.
\begin{lemma}[exceptional set]
    \label{exceptional-set}
    \leanok
    \lean{exceptional_set}
    \uses{first-exception,second-exception, third-exception}
    We have
    \begin{equation}
        \mu(G')\le 2^{-1}\mu(G)\, .
    \end{equation}
\end{lemma}

In \Cref{subsecforest}, we identify each set $\fC_5(k,n,j)$ outside $G'$ as forest and use Proposition
\ref{forest-operator} to prove the following lemma.

\begin{lemma}[forest union]
    \label{forest-union}
    \leanok
    \lean{forest_union}
    \uses{forest-operator,C-dens1,C6-forest, forest-geometry, forest-convex,forest-separation, forest-inner,forest-stacking}
    Let
    \begin{equation}
        \fP_1 =\bigcup_{k\ge 0}\bigcup_{n\ge k}
        \bigcup_{0\le j\le 2n+3}\fC_5(k,n,j)
    \end{equation}
    For all $f:X\to \C$ with $|f|\le \mathbf{1}_F$ we have
    \begin{equation}
        \label{disclesssim1}
        \int_{G \setminus G'} \left|\sum_{\fp \in \fP_1} T_{\fp} f \right|\, \mathrm{d}\mu \le \frac{2^{441a^3}}{(q-1)^4} \mu(G)^{1 - \frac{1}{q}} \mu(F)^{\frac{1}{q}}\,.
    \end{equation}
\end{lemma}

In \Cref{subsecantichain}, we decompose
the complement of the set of tiles in Lemma
\ref{forest-union} and apply the antichain estimate of
\Cref{antichain-operator} to prove the following lemma.

\begin{lemma}[forest complement]
    \label{forest-complement}
    \leanok
    \lean{forest_complement}
    \uses{antichain-operator,antichain-decomposition,L0-antichain,L2-antichain,L1-L3-antichain,C-dens1}
    Let
    \begin{equation}
        \fP_2 =\fP\setminus \fP_1\,.
    \end{equation}
    For all $f:X\to \C$ with $|f|\le \mathbf{1}_F$ we have
    \begin{equation}
        \label{disclesssim2}
        \int_{G \setminus G'} \left|\sum_{\fp \in \fP_2} T_{\fp} f\right| \, \mathrm{d}\mu \le \frac{2^{120a^3}}{(q-1)^5} \mu(G)^{1 - \frac{1}{q}} \mu(F)^{\frac{1}{q}}\,.
    \end{equation}
\end{lemma}

\begin{proof}[Proof of \Cref{discrete-Carleson}]
\leanok
\proves{discrete-Carleson}
\Cref{discrete-Carleson} follows by applying the
triangle inequality to \eqref{disclesssim}
according to the splitting in \Cref{forest-union}
and \Cref{forest-complement} and using both Lemmas as well
as the bound on the set $G'$ given by \Cref{exceptional-set}.
\end{proof}

\subsection{Proof of the Exceptional Sets Lemma}
\label{subsetexcset}

We prove separate bounds for $G_1$, $G_2$, and $G_3$
in Lemmas \ref{first-exception},
\ref{second-exception}, and \ref{third-exception}. Adding up these bounds proves \Cref{exceptional-set}.

The bound for $G_1$ follows from the Vitali covering lemma, \Cref{Hardy-Littlewood}.

\begin{lemma}[first exception]
    \label{first-exception}
    \leanok
    \lean{first_exception}
    \uses{Hardy-Littlewood}
    We have
    \begin{equation}
        \mu(G_1)\le 2^{-5}\mu(G)\, .
    \end{equation}
\end{lemma}
\begin{proof}
    \leanok
    Let
    $$
        K = 2^{2a+5}\frac{\mu(F)}{\mu(G)}\,.
    $$
    For each $\fp\in \fP_{F,G}$ pick a
    $r(\fp)>4D^{\ps(\fp)}$ with
    $$
    {\mu(F\cap B(\pc(\fp),r(\fp)))}\ge K{\mu(B(\pc(\fp),r(\fp)))}\, .
    $$
    This ball exists by definition of $\fP_{F,G}$
    and $\dens_2$. By applying \Cref{Hardy-Littlewood} to the collection of balls
    $$
        \mathcal{B} = \{B(\pc(\fp),r(\fp)) \ : \ \fp \in \fP_{F,G}\}
    $$
    and the function $u = \mathbf{1}_F$, we obtain
    $$
        \mu(\bigcup \mathcal{B}) \le 2^{2a} K^{-1} \mu(F)\,.
    $$
    We conclude with \eqref{eq-vol-sp-cube} and $r(\fp)>4D^{\ps(\fp)}$
    $$
        \mu(G_1)= \mu(\bigcup_{\fp\in \fP_{F,G}} \scI(\fp))
        \le \mu(\bigcup \mathcal{B})\le 2^{2a} K^{-1} \mu (F) = 2^{-5}\mu(G)\,.
    $$
\end{proof}

We turn to the bound of $G_2$, which relies on the Dyadic Covering \Cref{dense-cover} and the
John-Nirenberg \Cref{John-Nirenberg} below.

\begin{lemma}[dense cover]
\label{dense-cover}
\leanok
\lean{dense_cover}

For each $k\ge 0$, the union of all dyadic cubes
in $\mathcal{C}(G,k)$ has measure at most $2^{k+1} \mu(G)$ .
\end{lemma}
\begin{proof}
\leanok
The union of dyadic cubes in $\mathcal{C}(G,k)$
is contained the union of elements of the set $\mathcal{M}(k)$
of all dyadic cubes $J$ with
${\mu(G \cap J)} > 2^{-k-1}{\mu(J)}$.
The union of elements in the set $\mathcal{M}(k)$ is contained in the union of elements in
the set $\mathcal{M}^*(k)$ of maximal elements in
$\mathcal{M}(k)$ with respect to set inclusion. Hence
\begin{equation}\label{cbymstar}
\mu (\bigcup \mathcal{C}(G,k))\le \mu (\bigcup \mathcal{M}^*(k))\le
\sum_{J\in \mathcal{M}^*(k)}\mu(J)
\end{equation}
Using the definition of $\mathcal{M}(k)$ and then
the pairwise disjointedness of elements in
$\mathcal{M}^*(k)$,
we estimate \eqref{cbymstar} by
\begin{equation}
\le
2^{k+1}\sum_{J\in \mathcal{M}^*(k)}\mu(J\cap G)
\le 2^{k+1}\mu(G).
\end{equation}
This proves the lemma.
\end{proof}

\begin{lemma}[pairwise disjoint]
    \label{pairwise-disjoint}
    \leanok
    \lean{pairwiseDisjoint_E1}
    If $\fp, \fp' \in {\mathfrak{M}}(k,n)$ and
    \begin{equation}\label{eintersect}
        {E_1}(\fp)\cap {E_1}(\fp')\neq \emptyset,
    \end{equation}
    then $\fp=\fp'$.
\end{lemma}
\begin{proof}
\leanok
Let $\fp,\fp'$ be as in the lemma. By definition of $E_1$,
we have
$E_1(\fp)\subset \scI(\fp)$ and analogously for $\fp'$, we conclude from \eqref{eintersect} that $\scI(\fp)\cap \scI(\fp')\neq \emptyset$. Let without loss of generality $\scI(\fp)$ be maximal in
$\{\scI(\fp),\scI(\fp')\}$, then $\scI(\fp')\subset \scI(\fp)$.
By \eqref{eintersect}, we conclude by definition of $E_1$ that $\fc(\fp)\cap \fc(\fp')\neq \emptyset$. By
\eqref{eq-freq-dyadic} we conclude $\fc(\fp)\subset \fc(\fp')$. It follows that $\fp'\le \fp$. By maximality
\eqref{mnkmax}
of $\fp'$, we have $\fp'=\fp$. This proves the lemma.
\end{proof}

\begin{lemma}[dyadic union]
\label{dyadic-union}
\leanok
\lean{dyadic_union}
For each $x\in A(\lambda,k,n)$,
there is a dyadic cube $I$
that contains $x$ and is
a subset of
$A(\lambda,k,n)$.
\end{lemma}

\begin{proof}
\leanok
Fix $k,n,\lambda,x$ as in the lemma such that $x\in A(\lambda,k,n)$.
Let $\mathcal{M}$ be the set of dyadic cubes $\scI(\fp)$ with $\fp$ in $\mathfrak{M}(k,n)$ and $x\in \scI(\fp)$.
By definition of $A(\lambda,k,n)$, the cardinality of $\mathcal{M}$ is at least $1+\lambda 2^{n+1}$.
Let $I$ be a cube of smallest scale in $\mathcal{M}$.
Then $I$ is contained in all cubes of $\mathcal{M}$.
It follows that $I\subset A(\lambda,k,n)$.
\end{proof}

\begin{lemma}[John Nirenberg]
\label{John-Nirenberg}
\leanok
\lean{john_nirenberg}
\uses{dense-cover, pairwise-disjoint, dyadic-union}
    For all integers $k,n,\lambda\ge 0$, we have
    \begin{equation}\label{alambdameasure}
        \mu(A(\lambda,k,n)) \le 2^{k+1-\lambda}\mu(G)\,.
    \end{equation}

\end{lemma}
\begin{proof}
\leanok
Fix $k,n$ as in the lemma
and suppress notation to write
$A(\lambda)$ for $A(\lambda,k,n)$.
We prove the lemma by induction on $\lambda$.
For $\lambda=0$, we use that $A(\lambda)$ by definition of $\mathfrak{M}(k,n)$ is contained in the union of elements in $ \mathcal{C}(G,k)$. \Cref{dense-cover} then completes the base of the induction.

Now assume that the statement of \Cref{John-Nirenberg}
is proven for some integer $\lambda\ge 0$.
The set $A(\lambda+1)$ is contained in the set $A(\lambda)$.
Let $\mathcal{M}$ be the set of dyadic cubes which are a subset of $A(\lambda)$. By \Cref{dyadic-union}, the union of $\mathcal{M}$ is $A(\lambda)$.
Let $\mathcal{M}^*$ be the set of maximal dyadic cubes in $\mathcal{M}$.

Let $x\in A(\lambda+1)$ and $L\in \mathcal{M}^*$ such that $x\in L$. Then by the dyadic property \eqref{dyadicproperty}
\begin{equation}\label{suminout}
    \sum_{\fp \in {\mathfrak{M}}(k,n)} \mathbf{1}_{\scI(\fp)}(x) =
    \sum_{\fp \in {\mathfrak{M}}(k,n):\scI(\fp) \subset L} \mathbf{1}_{\scI(\fp)}(x) +
    \sum_{\fp \in {\mathfrak{M}}(k,n):L \subsetneq \scI(\fp)} \mathbf{1}_{\scI(\fp)}(x)\,.
\end{equation}
We now show
\begin{equation}\label{mnkonl}
    \sum_{\fp \in {\mathfrak{M}}(k,n):\scI(\fp) \subset L} \mathbf{1}_{\scI(\fp)}(x)\ge 2^{n+1}\,.
\end{equation}
The left-hand side of \eqref{suminout} is strictly greater than $(\lambda+1)2^{n+1}$.
If $L$ is the top cube the second sum $Q_2$ on the right-hand side of \eqref{suminout} is zero and
\eqref{mnkonl} follows immediately. Otherwise consider the inclusion-minimal cube $\hat{L}$ with $L\subsetneq\hat{L}$;
all tiles $\fp$ over which $Q_2$ is summed over satisfy $\hat{L}\subset\fp$, so $Q_2$ is constant for all $x\in\hat{L}$.
By maximality of $L$, $Q_2$ is at most $\lambda 2^{n+1}$ somewhere on $\hat{L}$,
thus on all of $\hat{L}$ and consequently also at $x$. Rearranging the inequality yields \eqref{mnkonl}.

By \Cref{pairwise-disjoint}, we have
\begin{equation}
\sum_{\fp \in {\mathfrak{M}}(k,n):\scI(\fp) \subset L} \mu({E_1}(\fp)) \leq \mu(L)\, .
\end{equation}
Multiplying by $2^n$ and applying \eqref{ebardense}, we obtain
\begin{equation}\label{mnkintl}
    \sum_{\fp \in {\mathfrak{M}}(k,n):\scI(\fp) \subset L} \mu(\scI(\fp)) \leq 2^n \mu(L)\, .
\end{equation}
We then have with \eqref{mnkonl} and \eqref{mnkintl}
\begin{equation}
2^{n+1}\mu(A(\lambda+1)\cap L) =
 \int_{A(\lambda+1)\cap L} 2^{n+1} d\mu
\end{equation}
\begin{equation}
\le
    \int \sum_{\fp \in {\mathfrak{M}}(k,n):\scI(\fp) \subset L} \mathbf{1}_{\scI(\fp)} d\mu
\le 2^n \mu(L)\, .
\end{equation}
Hence
\begin{equation}
    2\mu(A(\lambda+1))=2\sum_{L\in \mathcal{M}^*}
\mu(A(\lambda+1)\cap L)\le
\sum_{L\in \mathcal{M}^*}\mu( L)= \mu(A(\lambda))\, .
\end{equation}
Using the induction hypothesis, this proves
\eqref{alambdameasure} for $\lambda+1$ and completes the proof of the lemma.
\end{proof}

\begin{lemma}[second exception]
\label{second-exception}
\leanok
\lean{second_exception}
\uses{John-Nirenberg}
We have
\begin{equation}
    \mu(G_2)\le 2^{-2} \mu(G)\, .
\end{equation}
\end{lemma}
\begin{proof}
\leanok
We use \Cref{John-Nirenberg} and sum twice a geometric series
to obtain
\begin{equation}
    \sum_{0\le k}\sum_{k\le n}
\mu(A(2n+6,k,n))\le \sum_{0\le k}\sum_{k\le n} 2^{k-5-2n}\mu(G)
\end{equation}
\begin{equation}
   \le \sum_{0\le k} 2^{-k-3}\mu(G)\le 2^{-2}\mu(G)\, .
\end{equation}
This proves the lemma.
\end{proof}

We turn to the set $G_3$.

\begin{lemma}[top tiles]
\label{top-tiles}
\leanok
\lean{top_tiles}
 \uses{John-Nirenberg}
    We have
    \begin{equation}\label{eq-musum}
        \sum_{\mathfrak{m} \in \mathfrak{M}(k,n)} \mu(\scI(\mathfrak{m}))\le 2^{n+k+3}\mu(G).
    \end{equation}
\end{lemma}
\begin{proof}
\leanok
We write the left-hand side of \eqref{eq-musum}
\begin{equation}
    \int \sum_{\mathfrak{m} \in \mathfrak{M}(k,n)} \mathbf{1}_{\scI(\mathfrak{m})}(x) \, d\mu(x) \le
2^{n+1} \sum_{\lambda=0}^{|\mathfrak{M}|}\mu(A(\lambda, k,n))\,.
\end{equation}
Using \Cref{John-Nirenberg}
and then summing a geometric series, we estimate this by
\begin{equation}
    \le
2^{n+1}\sum_{\lambda=0}^{|\mathfrak{M}|}
2^{k+1-\lambda}\mu(G)
\le
2^{n+1}2^{k+2}\mu(G)\, .
\end{equation}
This proves the lemma.
\end{proof}

\begin{lemma}[tree count]
\label{tree-count}
\leanok
\lean{tree_count}
Let $k,n,j\ge 0$. We have for every $x\in X$
\begin{equation}
    \sum_{\fu\in \fU_1(k,n,j)} \mathbf{1}_{\scI(\fu)}(x)
    \le 2^{-j}
    2^{9a} \sum_{\mathfrak{m}\in \mathfrak{M}(k,n)}
     \mathbf{1}_{\scI(\mathfrak{m})}(x)
\end{equation}
\end{lemma}

\begin{proof}
\leanok
Let $x\in X$. For each
$\fu\in \fU_1(k,n,j)$ with $x\in \scI(\fu)$, as $\fu \in \fC_1(k,n,j)$,
there are at least $2^{j}$ elements $\mathfrak{m}\in \mathfrak{M}(k,n)$
with $100\fu \lesssim \mathfrak{m}$ and in particular
$x\in \scI(\mathfrak{m})$. Hence
\begin{equation}\label{ubymsum}
     \mathbf{1}_{\scI(\fu)}(x)
    \le 2^{-j}\sum_{\mathfrak{m} \in \mathfrak{M}(k,n): 100\fu\lesssim \mathfrak{m}} \mathbf{1}_{\scI(\mathfrak{m})}(x)\, .
\end{equation}
Conversely, for each $\mathfrak{m}\in \mathfrak{M}(k,n)$
with $x\in \scI(\mathfrak{m})$,
let $\fU(\mathfrak{m})$ be the set of
$\fu\in \fU_1(k,n,j)$ with $x\in \scI(\fu)$
and $100\fu \lesssim \mathfrak{m}$.
Summing \eqref{ubymsum} over $\fu$ and counting the pairs
$(\fu,\mathfrak{m})$ with $100\fu \lesssim \mathfrak{m}$
differently gives
\begin{equation}\label{usumbymsum}
     \sum_{\fu\in \fU_1(k,n,j)} \mathbf{1}_{\scI(\fu)}(x)
    \le 2^{-j}\sum_{\mathfrak{m} \in \mathfrak{M}(k,n)}
    \sum_{\fu \in \fU(\mathfrak{m})} \mathbf{1}_{\scI(\mathfrak{m})}(x)\, .
\end{equation}
We estimate the number of elements in $\fU(\mathfrak{m})$.
Let $\fu \in \fU(\mathfrak{m})$.
Then by definition of
$\fU(\mathfrak{m})$
\begin{equation}\label{dby2}
     d_{\fu}(\fcc(\fu),\fcc(\mathfrak{m}))\le 100\, .
\end{equation}
If $\fu'$ is a further element in $\fU(\mathfrak{m})$ with $\fu\neq \fu'$, then
\begin{equation}
    \fcc(\mathfrak{m})
    \in B_{\fu}(\fcc(\fu),100)\cap B_{\fu'}(\fcc(\fu'),100)\ .
\end{equation}
By the last display and definition of $\fU_1(k,n,j)$, none of $\scI(\fu)$, $\scI(\fu')$ is strictly contained in the other. As both contain $x$, we have $\scI(\fu)=\scI(\fu')$.
We then have $d_{\fu}=d_{\fu'}$.

By \eqref{eq-freq-comp-ball}, the balls
$B_{\fu}(\fcc(\fu),0.2)$ and
$B_{\fu}(\fcc(\fu'),0.2)$ are
contained respectively in $\fc(\fu)$
and $\fc(\fu')$ and thus are disjoint by \eqref{eq-dis-freq-cover}.
By \eqref{dby2} and the triangle inequality, both balls are contained in $B_{\fu}(\fcc(\mathfrak{m}), 100.2)$.

By \eqref{thirddb} applied nine times, there is a collection of at most
$2^{9a}$ balls of radius $0.2$ with respect to the metric $d_{\fu}$ which cover the ball $B_{\fu}(\fcc(\mathfrak{m}),100.2)$.
Let $B'$ be a ball in this cover.
As the center of $B'$ can be in at most one of the disjoint balls
$B_{\fu}(\fcc(\fu),0.2)$ and
$B_{\fu}(\fcc(\fu'),0.2)$,
the ball $B'$ can contain at most
one of the points $\fcc(\fu)$, $\fcc(\fu')$.

Hence the image of $\fU(\mathfrak{m})$ under $\fcc$ has at most
$2^{9a}$ elements; since $\fcc$ is injective on $\fU(\mathfrak{m})$,
the same is true of $\fU(\mathfrak{m})$.
Inserting this into \eqref{usumbymsum} proves the lemma.
\end{proof}

\begin{lemma}[boundary exception]
\label{boundary-exception}
\leanok
\lean{boundary_exception}

Let $\mathcal{L}(\fu)$ be as defined in \eqref{eq-L-def}. We have for each $\fu\in \fU_1(k,n,l)$,
\begin{equation}
\mu(\bigcup_{I\in \mathcal{L}(\fu)} I)
\le D^{1-\kappa Z(n+1)}
        \mu(\scI(\mathfrak{u})).
\end{equation}
\end{lemma}
\begin{proof}
\leanok
Let $\fu\in \fU_1(k,n,l)$.
Let $I \in \mathcal{L}(\fu)$. Then we have $s(I) = \ps(\fu) - Z(n+1) - 1$ and $I \subset \scI(\fu)$ and $B(c(I), 8D^{s(I)}) \not \subset \scI(\fu)$.
By \eqref{eq-vol-sp-cube}, the set $I$
is contained in $B(c(I), 4D^{s(I)})$.
By the triangle inequality, the set $I$
is contained in
\begin{equation}
    X(\fu):=\{x \in \scI(\fu) \, : \, \rho(x, X \setminus \scI(\fu)) \leq 12 D^{\ps(\fu) - Z(n+1)-1}\}\,.
\end{equation}
 By the small boundary property \eqref{eq-small-boundary}, noting that
 \begin{equation*}
     12D^{\ps(\fu) - Z(n+1) - 1} = 12D^{s(I)} > D^{-S}\ ,
 \end{equation*} we have
   $$
        \mu(X(\fu)) \le
        2\cdot(12 D^{-Z(n+1)-1})^\kappa
        \mu(\scI(\mathfrak{u})).
    $$
Using $\kappa<1$ and $D \ge 12$, this proves the lemma.
\end{proof}

\begin{lemma}[third exception]
    \label{third-exception}
    \leanok
    \lean{third_exception}
\uses{tree-count,boundary-exception, top-tiles}
       We have
\begin{equation}
    \mu(G_3)\le 2^{-4} \mu(G)\, .
\end{equation}
\end{lemma}
\begin{proof}
\leanok
As each $\fp\in \fL_4(k,n,j)$
is contained in $\cup\mathcal{L}(\fu)$ for some
$\fu\in \fU_1(k,n,l)$, we have
\begin{equation}
\mu(\bigcup_{\fp \in \fL_4 (k,n,j)}\scI(\fp))
\le \sum_{\fu\in \fU_1(k,n,j)}
\mu(\bigcup_{I \in \mathcal{L} (\fu)}I).
\end{equation}
Using \Cref{boundary-exception} and then \Cref{tree-count}, we estimate this further
 by
\begin{equation}
    \le \sum_{\fu\in \fU_1(k,n,j)}
    D^{1-\kappa Z(n+1)}
        \mu(\scI(\mathfrak{u}))
\end{equation}
\begin{equation}
    \le 2^{100a^2+9a+1-j} \sum_{\mathfrak{m}\in \mathfrak{M}(k,n)}
     D^{-\kappa Z(n+1)}
    \mu(\scI(\mathfrak{m}))\,.
\end{equation}
Using \Cref{top-tiles}, we estimate this by
  \begin{equation}
     \le
2^{100a^2 + 9a + 1-j} D^{-\kappa Z(n+1)}
     2^{n+k+3}\mu(G)\, .
\end{equation}
Now we estimate $G_3$ defined in \eqref{defineg3} by
\begin{equation}
    \mu(G_3)\le \sum_{k\ge 0}\, \sum_{n \geq k}\,
    \sum_{0\le j\le 2n+3}
    \mu(\bigcup_{\fp \in \fL_4 (k,n,j)}
    \scI(\fp))
\end{equation}
\begin{equation}
    \le \sum_{k\ge 0}\, \sum_{n \geq k}\,
    \sum_{0\le j\le 2n+3}
    2^{100a^2 + 9a + 3 + n + k -j} D^{-\kappa Z(n+1)}\mu(G)
\end{equation}
Summing geometric series, using that $D^{\kappa Z}\ge 8$ by \eqref{defineD}, \eqref{definekappa} and \eqref{defineZ}, we estimate this by
\begin{equation}
    \le \sum_{k\ge 0}\, \sum_{n \geq k}\,
    2^{100a^2 + 9a + 4 + n + k} D^{-\kappa Z(n+1)}\mu(G)
\end{equation}
\begin{equation}
    = \sum_{k\ge 0} 2^{100a^2 + 9a + 4 + 2k} D^{-\kappa Z(k+1)} \sum_{n \geq k}\,
    2^{n - k} D^{-\kappa Z(n-k)}\mu(G)
\end{equation}
\begin{equation}
    \le \sum_{k\ge 0} 2^{100a^2 + 9a + 5 + 2k} D^{-\kappa Z(k+1)}\mu(G)
\end{equation}
\begin{equation}
   \le 2^{100a^2 + 9a + 6} D^{-\kappa Z}\mu(G)
\end{equation}
Using $D = 2^{100a^2}$ and $a \ge 4$ and $\kappa Z \ge 2$ by \eqref{defineD} and \eqref{definekappa} proves the lemma.
\end{proof}

\begin{proof}[Proof of \Cref{exceptional-set}]
    \leanok
    \proves{exceptional-set}
Adding up the bounds in Lemmas \ref{first-exception}, \ref{second-exception}, and \ref{third-exception} proves \Cref{exceptional-set}.
\end{proof}

\subsection{Auxiliary lemmas}
\label{subsec-lessim-aux}
Before proving \Cref{forest-union} and \Cref{forest-complement}, we collect some useful properties of $\lesssim$.

\begin{lemma}[wiggle order 1]
    \label{wiggle-order-1}
    \leanok
    \lean{smul_mono}
    If $n\fp \lesssim m\fp'$ and
    $n' \ge n$ and $m \ge m'$ then $n'\fp \lesssim m'\fp'$.
\end{lemma}

\begin{proof}
    \leanok
    This follows immediately from the definition \eqref{wiggleorder} of $\lesssim$ and the two inclusions $B_{\fp}(\fcc(\fp), n) \subset B_{\fp}(\fcc(\fp), n')$ and $B_{\fp'}(\fcc(\fp'), m') \subset B_{\fp'}(\fcc(\fp'), m)$.
\end{proof}

\begin{lemma}[wiggle order 2]
    \label{wiggle-order-2}
    \leanok
    \lean{smul_C2_1_2}
    \uses{monotone-cube-metrics}

    Let $n, m \ge 1$ and $k > 0$.
    If $\fp, \fp' \in \fP$ with $\scI(\fp) \ne \scI(\fp')$ and
    \begin{equation}
        \label{eq-wiggle1}
        n \fp \lesssim k \fp'
    \end{equation}
    then
    \begin{equation}
        \label{eq-wiggle2}
        (n + 2^{-95 a} m) \fp \lesssim m\fp'\,.
    \end{equation}
\end{lemma}

\begin{proof}
    \leanok
    The assumption \eqref{eq-wiggle1} together with the definition \eqref{wiggleorder} of $\lesssim$ implies that $\scI(\fp) \subsetneq \scI(\fp')$. Let $\mfa \in B_{\fp'}(\fcc(\fp'), m)$. Then we have by the triangle inequality
    $$
        d_{\fp}(\fcc(\fp), \mfa) \le d_{\fp}(\fcc(\fp), \fcc(\fp')) + d_{\fp}(\fcc(\fp'), \mfa)
    $$
    The first summand is bounded by $n$ since
    $$
        \fcc(\fp') \in B_{\fp'}(\fcc(\fp'), k) \subset B_{\fp}(\fcc(\fp), n),
    $$
    using \eqref{wiggleorder}. For the second summand we use \Cref{monotone-cube-metrics}
    to show that the sum is estimated by
    $$
        n + 2^{-95a} d_{\fp'}(\fcc(\fp'), \mfa) < n + 2^{-95a} m\,.
    $$
    Thus $B_{\fp'}(\fcc(\fp'),m) \subset B_{\fp}(\fcc(\fp),n + 2^{-95a}m)$. Combined with $\scI(\fp) \subset \scI(\fp')$, this yields \eqref{eq-wiggle2}.
\end{proof}

\begin{lemma}[wiggle order 3]
\label{wiggle-order-3}
\leanok
\lean{wiggle_order_11_10, wiggle_order_100, wiggle_order_500}
\uses{wiggle-order-1, wiggle-order-2}
    The following implications hold for all $\fq, \fq' \in \fP$:
    \begin{equation}
        \label{eq-sc1}
        \fq \le \fq' \ \text{and} \ \lambda \ge 1.1 \implies \lambda \fq \lesssim \lambda \fq'\,,
    \end{equation}
    \begin{equation}
        \label{eq-sc2}
        10\fq \lesssim \fq' \ \text{and} \ \scI(\fq) \ne \scI(\fq') \implies 100 \fq \lesssim 100 \fq'\,,
    \end{equation}
    \begin{equation}
        \label{eq-sc3}
        2\fq \lesssim \fq' \ \text{and} \ \scI(\fq) \ne \scI(\fq') \implies 4 \fq \lesssim 500 \fq'\,.
    \end{equation}
\end{lemma}

\begin{proof}
    \leanok
    \eqref{eq-sc2} and \eqref{eq-sc3} are easy consequences of \Cref{wiggle-order-1}, \Cref{wiggle-order-2} and the fact that $a \ge 4$.
    For \eqref{eq-sc1}, if $\scI(\fq) = \scI(\fq')$ then
    we get $\fq = \fq'$ by \eqref{eq-dis-freq-cover} and \eqref{straightorder}.
    If $\scI(\fq) \ne \scI(\fq')$, then from \eqref{straightorder},
    \eqref{wiggleorder} and \eqref{eq-freq-comp-ball} it follows that
    $\fq \lesssim 0.2\fq'$, and \eqref{eq-sc1} follows from an easy calculation using
    \Cref{wiggle-order-2}.
\end{proof}

We call a collection $\mathfrak{A}$ of tiles convex if
\begin{equation}
    \label{eq-convexity}
    \fp \le \fp' \le \fp'' \ \text{and} \ \fp, \fp'' \in \mathfrak{A} \implies \fp' \in \mathfrak{A}\,.
\end{equation}

\begin{lemma}[P convex]
    \label{P-convex}
    \leanok
    \lean{ordConnected_tilesAt}

    For each $k$, the collection $\fP(k)$ is convex.
\end{lemma}

\begin{proof}
    \leanok
    Suppose that $\fp \le \fp' \le \fp''$ and $\fp, \fp'' \in \fP(k)$. By \eqref{eq-defPk} we have $\scI(\fp), \scI(\fp'') \in \mathcal{C}(G,k)$, so there exists by \eqref{muhj1} some $J \in \mathcal{D}$ with
    $$
        \scI(\fp') \subset \scI(\fp'') \subset J
    $$
    and $\mu(G \cap J) > 2^{-k-1} \mu(J)$. Thus \eqref{muhj1} holds for $\scI(\fp')$. On the other hand, by \eqref{muhj2}, there exists no $J \in \mathcal{D}$ with $\scI(\fp) \subset J$ and $\mu(G \cap J) > 2^{-k} \mu(J)$. Since $\scI(\fp) \subset \scI(\fp')$, this implies that \eqref{muhj2} holds for $\scI(\fp')$. Hence $\scI(\fp') \in \mathcal{C}(G,k)$, and therefore by \eqref{eq-defPk} $\fp' \in \fP(k)$.
\end{proof}

\begin{lemma}[C convex]
    \label{C-convex}
    \leanok
    \lean{ordConnected_C}
    \uses{P-convex}
    For each $k,n$, the collection $\fC(k,n)$ is convex.
\end{lemma}

\begin{proof}
    \leanok
    Let $\fp \le \fp' \le \fp''$ with $\fp, \fp'' \in \fC(k,n)$. Then, in particular, $\fp, \fp'' \in \fP(k)$, so, by \Cref{P-convex}, $\fp' \in \fP(k)$. Next, we show that if $\fq \le \fq' \in \fP(k)$ then $\dens'_k(\{\fq\}) \ge \dens_k'(\{\fq'\})$. If $\fp \in \fP(k)$ and $\lambda \ge 2$ with $\lambda \fq' \lesssim \lambda \fp$, then it follows from $\fq \le \fq'$, \eqref{eq-sc1} of \Cref{wiggle-order-3} and transitivity of $\lesssim$ that $\lambda \fq \lesssim \lambda \fp$. Thus the supremum in the definition \eqref{eq-densdef} of $\dens_k'(\{\fq\})$ is over a superset of the set the supremum in the definition of $\dens_k'(\{\fq'\})$ is taken over, which shows $\dens'_k(\{\fq\}) \ge \dens_k'(\{\fq'\})$. From $\fp' \le \fp''$, $\fp'' \in \fC(k,n)$ and \eqref{def-cnk} it then follows that
    $$
        2^{4a}2^{-n} < \dens_k'(\{\fp''\}) \le \dens_k'(\{\fp'\})\,.
    $$
    Similarly, it follows from $\fp \le \fp'$, $\fp \in \fC(k,n)$ and \eqref{def-cnk} that
    $$
        \dens_k'(\{\fp'\}) \le \dens_k'(\{\fp\}) \le 2^{4a}2^{-n+1}\,.
    $$
    Thus $\fp' \in \fC(k,n)$.
\end{proof}

\begin{lemma}[C1 convex]
    \label{C1-convex}
    \leanok
    \lean{ordConnected_C1}
    \uses{C-convex}
    For each $k,n,j$, the collection $\fC_1(k,n,j)$ is convex.
\end{lemma}

\begin{proof}
    \leanok
    Let $\fp\le\fp'\le\fp''$ with $\fp, \fp'' \in \fC_1(k,n,j)$. By \Cref{C-convex} and the inclusion $\fC_1(k,n,j) \subset \fC(k,n)$, which holds by definition \eqref{defcnkj}, we have $\fp' \in \fC(k,n)$. By \eqref{eq-sc1} and transitivity of $\lesssim$ we have that $\fq \le \fq'$ and $100 \fq' \lesssim \mathfrak{m}$ imply $100 \fq \lesssim \mathfrak{m}$. So, by \eqref{defbfp}, $\mathfrak{B}(\fp'') \subset \mathfrak{B}(\fp') \subset \mathfrak{B}(\fp)$. Consequently, by \eqref{defcnkj}
    $$
        2^j \le |\mathfrak{B}(\fp'')|\le |\mathfrak{B}(\fp')| \le |\mathfrak{B}(\fp)| < 2^{j+1}\,,
    $$
    thus $\fp' \in \fC_1(k,n,j)$.
\end{proof}

\begin{lemma}[C2 convex]
    \label{C2-convex}
    \leanok
    \lean{ordConnected_C2}
    \uses{C1-convex}
    For each $k,n,j$, the collection $\fC_2(k,n,j)$ is convex.
\end{lemma}

\begin{proof}
    \leanok
     Let $\fp \le \fp' \le \fp''$ with $\fp, \fp'' \in \fC_2(k,n,j)$. By \eqref{eq-C2-def}, we have
     \begin{equation*}
         \fC_2(k,n,j) \subset \fC_1(k,n,j)\ .
     \end{equation*} Combined with \Cref{C1-convex}, it follows that $\fp' \in \fC_1(k,n,j)$.
     If $\fp=\fp'$ the statement is trivially true, otherwise suppose that $\fp' \notin \fC_2(k,n,j)$.
     By \eqref{eq-C2-def}, this implies that there exists $0 \le l' \le Z(n+1)$ with $\fp' \in \fL_1(k,n,j,l')$.
     By the definition \eqref{eq-L1-def} of $\fL_1(k,n,j,l')$, this implies that $\fp'$ is minimal
     with respect to $\le$ in $\fC_1(k,n,j) \setminus \bigcup_{l < l'} \fL_1(k,n,j,l)$.
     Since $\fp\le\fp'$ and $\fp\in\fC_2(k,n,j)$, $\fp=\fp'$, a contradiction.
\end{proof}

\begin{lemma}[C3 convex]
    \label{C3-convex}
    \leanok
    \lean{ordConnected_C3}
    \uses{C2-convex}
    For each $k,n,j$, the collection $\fC_3(k,n,j)$ is convex.
\end{lemma}

\begin{proof}
    \leanok
    Let $\fp \le \fp' \le \fp''$ with $\fp, \fp'' \in \fC_3(k,n,j)$.
    By \eqref{eq-C3-def} and \Cref{C2-convex} it follows that $\fp' \in \fC_2(k,n,j)$.
    By \eqref{eq-C3-def} and \eqref{eq-L2-def}, there exists $\fu \in \fU_1(k,n,j)$ with
    $2\fp'' \lesssim \fu$ and $\scI(\fp'') \subsetneq \scI(\fu)$. From $\fp' \le \fp''$,
    \eqref{eq-sc1} and transitivity of $\lesssim$ we then have
    $2\fp' \lesssim \fu$ and $\scI(\fp') \subsetneq \scI(\fu)$, so $\fp' \in \fC_3(k,n,j)$.
\end{proof}

\begin{lemma}[C4 convex]
    \label{C4-convex}
    \leanok
    \lean{ordConnected_C4}
    \uses{C3-convex}
    For each $k,n,j$, the collection $\fC_4(k,n,j)$ is convex.
\end{lemma}

\begin{proof}
    \leanok
    The proof is entirely analogous to \Cref{C2-convex}, substituting $\fC_4$ for $\fC_2$,
    $\fC_3$ for $\fC_1$ and $\fp'\le\fp''$ for $\fp\le\fp'$.
\end{proof}

\begin{lemma}[C5 convex]
    \label{C5-convex}
    \leanok
    \lean{ordConnected_C5}
    \uses{C4-convex}
    For each $k,n,j$, the collection $\fC_5(k,n,j)$ is convex.
\end{lemma}

\begin{proof}
    \leanok
    Let $\fp \le \fp' \le\fp''$ with $\fp, \fp'' \in \fC_5(k,n,j)$.
    Then $\fp, \fp'' \in \fC_4(k,n,j)$ by \eqref{defc5}, and thus by \Cref{C4-convex} $\fp' \in \fC_4(k,n,j)$.
    It suffices to show that if $\fp' \in \fL_4(k,n,j)$ then $\fp \in \fL_4(k,n,j)$ by contraposition;
    this is true by \eqref{eq-L4-def} and $\fp\le\fp'$.
\end{proof}

\begin{lemma}[dens compare]
    \label{dens-compare}
    \leanok
    \lean{dens1_le_dens'}
    We have for every $k\ge 0$ and $\fP'\subset \fP(k)$
\begin{equation}
    \dens_1(\fP')\le \dens_k'(\fP')\, .
\end{equation}
\end{lemma}
\begin{proof}
\leanok
It suffices to show that for all $\fp'\in \fP'$
and $\lambda\ge 2$ and $\fp\in \fP(\fP')$ with $\lambda \fp' \lesssim \lambda \fp$ we have
\begin{equation}
    \frac{\mu({E}_2(\lambda, \fp))}{\mu(\scI(\fp))}
    \le \sup_{\fp'' \in \fP(k): \lambda \fp' \lesssim \lambda \fp''}
    \frac{\mu({E}_2(\lambda, \fp''))}{\mu(\scI(\fp''))}.
\end{equation}
    Let such $\fp'$, $\lambda$, $\fp$ be given.
    It suffices to show that $\fp\in \fP(k)$,
    that is, it satisfies \eqref{muhj1}
    and \eqref{muhj2}.

We show \eqref{muhj1}.
 As $\fp\in \fP(\fP')$, there exists
$\fp''\in \fP'$ with $\scI(\fp')\subset \scI(\fp'')$. By assumption on $\fP'$, we have $\fp''\in \fP(k)$ and there exists
$J\in \mathcal{D}$ with
   $\scI(\fp'')\subset J$ and
   \begin{equation}
       \mu(G\cap J)>2^{-k-1} \mu(J).
   \end{equation}
Then also $\scI(\fp')\subset J$, which proves
\eqref{muhj1} for $\fp$.

We show \eqref{muhj2}. Assume to get a contradiction that
there exists $J\in \mathcal{D}$ with
   $\scI(\fp)\subset J$ and
   \begin{equation}\label{mugj}
       \mu(G\cap J)>2^{-k} \mu(J).
   \end{equation}
   As $\lambda\fp'\lesssim \lambda\fp$, we have $\scI(\fp')\subset \scI(\fp)$, and therefore
    $\scI(\fp')\subset J$. This contradicts
   $\fp'\in \fP'\subset \fP(k)$. This proves
\eqref{muhj2} for $\fp$.
\end{proof}

\begin{lemma}[C dens1]
    \label{C-dens1}
    \leanok
    \lean{dens1_le}
    \uses{dens-compare}
    For each set $\mathfrak{A} \subset \mathfrak{C}(k,n)$, we have
    $$
        \dens_1(\mathfrak{A}) \le 2^{4a}2^{-n+1}\,.
    $$
\end{lemma}

\begin{proof}
    \leanok
    We have by \Cref{dens-compare} that
    $\dens_1(\mathfrak{A}) \le \dens_k'(\mathfrak{A})$. Since $\mathfrak{A} \subset \fC(k,n)$, it follows from monotonicity of suprema and the definition \eqref{eq-densdef} that
    $
        \dens_k'(\mathfrak{A}) \le \dens_k'(\fC(k,n))\,.
    $
    By \eqref{eq-densdef} and \eqref{def-cnk}, we have
    $$
        \dens_k'(\fC(k,n)) = \sup_{\fp \in \fC(k,n)} \dens_k'(\{\fp\}) \le 2^{4a}2^{-n+1}\,.
    $$
\end{proof}

\subsection{Proof of the Forest Union Lemma}
\label{subsecforest}

Fix $k,n,j\ge 0$.
Define
$$
    \fC_6(k,n,j)
$$
to be the set of all tiles $\fp \in \fC_5(k,n,j)$ such that $\scI(\fp) \not\subset G'$. The following chain of lemmas
establishes that the set $\fC_6(k,n,j)$ can be written as a union of a small number of $n$-forests.

For $\fu\in \fU_1(k,n,j)$, define
\begin{equation}
    \label{eq-T1-def}
    \mathfrak{T}_1(\fu):= \{\fp \in \fC_1(k,n,j) \ : \scI(\fp)\neq \scI(\fu), \ 2\fp \lesssim \fu\}\,.
\end{equation}
Define
\begin{equation}
    \label{eq-U2-def}
    \fU_2(k,n,j) := \{ \fu \in \fU_1(k,n,j) \, : \, \mathfrak{T}_1(\fu) \cap \fC_6(k,n,j) \ne \emptyset\}\,.
\end{equation}

Define a relation $\sim$ on $\fU_2(k,n,j)$
by setting $\fu\sim \fu'$
for $\fu,\fu'\in \fU_2(k,n,j)$
if $\fu=\fu'$ or there exists $\fp$ in $\mathfrak{T}_1(\fu)$
with $10 \fp\lesssim \fu'$.

\begin{lemma}[relation geometry]
    \label{relation-geometry}
    \leanok
    \lean{URel.eq, URel.not_disjoint}
    \uses{wiggle-order-3}
    If $\fu \sim \fu'$, then $\scI(u) = \scI(u')$ and
    \begin{equation*}
        B_{\fu}(\fcc(\fu), 100) \cap B_{\fu'}(\fcc(\fu'), 100) \neq \emptyset\ .
    \end{equation*}
\end{lemma}

\begin{proof}
    \leanok
    Let $\fu, \fu' \in \fU_2(k,n,j)$ with $\fu \sim \fu'$. If $\fu = \fu'$ then the conclusion of the Lemma clearly holds. Else, there exists $\fp \in \fC_1(k,n,j)$ such that $\scI(\fp) \ne \scI(\fu)$ and $2 \fp \lesssim \fu$ and $10 \fp \lesssim \fu'$.
    Using \Cref{wiggle-order-1} and \eqref{eq-sc2} of \Cref{wiggle-order-3}, we deduce that
    \begin{equation}
        \label{eq-Fefferman-trick0}
        100 \fp\lesssim 100 \fu\,, \qquad 100 \fp \lesssim 100\fu'\,.
    \end{equation}
    Now suppose that $B_{\fu}(\fcc(\fu), 100) \cap B_{\fu'}(\fcc(\fu'), 100) = \emptyset$. Then we have $\mathfrak{B}(\fu) \cap \mathfrak{B}(\fu') = \emptyset$, by the definition \eqref{defbfp} of $\mathfrak{B}$ and the definition \eqref{wiggleorder} of $\lesssim$, but also $\mathfrak{B}(\fu) \subset \mathfrak{B}(\fp)$ and $\mathfrak{B}(\fu') \subset \mathfrak{B}(\fp)$, by \eqref{defbfp}, \eqref{wiggleorder} and \eqref{eq-Fefferman-trick0}.
    Hence,
    $$
        |\mathfrak{B}(\fp)| \geq |\mathfrak{B}(\fu)| + |\mathfrak{B}(\fu')| \geq 2^{j} + 2^j = 2^{j+1}\,,
    $$
    which contradicts $\fp \in \fC_1(k,n,j)$. Therefore we must have
    \begin{equation*}
        B_{\fu}(\fcc(\fu), 100) \cap B_{\fu'}(\fcc(\fu'), 100) \ne \emptyset\, .
    \end{equation*}

    It follows from $2\fp \lesssim \fu$ and $10\fp \lesssim \fu'$ that $\scI(\fp) \subset \scI(\fu)$ and $\scI(\fp) \subset \scI(\fu')$. By \eqref{dyadicproperty}, it follows that $\scI(\fu)$ and $\scI(\fu')$ are nested.
    Combining this with the conclusion of the last paragraph and definition \eqref{defunkj} of $\fU_1(k,n,j)$, we obtain that $\scI(\fu) = \scI(\fu')$.
\end{proof}

\begin{lemma}[equivalence relation]
\label{equivalence-relation}
\uses{relation-geometry}
\leanok
\lean{equivalenceOn_urel}
For each $k,n,j$, the relation $\sim$ on
$\fU_2(k,n,j)$ is an equivalence relation.
\end{lemma}

\begin{proof}
    \leanok
    Reflexivity holds by definition.
    For transitivity, suppose that
    \begin{equation*}
        \fu, \fu', \fu'' \in \fU_1(k,n,j)
    \end{equation*}
    and $\fu \sim \fu'$, $\fu' \sim \fu''$.
    By \Cref{relation-geometry}, it follows that $\scI(\fu) =\scI(\fu') = \scI(\fu'')$, that there exists
    \begin{equation*}
        \mfa \in B_{\fu}(\fcc(\fu), 100) \cap B_{\fu'}(\fcc(\fu'), 100)
    \end{equation*}
    and that there exists
    \begin{equation*}
        \mfb \in B_{\fu'}(\fcc(\fu'), 100) \cap B_{\fu''}(\fcc(\fu''), 100)\, .
    \end{equation*}
    If $\fu = \fu'$, then $\fu \sim \fu''$ holds by assumption. Else, there exists by the definition of $\sim$ some $\fp \in \mathfrak{T}_1(\fu)$ with $10\fp\lesssim \fu'$.
    Then we have $2\fp \lesssim \fu$ and $\fp \ne \fu$ by definition of $\mathfrak{T}_1(\fu)$, so $4 \fp \lesssim 500 \fu$ by \eqref{eq-sc3}. For $q \in B_{\fu''}(\fcc(\fu''), 1)$ it follows by the triangle inequality that
    \begin{align*}
        d_{\fu}(\fcc(\fu), q) &\le d_{\fu}(\fcc(\fu), \mfa) + d_{\fu}(\mfa, \fcc(\fu'))\\
        &\quad+ d_{\fu}(\fcc(\fu'), \mfb) + d_{\fu}(\mfb, \fcc(\fu'')) +
        d_{\fu}(\fcc(\fu''), q)\,.
    \end{align*}
    Using \eqref{defdp} and the fact that $\scI(\fu) = \scI(\fu') = \scI(\fu'')$ this equals
    \begin{align*}
        &\quad d_{\fu}(\fcc(\fu), \mfa) + d_{\fu'}(\mfa, \fcc(\fu'))\\
        &\quad+ d_{\fu'}(\fcc(\fu'), \mfb) + d_{\fu''}(\mfb, \fcc(\fu'')) +
        d_{\fu''}(\fcc(\fu''), q)\\
        &< 100 + 100 + 100 + 100 + 1 < 500\,.
    \end{align*}
    Since $4\fp \lesssim 500 \fu$, it follows that $d_{\fp}(\fcc(\fp), q) < 4 < 10$. We have shown that $B_{\fu''}(\fcc(\fu''), 1) \subset B_{\fp}(\fcc(\fp), 10)$, combining this with $\scI(\fu'') = \scI(\fu)$ gives $\fu \sim \fu''$.

    For symmetry suppose that $\fu \sim \fu'$. By Lemma \eqref{relation-geometry}, it follows that $\scI(\fu) = \scI(\fu')$ and that there exists $\mfa \in B_{\fu}(\fcc(\fu), 100) \cap B_{\fu'}(\fcc(\fu'), 100)$. Again, for $\fu = \fu'$ symmetry is obvious, so suppose that $\fu \ne \fu'$. There exists $\fp \in \mathfrak{T}_1(\fu')$, which then satisfies $2\fp\lesssim \fu'$ and $\scI(\fp) \neq \scI(\fu')$. By \Cref{wiggle-order-1} and \eqref{eq-sc3}, it follows that
    \begin{equation}
        \label{eq-rel1}
        10\fp \lesssim 4\fp \lesssim 500 \fu'\,.
    \end{equation}
    If $q \in B_{\fu}(\fcc(\fu),1)$ then we have from the triangle inequality and the fact that $\scI(\fu) = \scI(\fu')$:
    \begin{align*}
        d_{\fu'}(\fcc(\fu'), q) &\le d_{\fu'}(\fcc(\fu'), \mfa) + d_{\fu'}(\mfa, \fcc(\fu)) + d_{\fu'}(\fcc(\fu), q)\\
        &= d_{\fu'}(\fcc(\fu'), \mfa) + d_{\fu}(\mfa, \fcc(\fu)) + d_{\fu}(\fcc(\fu), q)\\
        &< 100 + 100 + 1 < 500\,.
    \end{align*}
    Combining this with \eqref{eq-rel1} and \eqref{wiggleorder}, we get
    \begin{equation*}
     B_{\fu}(\fcc(\fu), 1) \subset B_{\fp}(\fcc(\fp), 10)\, .
    \end{equation*}
    Since $2\fp \lesssim \fu'$, we have $\scI(\fp) \subset \scI(\fu') = \scI(\fu)$. Thus, $10\fp \lesssim \fu$ which completes the proof of $\fu' \sim \fu$.
\end{proof}

Choose a set $\fU_3(k,n,j)$ of representatives for the equivalence
classes of $\sim$ in $\fU_2(k,n,j)$.
Define for each $\fu\in \fU_3(k,n,j)$
\begin{equation}\label{definesv}
\fT_2(\fu):=
   \bigcup_{\fu\sim \fu'}\mathfrak{T}_1(\fu')\cap \fC_6(k,n,j)\, .
\end{equation}

\begin{lemma}[C6 forest]
\label{C6-forest}
\leanok
\lean{C6_forest}
\uses{equivalence-relation}
We have
\begin{equation}
    \fC_6(k,n,j)=\bigcup_{\fu\in \fU_3(k,n,j)}\mathfrak{T}_2(\fu)\, .
\end{equation}
\end{lemma}
\begin{proof}
    \leanok
    Let $\fp \in \fC_6(k,n,j)$.
    By \eqref{eq-C4-def} and \eqref{defc5}, we have $\fp \in \fC_3(k,n,j)$. By \eqref{eq-L2-def} and \eqref{eq-C3-def}, there exists $\fu \in \fU_1(k,n,j)$ with $2\fp \lesssim \fu$ and $\scI(\fp) \ne \scI(\fu)$, that is, with $\fp \in \mathfrak{T}_1(\fu)$. Then $\mathfrak{T}_1(\fu)$ is clearly nonempty, so $\fu \in \fU_2(k,n,j)$. By the definition of $\fU_3(k,n,j)$, there exists $\fu' \in \fU_3(k,n,j)$ with $\fu \sim \fu'$. By \eqref{definesv}, we have $\fp \in \mathfrak{T}_2(\fu')$.
\end{proof}

\begin{lemma}[forest geometry]
    \label{forest-geometry}
    \uses{relation-geometry}
    \leanok
    \lean{forest_geometry}
    For each $\fu\in \fU_3(k,n,j)$,
    the set $\mathfrak{T}_2(\fu)$
    satisfies \eqref{forest1}.
\end{lemma}
\begin{proof}
    \leanok
    Let $\fp \in \mathfrak{T}_2(\fu)$. By \eqref{definesv}, there exists $\fu' \sim \fu$ with $\fp \in \mathfrak{T}_1(\fu')$. Then we have $2\fp \lesssim \fu'$ and $\scI(\fp) \ne \scI(\fu')$, so by \eqref{eq-sc3} $4\fp \lesssim 500\fu'$.
    Further, by \Cref{relation-geometry}, we have that $\scI(\fu') = \scI(\fu)$ and there exists $\mfa \in B_{\fu'}(\fcc(\fu'),100) \cap B_{\fu}(\fcc(\fu),100)$.
    Let $\mfb \in B_{\fu}(\fcc(\fu), 1)$.
    Using the triangle inequality and the fact that $\scI(\fu') =\scI(\fu)$, we obtain
    \begin{align*}
        d_{\fu'}(\fcc(\fu'), \mfb) &\le d_{\fu'}(\fcc(\fu'), \mfa) + d_{\fu'}(\fcc(\fu), \mfa) + d_{\fu'}(\fcc(\fu), \mfb)\\
        &= d_{\fu'}(\fcc(\fu'), \mfa) + d_{\fu}(\fcc(\fu), \mfa) + d_{\fu}(\fcc(\fu), \mfb)\\
        &< 100 + 100 + 1 < 500\,.
    \end{align*}
    Combining this with $4\fp \lesssim 500\fu'$, we obtain
    $$
        B_{\fu}(\fcc(\fu), 1) \subset B_{\fu'}(\fcc(\fu'), 500) \subset B_{\fp}(\fcc(\fp), 4)\,.
    $$
    Together with $\scI(\fp) \subset \scI(\fu') = \scI(\fu)$, this gives $4\fp \lesssim \fu$, which is \eqref{forest1}.
\end{proof}

\begin{lemma}[forest convex]
    \label{forest-convex}
    \uses{C5-convex}
    \leanok
    \lean{forest_convex}
    For each $\fu\in \fU_3(k,n,j)$,
    the set $\mathfrak{T}_2(\fu)$
    satisfies the convexity condition \eqref{forest2}.
\end{lemma}

\begin{proof}
    \leanok
    Let $\fp, \fp'' \in \mathfrak{T}_2(\fu)$ and $\fp' \in \fP$ with $\fp \le \fp' \le \fp''$. By \eqref{definesv} we have $\fp, \fp'' \in \fC_6(k,n,j) \subset \fC_5(k,n,j)$. By \Cref{C5-convex}, we have $\fp' \in \fC_5(k,n,j)$. Since $\fp \in \fC_6(k,n,j)$ we have $\scI(\fp) \not \subset G'$, so $\scI(\fp') \not \subset G'$ and therefore also $\fp' \in \fC_6(k,n,j)$.

    By \eqref{definesv} there exists $\fu' \in \fU_2(k,n,j)$ with $\fp'' \in \mathfrak{T}_1(\fu')$ and hence $2\fp'' \lesssim \fu'$ and $\scI(\fp'') \ne \scI(\fu')$. Together this implies $\scI(\fp'') \subsetneq \scI(\fu')$. With the inclusion $\scI(\fp') \subset \scI(\fp'')$ from $\fp' \le \fp''$, it follows that $\scI(\fp') \subsetneq \scI(\fu')$ and hence $\scI(\fp') \ne \scI(\fu')$.
    By \eqref{eq-sc1} and transitivity of $\lesssim$ we further have $2\fp' \lesssim \fu'$, so $\fp' \in \mathfrak{T}_1(\fu')$.
    It follows that $\fp' \in \mathfrak{T}_2(\fu)$, which shows \eqref{forest2}.
\end{proof}

\begin{lemma}[forest separation]
    \label{forest-separation}
    \uses{monotone-cube-metrics}
    \leanok
    \lean{forest_separation}
    For each $\fu,\fu'\in \fU_3(k,n,j)$ with $\fu\neq \fu'$ and each $\fp \in \fT_2(\fu)$
    with $\scI(\fp)\subset \scI(\fu')$ we have
    \begin{equation}
    d_{\fp}(\fcc(\fp), \fcc(\fu')) > 2^{Z(n+1)}\,.
    \end{equation}
\end{lemma}

\begin{proof}
    \leanok
    By the definition \eqref{eq-C2-def} of $\fC_2(k,n,j)$, there exists a tile $\fp' \in \fC_1(k,n,j)$ with $\fp' \le \fp$ and $\ps(\fp') \le \ps(\fp)- Z(n+1)$.
    By \Cref{monotone-cube-metrics} we have
    $$
        d_{\fp}(\fcc(\fp), \fcc(\fu')) \ge 2^{95a Z(n+1)} d_{\fp'}(\fcc(\fp), \fcc(\fu'))\,.
    $$
    By \eqref{eq-sc1} we have $2\fp' \lesssim 2\fp$, so by transitivity of $\lesssim$ there exists $\mathfrak{v} \sim \fu$ with $2\fp' \lesssim \mathfrak{v}$ and $\scI(\fp') \ne \scI(\mathfrak{v})$. Since $\fu, \fu'$ are not equivalent under $\sim$, we have $\mathfrak{v} \not \sim \fu'$, thus $10\fp' \not\lesssim \fu'$. This implies that there exists $q \in B_{\fu'}(\fcc(\fu'), 1) \setminus B_{\fp'}(\fcc(\fp'), 10)$.

    From $\fp' \le \fp$, $\scI(\fp') \subset \scI(\fp) \subset \scI(\fu')$ and \Cref{monotone-cube-metrics} it then follows that
    \begin{align*}
        &\quad d_{\fp'}(\fcc(\fp), \fcc(\fu'))\\
        &\ge -d_{\fp'}(\fcc(\fp), \fcc(\fp')) + d_{\fp'}(\fcc(\fp'), q) - d_{\fp'}(q, \fcc(\fu'))\\
        &\ge -d_{\fp'}(\fcc(\fp), \fcc(\fp')) + d_{\fp'}(\fcc(\fp'), q) - d_{\fu'}(q, \fcc(\fu'))\\
        &> -1 + 10 - 1 = 8\,.
    \end{align*}
    The lemma follows by combining the two displays with the fact that $95 a \ge 1$.
\end{proof}

\begin{lemma}[forest inner]
    \label{forest-inner}
    \uses{relation-geometry}
    \leanok
    \lean{forest_inner}
    For each $\fu\in \fU_3(k,n,j)$
    and each $\fp \in \mathfrak{T}_2(\fu)$
    we have
    \begin{equation}
        B(\pc(\fp), 8 D^{\ps(\fp)}) \subset \scI(\fu).
    \end{equation}
\end{lemma}

\begin{proof}
    \leanok
    Let $\fp \in \mathfrak{T}_2(\fu)$. Then $\fp \in \fC_4(k,n,j)$, hence there exists a chain
    $$
        \fp \le \fp_{Z(n+1)} \le \dotsb \le \fp_0
    $$
    of distinct tiles in $\fC_3(n,k,j)$. We pick such a chain and set $\fq = \fp_0$.
    Then we have from distinctness of the tiles in the chain that
    $\ps(\fp) \le \ps(\fq) - Z(n+1)$.
    By \eqref{eq-C3-def} there exists $\fu'' \in \fU_1(k,n,j)$ with $2\fq \lesssim \fu''$ and $\ps(\fq) < \ps(\fu'')$.
    Then we have in particular by \Cref{wiggle-order-1} that $10 \fp \lesssim \fu''$.
    Let $\fu' \sim \fu$ be such that $\fp \in \mathfrak{T}_1(\fu')$.
    By the definition of $\sim$, it follows that $\fu' \sim \fu''$.
    By transitivity of $\sim$, we have $\fu \sim \fu''$.
    By \Cref{relation-geometry}, we have $\ps(\fu'') = \ps(\fu)$, hence $\ps(\fq) < \ps(\fu)$ and $\ps(\fp) \le \ps(\fq) - Z(n+1) \le \ps(\fu) - Z(n+1) - 1$.

    Thus, there exists some cube $I \in \mathcal{D}$ with $s(I) = \ps(\fu) - Z(n+1) - 1$ and $I \subset \scI(\fu)$ and $\scI(\fp) \subset I$.
    Since $\fp \in \fC_5(k,n,j)$, we have that $I \notin \mathcal{L}(\fu)$, so $B(c(I), 8D^{s(I)}) \subset \scI(\fu)$.
    By the triangle inequality, \eqref{defineD} and $a \ge 4$, the same then holds for the subcube $\scI(\fp) \subset I$.
\end{proof}

\begin{lemma}[forest stacking]
    \label{forest-stacking}
    \leanok
    \lean{forest_stacking}
    It holds for $k\le n$ that
    \begin{equation}
        \sum_{\fu \in \fU_3(k,n,j)} \mathbf{1}_{\scI(\fu)} \le (4n+12)2^{n}\,.
    \end{equation}
\end{lemma}

\begin{proof}
    \leanok
    Suppose that a point $x$ is contained in more than $(4n + 12)2^n$ cubes $\scI(\fu)$ with $\fu \in \fU_3(k,n,j)$.
    Since $\fU_3(k,n,j) \subset \fC_1(k,n,j)$ for each such $\fu$,
    there exists $\mathfrak{m} \in \mathfrak{M}(k,n)$ such that $100\fu \lesssim \mathfrak{m}$.
    We fix such an $\mathfrak{m}(\fu) := \mathfrak{m}$ for each $\fu$, and claim that the map $\fu \mapsto\mathfrak{m}(\fu)$ is injective.
    Indeed, assume for $\fu\neq \fu'$ there is $\mathfrak{m} \in \mathfrak{M}(k,n)$ such that
    $100\fu \lesssim \mathfrak{m}$ and $100\fu' \lesssim \mathfrak{m}$. By \eqref{dyadicproperty},
    either $\scI(\fu) \subset \scI(\fu')$ or $\scI(\fu') \subset \scI(\fu)$.
    By \eqref{defunkj}, $B_{\fu}(\fcc(\fu),100) \cap B_{\fu'}(\fcc(\fu'), 100) = \emptyset$.
    This contradicts $\Omega(\mathfrak{m})$ being contained in both sets by \eqref{eq-freq-comp-ball}.
    Thus $x$ is contained in more than $(4n + 12)2^n$ cubes $\scI(\mathfrak{m})$, $\mathfrak{m} \in \mathfrak{M}(k,n)$.
    Consequently, we have by \eqref{eq-Aoverlap-def} that $x \in A(2n + 6, k,n) \subset G_2$.
    Let $\scI(\fu)$ be an inclusion minimal cube among the $\scI(\fu'), \fu' \in \fU_3(k,n,j)$ with $x \in \scI(\fu)$.
    By the dyadic property \eqref{dyadicproperty}, we have $\scI(\fu) \subset \scI(\fu')$ for all cubes $\scI(\fu')$ containing $x$. Thus
    $$
        \scI(\fu) \subset \{y \ : \ \sum_{\fu \in \fU_3(k,n,j)} \mathbf{1}_{\scI(\fu)}(y) > 1 + (4n+12)2^{n}\} \subset G_2\,.
    $$
    Thus $\mathfrak{T}_1(\fu) \cap \fC_6(k,n,j) = \emptyset$.
    This contradicts $\fu \in \fU_2(k,n,j)$.
\end{proof}
We now turn to the proof of \Cref{forest-union}.
\begin{proof}[Proof of \Cref{forest-union}]
    \proves{forest-union}\leanok

    We first fix $k,n, j$.
    By \eqref{definetp} and \eqref{defineep}, we have that
    $\mathbf{1}_{\scI(\fp)} T_{\fp}f(x) = T_{\fp}f(x)$ and hence $\mathbf{1}_{G \setminus G'} T_{\fp}f(x)= 0$ for all $\fp \in \fC_5(k,n,j) \setminus \fC_6(k,n,j)$.
    Thus it suffices to estimate the contribution of the sets $\fC_6(k,n,j)$. By \Cref{forest-stacking}, we can decompose $\fU_3(k,n,j)$ as a disjoint union of at most $4n + 12$ collections $\fU_4(k,n,j,l)$, $1 \le l \le 4n+12$, each satisfying
    $$
        \sum_{\fu \in \fU_4(k,n,j,l)} \mathbf{1}_{\scI(\fu)} \le 2^n\,.
    $$
    By Lemmas \ref{forest-geometry}, \ref{forest-convex}, \ref{forest-separation}, \ref{forest-inner} and \ref{C-dens1}, the pairs
    $$
        (\fU_4(k,n,j,l), \mathfrak{T}_2|_{\fU_4(k,n,j,l)})
    $$
    are $n$-forests for each $k,n,j,l$, and by \Cref{C6-forest}, we have
    $$
        \fC_6(k,n,j) = \bigcup_{l = 1}^{4n + 12} \bigcup_{\fu \in \fU_4(k,n,j,l)} \mathfrak{T}_2(\fu)\,.
    $$
    Since $\scI(\fp) \not\subset G_1$ for all $\fp \in \fC_6(k,n,j)$, we have $\fC_6(k,n,j) \cap \fP_{F,G} = \emptyset$ and hence
    $$
        \dens_2(\bigcup_{\fu \in \fU_4(k,n,j,l)} \mathfrak{T}_2(\fu)) \le 2^{2a + 5} \frac{\mu(F)}{\mu(G)}\,.
    $$
    Using the triangle inequality according to the splitting by $k,n,j$ and $l$ in \eqref{disclesssim1} and applying \Cref{forest-operator} to each term, we obtain the estimate
    $$
        \sum_{k \ge 0}\sum_{n \ge k} (2n+3)(4n+12) 2^{440a^3}2^{-(1-\frac{1}{q})n}(2^{2a+5} \frac{\mu(F)}{\mu(G)})^{\frac{1}{q} - \frac{1}{2}} \|f\|_2 \|\mathbf{1}_{G\setminus G'}\|_2
    $$
    for the left hand side of \eqref{disclesssim1}. Since $|f| \le \mathbf{1}_F$, we have $\|f\|_2 \le \mu(F)^{1/2}$, and we have $\|\mathbf{1}_{G\setminus G'}\|_2 \le \mu(G)^{1/2}$.
    We get a bound
    $$
        2^{440a^3} \mu(F)^{\frac{1}{q}} \mu(G)^{1
        -\frac{1}{q}} 2^{2a+5}\sum_{k \ge 0}\sum_{n \ge k}(2n+3)(4n+12) 2^{-(1-\frac{1}{q})n}\,.
    $$ 
    Interchanging the order of summation, the last factor equals
    $$
        2^{2a+5} \sum_{n \ge 0} (2n+3)(4n+12) (n+1) 2^{-\frac{q-1}{q}n}\,.
    $$
    Up to an explicit constant, the sum is bounded by $\sum n^3 2^{-\frac{q-1}{q}n}$, which is at most
    some constant times $1/(q-1)^4$ by comparing to an integral. Since $a \ge 4$, this is overall bounded by $2^{a^3}/(q-1)^4$,
    which completes the proof of the lemma.
\end{proof}

\subsection{Proof of the Forest Complement Lemma}
\label{subsecantichain}

Define $\fP_{G \setminus G'}$ to be the set of all $\fp \in \fP$ such that $\mu(\scI(\fp) \cap (G \setminus  G')) > 0$.
\begin{lemma}[antichain decomposition]
    \label{antichain-decomposition}
    \leanok
    \lean{antichain_decomposition}
    We have that
    \begin{align}
        \label{eq-fp'-decomposition}
        &\quad \fP_2 \cap \fP_{G \setminus G'}\\
        &= \bigcup_{k \ge 0} \bigcup_{n \ge k} \fL_0(k,n) \cap \fP_{G \setminus G'} \\
        &\quad\cup \bigcup_{k \ge 0} \bigcup_{n \ge k}\bigcup_{0 \le j \le 2n+3} \fL_2(k,n,j) \cap \fP_{G \setminus G'}\\
        &\quad\cup \bigcup_{k \ge 0} \bigcup_{n \ge k}\bigcup_{0 \le j \le 2n+3} \bigcup_{0 \le l \le Z(n+1)} \fL_1(k,n,j,l) \cap \fP_{G \setminus G'}\\
        &\quad\cup \bigcup_{k \ge 0} \bigcup_{n \ge k}\bigcup_{0 \le j \le 2n+3} \bigcup_{0 \le l \le Z(n+1)} \fL_3(k,n,j,l)\cap \fP_{G \setminus G'}\,.
    \end{align}
\end{lemma}

\begin{proof}
    \leanok
    Let $\fp \in \fP_2 \cap \fP_{G \setminus G'}$. Clearly, for every cube $J = \scI(\fp)$ with $\fp \in \fP_{G \setminus G'}$ there exists some $k \ge 0$ such that \eqref{muhj1} holds, and for no cube $J \in \mathcal{D}$ and no $k < 0$ does \eqref{muhj2} hold. Thus $\fp \in \fP(k)$ for some $k \ge 0$.

    Next, since $E_2(\lambda, \fp') \subset \scI(\fp')\cap G$ for every $\lambda \ge 2$ and every tile $\fp' \in \fP(k)$ with $\lambda\fp \lesssim \lambda \fp'$, it follows from \eqref{muhj2} that $\mu(E_2(\lambda, \fp')) \le 2^{-k} \mu(\scI(\fp'))$ for every such $\fp'$, so $\dens_k'(\{\fp\}) \le 2^{-k}$. Combining this with $a \ge 0$, it follows from \eqref{def-cnk} that there exists $n\ge k$ with $\fp \in \fC(k,n)$.

    Since $\fp \in \fP_{G \setminus G'}$, we have in particular $\scI(\fp) \not \subset A(2n + 6, k, n)$, so there exist at most $1 + (4n + 12)2^n < 2^{2n+4}$ tiles $\mathfrak{m} \in \mathfrak{M}(k,n)$ with $\fp \le \mathfrak{m}$. It follows that $\fp \in \fL_0(k,n)$ or $\fp \in \fC_1(k,n,j)$ for some $1 \le j \le 2n + 3$. In the former case we are done, in the latter case the equality to be shown follows from the definitions of the collections $\fC_i$ and $\fL_i$.
\end{proof}

\begin{lemma}[L0 antichain]
\label{L0-antichain}
\uses{monotone-cube-metrics}
\leanok
\lean{iUnion_L0', pairwiseDisjoint_L0', antichain_L0'}
    We have that
    $$
        \fL_0(k,n) = \dot{\bigcup_{0 \le l < n}} \fL_0(k,n,l)\,,
    $$
    where each $\fL_0(k,n,l)$ is an antichain.
\end{lemma}

\begin{proof}
    \leanok
    It suffices to show that $\fL_0(k,n)$ contains no chain of length $n + 1$. Suppose that we had such a chain $\fp_0 \le \fp_1 \le \dotsb \le \fp_{n}$ with $\fp_i \ne \fp_{i+1}$ for $i =0, \dotsc, n-1$. By \eqref{def-cnk}, we have that $\dens_k'(\{\fp_n\}) > 2^{-n}$. Thus, by \eqref{eq-densdef}, there exists $\fp' \in \fP(k)$ and $\lambda \ge 2$ with $\lambda \fp_n \le \lambda \fp'$ and
    \begin{equation}
        \label{eq-p'}
        \frac{\mu(E_2(\lambda, \fp'))}{\mu(\scI(\fp'))} > \lambda^{a} 2^{4a} 2^{-n}\,.
    \end{equation}
    Let $\mathfrak{O}$ be the set of all $\fp'' \in \fP(k)$ such that we have $ \scI(\fp'') = \scI(\fp')$ and $B_{\fp'}(\fcc(\fp'), \lambda) \cap \Omega(\fp'') \neq \emptyset$.
    We now show that
    \begin{equation}
        \label{eq-O-bound}
        |\mathfrak{O}| \le 2^{4a}\lambda^a\,.
    \end{equation}
    The balls $B_{\fp'}(\fcc(\fp''), 0.2)$, $\fp'' \in \mathfrak{O}$ are disjoint by \eqref{eq-freq-comp-ball},
    and by the triangle inequality contained in $B_{\fp'}(\fcc(\fp'), \lambda+1.2)$.
    By assumption \eqref{thirddb} on $\Theta$, this ball can be covered with
    $$
        2^{a\lceil \log_2(\lambda+1.2) + \log_2(5)\rceil} \le 2^{a(\log_2(\lambda) + 4)} = 2^{4a}\lambda^a
    $$
    many $d_{\fp'}$-balls of radius $0.2$. Here we have used that for $\lambda \ge 2$
    $$
        \lceil \log_2(\lambda + 1.2)  + \log_2(5) \rceil \le 1+ \log_2(1.6  \lambda) + \log_2(5) = 4 + \log_2(\lambda)\,.
    $$
    By the triangle inequality, each such ball contains at most one $\fcc(\fp'')$, and each $\fcc(\fp'')$ is contained in one of the balls. Thus we get \eqref{eq-O-bound}.

    By \eqref{definee1} and \eqref{definee2} we have $E_2(\lambda, \fp') \subset \bigcup_{\fp'' \in \mathfrak{O}} E_1(\fp'')$, thus
    $$
        2^{4a}\lambda^a 2^{-n} < \sum_{\fp'' \in \mathfrak{O}} \frac{\mu(E_1(\fp''))}{\mu(\scI(\fp''))}\,.
    $$
    Hence there exists a tile $\fp'' \in \mathfrak{O}$ with
    \begin{equation*}
        \mu(E_1(\fp'')) \ge 2^{-n} \mu(\scI(\fp'))\,.
    \end{equation*}
    By the definition \eqref{mnkmax} of $\mathfrak{M}(k,n)$, there exists a tile $\mathfrak{m} \in \mathfrak{M}(k,n)$ with $\fp' \leq \mathfrak{m}$. From \eqref{eq-p'}, the inclusion $E_2(\lambda, \fp') \subset \scI(\fp')$ and $a\ge 1$ we obtain
    $$
        2^n \geq 2^{4a} \lambda^{a} \geq \lambda\,.
    $$
    From the triangle inequality, \Cref{monotone-cube-metrics} and $a \ge 1$, we now obtain for all $\mfa \in B_{\mathfrak{m}}(\fcc(\mathfrak{m}), 1)$ that
    \begin{align*}
        &\quad d_{\fp_0}(\fcc(\fp_0), \mfa)\\
        &\leq d_{\fp_0}(\fcc(\fp_0), \fcc(\fp_{n})) + d_{\fp_0}(\fcc(\fp_{n}), \fcc(\fp')) + d_{\fp_0}(\fcc(\fp'), \fcc(\fp''))\\
        &\quad+ d_{\fp_0}(\fcc(\fp''), \fcc(\mathfrak{m})) +
        d_{\fp_0}(\fcc(\mathfrak{m}), \mfa)\\
        &\leq 1 + 2^{-95an} (d_{\fp_{n}}(\fcc(\fp_n), \fcc(\fp')) + d_{\fp'}(\fcc(\fp'), \fcc(\fp''))\\
        &\quad+ d_{\fp''}(\fcc(\fp''), \fcc(\mathfrak{m})) +
        d_{\mathfrak{m}}(\fcc(\mathfrak{m}), \mfa))\\
        &\leq 1 + 2^{-95an}(\lambda + (\lambda + 1) + 1 + 1) \leq 100\,.
    \end{align*}
    Thus, by \eqref{straightorder}, $100\fp_0 \lesssim \mathfrak{m}$, a contradiction to $\fp_0 \notin \fC(k,n)$.
\end{proof}

\begin{lemma}[L2 antichain]
\label{L2-antichain}
\uses{monotone-cube-metrics}
\leanok
\lean{antichain_L2}
    Each of the sets $\fL_2(k,n,j)$ is an antichain.
\end{lemma}

\begin{proof}
    \leanok
    Suppose that there are $\fp_0, \fp_1 \in \fL_2(k,n,j)$ with $\fp_0 \ne \fp_1$ and $\fp_0 \le \fp_1$. By \Cref{wiggle-order-1} and \Cref{wiggle-order-2}, it follows that $2\fp_0 \lesssim 200\fp_1$. Since $\fL_2(k,n,j)$ is finite, there exists a maximal $l \ge 1$ such that there exists a chain $2\fp_0 \lesssim 200 \fp_1 \lesssim \dotsb \lesssim 200 \fp_l$ with all $\fp_i$ in $\fC_1(k,n,j)$ and $\fp_i \ne \fp_{i+1}$ for $i = 0, \dotsc, l-1$.
    If we have $\fp_l \in \fU_1(k,n,j)$, then it follows from $2\fp_0 \lesssim 200 \fp_l \lesssim \fp_l$ and \eqref{eq-L2-def} that $\fp_0 \not\in \fL_2(k,n,j)$, a contradiction. Thus, by the definition \eqref{defunkj} of $\fU_1(k,n,j)$, there exists $\fp_{l+1} \in \fC_1(k,n,j)$ with $\scI(\fp_l) \subsetneq \scI(\fp_{l+1}) $ and $\mfa \in B_{\fp_l}(\fcc(\fp_l), 100) \cap B_{\fp_{l+1}}(\fcc(\fp_{l+1}), 100)$. Using the triangle inequality and \Cref{monotone-cube-metrics}, one deduces that $200 \fp_l \lesssim 200\fp_{l+1}$. This contradicts maximality of $l$.
\end{proof}

\begin{lemma}[L1 L3 antichain]
\label{L1-L3-antichain}
\leanok
\lean{antichain_L1, antichain_L3}
    Each of the sets $\fL_1(k,n,j,l)$ and $\fL_3(k,n,j,l)$ is an antichain.
\end{lemma}

\begin{proof}
    \leanok
    By its definition \eqref{eq-L1-def}, each set $\fL_1(k,n,j,l)$ is a set of minimal elements in some set of tiles with respect to $\le$. If there were distinct $\fp, \fq \in \fL_1(k,n,j,l)$ with $\fp \le \fq$, then $\fq$ would not be minimal. Hence such $\fp, \fq$ do not exist. Similarly, by \eqref{eq-L3-def}, each set $\fL_3(k,n,j,l)$ is a set of maximal elements in some set of tiles with respect to $\le$. If there were distinct $\fp, \fq \in \fL_3(k,n,j,l)$ with $\fp \le \fq$, then $\fp$ would not be maximal.
\end{proof}

We now turn to the proof of \Cref{forest-complement}.
\begin{proof}[Proof of \Cref{forest-complement}]
    \proves{forest-complement}\leanok
    If $\fp \not\in \fP_{G \setminus G'}$, then $\mu(\scI(\fp) \cap (G \setminus G')) = 0$. By \eqref{definetp} and \eqref{definee1}, it follows that
    $\mathbf{1}_{G \setminus G'} T_{\fp}f(x) = 0$ for almost every $x$. We thus have, almost everywhere,
    $$
        \mathbf{1}_{G\setminus G'} \sum_{\fp \in \fP_2} T_{\fp}f(x) = \mathbf{1}_{G\setminus G'} \sum_{\fp \in \fP_2 \cap \fP_{G \setminus G'}} T_{\fp}f(x)\,.
    $$
    Let $\fL(k,n)$ denote any of the terms $\fL_i(k,n,j,l) \cap \fP_{G \setminus G'}$ on the right hand side of \eqref{eq-fp'-decomposition}, where the indices $j, l$ may be void. Then $\fL(k,n)$ is an antichain, by Lemmas \ref{L0-antichain},\ref{L2-antichain}, \ref{L1-L3-antichain}. Further, we have
    \begin{equation*}
    \dens_1(\fL(k,n)) \le 2^{4a+1 - n}
    \end{equation*}
    by \Cref{C-dens1}, and we have
    \begin{equation*}
     \dens_2(\fL(k,n)) \le 2^{2a+5} \frac{\mu(F)}{\mu(G)},
     \end{equation*}
     since
     \begin{equation*}
     \fL(k,n) \cap \fP_{F,G} \subset \fP_{G \setminus G'} \cap \fP_{F, G} = \emptyset.
     \end{equation*}

    Applying now the triangle inequality according to the decomposition coming from \Cref{antichain-decomposition}, and then applying \Cref{antichain-operator} to each term, we obtain the estimate
    \begin{multline*}
        \le \sum_{k \ge 0} \sum_{n \ge k} (n + (2n+4) + 2(2n+4) (1+Z(n+1))) \\
        \times 2^{117a^3}(q-1)^{-1} (2^{4a+1-n})^{\frac{q-1}{8a^4}} (2^{2a+5} \frac{\mu(F)}{\mu(G)})^{\frac{1}{q} - \frac{1}{2}} \|f\|_2\|\mathbf{1}_{G\setminus G'}\|_2\,.
    \end{multline*}
    Because $|f| \le \mathbf{1}_F$, we have $\|f\|_2 \le \mu(F)^{1/2}$, and we have $\|\mathbf{1}_{G\setminus G'}\|_2 \le \mu(G)^{1/2}$. Using this and \eqref{defineZ}, we bound
    $$
        \le 2^{118a^3} (q - 1)^{-1} \mu(F)^{\frac{1}{q}} \mu(G)^{\frac{1}{q'}} \sum_{k \ge 0} \sum_{n \ge k} n^2 2^{-n\frac{q-1}{8a^4}}\,.
    $$
    The last sum equals, by changing the order of summation,
    $$
        \sum_{n \ge 0} n^2(n+1) 2^{-n\frac{q-1}{8a^4}} \le \frac{2^{2a^3}}{(q-1)^4}\,.
    $$
    This completes the proof.
\end{proof}

\section{Proof of the Antichain Operator Proposition}

\label{antichainboundary}

Let an antichain $\mathfrak{A}$
and functions $f$, $g$ as in \Cref{antichain-operator} be given.
We prove \eqref{eq-antiprop}
in \Cref{sec-TT*-T*T}
as the geometric mean of two inequalities,
each involving one of the two densities.
One of these two inequalities will need a careful estimate formulated in
\Cref{tile-correlation} of
the $TT^*$ correlation between two tile operators.
\Cref{tile-correlation} will be proven in
\Cref{sec-tile-operator}.

The summation of the contributions of these individual correlations will require a
geometric \Cref{antichain-tile-count} counting the relevant tile pairs.
\Cref{antichain-tile-count} will be proven in Subsection
\ref{subsec-geolem}.

\subsection{The density arguments}\label{sec-TT*-T*T}

We begin with the following crucial disjointedness property of the sets $E(\fp)$ with $\fp \in \mathfrak{A}$.
\begin{lemma}[tile disjointness]
\label{tile-disjointness}
\leanok
\lean{tile_disjointness}
Let $\fp,\fp'\in \mathfrak{A}$.
If there exists an $x\in X$ with $x\in E(\fp)\cap E(\fp')$,
then $\fp= \fp'$.
\end{lemma}
\begin{proof}
\leanok
Let $\fp,\fp'$ and $x$ be given.
Assume without loss of generality that $\ps(\fp)\le \ps(\fp')$.
As we have $x\in E(\fp)\subset \scI(\fp)$ and $x\in E(\fp')\subset \scI(\fp')$ by Definition \eqref{defineep}, we conclude
for $i=1,2$ that
$\tQ(x)\in\fc(\fp)$ and $\tQ(x)\in\fc(\fp')$. By \eqref{eq-freq-dyadic} we have $\fc(\fp')\subset \fc(\fp)$. By Definition
\eqref{straightorder}, we conclude $\fp\le \fp'$. As $\mathfrak{A}$ is an antichain, we conclude $\fp=\fp'$.
This proves the lemma.
\end{proof}

Let $\mathcal{B}$ be the collection of balls
\begin{equation}
    B(\pc(\fp), 8D^{\ps(\fp)})
\end{equation}
with $\fp\in \mathfrak{A}$ and recall the definition of
$M_{\mathcal{B}}$ from Definition \ref{def-hlm}.

\begin{lemma}[maximal bound antichain]
    \label{maximal-bound-antichain}
    \uses{tile-disjointness}
    \leanok
    \lean{maximal_bound_antichain}
    Let $x\in X$.
    Then
    \begin{equation}\label{hlmbound}
    | \sum_{\fp \in \mathfrak{A}}T_{\fp} f(x)|\le 2^{102 a^3} M_{\mathcal{B}} f (x) \, .
    \end{equation}
\end{lemma}

\begin{proof}
\leanok
Fix $x\in X$. By \Cref{tile-disjointness}, there is at most one $\fp \in \mathfrak{A}$
such that
 $T_{\fp} f(x)$ is not zero.
 If there is no such $\fp$, the estimate \eqref{hlmbound} follows.

 Assume there is such a $\fp$.
 By definition of $T_{\fp}$ we have $x\in E(\fp)\subset \scI(\fp)$ and by the squeezing property \eqref{eq-vol-sp-cube}
\begin{equation}\label{eqtttt0}
    \rho(x, \pc(\fp))\le 4D^{\ps(\fp)}\, .
\end{equation}

Let $y\in X$ with $K_{\ps(\fp)}(x,y)\neq 0$. By Definition \eqref{defks} of $K_{\ps(\fp)}$
we have
\begin{equation}\label{supp-Ks1}
   \frac{1}{4} D^{\ps(\fp)-1}
   \leq \rho(x,y) \leq \frac{1}{2} D^{\ps(\fp)}\, .
\end{equation}
The triangle inequality with \eqref{eqtttt0} and \eqref{supp-Ks1} implies
\begin{equation}
    \rho(\pc(\fp),y) < 8D^{\ps(\fp)}\, .
\end{equation}
Using the kernel bound \eqref{eqkernel-size} and the lower bound in \eqref{supp-Ks}
we obtain
\begin{equation}
|K_{\ps(\fp)}(x,y)|\le \frac{2^{a^3}}{\mu(B(x,\frac 14 D^{{\ps(\fp)}-1}))}\, .
\end{equation}
Using $D=2^{100a^2}$
and the doubling property \eqref{doublingx} $5 +100a^2$ times estimates
the last display by
\begin{equation}
\le \frac{2^{5a+101a^3}}{\mu(B(x, 8D^{\ps(\fp)}))}\,
\end{equation}
which, thanks to the closeness of the points $x$ and $\pc(\fp)$ shown in \eqref{eqtttt0}, is in turn bounded by
\begin{equation}
\le \frac{2^{6a+101a^3}}{\mu(B(\pc(\fp), 8D^{\ps(\fp)}))}\, .
\end{equation}
 Using that {$|e(\mfa)|$} is bounded by $1$
for every $\mfa\in \Mf$, we estimate with the triangle inequality and the above information
 \begin{equation}
  | T_{\fp} f(x)|
    \le \frac{2^{6a+101 a^3}}{\mu(B(\pc(\fp), 8D^{\ps(\fp)}))} \int _{\mu(B(\pc(\fp), 8D^{\ps(\fp)}))} |f(y)|\, dy
  \end{equation}
This together with $a\ge 4$ proves the Lemma.
\end{proof}

Set
\begin{equation}
    \tilde{q}=\frac {2q}{1+q}\,.
\end{equation}
Since $1< q\le 2$, we have $1<\tilde{q}<q\le 2$.
\begin{lemma}[dens2 antichain]
\label{dens2-antichain}
\uses{Hardy-Littlewood,maximal-bound-antichain}
\leanok
\lean{dens2_antichain}
We have that
\begin{equation}\label{eqttt9}
  \left|\int \overline{g(x)} \sum_{\fp \in \mathfrak{A}} T_{\fp} f(x)\, d\mu(x)\right|\le
  2^{103a^3}({q}-1)^{-1} \dens_2(\mathfrak{A})^{\frac 1{\tilde{q}}-\frac 12} \|f\|_2\|g\|_2\, .
\end{equation}
\end{lemma}
\begin{proof}
\leanok
We have $f=\mathbf{1}_Ff$. Using H\"older's inequality, we obtain for
each $x\in B'$ and each $B'\in \mathcal{B}$ using $1<\tilde{q}\le 2$
\begin{equation}
    \frac 1{\mu(B')}\int_{B'} |f(y)|\, d\mu(y)
\end{equation}
\begin{equation}
    \le
    \left(\frac 1{\mu(B')}\int_{B'} |f(y)|^{\frac {2{\tilde{q}}}{3\tilde{q}-2}}\, d\mu(y)\right)^{\frac 32-\frac 1{\tilde{q}}}
    \left(\frac 1{\mu(B')}\int_{B'} \mathbf{1}_F(y)\, d\mu(y)\right)^{\frac 1{\tilde{q}}-\frac 12}
\end{equation}
\begin{equation}
    \le \left(M_{\mathcal{B}} (|f|^{\frac {2{\tilde{q}}}{3{\tilde{q}}-2}})(x)\right)^{\frac 32-\frac 1{\tilde{q}}}
\dens_2(\mathfrak{A})^{\frac 1{\tilde{q}}-\frac 12}\, .
\end{equation}
Taking the maximum over all $B'$ containing $x$, we obtain
\begin{equation} \label{eqttt1}
    M_{\mathcal{B}}|f|\le
    M_{\mathcal{B},\frac {2{\tilde{q}}}{3{\tilde{q}}-2} } |f|
    \dens_2(\mathfrak{A})^{\frac 1{\tilde{q}}-\frac 12}\, .
\end{equation}
We have with \Cref{Hardy-Littlewood}
\begin{equation}
\left\|M_{\mathcal{B}, \frac {2{\tilde{q}}}{3{\tilde{q}}-2}} f\right\|_2\le 2^{2a}(3\tilde{q}-2)(2\tilde{q}-2)^{-1}\|f\|_2\, .
\end{equation}
Using $1<\tilde{q}\le 2$ estimates the last display by
\begin{equation}\label{eqttt2}
 2^{2a+2} (\tilde{q}-1)^{-1} \|f\|_2\, .
\end{equation}
We obtain with Cauchy-Schwarz
and then \Cref{maximal-bound-antichain}
 \begin{equation}
     |\int \overline{g(x)} \sum_{\fp \in \mathfrak{A}} T_{\fp} f(x)\, d\mu(x)|
\end{equation}
 \begin{equation}
     \le \|g\|_2 \Big\| \sum_{\fp \in \mathfrak{A}} T_{\fp} f \Big\|_2
\end{equation}
 \begin{equation}
     \le 2^{102a^3}\|g\|_2 \| M_{\mathcal{B}}f \|_2
\end{equation}
With \eqref{eqttt1} and
\eqref{eqttt2} we can estimate the last display by
\begin{equation}
    \le 2^{102a^3+2a+2}(\tilde{q}-1)^{-1} \|g\|_2 \|f\|_2\dens_2(\mathfrak{A})^{\frac 1{\tilde{q}}-\frac 12}
\end{equation}
Using $a\ge 4$ and
$(\tilde q - 1)^{-1} = (q+1)/(q-1) \le 3(q-1)^{-1}$
proves the lemma.
\end{proof}

\begin{lemma}[dens1 antichain]\label{dens1-antichain}
    \leanok
    \lean{dens1_antichain}
    \uses{Hardy-Littlewood, tile-correlation,antichain-tile-count}
Set $p:=4a^4$. We have
    \begin{equation}\label{eqttt3}
  \left|\int \overline{g(x)} \sum_{\fp \in \mathfrak{A}} T_{\fp} f(x)\, d\mu(x)\right|\le
   2^{117a^3}\dens_1(\mathfrak{A})^{\frac 1{2p}} \|f\|_2\|g\|_2\,.
\end{equation}
\end{lemma}

\begin{proof}
\leanok
We write for the expression inside the absolute values on the left-hand side of \eqref{eqttt3}
\begin{equation}
  \sum_{\fp \in \mathfrak{A}}\iint \overline{g(x)} \mathbf{1}_{E(\fp)}(x)
  {K_{\ps(\fp)}(x,y)}e(\tQ(x)(y) -
   \tQ(x)(x))
   f(y)\, d\mu(y)\,d\mu(x)
\end{equation}
\begin{equation}
  =\int \sum_{\fp \in \mathfrak{A}} \overline{T_{\fp} ^*g(y)} f(y)\, d\mu(y)
\end{equation}
with the adjoint operator
\begin{equation}\label{eq-tstarwritten}
    T_{\fp}^*g(y)=\int_{E(\fp)} \overline{K_{\ps(\fp)}(x,y)}e(-\tQ(x)(y)+
    \tQ(x)(x))g(x)\, d\mu(x)\, .
\end{equation}
We have by expanding the square
\begin{equation}
    \int \Big|\sum_{\fp\in \mathfrak{A}}T^*_{\fp}g(y)\Big|^2\, d\mu(y)=
    \int \left(\sum_{\fp\in \mathfrak{A}} T^*_{\fp}g(y)\right)
    \left(\sum_{\fp'\in \mathfrak{A}}\overline{T^*_{\fp'}g(y)}\right)\, d\mu(y)
\end{equation}
\begin{equation}\label{eqtts1}
    \le \sum_{\fp\in \mathfrak{A}} \sum_{\fp'\in \mathfrak{A}}
    \Big|\int T^*_{\fp}g(y)\overline{T^*_{\fp'}g(y)}\, d\mu(y)\Big|\,.
\end{equation}
We split the sum into the terms with $\ps(\fp')\le \ps(\fp)$
and $\ps(\fp)< \ps(\fp')$. Using the symmetry of each summand,
we may switch $\fp$ and $\fp'$ in the second sum. Using further positivity
of each summand to replace the condition $\ps(\fp')< \ps(\fp)$
by $\ps(\fp')\le \ps(\fp)$ in the second sum, we estimate \eqref{eqtts1} by
\begin{equation}\label{eqtts2}
    \le2 \sum_{\fp\in \mathfrak{A}} \sum_{\fp'\in \mathfrak{A}: \ps(\fp')\le \ps(\fp)}
    \Big|\int T^*_{\fp}g(y)\overline{T^*_{\fp'}g(y)}\, d\mu(y)\Big|\,.
\end{equation}

Define for $\fp\in \fP$
\begin{equation}
    B(\fp):=B(\pc(\fp), 14D^{\ps(\fp)})
\end{equation}
and define
\begin{equation}
    \label{eq-Dp-definition}
    \mathfrak{A}(\fp):=\{\fp'\in\mathfrak{A}: \ps(\fp')\leq \ps(\fp) \land \scI(\fp') \subset B(\fp)\}.
\end{equation}
Note that by the squeezing property \eqref{eq-vol-sp-cube}
and the doubling property \eqref{doublingx} applied
$6$ times we have
\begin{equation}\label{eqttt4}
    \mu(B(\fp))\le 2^{6a} \mu(\scI(\fp))\, .
\end{equation}
Using \Cref{tile-correlation} and \eqref{eqttt4}, we estimate \eqref{eqtts2} by
\begin{equation}\label{eqtts3}
     \le 2^{232a^3+6a+1} \sum_{\fp\in \mathfrak{A}}
    \int_{E(\fp)}|g|(y) h(\fp)\, d\mu(y)
\end{equation}
with $h(\fp)$ defined as
\begin{equation}\label{def-hp}
    \frac 1{\mu(B(\fp))}\int \sum_{\fp'\in \mathfrak{A}(\fp)}
    {(1+d_{\fp'}(\fcc(\fp'), \fcc(\fp))^{-1/(2a^2+a^3)}}(\mathbf{1}_{E(\fp')}|g|)(y')\, d\mu(y')\,.
\end{equation}

Note that $p>4$ since $a\ge 4$. Let $p'$ be the dual exponent of $p$, satisfying $1/p+1/p'=1$.
We estimate $h(\fp)$ as defined in \eqref{def-hp} with H\"older using $|g|\le \mathbf{1}_G$ and $E(\fp')\subset B(\fp)$ by

\begin{equation}
    \frac{\|g\mathbf{1}_{B(\fp)}\|_{p'}}{\mu(B(\fp))}
    \Big\|\sum_{\fp\in\mathfrak{A}(\fp)}(1+d_{\fp}(\fcc(\fp), \fcc(\fp'))^{-1/(2a^2+a^3)}\mathbf{1}_{E(\fp)}\mathbf{1}_G\Big\|_{p}\, .
\end{equation}
Then we apply \Cref{antichain-tile-count} to estimate this by
\begin{equation}\label{eqttt5}
    \le 2^{5a}
    \frac{\|g\mathbf{1}_{B(\fp)}\|_{p'}}{\mu(B(\fp))}
    \dens_1(\mathfrak{A})^{\frac 1p}\mu(B(\fp))^{\frac 1p}\,.
\end{equation}
Let $\mathcal{B}'$ be the collection of all balls
$B(\fp)$ with $\fp\in \mathfrak{A}$. Then
for each $\fp\in \mathfrak{A}$ and $x\in B(\fp)$ we have by
definition \eqref{def-hlm} of $M_{\mathcal{B}',p'}$
\begin{equation}
    \|g\mathbf{1}_{B(\fp)}\|_{p'}\le
    \mu(B(\fp))^{\frac 1{p'}} M_{\mathcal{B}',p'}g(x) \, .
\end{equation}
Hence we can estimate \eqref{eqttt5} by
\begin{equation}
\label{eqttt5b}
    \le
    2^{5a}
    (M_{\mathcal{B}', p'}g(x))
   \dens_1(\mathfrak{A})^{\frac 1p}\, .
\end{equation}
With this estimate of $h(\fp)$, using $E(\fp)\subset B(\fp)$ by construction of $B(\fp)$, we estimate \eqref{eqtts3} by
 \begin{equation}\label{eqtts4}
 \le 2^{232a^3+11a + 1} { \dens_1(\mathfrak{A})^{\frac 1p}}\sum_{\fp\in \mathfrak{A}}
 \int_{E(\fp)}|g|(y)M_{\mathcal{B}', p'}g(y) \, dy\,.
         \end{equation}
Using \Cref{tile-disjointness}, the last display is observed to be
\begin{equation}\label{eqtts4a}
\le 2^{232a^3+11a + 1}
 {\dens_1(\mathfrak{A})^{\frac 1p}} \int |g|(y)(M_{\mathcal{B}', p'}g)(y) \, dy\,.
         \end{equation}
Applying Cauchy-Schwarz and using \Cref{Hardy-Littlewood} and $1<p'<\frac 32$ estimates the last display by
\begin{equation}
    \le 2^{232a^3+11a + 1} \dens_1(\mathfrak{A})^{\frac 1p}
    \|g\|_2 \|M_{\mathcal{B}', p'} g\|_2
\end{equation}
\begin{equation}
    \le 2^{232a^3+12a+3} \dens_1(\mathfrak{A})^{\frac 1p}\|g\|_2 ^2\,.
\end{equation}
Now \Cref{dens1-antichain} follows by applying Cauchy-Schwarz on the left-hand side and using $a\ge 4$.
\end{proof}
The following basic $TT^*$ estimate will be proved in \Cref{sec-tile-operator}.
\begin{lemma}[tile correlation]
    \label{tile-correlation}
    \leanok
    \lean{Tile.correlation_le, Tile.correlation_zero_of_ne_subset}
    \uses{Holder-van-der-Corput,correlation-kernel-bound,tile-uncertainty,tile-range-support}
    Let $\fp, \fp'\in \fP$ with
    $\ps({\fp'})\leq \ps({\fp})$.
    Then
    \begin{equation}
        \label{eq-basic-TT*-est}
        \left|\int T^*_{\fp'}g\overline{T^*_{\fp}g}\right|
    \end{equation}
    \begin{equation}
        \le 2^{232a^3}\frac{(1+d_{\fp'}(\fcc(\fp'), \fcc(\fp))^{-1/(2a^2+a^3)}}{\mu(\scI(\fp))}\int_{E(\fp')}|g|\int_{E(\fp)}|g|\,.
    \end{equation}
    Moreover, the term \eqref{eq-basic-TT*-est} vanishes unless
    \begin{equation}
        \scI(\fp') \subset B(\pc(\fp), 14D^{\ps(\fp)})\, .
    \end{equation}
\end{lemma}

The following lemma will be proved in \Cref{subsec-geolem}.
\begin{lemma}[antichain tile count]
    \label{antichain-tile-count}
\leanok
\lean{Antichain.tile_count}
    \uses{global-antichain-density}
    Set $p:=4a^4$. For every $\mfa\in\Mf$ and every antichain $\mathfrak{A}$ we have
    \begin{equation}
        \label{eq-antichain-Lp}
        \Big\|\sum_{\fp\in\mathfrak{A}}(1+d_{\fp}(\fcc(\fp), \mfa))^{-1/(2a^2+a^3)}\mathbf{1}_{E(\fp)}\mathbf{1}_G\Big\|_{p}
    \end{equation}
    \begin{equation}
        \le
        2^{5a}\dens_1(\mathfrak{A})^{\frac 1p}\mu\left(\cup_{\fp\in\mathfrak{A}}I_{\fp}\right)^{\frac 1p}\, .
    \end{equation}
\end{lemma}

From these lemmas it is easy to prove \Cref{antichain-operator}.
\begin{proof}[Proof of \Cref{antichain-operator}]
\proves{antichain-operator}
\leanok
We have
\begin{equation}
    \left(\frac 1{\tilde{q}} -\frac 12\right) (2-q)= \frac 1q -\frac 12\,.
\end{equation}
Multiplying the $(2-q)$-th power of \eqref{eqttt9} and the $(q-1)$-th power of \eqref{eqttt3}
and estimating gives after simplification of some factors
\begin{equation}\label{eqttt8}
    \Big|\int \overline{g(x)} \sum_{\fp \in \mathfrak{A}} T_{\fp} f(x)\, d\mu(x)\Big|
\end{equation}
 \begin{equation}
    \le 2^{117a^3}(q-1)^{-1} \dens_1(\mathfrak{A})^{\frac {q-1}{2p}}\dens_2(\mathfrak{A})^{\frac 1{q}-\frac 12} \|f\|_2\|g\|_2\, .
\end{equation}
With the definition of $p$, this implies
\Cref{antichain-operator}.
\end{proof}

\subsection{Proof of the Tile Correlation Lemma}\label{sec-tile-operator}

The next lemma prepares an application of
\Cref{Holder-van-der-Corput}.
\begin{lemma}[correlation kernel bound]\label{correlation-kernel-bound}
\leanok
\lean{Tile.correlation, Tile.mem_ball_of_correlation_ne_zero, Tile.correlation_kernel_bound}
Let $-S\le s_1\le s_2\le S$ and let $x_1,x_2\in X$.
Define \begin{equation}
 \varphi(y) := \overline{K_{s_1}(x_1, y)}
 K_{s_2}(x_2, y) \, .
\end{equation}
If $\varphi(y)\neq 0$, then
\begin{equation}\label{eqt10}
    y\in B(x_1, D^{s_1})\, .
\end{equation}
Moreover, we have with $\tau = 1/a$
\begin{equation}\label{eqt11}
  \|\varphi\|_{C^\tau(B(x_1, 2 D^{s_1}))}\le
\frac{2^{231 a^3}}{\mu(B(x_1, D^{s_1}))\mu(B(x_2, D^{s_2}))}
      \, .
\end{equation}

\end{lemma}
\begin{proof}
\leanok

If $\varphi(y)$ is not zero, then $K_{s_1}(x_1, y)$ is not zero and thus
\eqref{supp-Ks} gives \eqref{eqt10}.

We next have for $y$ with \eqref{eq-Ks-size}
\begin{equation}\label{suppart}
    |\varphi(y)|\le
    \frac{2^{204 a^3}}{\mu(B(x_1, D^{s_1}))\mu(B(x_2, D^{s_2}))}
\end{equation}
and for $y'\neq y$ additionally with \eqref{eq-Ks-smooth}
\begin{equation}
    |\varphi(y)-\varphi(y')|
 \end{equation}
 \begin{equation}
 \le
 |K_{s_1}(x_1,y)-K_{s_1}(x_1,y'))||
 K_{s_2}(x_2, y)|
\end{equation}
 \begin{equation}+|K_{s_1}(x_1, y')|
 |K_{s_2}(x_2, y) - K_{s_2}(x_2, y'))|
\end{equation}
\begin{equation}
      \le \frac{2^{229 a^3}}{\mu(B(x_1, D^{s_1}))\mu(B(x_2, D^{s_2}))}
       \left(\left(\frac{ \rho(y,y')}{D^{s_1}}\right)^{1/a}+
       \left(\frac{ \rho(y,y')}{D^{s_2}}\right)^{1/a}\right)
\end{equation}
\begin{equation}\label{holderpart}
      \le \frac{2^{230 a^3}}{\mu(B(x_1, D^{s_1}))\mu(B(x_2, D^{s_2}))}
       \left(\frac{ \rho(y,y')}{D^{s_1}}\right)^{1/a}\,.
\end{equation}
Adding the estimates \eqref{suppart} and \eqref{holderpart} gives \eqref{eqt11}.
This proves the lemma.
\end{proof}

The following auxiliary statement about the support of $T_\fp^*g$ will be
used repeatedly.

\begin{lemma}[tile range support]\label{tile-range-support}
    \leanok
    \lean{Tile.range_support}
        For each $\fp\in \fP$, and each $y\in X$, we have that
    \begin{equation}\label{tstargnot0}
         T_{\fp}^* g(y)\neq 0
    \end{equation}
       implies
    \begin{equation}\label{ynotfar}
        y\in B(\pc(\fp),5D^{\ps(\fp)})\, .
    \end{equation}
    \end{lemma}
    \begin{proof}
    \leanok
    Fix $\fp$ and $y$ with \eqref{tstargnot0}.
    Then there exists $x\in E(\fp)$ with
    \begin{equation}
       \overline{K_{\ps(\fp)}(x,y)}e(-\tQ(x)(y)
        +\tQ(x)(x))g(x) \neq 0\, .
    \end{equation}
    As $E(\fp)\subset \scI(\fp)$ and by the squeezing property
    \eqref{eq-vol-sp-cube}, we have
    \begin{equation}
        \rho(x,\pc(\fp)) < 4D^{\ps(\fp)}\, .
    \end{equation}
    As $K_{\ps(\fp)}(x,y)\neq 0$, we have by \eqref{supp-Ks}
    that
    \begin{equation}
    \rho(x,y)\le \frac 12 D^{\ps(\fp)}\, .
    \end{equation}
    Now \eqref{ynotfar} follows by the triangle inequality.
\end{proof}

The next lemma is a geometric estimate for two tiles.
\begin{lemma}\label{tile-uncertainty}
\leanok
\lean{Tile.uncertainty}
\uses{monotone-cube-metrics,tile-range-support}
    Let $\fp_1, \fp_2\in \fP$ with
    $B(\pc(\fp_1),5D^{\ps(\fp_1)}) \cap B(\pc(\fp_2),5D^{\ps(\fp_2)}) \ne \emptyset$ and
$\ps({\fp_1})\leq \ps({\fp_2})$. For each $x_1\in E(\fp_1)$ and
$x_2\in E(\fp_2)$ we have
\begin{equation}\label{tgeo}
  1+d_{\fp_1}(\fcc(\fp_1), \fcc(\fp_2))\le
    2^{8a}(1 + d_{B(x_1, D^{\ps(\fp_1)})}(\tQ(x_1),\tQ(x_2)))\, .
\end{equation}
\end{lemma}
\begin{proof}
\leanok
Let $i\in \{1,2\}$.
By Definition \eqref{defineep} of $E$,
we have $\tQ(x_i)\in \fc(\fp_i)$
With \eqref{eq-freq-comp-ball} we then conclude
\begin{equation}\label{dponetwo}
    d_{\fp_i}(\tQ(x_i),\fcc(\fp_i)) < 1\, .
\end{equation}
We have by the triangle inequality and \eqref{eq-vol-sp-cube} that $\scI(\fp_1) \subset B(\pc(\fp_2),14D^{\ps(\fp_2)})$.
Thus, using again \eqref{eq-vol-sp-cube} and the doubling property \eqref{firstdb}
\begin{equation}\label{tgeo0.5}
    d_{\fp_1}(\tQ(x_2), \fcc(\fp_2)) \le 2^{6a} d_{\fp_2}(\tQ(x_2), \fcc(\fp_2)) \le 2^{6a}\,.
\end{equation}
By the triangle inequality, we obtain from \eqref{dponetwo} and
\eqref{tgeo0.5}
\begin{equation}\label{tgeo1}
     1+d_{\fp_1}(\fcc(\fp_1), \fcc(\fp_2))\le 2 + 2^{6a} +d_{\fp_1}(\tQ(x_1), \tQ(x_2))\, .
\end{equation}
As $x_1\in \scI(\fp_1)$ by Definition \eqref{defineep} of $E$, we have by the squeezing property \eqref{eq-vol-sp-cube}
\begin{equation}
    d(x_1,\pc(\fp_1))\le 4D^{\ps(\fp_1)}
\end{equation}
and thus by \eqref{eq-vol-sp-cube} again and the triangle inequality
\begin{equation}
    \scI(\fp_1)\subset B(x_1,8D^{\ps(\fp_1)})\, .
\end{equation}
We thus estimate the right-hand side of \eqref{tgeo1} with monotonicity \eqref{monotonedb} of the metrics $d_B$ by
\begin{equation}\label{tgeo1.5}
    \le 2 + 2^{6a} + d_{B(x_1,8D^{\ps(\fp_1)})}(\tQ(x_1), \tQ(x_2))\, .
\end{equation}
This is further estimated by applying the doubling property \eqref{firstdb} three times by
\begin{equation}\label{tgeo2}
    \le 2 + 2^{6a} +2^{3a}d_{B_1(x_1, D^{\ps(\fp_1)})}(\tQ(x_1), \tQ(x_2))\, .
\end{equation}
Now \eqref{tgeo} follows with $a\ge 4$.
\end{proof}

We now prove \Cref{tile-correlation}.
\begin{proof}[Proof of \Cref{tile-correlation}]
    \leanok
    \proves{tile-correlation}
We begin with \eqref{eq-basic-TT*-est}.
By Lemma \ref{tile-range-support}, the left-hand side of \eqref{eq-basic-TT*-est} vanishes if
$B(\pc(\fp'),5D^{\ps(\fp')}) \cap B(\pc(\fp),5D^{\ps(\fp)}) = \emptyset$.
Thus we can assume for the remainder of the proof that
\begin{equation}
        \label{intersec5B}
        B(\pc(\fp'),5D^{\ps(\fp')}) \cap B(\pc(\fp),5D^{\ps(\fp)}) \neq \emptyset\,.
\end{equation}
We expand the left-hand side of \eqref{eq-basic-TT*-est} as
\begin{equation}\label{tstartstar}
\left|\int \int_{E(\fp')} \overline{K_{\ps(\fp')}(x_1,y)}e(-\tQ(x_1)(y)+
    \tQ(x_1)(x_1))g(x_1)\, d\mu(x_1) \right.
\end{equation}
\begin{equation}\label{tstartstar'}
 \times \left.\int_{E(\fp)} {K_{\ps(\fp)}(x_2,y)}e(\tQ(x_2)(y)
    -\tQ(x_2)(x_2))\overline{g(x_2)}\, d\mu(x_2)\, d\mu(y)\right|\, .
\end{equation}
By Fubini and the triangle inequality and
the fact $|e(\tQ(x_i)(x_i))|=1$ for $i=1,2$, we can estimate
\eqref{tstartstar} and \eqref{tstartstar'} from above by
\begin{equation}\label{eqa1}
    \int_{E(\fp')} \int_{E(\fp)} {\bf I}(x_1, x_2)\, d\mu(x_1)d\mu(x_2)\,.
\end{equation}
with
\begin{equation}
    {\bf I}(x_1, x_2):=
    \left|\int
    e(-\tQ(x_1)(y)+\tQ(x_2)(y))\varphi_{x_1,x_2}(y)
    d\mu(y) \, g(x_1)g(x_2)\right|
\end{equation}

We estimate for fixed $x_1\in E(\fp')$ and
$x_2\in E(\fp)$ the inner integral of \eqref{eqa1} with
\Cref{Holder-van-der-Corput}. The function
$\varphi:=\varphi_{x_1,x_2}$ satisfies the assumptions of
\Cref{Holder-van-der-Corput} with $z = x_1$ and $R = D^{s_1}$ by \Cref{correlation-kernel-bound}.
We obtain with $B':= B(x_1, D^{\ps(\fp')})$,
\begin{equation*}
 {\bf I}(x_1, x_2) \le 2^{8a} \mu(B') \|{\varphi}\|_{C^\tau(B')}
       (1 + d_{B'}(\tQ(x_1),\tQ(x_2)))^{-1/(2a^2+a^3)}|g(x_1)g(x_2)|
\end{equation*}
\begin{equation}
\label{eqa1.5}
 \le \frac{2^{231a^3+8a}}
 {\mu(B(x_2, D^{\ps(\fp)}))}
       (1 + d_{B'}(\tQ(x_1),\tQ(x_2)))^{-1/(2a^2+a^3)}\,.
\end{equation}
Using \eqref{intersec5B}, \Cref{tile-uncertainty} and $a\ge 1$ estimates \eqref{eqa1.5} by
\begin{equation}\label{eqa2}
 \le \frac{2^{231a^3 + 8a + 1}}
 {\mu(B(x_2, D^{\ps(\fp)}))}
       (1+d_{\fp'}(\fcc(\fp'), \fcc(\fp)))^{-1/(2a^2+a^3)}|g(x_1)g(x_2)|\,.
\end{equation}
As $x_2\in \scI(\fp)$ by Definition \eqref{defineep} of $E$, we have by \eqref{eq-vol-sp-cube}
\begin{equation}
    \rho(x_2,\pc(\fp)) < 4D^{\ps(\fp)}
\end{equation}
and thus by \eqref{eq-vol-sp-cube} again and the triangle inequality
\begin{equation}
    \scI(\fp)\subset B(x_2,8D^{\ps(\fp)})\, .
\end{equation}
Using three iterations of the doubling property \eqref{doublingx} give
\begin{equation}
    \mu(\scI(\fp))\le 2^{3a}\mu(B(x_2,D^{\ps(\fp)}))\, .
\end{equation}
With $a\ge 4$ and \eqref{eqa2} we conclude \eqref{eq-basic-TT*-est}.

Now assume the left-hand side of \eqref{eq-basic-TT*-est} is not zero.
There is a $y\in X$ with
\begin{equation}
    T^*_{\fp'}g(y)\overline{T^*_{\fp}g(y)}\neq 0
\end{equation}
By the triangle inequality and \Cref{tile-range-support}, we conclude
\begin{equation}
   \rho(\pc(\fp),\pc(\fp'))\le \rho(\pc(\fp),y) +\rho(\pc(\fp'),y)
   \le 5D^{\ps(\fp)}+5D^{\ps(\fp')}\le 10 D^{\ps(\fp)}\, .
\end{equation}
By the squeezing property \eqref{eq-vol-sp-cube} and the triangle inequality,
we conclude
\begin{equation}
    \scI(\fp') \subset B(\pc(\fp), 14D^{\ps(\fp)})\, .
\end{equation}
   This completes the proof of Lemma \ref{tile-correlation}.
\end{proof}

\subsection{Proof of the Antichain Tile Count Lemma}
\label{subsec-geolem}

\begin{lemma}[tile reach]\label{tile-reach}
\leanok
\lean{Antichain.tile_reach}
\uses{monotone-cube-metrics}
Let $\mfa\in \Mf$ and $N\ge0$ be an integer.
Let $\fp, \fp'\in \fP$ with
\begin{equation}\label{eqassumedismfa}
    d_{\fp}(\fcc(\fp), \mfa))\le 2^N\,
\end{equation}
\begin{equation}\label{eqassumedismfap}
    d_{\fp'}(\fcc(\fp'), \mfa))\le 2^N\, .
\end{equation}
Assume $\scI(\fp)\subset \scI(\fp')$ and $\ps(\fp)<\ps(\fp')$.
Then
\begin{equation}\label{lp'lp''}2^{N+2}\fp\lesssim 2^{N+2} \fp'\, .
\end{equation}
\end{lemma}

\begin{proof}
\leanok
By \Cref{monotone-cube-metrics}, we have
\begin{equation}
     d_{\fp}(\fcc(\fp'),\mfa)
     \le d_{\fp'}(\fcc(\fp'),\mfa)
     \le 2^{N} \, .
\end{equation}
Together with \eqref{eqassumedismfa} and the triangle inequality, we obtain
\begin{equation}\label{eqdistqpqp}
    d_{\fp}(\fcc(\fp'),\fcc(\fp))\le 2^{N+1} \, .
\end{equation}
Now assume
\begin{equation}
    \mfa'\in B_{\fp'}(\fcc(\fp'),2^{N+2}).
\end{equation}
By the doubling property \eqref{firstdb}, applied five times, we have
\begin{equation}\label{ageo1} d_{B(\pc(\fp'),8D^{\ps(\fp')})}(\fcc(\fp'),\mfa') < 2^{5a+N+2}\, .
\end{equation}
We have by the squeezing property \eqref{eq-vol-sp-cube}
\begin{equation}
 \pc(\fp)\in
B(\pc(\fp'),4D^{\ps(\fp')})\, .
\end{equation}
Hence by the triangle inequality
\begin{equation}
 B(\pc(\fp), 4D^{\ps(\fp')})
 \subseteq
B(\pc(\fp'),8D^{\ps(\fp')})\, .
\end{equation}
Together with \eqref{ageo1} and monotonicity \eqref{monotonedb} of $d$
\begin{equation}
    d_{B(\pc(\fp),4D^{\ps(\fp')})}(\fcc(\fp'),\mfa') < 2^{5a+N+2}\, .
\end{equation}
Using the doubling property \eqref{seconddb} $5a+2$ times gives
\begin{equation}
    d_{B(\pc(\fp),2^{2-5a^2-2a}D^{\ps(\fp')})}(\fcc(\fp'),\mfa') < 2^{N}\, .
\end{equation}
Using $\ps(\fp)<\ps(\fp')$ and $D=2^{100a^2}$ and $a\ge 4$ gives
\begin{equation}
    d_{\fp}(\fcc(\fp'),\mfa') < 2^{N}\, .
\end{equation}
With the triangle inequality and \eqref{eqdistqpqp},
\begin{equation}
    d_{\fp}(\fcc(\fp),\mfa') < 2^{N+2}\, .
\end{equation}
This shows
\begin{equation}
B_{\fp'}(\fcc(\fp'),2^{N+2})\subset B_{\fp}(\fcc(\fp),2^{N+2})\, .
\end{equation}
This implies \eqref{lp'lp''} and completes the proof of the lemma.
\end{proof}

For $\mfa \in \Mf$ and $N\ge 0$ define
\begin{equation}\label{eqantidefap}
    \mathfrak{A}_{\mfa,N}:=\{\fp\in\mathfrak{A}: 2^{N}\le 1+d_{\fp}(\fcc(\fp), \mfa) < 2^{N+1}\} \, .
\end{equation}

\begin{lemma}[stack density]
\label{stack-density}
\leanok
\lean{Antichain.stack_density}
Let $\mfa \in \Mf$, $N\ge 0$ and
$L\in \mathcal{D}$. Then
\begin{equation}\label{eqanti-1}
    \sum_{\fp\in\mathfrak{A}_{\mfa,N}:\scI(\fp)=L}\mu(E(\fp)\cap G)\le 2^{a(N+5)}\dens_1(\mathfrak{A})\mu(L)\, .
\end{equation}
\end{lemma}
\begin{proof}
\leanok
Let $\mfa,N,L$ be given and set
\begin{equation}
\mathfrak{A}':=\{\fp\in\mathfrak{A}_{\mfa,N}:\scI(\fp)=L\}\, .
\end{equation}
Let
$\fp\in\mathfrak{A}'$.
We have
by Definition \eqref{definedens1}
using $\lambda=2$ and the squeezing property \eqref{eq-freq-comp-ball}
\begin{equation}\label{eqanti-3}
\mu(E(\fp)\cap G)\le \mu(E_2(2, \fp))\le 2^{a}\dens_1(\mathfrak{A}')\mu(L)\, .
\end{equation}
By the covering property \eqref{thirddb}, applied $N+4$ times, there is a collection $\Mf'$ of at most $2^{a(N+4)}$
elements such that
\begin{equation}\label{eqanti-4}
    B_{\fp}(\mfa, 2^{N+1})\subset \bigcup_{\mfa'\in \Mf'}
    B_{\fp}(\mfa', 0.2)\, .
\end{equation}
As each $\fcc(\fp')$ with $\fp'\in \mathfrak{A}'$
is contained in the left-hand-side
of \eqref{eqanti-4}
by definition (because $\scI(\fp') = \scI(\fp))$, it is in at least one $B_{\fp}(\mfa', 0.2)$
with $\mfa'\in \Mf'$.

For two different $\fp',\fp''\in \mathfrak{A}'$, we have by
\eqref{eq-dis-freq-cover} that
$\fc(\fp')$ and $\fc(\fp'')$ are disjoint and thus by the squeezing property \eqref{eq-freq-comp-ball} we have for every $\mfa'\in \Mf'$
\begin{equation}
    \mfa'\not\in B_{\fp}(\fcc(\fp'), 0.2)\cap
B_{\fp}(\fcc(\fp''), 0.2)\, .
\end{equation}
Hence at most one of $\fcc(\fp')$
and $\fcc(\fp'')$ is in
$B_{\fp}(\mfa', 0.2)$.
It follows that there are at most $2^{a(N+4)}$ elements in
$\mathfrak{A}'$. Adding \eqref{eqanti-3} over $\mathfrak{A}'$ proves
\eqref{eqanti-1}.
\end{proof}

\begin{lemma}[local antichain density]\label{local-antichain-density}
\leanok
\lean{Antichain.local_antichain_density, Antichain.Ep_inter_G_inter_Ip'_subset_E2}
\uses{tile-disjointness,tile-reach}
Let $\mfa\in\Mf$ and {$N$} be
an integer. Let $\fp_{\mfa}$ be a tile with $\mfa\in B_{\fp_{\mfa}}(\fcc(\fp_{\mfa}), 2^{N+1})$.
Then we have
\begin{equation}\label{eqanti-0.5}
    \sum_{\fp\in\mathfrak{A}_{\mfa,N}: \ps(\fp_{\mfa})<\ps(\fp)}\mu(E(\fp)\cap G \cap \scI(\fp_{\mfa}))
    \le \mu (E_2(2^{N+3},\fp_{\mfa}))
 \, .
\end{equation}

\end{lemma}

\begin{proof}
\leanok

Let $\fp$ be any tile in $\mathfrak{A}_{\mfa,N}$ with $\ps(\fp_{\mfa})<\ps(\fp)$. By definition of
$E$, the tile contributes zero to the sum on the left-hand side of \eqref{eqanti-0.5} unless
 $\scI(\fp)\cap \scI(\fp_{\mfa}) \neq \emptyset$, which we may assume. With $\ps(\fp_{\mfa})<\ps(\fp)$
and the dyadic property
\eqref{dyadicproperty} we conclude $\scI(\fp_{\mfa})\subset \scI(\fp)$.
We conclude from $\fp \in \mathfrak{A}_{\mfa,N}$ that
\begin{equation}
    \mfa \in B_\fp(\fcc(\fp), 2^{N+1})\, .
\end{equation}
With \Cref{tile-reach} and the assumption on $\fp_\mfa$, we conclude
\begin{equation}
    2^{N+3}\fp_{\mfa} \lesssim 2^{N+3}\fp \, .
\end{equation}
By Definition \eqref{definee2} of $E_2$, we conclude
\begin{equation}
    E(\fp)\cap G \cap \scI(\fp_{\mfa}) \subset E_2(2^{N+3},\fp_{\mfa})\, .
\end{equation}
Using disjointedness of the various $E(\fp)$ with $\fp\in \mathfrak{A}$ by \Cref{tile-disjointness}, we obtain \eqref{eqanti-0.5}.
This proves the lemma.
\end{proof}
\begin{lemma}[global antichain density]
\label{global-antichain-density}
\leanok
\lean{Antichain.global_antichain_density}
\uses{stack-density,local-antichain-density}
Let $\mfa\in Q(X)$ and let $N\ge 0$ be
an integer. Then we have
\begin{equation}\label{eqanti00}
    \sum_{\fp\in\mathfrak{A}_{\mfa,N}}\mu(E(\fp)\cap G)
    \le
 2^{101a^3+Na}\dens_1(\mathfrak{A})\mu\left(\cup_{\fp\in\mathfrak{A}}I_{\fp}\right)\, .
\end{equation}
\end{lemma}

\begin{proof}
\leanok
{Fix $\mfa$ and $N$. Let
$\mathfrak{A}'$ be the set of $\fp\in\mathfrak{A}_{\mfa,N}$ such that $\scI(\fp)\cap G$ is not empty and $s(\fp) > -S$.
Let $\mathfrak{A}_{-S}$ be the set of $\fp\in\mathfrak{A}_{\mfa,N}$ such that $\scI(\fp)\cap G$ is not empty and $s(\fp) = -S$ }
Then we have
\begin{equation*}
    \sum_{\fp\in\mathfrak{A}_{\mfa,N}}\mu(E(\fp)\cap G) =
    \sum_{\fp\in\mathfrak{A}'}\mu(E(\fp)\cap G) + \sum_{\fp\in\mathfrak{A}_{-S}}\mu(E(\fp)\cap G)
\end{equation*}

We start by estimating the contribution of $\mathfrak{A}_{-S}$. Let $\mathcal{L}_{-S}$ be the
collection of dyadic cubes $\scI(\fp)$ with $\fp \in \mathfrak{A}_{-S}$. They all have scale $-S$,
by definition of $\mathcal{L}_{-S}$, and hence they are pairwise disjoint by the dyadic property
\eqref{dyadicproperty}. We write
\[
    \sum_{\fp\in\mathfrak{A}_{-S}}\mu(E(\fp)\cap G)
    = \sum_{L \in \mathcal{L}_{-S}} \sum_{\fp\in\mathfrak{A}_{-S}, \scI(\fp) = L} \mu(E(\fp)\cap G),
\]
and using \Cref{stack-density}, we estimate
\begin{equation}\label{eq_minus_S}
    \le 2^{a(N+5)} \dens_1(\mathfrak{A}) \sum_{L \in \mathcal{L}_{-S}} \mu(L)
    \le 2^{a(N+5)} \dens_1(\mathfrak{A}) \mu\left(\cup_{\fp\in\mathfrak{A}}I_{\fp}\right).
\end{equation}

We turn to $\mathfrak{A}'$.
Let $\mathcal{L}$ be the collection of dyadic cubes $I\in\mathcal{D}$ such that $I \le \scI(\fp)$
for some $\fp\in\mathfrak{A}'$ and such that $\scI(\fp)\not \subset I$ for all $\fp\in\mathfrak{A}'$.
By \eqref{coverdyadic}, for each $\fp \in \mathfrak{A}'$
and each $x\in \scI(\fp)\cap G$, there is $I\in \mathcal{D}$ with $s(I)=-S$ and $x\in I$.
By \eqref{dyadicproperty}, we have $I\subset \scI(\fp)$. Since for all $\fp' \in \mathfrak{A}'$ we
have $s(\fp') > -S$, we have $I \in \mathcal{L}$. Hence
\begin{equation}
    \scI(\fp)\subset \bigcup\{I\in \mathcal{D}: s(I)=-S, I\subset \scI(\fp)\}\subset \bigcup \mathcal{L}\, .
\end{equation}
As each $I\in \mathcal{L}$ satisfies $I\subset \scI(\fp)$ for some $\fp$ in $\mathfrak{A'}$, we conclude
\begin{equation}
    \bigcup\mathcal{L}=\bigcup_{\fp \in \mathfrak{A}'}\scI(\fp)\, .
\end{equation}
Let $\mathcal{L}^*$ be the set of maximal elements in $\mathcal{L}$ with respect to set inclusion.
By \eqref{dyadicproperty}, the elements in $\mathcal{L}^*$ are pairwise disjoint and we have
 \begin{equation}\label{eqdecAprime}
\bigcup\mathcal{L}^*=\bigcup_{\fp \in \mathfrak{A}'}\scI(\fp)\, .
   \end{equation}
Using the partition \eqref{eqdecAprime} into elements of $\mathcal{L}$ in \eqref{eqanti0}, it suffices to show for each $L\in \mathcal{L}^*$
\begin{equation}\label{eqanti0}
    \sum_{\fp\in\mathfrak{A}'}\mu(E(\fp)\cap G \cap L)
    \le
    2^{101a^3+aN}
    \dens_1(\mathfrak{A})\mu(L)\,.
\end{equation}
Fix $L\in \mathcal{L}^*$. By definition of $\mathcal{L}^*$, there exists an element $\fp'\in \mathfrak{A}'$
such that $L\subset \scI(\fp')$. Pick such an element $\fp'$ in $\mathfrak{A}$ with
minimal $\ps(\fp')$. As $\scI(\fp')\not \subset L$ by definition of $L$,
we have with \eqref{dyadicproperty} that $s(L)< \ps(\fp')$.
In particular $s(L)<S$, thus $L \ne I_0$ and hence by \eqref{subsetmaxcube}
there exists a cube $J \in \mathcal{D}$ with $L \subsetneq J$.
By \eqref{coverdyadic}, there is an $L'\in \mathcal{D}$ with $s(L')=s(L)+1$ and $L \le L'$.
By \eqref{dyadicproperty}, we have $L\subset L'$.

We split the left-hand side of \eqref{eqanti0} as
\begin{equation}\label{eqanti1}
    \sum_{\fp\in\mathfrak{A}':\scI(\fp)=L'}\mu(E(\fp)\cap G\cap L)
\end{equation}
\begin{equation}\label{eqanti2}
    +
     \sum_{\fp\in\mathfrak{A}':\scI(\fp)\neq L'}\mu(E(\fp)\cap G\cap L)\, ,
\end{equation}

We first estimate \eqref{eqanti1}
with \Cref{stack-density} by
\begin{equation}\label{equanti1.5}
    \le \sum_{\fp\in\mathfrak{A}':\scI(\fp)=L'}\mu(E(\fp)\cap G\cap L')\le 2^{a(N+5)}\dens_1(\mathfrak{A})\mu(L')\, .
\end{equation}

We turn to \eqref{eqanti2}.
Consider the element $\fp'\in \mathfrak{A}'$ as above
with $L\subset \scI(\fp')$ and $s(L)<\ps(\fp')$.
As $L\subset L'$ and $s(L')=s(L)+1$, we conclude with the dyadic property that $L'\subset \scI(\fp')$.
By maximality of $L$, we have
$L'\not\in \mathcal{L}$.
This together with the existence of the given $\fp'\in \mathfrak{A}$
with $L'\subset \scI(\fp')$
shows by definition of $\mathcal{L}$ that there exists $\fp''\in \mathfrak{A}'$ with
$\scI(\fp'')\subset L'$.

If $\scI(\fp'') = L'$, then we set $\fp_{\mfa} = \fp''$ and note that as $\fp'' \in \mathfrak{A}_{\mfa,N}$
\begin{equation}
    \mfa \in B(\fcc(\fp''), 2^{N+1})\, .
\end{equation}

If $\scI(\fp'') \ne L'$, then it follows that $s(\fp'') < s(L')$.
By the covering property \eqref{eq-dis-freq-cover}, there exists a unique $\fp_{\mfa}$ with
\begin{equation*}
    \scI(\fp_{\mfa})=L'
\end{equation*}
such that $\mfa\in \fc(\fp_{\mfa})$. We take this as the definition of $\fp_\mfa$ in this case.
Note that
\begin{equation}
    \mfa\in B(\fcc(\fp_{\mfa}), 1)
\end{equation}
so by \Cref{tile-reach}, we conclude
\begin{equation}
    2^{N+3}\fp'' \lesssim 2^{N+3}\fp_{\mfa} \, .
\end{equation}
This clearly also holds in the case $\scI(\fp'') = L'$, since then $\fp'' = \fp_\mfa$. Furthermore,
in both cases it also holds that
\begin{equation}
    \mfa\in B_{\fp_{\mfa}}(\fcc(\fp_{\mfa}), 2^{N+1}).
\end{equation}

As $\fp''\in \mathfrak{A}'$, we have by Definition
\eqref{definedens1} of $\dens_1$ that
\begin{equation}\label{pmfadens}
   \mu(E_2(2^{N+3}, \fp_{\mfa}))\le 2^{Na+3a}\dens_1(\mathfrak{A}) {\mu(L')}\, .
\end{equation}
Now let $\fp$ be any tile in the summation set in \eqref{eqanti2}, that is, $\fp\in \mathfrak{A}'$ and $\scI(\fp)\neq L'$.
Then $\scI(\fp)\cap L\neq \emptyset$. It follows by the dyadic property \eqref{dyadicproperty}
and the definition of $L$ that
$L\subset \scI(\fp)$ and $L\neq \scI(\fp)$. By the dyadic property \eqref{dyadicproperty}, we have
$s(L)<\ps(\fp)$ and thus $s(L')\le \ps(\fp)$. By the dyadic property
   \eqref{dyadicproperty} again, we have $L'\subset \scI(\fp)$.
As $L'\neq \scI(\fp)$, we conclude $s(L)<\ps(\fp)$.
By \Cref{local-antichain-density}, we can thus estimate \eqref{eqanti2} by
\begin{equation}\label{eqanti0.5}
    \sum_{\fp\in\mathfrak{A}':\scI(\fp)\neq L'}\mu(E(\fp)\cap G\cap L')
    \le \mu (E_2(2^{N+3},\fp_{\mfa}))\, .
\end{equation}
With the decomposition in \eqref{eqanti1} and \eqref{eqanti2} and the
estimates \eqref{equanti1.5}, \eqref{eqanti-0.5}, \eqref{pmfadens} we obtain
the estimate
\begin{equation}\label{eqanti3.14}
\sum_{\fp\in\mathfrak{A}'}\mu(E(\fp)\cap G \cap L)
    \le (2^{a(N+5)}+2^{Na+3a})\dens_1(\mathfrak{A})\mu(L')\,.
\end{equation}

Using $s(L')=s(L)+1$ and $D=2^{100a^2}$ and the
squeezing property \eqref{eq-vol-sp-cube}
and the doubling property \eqref{doublingx} $100a^2+4$ times , we obtain
\begin{equation}
    \mu(L')\le 2^{100a^3+4a}\mu(L)\, .
\end{equation}
Inserting in \eqref{eqanti3.14}, adding the estimate \eqref{eq_minus_S} and using $a\ge 4$ gives \eqref{eqanti0}.
This completes the proof of the lemma.
\end{proof}

We turn to the proof of \Cref{antichain-tile-count}.

\begin{proof}[Proof of \Cref{antichain-tile-count}]
\leanok
\proves{antichain-tile-count}

Using that $\mathfrak{A}$ is the disjoint union of the $\mathfrak{A}_{\mfa,N}$ with $N\ge 0$,
we estimate the left-hand side \eqref{eq-antichain-Lp}
with the triangle inequality by
\begin{equation}\label{eqanti23}
\le \sum_{N\ge 0} \left\|\sum_{\fp\in \mathfrak{A}_{\mfa,N}} 2^{-N/(2a^2+a^3)}\mathbf{1}_{E(\fp)} \mathbf{1}_G\right\|_{p}
\end{equation}
We consider each individual term in this sum and estimate its $p$-th power.
Using that for each $x\in X$ by \Cref{tile-disjointness} there is at most one $\fp\in \mathfrak{A}$ with $x\in E(\fp)$, we have
\begin{equation}
    \left\|\sum_{\fp\in \mathfrak{A}_{\mfa,N}} 2^{-N/(2a^2+a^3)}\mathbf{1}_{E(\fp)} \mathbf{1}_G\right\|_{p}^p
\end{equation}
\begin{equation}
    =\int\Big(\sum_{\fp\in \mathfrak{A}_{\mfa,N}}2^{-N/(2a^2+a^3)}\mathbf{1}_{E(\fp)\cap G}(x)\Big)^p\, d\mu(x)
\end{equation}
\begin{equation}
    =\int\sum_{\fp\in\mathfrak{A}_{\mfa,N}}2^{-pN/(2a^2+a^3)}\mathbf{1}_{E(\fp)\cap G}(x)\, d\mu(x)
\end{equation}
\begin{equation}
    = 2^{-pN/(2a^2+a^3)} \sum_{\fp\in\mathfrak{A}_{\mfa,N}}\mu(E(\fp)\cap G)
\end{equation}

Using \Cref{global-antichain-density}, we estimate the last display by
\begin{equation}\label{eqanti21}
    \le 2^{-pN/(2a^2+a^3)+101a^3+Na}\dens_1(\mathfrak{A})\mu\left(\cup_{\fp\in\mathfrak{A}}\scI(\fp)\right)
\end{equation}
Using that $a\ge 4$ and since $p = 4a^4$, we have
$$
    pN/(2a^2+a^3)\ge 4a^4N/(3a^3) \ge Na+N\, .
$$
Hence we have for \eqref{eqanti21} the upper bound
$$
    \le 2^{101a^3-N}\dens_1(\mathfrak{A})\mu\left(\cup_{\fp\in\mathfrak{A}}\scI(\fp)\right)\, .
$$
Taking the $p$-th root and summing over $N\ge 0$ gives for \eqref{eqanti23} the upper bound
$$
    \le \left(\sum_{N\ge 0} 2^{-N/p}\right)2^{101a^3/p}\dens_1(\mathfrak{A})^{{\frac{1}{p}}}\mu\left(\cup_{\fp\in\mathfrak{A}}\scI(\fp)\right)^{{\frac{1}{p}}}
$$
$$
    = \left(1-2^{-1/p}\right)^{-1} 2^{101a^3/p} \dens_1(\mathfrak{A})^{\frac 1p}\mu\left(\cup_{\fp\in\mathfrak{A}}\scI(\fp)\right)^{\frac 1p}\, .
$$
Using that $p = 4a^4$ and $a \ge 4$, this proves the lemma.
\end{proof}

\section{Proof of the Forest Operator Proposition}

\label{treesection}

\subsection{The pointwise tree estimate}
Fix a forest $(\fU, \fT)$. The main result of this subsection is \Cref{pointwise-tree-estimate}, we begin this section with some definitions necessary to state the lemma.

For $\fu \in \fU$ and $x\in X$, we define
$$
    \sigma (\fu, x):=\{\ps(\fp):\fp\in \fT(\fu), x\in E(\fp)\}\,.
$$
This is a subset of $\mathbb{Z} \cap [-S, S]$, so has a minimum and a maximum. We set
$$
    \overline{\sigma} (\fu, x) := \max \sigma(\fT(\fu), x)
$$
$$
    \underline{\sigma} (\fu, x) := \min\sigma(\fT(\fu), x)\,.
$$
\begin{lemma}[convex scales]
\label{convex-scales}
\leanok
\lean{TileStructure.Forest.convex_scales}
    For each $\fu \in \fU$, we have
    $$
        \sigma(\fu, x) = \mathbb{Z} \cap [\underline{\sigma} (\fu, x), \overline{\sigma} (\fu, x)]\,.
    $$
\end{lemma}

\begin{proof}
    \leanok
    Let $s \in \mathbb{Z}$ with $\underline{\sigma} (\fu, x) \le s \le \overline{\sigma} (\fu, x)$. By definition of $\sigma$, there exists $\fp \in \fT(\fu)$ with $\ps(\fp) = \underline{\sigma} (\fu, x)$ and $x \in E(\fp)$, and there exists $\fp'' \in \fT(\fu)$ with $\ps(\fp'') = \overline{\sigma} (\fu, x)$ and $x \in E(\fp'') \subset \scI(\fp'')$. By property \eqref{coverdyadic} of the dyadic grid, there exists a cube $I \in \mathcal{D}$ of scale $s$ with $x \in I$. By property \eqref{eq-dis-freq-cover}, there exists a tile $\fp' \in \fP(I)$ with $\tQ(x) \in \fc(\fp')$. By the dyadic property \eqref{dyadicproperty} we have $\scI(\fp) \subset \scI(\fp') \subset \scI(\fp'')$, and by \eqref{eq-freq-dyadic}, we have $\fc(\fp'') \subset \fc(\fp') \subset \fc(\fp)$. Thus $\fp \le \fp' \le\fp''$, which gives with the convexity property \eqref{forest2} of $\fT(\fu)$ that $\fp' \in \fT(\fu)$, so $s \in \sigma(\fu, x)$.
\end{proof}

For a nonempty collection of tiles $\mathfrak{S} \subset \fP$ we define
$$
    \mathcal{J}_0(\mathfrak{S})
$$
to be the collection of all dyadic cubes $J \in \mathcal{D}$ such that $s(J) = -S$ or
$$
    \scI(\fp) \not\subset B(c(J), 100D^{s(J) + 1})
$$
for all $\fp \in \mathfrak{S}$. We define $\mathcal{J}(\mathfrak{S})$ to be the collection of inclusion maximal cubes in $\mathcal{J}_0(\mathfrak{S})$.

We further define
$$
    \mathcal{L}_0(\mathfrak{S})
$$
to be the collection of dyadic cubes $L \in \mathcal{D}$ such that $s(L) = -S$, or there exists $\fp \in \mathfrak{S}$ with $L \subset \scI(\fp)$ and there exists no $\fp \in \mathfrak{S}$ with $\scI(\fp) \subset L$. We define $\mathcal{L}(\mathfrak{S})$ to be the collection of inclusion maximal cubes in $\mathcal{L}_0(\mathfrak{S})$.

\begin{lemma}[dyadic partitions]
    \label{dyadic-partitions}
    \leanok
    \lean{TileStructure.Forest.biUnion_𝓙, TileStructure.Forest.pairwiseDisjoint_𝓙,
    TileStructure.Forest.biUnion_𝓛, TileStructure.Forest.pairwiseDisjoint_𝓛}
    For each $\mathfrak{S} \subset \fP$, we have
    \begin{equation}
        \label{eq-J-partition}
        \bigcup_{I \in \mathcal{D}} I = \dot{\bigcup_{J \in \mathcal{J}(\mathfrak{S})}} J
    \end{equation}
    and
    \begin{equation}
        \label{eq-L-partition}
        \bigcup_{I \in \mathcal{D}} I = \dot{\bigcup_{L \in \mathcal{L}(\mathfrak{S})}} L\,.
    \end{equation}
\end{lemma}

\begin{proof}
    \leanok
    Since $\mathcal{J}(\mathfrak{S})$ is the set of inclusion maximal cubes in $\mathcal{J}_0(\mathfrak{S})$, cubes in $\mathcal{J}(\mathfrak{S})$ are pairwise disjoint by \eqref{dyadicproperty}. The same applies to $\mathcal{L}(\mathfrak{S})$.

    If $x \in \bigcup_{I \in \mathcal{D}} I$, then there exists by \eqref{coverdyadic} a cube $I \in \mathcal{D}$ with $x \in I$ and $s(I) = -S$. Then $I \in \mathcal{J}_0(\mathfrak{S})$. There exists an inclusion maximal cube in $\mathcal{J}_0(\mathfrak{S})$ containing $I$. This cube contains $x$ and is contained in $\mathcal{J}(\mathfrak{S})$. This shows one inclusion in \eqref{eq-J-partition}, the other one follows from $\mathcal{J}(\mathfrak{S}) \subset \mathcal{D}$.

    The proof of the two inclusions in \eqref{eq-L-partition} is similar.
\end{proof}

For a finite collection of pairwise disjoint cubes $\mathcal{C}$, define the projection operator
$$
    P_{\mathcal{C}}f(x) :=\sum_{J\in\mathcal{C}}\mathbf{1}_J(x) \frac{1}{\mu(J)}\int_J f(y) \, \mathrm{d}\mu(y)\,.
$$
Given a scale $-S \le s\le S$ and a point $x \in \bigcup_{I\in \mathcal{D}, s(I) = s} I$, there exists a unique cube in $\mathcal{D}$ of scale $s$ containing $x$ by \eqref{coverdyadic}. We denote it by $I_s(x)$. Define for $\mfa \in \Mf$ the nontangential maximal operator
\begin{equation}
    \label{eq-TN-def}
    T_{\mathcal{N}}^\mfa f(x) := \sup_{-S \le s_1 < S} \sup_{x' \in I_{s_1}(x)} \sup_{\substack{s_1 \le s_2 \le S\\ D^{s_2-1} \le R_Q(\mfa, x')}} \left| \sum_{s = s_1}^{s_2} \int K_s(x',y) f(y) \, \mathrm{d}\mu(y) \right|\,.
\end{equation}
Define for each $\fu \in \fU$ the auxiliary operator
$$
    S_{1,\fu}f(x)
$$
\begin{equation}
    \label{eq-def-S-op}
    :=\sum_{I\in\mathcal{D}} \mathbf{1}_{I}(x) \sum_{\substack{J\in \mathcal{J}(\fT(\fu))\\
    J\subset B(c(I), 16 D^{s(I)})\\ s(J) \le s(I)}} \frac{D^{(s(J) - s(I))/a}}{\mu(B(c(I), 16D^{s(I)}))}\int_J |f(y)| \, \mathrm{d}\mu(y)\,.
\end{equation}
Define also the collection of balls
$$
    \mathcal{B} = \{B(c(I), 2^s D^{s(I)+t}) \ : \ I \in \mathcal{D}\,, 0 \le s \le S + 5\,, 0 \le t \le 2S+3\}\,.
$$

The following pointwise estimate for operators associated to sets $\fT(\fu)$ is the main result of this subsection.

\begin{lemma}[pointwise tree estimate]
    \label{pointwise-tree-estimate}
    \leanok
    \lean{TileStructure.Forest.pointwise_tree_estimate}
    \uses{first-tree-pointwise,second-tree-pointwise,third-tree-pointwise}
    Let $\fu \in \fU$ and $L \in \mathcal{L}(\fT(\fu))$. Let $x, x' \in L$.
    Then for all bounded functions $f$ with bounded support
    $$
        \left|\sum_{\fp \in \fT(\fu)} T_{\fp}[ e(-\fcc(\fu))f](x)\right|
    $$
    \begin{equation}
        \label{eq-LJ-ptwise}
        \leq 2^{129a^3}(M_{\mathcal{B},1}+S_{1,\fu})P_{\mathcal{J}(\fT(\fu))}|f|(x')+|T_{\mathcal{N}}^{\fcc(\fu)} P_{\mathcal{J}(\fT(\fu))}f(x')|,
    \end{equation}
\end{lemma}

\begin{proof}
    \leanok
    By \eqref{definetp}, if $T_{\fp}[ e(-\fcc(\fu))f](x) \ne 0$, then $x \in E(\fp)$. Combining this with $|e(\fcc(\fu)(x))| = 1$, we obtain
    $$
        |\sum_{\fp \in \fT(\fu)} T_{\fp}[ e(-\fcc(\fu))f](x)|
    $$
    \begin{multline*}
        = \Bigg| \sum_{s \in \sigma(\fu, x)} \int e(-\fcc(\fu)(y) + \tQ(x)(y) + \fcc(\fu)(x) -\tQ(x)(x))\times\\
        K_s(x,y)f(y) \, \mathrm{d}\mu(y) \Bigg|\,.
    \end{multline*}
    Using the triangle inequality, we bound this by the sum of three terms:
    \begin{multline}
        \label{eq-term-A}
        \le \Bigg| \sum_{s \in \sigma(\fu, x)} \int (e(-\fcc(\fu)(y) + \tQ(x)(y) + \fcc(\fu)(x) -\tQ(x)(x))-1)\times\\
        K_s(x,y)f(y) \, \mathrm{d}\mu(y) \Bigg|
    \end{multline}
    \begin{equation}
        \label{eq-term-B}
        + \Bigg| \sum_{s \in \sigma(\fu, x)} \int K_s(x,y) P_{\mathcal{J}(\fT(\fu))} f(y) \, \mathrm{d}\mu(y) \Bigg|
    \end{equation}
    \begin{equation}
        \label{eq-term-C}
        + \Bigg| \sum_{s \in \sigma(\fu, x)} \int K_s(x,y) (f(y) - P_{\mathcal{J}(\fT(\fu))} f(y)) \, \mathrm{d}\mu(y) \Bigg|\,.
    \end{equation}
    The proof is completed using the bounds for these three terms proven respectively in \Cref{first-tree-pointwise}, \Cref{second-tree-pointwise} and \Cref{third-tree-pointwise}.
\end{proof}

\begin{lemma}[first tree pointwise]
    \label{first-tree-pointwise}
    \leanok
    \lean{TileStructure.Forest.first_tree_pointwise}
    \uses{convex-scales, kernel-summand}
    For all $\fu \in \fU$, all $L \in \mathcal{L}(\fT(\fu))$, all $x, x' \in L$ and all bounded $f$ with bounded support, we have
    $$
        \eqref{eq-term-A} \le 10 \cdot 2^{104a^3} M_{\mathcal{B}, 1}P_{\mathcal{J}(\fT(\fu))}|f|(x')\,.
    $$
\end{lemma}

\begin{proof}
    \leanok
    Let $s \in \sigma(\fu,x)$.
    If $x, y \in X$ are such that $K_s(x,y)\neq 0$, then, by \eqref{supp-Ks}, we have $\rho(x,y)\leq 1/2 D^s$. By $1$-Lipschitz continuity of the function $t \mapsto \exp(it) = e(t)$ and the property \eqref{osccontrol} of the metrics $d_B$, it follows that
    \begin{multline*}
        |e(-\fcc(\fu)(y)+\tQ(x)(y)+\fcc(\fu)(x)-\tQ(x)(x))-1|\\
        \leq d_{B(x, 1/2 D^{s})}(\fcc(\fu), \tQ(x))\,.
    \end{multline*}
    Let $\fp_s \in \fT(\fu)$ be a tile with $\ps(\fp_s) = s$ and $x \in E(\fp_s)$, and let $\fp'$ be a tile with $\ps(\fp') = \overline{\sigma}(\fu, x)$ and $x \in E(\fp')$.
    Using the monotonicity property \eqref{monotonedb}, the doubling property \eqref{firstdb} repeatedly, the definition of $d_{\fp}$ and \Cref{monotone-cube-metrics}, we can bound the previous display by
    $$
        d_{B(x, 4 D^{s})}(\fcc(\fu), \tQ(x)) \leq 2^{4a} d_{\fp_s}(\fcc(\fu), \tQ(x)) \le 2^{4a} 2^{s - \overline{\sigma}(\fu, x)} d_{\fp'}(\fcc(\fu), \tQ(x))\,.
    $$
    Since $\fcc(\fu) \in B_{\fp'}(\fcc(\fp'), 4)$ by \eqref{forest1} and $\tQ(x) \in \Omega(\fp') \subset B_{\fp'}(\fcc(\fp'), 1)$ by \eqref{eq-freq-comp-ball}, this is estimated by
    $$
        \le 5 \cdot 2^{4a} 2^{s - \overline{\sigma}(\fu, x)} \,.
    $$
    Using \eqref{eq-Ks-size}, it follows that
    $$
        \eqref{eq-term-A} \le 5\cdot 2^{103a^3} \sum_{s\in\sigma(x)}2^{s - \overline{\sigma}(\fu, x)} \frac{1}{\mu(B(x,D^s))}\int_{B(x,0.5D^{s})}|f(y)|\,\mathrm{d}\mu(y)\,.
    $$
    By \eqref{eq-J-partition}, the collection $\mathcal{J}$ is a partition of $\bigcup_{I \in \mathcal{D}} I$, so this is estimated by
    $$
         5\cdot 2^{103a^3} \sum_{s\in\sigma(x)}2^{s - \overline{\sigma}(\fu, x)} \frac{1}{\mu(B(x,D^s))}\sum_{\substack{J \in \mathcal{J}(\fT(\fu))\\J \cap B(x, 0.5D^s) \ne \emptyset} }\int_{J}|f(y)|\,\mathrm{d}\mu(y)\,.
    $$
    This expression does not change if we replace $|f|$ by $P_{\mathcal{J}(\fT(\fu))}|f|$.

    Let $J \in \mathcal{J}(\fT(\fu))$ with $B(x, 0.5 D^s) \cap J \ne \emptyset$. By the triangle inequality and since $x \in E(\fp_s) \subset B(\pc(\fp_s), 4D^{s})$, it follows that $B(\pc(\fp_s), 4.5D^s) \cap J \ne \emptyset$. If $s(J) \ge s$ and $s(J) > -S$, then it follows from the triangle inequality, \eqref{eq-vol-sp-cube} and \eqref{defineD} that $\scI(\fp_s) \subset B(c(J), 100 D^{s(J)+1})$, contradicting $J \in \mathcal{J}(\mathfrak{T}(\fu))$. Thus $s(J) \le s - 1$ or $s(J) = -S$. If $s(J) = -S$ and $s(J) > s - 1$, then $s = -S$. Thus we always have $s(J) \le s$. It then follows from the triangle inequality and \eqref{eq-vol-sp-cube} that $J \subset B(\pc(\fp_s), 16 D^s)$.

    Thus we can continue our chain of estimates with
    $$
        5\cdot 2^{103a^3} \sum_{s\in\sigma(x)}2^{s - \overline{\sigma}(\fu, x)} \frac{1}{\mu(B(x,D^s))}\int_{B(\pc(\fp_s),16 D^s)}P_{\mathcal{J}(\fT(\fu))}|f(y)|\,\mathrm{d}\mu(y)\,.
    $$
    We have $B(\pc(\fp_s), 16D^s)) \subset B(x, 32D^s)$, by \eqref{eq-vol-sp-cube} and the triangle inequality, since $x \in \scI(\fp_s)$. Combining this with the doubling property \eqref{doublingx}, we obtain
    $$
        \mu(B(\pc(\fp_s), 16D^s)) \le 2^{5a} \mu(B(x, D^s))\,.
    $$
    Since $a \ge 4$, it follows that \eqref{eq-term-A} is bounded by
    $$
        5\cdot 2^{103a^3} \sum_{s\in\sigma(x)}2^{s - \overline{\sigma}(\fu, x)} \frac{1}{\mu(B(\pc(\fp_s),16D^s))}\int_{B(\pc(\fp_s),16D^s)}P_{\mathcal{J}(\fT(\fu))}|f(y)|\,\mathrm{d}\mu(y)\,.
    $$
    Since $L \in \mathcal{L}(\fT(\fu))$ and $x\in L \cap \scI(\fp_s)$, we have $s(L) \le \ps(\fp_s)$. It follows by \eqref{dyadicproperty} that $L \subset \scI(\fp_s)$, in particular $x' \in \scI(\fp_s) \subset B(\pc(\fp_s), 16D^s)$. Thus
    $$
        \le 5\cdot 2^{104a^3} \sum_{s\in\sigma(x)}2^{s - \overline{\sigma}(\fu, x)} M_{\mathcal{B}, 1}P_{\mathcal{J}(\fT(\fu))}|f|(x')
    $$
    $$
        \le 10\cdot 2^{104a^3} M_{\mathcal{B}, 1}P_{\mathcal{J}(\fT(\fu))}|f|(x')\,.
    $$
    This completes the estimate for term \eqref{eq-term-A}.
\end{proof}

\begin{lemma}[second tree pointwise]
    \label{second-tree-pointwise}
    \leanok
    \lean{TileStructure.Forest.second_tree_pointwise}
    For all $\fu \in \fU$, all $L \in \mathcal{L}(\fT(\fu))$, all $x, x' \in L$ and all bounded $f$ with bounded support, we have
    $$
         \Bigg| \sum_{s \in \sigma(\fu, x)} \int K_s(x,y) P_{\mathcal{J}(\fT(\fu))} f(y) \, \mathrm{d}\mu(y) \Bigg| \le T_{\mathcal{N}}^{\fcc(\fu)} P_{\mathcal{J}(\fT(\fu))} f(x')\,.
    $$
\end{lemma}

\begin{proof}
    \leanok
    Let $s_1 = \underline{\sigma}(\fu, x)$. By definition, there exists a tile $\fp \in \fT(\fu)$ with $\ps(\fp) = s_1$ and $x \in E(\fp)$. Then $x \in \scI(\fp) \cap L$. By \eqref{dyadicproperty} and the definition of $\mathcal{L}(\fT(\fu))$, it follows that $L \subset \scI(\fp)$, in particular $x' \in \scI(\fp)$, so $x \in I_{s_1}(x')$.
    Next, let $s_2 = \overline{\sigma}(\fu, x)$ and let $\fp' \in \fT(\fu)$ with $\ps(\fp') = s_2$ and $x \in E(\fp')$. Since $\fp' \in \fT(\fu)$, we have $4\fp' \lesssim \fu$. Since $\tQ(x) \in \fc(\fp')$, it follows that
    $$
        d_{\fp}(\fcc(\fu), \tQ(x)) \le 5\,.
    $$
    Applying the doubling property \eqref{firstdb} five times, we obtain
    $$
        d_{B(c(\fp), 8D^{s_2})}(\fcc(\fu), \tQ(x)) \le 5 \cdot 2^{5a}\,.
    $$
    By the triangle inequality, we have $B(x, D^{s_2}) \subset B(c(\fp), 8 D^{s_2})$, so by \eqref{monotonedb}
    $$
        d_{B(x, D^{s_2})}(\fcc(\fu), \tQ(x)) \le 5 \cdot 2^{5a}\,.
    $$
    Finally, by applying \eqref{seconddb} $100a$ times, we obtain
    $$
        d_{B(x, D^{s_2-1})}(\fcc(\fu), \tQ(x)) \le 5 \cdot 2^{-95a} < 1\,.
    $$
    Consequently, $D^{s_2 - 1} < R_Q(\fcc(\fu), x)$.
    The lemma now follows from the definition of $T_{\mathcal{N}}$.
\end{proof}

\begin{lemma}[third tree pointwise]
    \label{third-tree-pointwise}
    \leanok
    \lean{TileStructure.Forest.third_tree_pointwise}
    For all $\fu \in \fU$, all $L \in \mathcal{L}(\fT(\fu))$, all $x, x' \in L$ and all bounded $f$ with bounded support, we have
    \begin{equation*}
        \Bigg| \sum_{s \in \sigma(\fu, x)} \int K_s(x,y) (f(y) - P_{\mathcal{J}(\fT(\fu))} f(y)) \, \mathrm{d}\mu(y) \Bigg|
    \end{equation*}
    \begin{equation*}
          \le 2^{128a^3} S_{1,\fu} P_{\mathcal{J}(\fT(\fu))}|f|(x')\,.
    \end{equation*}
\end{lemma}

\begin{proof}
    \leanok
    We have for $J \in \mathcal{J}(\fT(\fu))$:
    $$
        \int_J K_{s}(x,y)(1 - P_{\mathcal{J}(\fT(\fu))})f(y) \, \mathrm{d}\mu(y)
    $$
    \begin{equation}
    \label{eq-canc-comp}
        = \int_J \frac{1}{\mu(J)} \int_J K_s(x,y) - K_s(x,z) \, \mathrm{d}\mu(z) \,f(y) \, \mathrm{d}\mu(y)\,.
    \end{equation}
    By \eqref{eq-Ks-smooth} and \eqref{eq-vol-sp-cube}, we have for $y, z \in J$
    $$
        |K_s(x,y) - K_s(x,z)| \le \frac{2^{127a^3}}{\mu(B(x, D^s))} \left(\frac{8 D^{s(J)}}{D^s}\right)^{1/a}\,.
    $$
    Suppose that $s \in \sigma(\fu, x)$.
    If $K_s(x,y) \ne 0$ for some $y \in J \in \mathcal{J}(\fT(\fu))$ then, by \eqref{supp-Ks}, $y \in B(x, 0.5 D^s) \cap J \ne \emptyset$. Let $\fp \in \fT(\fu)$ with $\ps(\fp) = s$ and $x \in E(\fp)$. Then $B(\pc(\fp_s), 4.5D^s) \cap J \ne \emptyset$ by the triangle inequality. If $s(J) \ge s$ and $s(J) > -S$, then it follows from the triangle inequality, \eqref{eq-vol-sp-cube} and \eqref{defineD} that $\scI(\fp) \subset B(c(J), 100 D^{s(J)+1})$, contradicting $J \in \mathcal{J}(\mathfrak{T}(\fu))$. Thus $s(J) \le s - 1$ or $s(J) = -S$. If $s(J) = -S$ and $s(J) > s - 1$, then $s = -S$. So in both cases, $s(J) \le s$. It then follows from the triangle inequality and \eqref{eq-vol-sp-cube} that $J \subset B(x, 16 D^s)$.

    Thus, we can estimate \eqref{eq-term-C} by
    $$
        2^{127a^3 + 3/a}\sum_{\fp\in \mathfrak{T}}\frac{\mathbf{1}_{E(\fp)}(x)}{\mu(B(x,D^{\ps(\fp)}))}\sum_{\substack{J\in \mathcal{J}(\fT(\fu))\\J\subset B(x, 16D^{\ps(\fp)})\\ s(J) \le \ps(\fp)}} D^{(s(J) - \ps(\fp))/a} \int_J |f|\,.
    $$
    $$
        = 2^{127a^3 + 3/a}\sum_{I \in \mathcal{D}} \sum_{\substack{\fp\in \mathfrak{T}\\ \scI(\fp) = I}}\frac{\mathbf{1}_{E(\fp)}(x)}{\mu(B(x, D^{s(I)}))}\sum_{\substack{J\in \mathcal{J}(\fT(\fu))\\J\subset B(x, 16 D^{s(I)})\\ s(J) \le s(I)}} D^{(s(J) - s(I))/a} \int_J |f|\,.
    $$
    By \eqref{eq-dis-freq-cover} and \eqref{defineep}, the sets $E(\fp)$ for tiles $\fp$ with $\scI(\fp) = I$ are pairwise disjoint. It follows from the definition of $\mathcal{L}(\fT(\fu))$ that $x \in \scI(\fp)$ if and only if $x' \in \scI(\fp)$, thus we can estimate the sum over $\mathbf{1}_{E(\fp)}(x)$ by $\mathbf{1}_{I}(x')$.
    If $x \in E(\fp)$ then in particular $x \in \scI(\fp)$, so by \eqref{eq-vol-sp-cube} $B(c(I),16D^{s(I)}) \subset B(x, 32D^{s(I)})$. By the doubling property \eqref{doublingx}
    $$
        \mu(B(c(I), 16D^{s(I)})) \le 2^{5a} \mu(B(x, D^{s(I)}))\,.
    $$
    Since $a \ge 4$ we can continue our estimate with
    $$
        \le 2^{128a^3}\sum_{I \in \mathcal{D}} \frac{\mathbf{1}_{I}(x')}{\mu(B(c(I), 16D^{s(I)}))}\sum_{\substack{J\in \mathcal{J}(\fT(\fu))\\J\subset B(x, 16 D^{s(I)})\\ s(J) \le s(I)}} D^{(s(J) - s(I))/a} \int_J |f|
    $$
    $$
         = 2^{128a^3} S_{1,\fu} P_{\mathcal{J}(\fT(\fu))}|f|(x')\,.
    $$
    This completes the proof.
\end{proof}

\subsection{An auxiliary \texorpdfstring{$L^2$}{L2} tree estimate}

In this subsection we prove the following estimate on $L^2$ for operators associated to trees.

\begin{lemma}[tree projection estimate]
    \label{tree-projection-estimate}
    \leanok
    \lean{TileStructure.Forest.tree_projection_estimate}
    \uses{dyadic-partitions,pointwise-tree-estimate,nontangential-operator-bound,boundary-operator-bound}
    Let $\fu \in \fU$.
    Then we have for all $f, g$ bounded with bounded support
    $$
        \Bigg|\int_X \sum_{\fp \in \fT(\fu)} \bar g(y) T_{\fp}f(y) \, \mathrm{d}\mu(y) \Bigg|
    $$
    \begin{equation}
        \label{eq-tree-est}
         \le 2^{130a^3}\|P_{\mathcal{J}(\fT(\fu))}|f|\|_{2}\|P_{\mathcal{L}(\fT(\fu))}|g|\|_{2}.
    \end{equation}
\end{lemma}

Below, we deduce \Cref{tree-projection-estimate} from \Cref{pointwise-tree-estimate} and the following estimates for the operators in \Cref{pointwise-tree-estimate}.

\begin{lemma}[nontangential operator bound]
    \label{nontangential-operator-bound}
    \leanok
    \lean{TileStructure.Forest.nontangential_operator_bound}
    \uses{Hardy-Littlewood}
    For all bounded $f$ with bounded support and all $\mfa \in \Mf$
    $$
        \|T_{\mathcal{N}}^{\mfa} f\|_2 \le 2^{102a^3} \|f\|_2\,.
    $$
\end{lemma}

\begin{lemma}[boundary operator bound]
    \label{boundary-operator-bound}
    \leanok
    \lean{TileStructure.Forest.boundary_operator_bound}
    \uses{Hardy-Littlewood,boundary-overlap}
    For all $\fu \in \fU$ and all bounded functions $f$ with bounded support
    \begin{equation}
        \label{eq-S-bound}
        \|S_{1,\fu}f\|_2 \le 2^{12a} \|f\|_2\,.
    \end{equation}
\end{lemma}

\begin{proof}[Proof of \Cref{tree-projection-estimate}]
    \proves{tree-projection-estimate}\leanok
    Let $L \in \mathcal{L}(\fT(\fu))$.
    Let $b(x')$ denote the right-hand side of \Cref{pointwise-tree-estimate}. Apply this lemma to $e(\fcc(\fu)) f$, to obtain for all $y, x' \in L$
    $$
        \Bigg| \sum_{\fp \in \fT(\fu)} T_{\fp} f(y) \Bigg| \le b(x').
    $$
    Hence, by taking an infimum, we have for $y \in L$
    $$
        \Bigg| \sum_{\fp \in \fT(\fu)} T_{\fp} f(y) \Bigg| \le \inf_{x' \in L} b(x').
    $$
    Integrating this estimate yields
    $$
        \int_L |g(y)| \Bigg| \sum_{\fp \in \fT(\fu)} T_{\fp} f(y) \Bigg| \, \mathrm{d}\mu(y)
    $$
    $$
        \le \int_L |g(y)| \times \inf_{x' \in L} b(x') \, \mathrm{d}\mu(y)
    $$
    $$
        =   \int_L P_{\mathcal{L}(\fT(\fu))}|g|(y) \times \inf_{x' \in L} b(x') \, \mathrm{d}\mu(y)
    $$
    $$
        \le \int_L P_{\mathcal{L}(\fT(\fu))}|g|(y) \times  b(y) \, \mathrm{d}\mu(y)
    $$
    By \eqref{definetp}, we have $T_{\fp} f = \mathbf{1}_{\scI(\fp)} T_{\fp} f$ for all $\fp \in \fP$, so
    $$
        \Bigg| \int \bar g(y) \sum_{\fp \in \fT(\fu)} T_{\fp} f(y) \, \mathrm{d}\mu(y) \Bigg| = \Bigg| \int_{\bigcup_{\fp \in \fT(\fu)} \scI(\fp)} \bar g(y) \sum_{\fp \in \fT(\fu)} T_{\fp} f(y) \, \mathrm{d}\mu(y) \Bigg|\,.
    $$
    Since $\mathcal{L}(\fT(\fu))$ partitions $\bigcup_{\fp \in \fT(\fu)} \scI(\fp)$ by \Cref{dyadic-partitions},
    we get from the triangle inequality
    $$
        \le \sum_{L \in \mathcal{L}(\fT(\fu))} \int_L |g(y)| \Bigg| \sum_{\fp \in \fT(\fu)} T_{\fp} f(y) \Bigg| \, \mathrm{d}\mu(y)
    $$
    which by the above computation is bounded by
    $$
        \sum_{L \in \mathcal{L}(\fT(\fu))} \int_L P_{\mathcal{L}(\fT(\fu))}|g|(y) \times  b(y) \, \mathrm{d}\mu(y)
    $$
    $$
        = \int_X P_{\mathcal{L}(\fT(\fu))}|g|(y) \times  b(y) \, \mathrm{d}\mu(y)
    $$
    Applying Cauchy-Schwarz, this is bounded by $\|P_{\mathcal{L}(\fT(\fu))}|g|\|_2 \times \|b\|_2$.
    By Minkowski's inequality, \Cref{Hardy-Littlewood}, \Cref{nontangential-operator-bound} and \Cref{boundary-operator-bound}, $\|b\|_2$ is at most
    $$
        2^{129a^3} (2^{2a+1} + 2^{12a})\|P_{\mathcal{J}(\fT(\fu))}|f|\|_2 + 2^{103a^3} \|P_{\mathcal{J}(\fT(\fu))}[e(\fcc(\fu))f]\|_2\,.
    $$
    By the triangle inequality we have for all $x \in X$ that $|P_{\mathcal{J}(\fT(\fu))}[e(\fcc(\fu))f]|(x) \le P_{\mathcal{J}(\fT(\fu))}|f|(x)$, thus we can further estimate the above by
    $$
        (2^{129a^3} (2^{2a+1} + 2^{12a}) + 2^{103a^3}) \|P_{\mathcal{J}(\fT(\fu))}|f|\|_2\,.
    $$
    This completes the proof since $a \ge 4$.
\end{proof}

Now we prove the two auxiliary lemmas. We begin with the nontangential maximal operator $T_{\mathcal{N}}$.

\begin{proof}[Proof of \Cref{nontangential-operator-bound}]
    \leanok
    \proves{nontangential-operator-bound}
    Fix $s_1, s_2$. By \eqref{eq-psisum} we have for all $x \in (0, \infty)$
    $$
        \sum_{s = s_1}^{s_2} \psi(D^{-s}x) = 1 - \sum_{s < s_1} \psi(D^{-s}x) - \sum_{s > s_1} \psi(D^{-s}x)\,.
    $$
    Since $\psi$ is supported in $[\frac{1}{4D}, \frac{1}{2}]$, the two sums on the right hand side are zero for all $x \in [\frac{1}{2}D^{s_1-1}, \frac{1}{4} D^{s_2 - 1}]$, hence
    $$
        x \in [\frac{1}{2}D^{s_1-1}, \frac{1}{4} D^{s_2}] \implies \sum_{s = s_1}^{s_2} \psi(D^{-s}x) = 1\,.
    $$
    Since $\psi$ is supported in $[\frac{1}{4D}, \frac{1}{2}]$, we further have
    $$
        x \notin [\frac{1}{4}D^{s_1 - 1}, \frac{1}{2}D^{s_2}] \implies \sum_{s = s_1}^{s_2} \psi(D^{-s}x) = 0\,.
    $$
    Finally, since $\psi \ge 0$ and $\sum_{s \in \mathbb{Z}} \psi(D^{-s}x) = 1$, we have for all $x$
    $$
        0 \le \sum_{s = s_1}^{s_2} \psi(D^{-s}x) \le 1\,.
    $$
    Let $x' \in I_{s_1}(x)$ and suppose that $D^{s_2 - 1} \le R_Q(\mfa, x')$. By the triangle inequality and \eqref{eq-vol-sp-cube}, it holds that $\rho(x,x') \le 8D^{s_1}$. We have
    $$
        \Bigg|\sum_{s = s_1}^{s_2} \int K_s(x',y) f(y) \, \mathrm{d}\mu(y)\Bigg|
    $$
    $$
        = \Bigg|\int \sum_{s = s_1}^{s_2} \psi(D^{-s}\rho(x',y)) K(x',y) f(y) \, \mathrm{d}\mu(y)\Bigg|
    $$
    \begin{equation}
        \label{eq-sharp-trunc-term}
        \le \Bigg| \int_{8D^{s_1} < \rho(x',y) \le \frac{1}{4}D^{s_2-1}} K(x',y) f(y) \, \mathrm{d}\mu(y) \Bigg|
    \end{equation}
    \begin{equation}
        \label{eq-lower-bound-term}
        + \int_{\frac{1}{4}D^{s_1-1} \le \rho(x',y) \le 8D^{s_1}} |K(x', y)| |f(y)| \, \mathrm{d}\mu(y)
    \end{equation}
    \begin{equation}
        \label{eq-upper-bound-term}
        + \int_{\frac{1}{4}D^{s_2-1} \le \rho(x',y) \le \frac{1}{2}D^{s_2}} |K(x', y)| |f(y)| \, \mathrm{d}\mu(y)\,.
    \end{equation}

    The first term \eqref{eq-sharp-trunc-term} is at most $2T_{\tQ}^\mfa f(x)$, using with $R_1 := 8D^{s_1}$, $R_2 := \frac{1}{4}D^{s_2-1}$ and $R_1 < R_2 < R_{\tQ}(\mfa,x')$ the triangle inequality in the form
    \begin{equation}
        \left|\int_{R_1 < \rho(x',y) \le R_2} K(x',y) f(y) \, \mathrm{d}\mu(y) \right|
    \end{equation}
    \begin{equation}
        \le \left|\int_{R_1 < \rho(x',y) < R_{\tQ}(\mfa,x')} K(x',y) f(y) \, \mathrm{d}\mu(y) \right|
    \end{equation}
    \begin{equation}
        + \left|\int_{R_2 < \rho(x',y) < R_{\tQ}(\mfa,x')} K(x',y) f(y) \, \mathrm{d}\mu(y) \right|.
    \end{equation}

    The other two terms will be estimated by the finitary maximal function from \Cref{Hardy-Littlewood}.
    For the second term \eqref{eq-lower-bound-term} we use \eqref{eqkernel-size} which implies that for all $y$ with $\rho(x', y) \ge \frac{1}{4}D^{s_1 - 1}$, we have
    $$
        |K(x', y)| \le \frac{2^{a^3}}{\mu(B(x', \frac{1}{4}D^{s_1 - 1}))}\,.
    $$
    Using $D=2^{100a^2}$
    and the doubling property \eqref{doublingx} $7 +100a^2$ times estimates
    the last display by
    \begin{equation}
    \label{pf-nontangential-operator-bound-imeq}
        \le \frac{2^{7a+101a^3}}{\mu(B(x', 32D^{s_1}))}\, .
    \end{equation}
    By the triangle inequality and \eqref{eq-vol-sp-cube}, we have
    $$
        B(x', 8D^{s_1}) \subset B(c(I_{s_1}(x)), 16D^{s(I_{s_1}(x))}) \subset B(x', 32D^{s_1})\,.
    $$
    Combining this with \eqref{pf-nontangential-operator-bound-imeq}, we conclude that \eqref{eq-lower-bound-term} is at most
    $$
        2^{7a + 101a^3} M_{\mathcal{B},1} f(x)\,.
    $$

    For \eqref{eq-upper-bound-term} we argue similarly. We have for all $y$ with $\rho(x', y) \ge \frac{1}{4}D^{s_2-1}$
    $$
        |K(x', y)| \le \frac{2^{a^3}}{\mu(B(x', \frac{1}{4}D^{s_2-1}))}\,.
    $$
    Using the doubling property \eqref{doublingx} $7 + 100a^2$ times estimates
    the last display by
    \begin{equation}
        \le \frac{2^{7a + 101a^3}}{\mu(B(x', 32 D^{s_2}))}\, .
    \end{equation}
    Note that by \eqref{dyadicproperty} we have $I_{s_1}(x) \subset I_{s_2}(x)$, in particular $x' \in I_{s_2}(x)$.
    By the triangle inequality and \eqref{eq-vol-sp-cube}, we have
    $$
        B(x', 8D^{s_2}) \subset B(c(I_{s_2}(x)), 16D^{s(I_{s_2}(x))}) \subset B(x', 32D^{s_2})\,.
    $$
    Combining this, \eqref{eq-upper-bound-term} is at most
    $$
        2^{7a+101a^3} M_{\mathcal{B},1} f(x)\,.
    $$

    Using $a \ge 4$, taking a supremum over all $x' \in I_{s_1}(x)$ and then a supremum over all $-S \le s_1 < s_2 \le S$, we obtain
    $$
        T_{\mathcal{N}} f(x) \le 2T_{\tQ}^\mfa f(x) + 2^{102a^3} M_{\mathcal{B},1} f(x)\,.
    $$
    The lemma now follows from assumption \eqref{linnontanbound}, \Cref{Hardy-Littlewood} and $a \ge 4$.
\end{proof}

We need the following lemma to prepare the $L^2$-estimate for the auxiliary operators $S_{1, \fu}$.

\begin{lemma}[boundary overlap]
    \label{boundary-overlap}
    \leanok
    \lean{TileStructure.Forest.boundary_overlap}
    For every cube $I \in \mathcal{D}$, there exist at most $2^{9a}$ cubes $J \in \mathcal{D}$ with $s(J) = s(I)$ and $B(c(I), 16D^{s(I)}) \cap B(c(J), 16 D^{s(J)}) \ne \emptyset$.
\end{lemma}

\begin{proof}
    \leanok
    Suppose that $B(c(I), 16 D^{s(I)}) \cap B(c(J), 16 D^{s(J)}) \ne \emptyset$ and $s(I) = s(J)$. Then $B(c(I), 33 D^{s(I)}) \subset B(c(J), 128 D^{s(J)})$. Hence by the doubling property \eqref{doublingx}
    $$
        2^{9a}\mu(B(c(J), \frac{1}{4}D^{s(J)})) \ge \mu(B(c(I), 33 D^{s(I)}))\,,
    $$
    and by the triangle inequality, $B(c(J), \frac{1}{4}D^{s(J)})$ is contained in $B(c(I), 33 D^{s(I)})$.

    If $\mathcal{C}$ is any finite collection of cubes $J \in \mathcal{D}$ satisfying $s(J) = s(I)$ and
    \begin{equation*}
        B(c(I), 16 D^{s(I)}) \cap B(c(J), 16 D^{s(J)}) \ne\emptyset\ ,
    \end{equation*} then it follows from \eqref{eq-vol-sp-cube} and pairwise disjointedness of cubes of the same scale \eqref{dyadicproperty} that the balls $B(c(J), \frac{1}{4} D^{s(J)})$ are pairwise disjoint. Hence
    \begin{align*}
        \mu(B(c(I), 33 D^{s(I)})) &\ge \sum_{J \in \mathcal{C}} \mu(B(c(J), \frac{1}{4}D^{s(J)}))\\
        &\ge |\mathcal{C}| 2^{-9a} \mu(B(c(I), 33 D^{s(I)}))\,.
    \end{align*}
    Since $\mu$ is doubling and $\mu \ne 0$, we have $\mu(B(c(I), 33D^{s(I)})) > 0$. The lemma follows after dividing by $2^{-9a}\mu(B(c(I), 33D^{s(I)}))$.
\end{proof}

Now we can bound the operators $S_{1, \fu}$.

\begin{proof}[Proof of \Cref{boundary-operator-bound}]
    \leanok
    \proves{boundary-operator-bound}
    Note that by definition, $S_{1,\fu}f$ is a finite sum of indicator functions of cubes $I \in \mathcal{D}$ for each locally integrable $f$, and hence is bounded, has bounded support and is integrable. Let $g$ be another function with the same three properties. Then $\bar g S_{1,\fu}f$ is integrable, and we have
    $$
        \Bigg|\int \bar g(y) S_{1,\fu}f(y) \, \mathrm{d}\mu(y)\Bigg|
    $$
    \begin{multline*}
        = \Bigg|\sum_{I\in\mathcal{D}} \frac{1}{\mu(B(c(I), 16 D^{s(I)}))} \int_I \bar g(y) \, \mathrm{d}\mu(y)\\
        \times \sum_{\substack{J\in \mathcal{J}(\fT(\fu))\,:\,J\subseteq B(c(I), 16 D^{s(I)})\\s(J)\le s(I)}} D^{(s(J)-s(I))/a}\int_J |f(y)| \,\mathrm{d}\mu(y)\Bigg|
    \end{multline*}
    \begin{multline*}
        \le \sum_{I\in\mathcal{D}} \frac{1}{\mu(B(c(I), 16D^{s(I)}))} \int_{B(c(I), 16D^{s(I)})} | g(y)| \, \mathrm{d}\mu(y)\\ \times \sum_{\substack{J\in \mathcal{J}(\fT(\fu))\,:\,J\subseteq B(c(I), 16 D^{s(I)})\\ s(J) \le s(I)}} D^{(s(J)-s(I))/a}\int_J |f(y)| \,\mathrm{d}\mu(y)\,.
    \end{multline*}
    Changing the order of summation and using $J \subset B(c(I), 16 D^{s(I)})$ to bound the first average integral by $M_{\mathcal{B},1}|g|(y)$ for any $y \in J$, we obtain
    \begin{align}
    \label{eq-boundary-operator-bound-1}
        \le \sum_{J\in\mathcal{J}(\fT(\fu))}\int_J|f(y)| M_{\mathcal{B},1}|g|(y) \, \mathrm{d}\mu(y) \sum_{\substack{I \in \mathcal{D} \, : \, J\subset B(c(I),16 D^{s(I)})\\ s(I) \ge s(J)}} D^{(s(J)-s(I))/a}.
    \end{align}
    By \Cref{boundary-overlap}, there are at most $2^{9a}$ cubes $I$ at each scale satisfying the inclusion $J \subset B(c(I),16D^{s(I)})$.
    Since $D^{-1/a}\le\frac12$, $(1 - D^{-1/a})^{-1} \le 2$.
    Using this estimate for the sum of the geometric series, we conclude that \eqref{eq-boundary-operator-bound-1} is at most
    $$
        2^{9a+1} \sum_{J\in\mathcal{J}(\fT(\fu))}\int_J|f(y)| M_{\mathcal{B},1}|g|(y) \, \mathrm{d}\mu(y)\,.
    $$
    The collection $\mathcal{J}$ is a partition of $\bigcup_{I \in \mathcal{D}} I$, so this equals
    $$
        2^{9a+1} \int_X|f(y)| M_{\mathcal{B},1}|g|(y) \, \mathrm{d}\mu(y)\,.
    $$
    Using Cauchy-Schwarz and \Cref{Hardy-Littlewood}, we conclude
    $$
        \left|\int \bar g S_{1,\fu}f \, \mathrm{d}\mu \right| \le 2^{12a} \|g\|_2\|f\|_2\,.
    $$
    The lemma now follows by choosing $g = S_{1,\fu}f$ and dividing on both sides by the finite $\|S_{1,\fu}f\|_2$.
\end{proof}

\subsection{The quantitative \texorpdfstring{$L^2$}{L2} tree estimate}

The main result of this subsection is the following quantitative bound for operators associated to trees, with decay in the densities $\dens_1$ and $\dens_2$.

\begin{lemma}[densities tree bound]
    \label{densities-tree-bound}
    \leanok
    \lean{TileStructure.Forest.density_tree_bound1, TileStructure.Forest.density_tree_bound2}
    \uses{tree-projection-estimate,local-dens1-tree-bound,local-dens2-tree-bound}
    Let $\fu \in \fU$. Then for all bounded $f$ with bounded support and $g$ with $\text{support}(g) \subseteq G$
    we have
    \begin{equation}
        \label{eq-cor-tree-est}
        \left|\int_X \bar g \sum_{\fp \in \fT(\fu)} T_{\fp }f \, \mathrm{d}\mu \right| \le 2^{181a^3} \dens_1(\fT(\fu))^{1/2} \|f\|_2\|g\|_2\,.
    \end{equation}
    If additionally $\text{support}(f) \subseteq F$, then we have
    \begin{equation}
        \label{eq-cor-tree-est-F}
        \left| \int_X \bar g \sum_{\fp \in \fT(\fu)} T_{\fp }f\, \mathrm{d}\mu \right| \le 2^{282a^3} \dens_1(\fT(\fu))^{1/2} \dens_2(\fT(\fu))^{1/2} \|f\|_2\|g\|_2\,.
    \end{equation}
\end{lemma}

Below, we deduce this lemma from \Cref{tree-projection-estimate} and the following two estimates controlling the size of support of the operator and its adjoint.

\begin{lemma}[local dens1 tree bound]
    \label{local-dens1-tree-bound}
    \leanok
    \lean{TileStructure.Forest.local_dens1_tree_bound}
    \uses{monotone-cube-metrics}
    Let $\fu \in \fU$ and $L \in \mathcal{L}(\fT(\fu))$. Then
    \begin{equation}
    \label{eq-1density-estimate-tree}
        \mu(L \cap G \cap \bigcup_{\fp \in \fT(\fu)} E(\fp)) \le 2^{101a^3} \dens_1(\fT(\fu)) \mu(L)\,.
    \end{equation}
\end{lemma}

\begin{lemma}[local dens2 tree bound]
    \label{local-dens2-tree-bound}
    \leanok
    \lean{TileStructure.Forest.local_dens2_tree_bound}
    Let $\fu \in \fU$ and $J \in \mathcal{J}(\fT(\fu))$. Then
    $$
        \mu(F \cap J) \le 2^{201a^3} \dens_2(\fT(\fu)) \mu(J)\,.
    $$
\end{lemma}

\begin{proof}[Proof of \Cref{densities-tree-bound}]
    \proves{densities-tree-bound}
    \leanok

    Denote
    $$
        \mathcal{E}(\fu) = \bigcup_{\fp \in \fT(\fu)} E(\fp)\,.
    $$
    Then we have
    $$
        \left| \int_X \bar g \sum_{\fp \in \fT(\fu)} T_{\fp} f \, \mathrm{d}\mu \right| = \left| \int_X \overline{ g\mathbf{1}_{\mathcal{E}(\fu)}} \sum_{\fp \in \fT(\fu)} T_{\fp} (\mathbf{1}_{\scI(\fu)}f) \, \mathrm{d}\mu \right|\,.
    $$
    By \Cref{tree-projection-estimate}, this is bounded by
    \begin{equation}
        \label{eq-both-factors-tree}
        \le 2^{130a^3}\|P_{\mathcal{J}(\fT(\fu))}|f|\|_2 \|P_{\mathcal{L}(\fT(\fu))} |\mathbf{1}_{\mathcal{E}(\fu)}g|\|_2\,.
    \end{equation}
    We bound the two factors separately.
    We have
    $$
        \|P_{\mathcal{L}(\fT(\fu))} |\mathbf{1}_{\mathcal{E}(\fu)}g|\|_2 = \left( \sum_{L \in \mathcal{L}(\fT(\fu))} \frac{1}{\mu(L)} \left(\int_{L \cap \mathcal{E}(\fu)} |g(y)| \, \mathrm{d}\mu(y)\right)^2 \right)^{1/2}\,.
    $$
    By Cauchy-Schwarz and \Cref{local-dens1-tree-bound} this is at most
    $$
        \le \left( \sum_{L \in \mathcal{L}(\fT(\fu))} 2^{101a^3} \dens_1(\fT(\fu)) \int_{L \cap \mathcal{E}(\fu)} |g(y)|^2 \, \mathrm{d}\mu(y) \right)^{1/2}\,.
    $$
    Since cubes $L \in \mathcal{L}(\fT(\fu))$ are pairwise disjoint by \Cref{dyadic-partitions}, this is
    \begin{equation}
        \label{eq-factor-L-tree}
         \le 2^{51 a^3} \dens_1(\fT(\fu))^{1/2} \|g\|_2\,.
    \end{equation}
    Similarly, we have
    \begin{equation}
        \label{eq-cor-tree-proof}
        \|P_{\mathcal{J}(\fT(\fu))}|f|\|_2 = \left( \sum_{J \in \mathcal{J}(\fT(\fu))} \frac{1}{\mu(J)} \left(\int_J |f(y)| \, \mathrm{d}\mu(y)\right)^2 \right)^{1/2}\,.
    \end{equation}
    By Cauchy-Schwarz, this is
    $$
        \le \left( \sum_{J \in \mathcal{J}(\fT(\fu))} \int_J |f(y)|^2 \, \mathrm{d}\mu(y) \right)^{1/2}\,.
    $$
    Since cubes in $\mathcal{J}(\fT(\fu))$ are pairwise disjoint by \Cref{dyadic-partitions}, this is at most
    \begin{equation}
        \label{eq-factor-J-tree}
        \|f\|_2\,.
    \end{equation}
    Combining \eqref{eq-both-factors-tree}, \eqref{eq-factor-L-tree} and \eqref{eq-factor-J-tree} gives \eqref{eq-cor-tree-est}.

    If $f \le \mathbf{1}_F$ then $f = f\mathbf{1}_F$, so
    $$
        \left( \sum_{J \in \mathcal{J}(\fT(\fu))} \frac{1}{\mu(J)} \left(\int_J |f(y)| \, \mathrm{d}\mu(y)\right)^2 \right)^{1/2}
    $$
    $$
        =\left(\sum_{J \in \mathcal{J}(\fT(\fu))} \frac{1}{\mu(J)} \left(\int_{J \cap F} |f(y)| \, \mathrm{d}\mu(y)\right)^2 \right)^{1/2}\,.
    $$
    We estimate as before, using now \Cref{local-dens2-tree-bound} and Cauchy-Schwarz, and obtain that this is
    $$
        \le 2^{101a^3} \dens_2(\fT(\fu))^{1/2} \|f\|_2\,.
    $$
    Combining this with \eqref{eq-both-factors-tree} and \eqref{eq-factor-L-tree} gives \eqref{eq-cor-tree-est-F}.
\end{proof}

Now we prove the two auxiliary estimates.

\begin{proof}[Proof of \Cref{local-dens1-tree-bound}]
    \leanok
    \proves{local-dens1-tree-bound}
    If the set on the right hand side is empty, then \eqref{eq-1density-estimate-tree} holds. If not, then there exists $\fp \in \fT(\fu)$ with $L \cap \scI(\fp) \ne \emptyset$.

    Suppose first that there exists such $\fp$ with $\ps(\fp) \le s(L)$. Then by \eqref{dyadicproperty} $\scI(\fp) \subset L$, which gives by the definition of $\mathcal{L}(\fT(\fu))$ that $s(L) = -S$ and hence $L = \scI(\fp)$. Let $\fq \in \fT(\fu)$ with $E(\fq) \cap L \ne \emptyset$. Since $s(L) = -S \le \ps(\fq)$ it follows from \eqref{dyadicproperty} that $\scI(\fp) = L \subset \scI(\fq)$. We have then by \Cref{monotone-cube-metrics}
    \begin{align*}
        d_{\fp}(\fcc(\fp), \fcc(\fq)) &\le d_{\fp}(\fcc(\fp), \fcc(\fu)) + d_{\fp}(\fcc(\fq), \fcc(\fu))\\
        &\le d_{\fp}(\fcc(\fp), \fcc(\fu)) + d_{\fq}(\fcc(\fq), \fcc(\fu))\,.
    \end{align*}
    Using that $\fp, \fq \in \fT(\fu)$ and \eqref{forest1}, this is at most $8$. Using again the triangle inequality and \Cref{monotone-cube-metrics}, we obtain that for each $q \in B_{\fq}(\fcc(\fq), 1)$
    $$
        d_{\fp}(\fcc(\fp), q) \le d_{\fp}(\fcc(\fp), \fcc(\fq)) + d_{\fq}(\fcc(\fq), q) \le 9\,.
    $$
    Thus $L \cap G \cap E(\fq) \subset E_2(9, \fp)$. We obtain
    $$
        \mu(L \cap G \cap \bigcup_{\fq \in \fT(\fu)} E(\fq)) \le \mu(E_2(9, \fp))\,.
    $$
    By the definition of $\dens_1$, this is bounded by
    $$
        9^a \dens_1(\fT(\fu)) \mu(\scI(\fp)) =9^a \dens_1(\fT(\fu)) \mu(L)\,.
    $$
    Since $a \ge 4$, \eqref{eq-1density-estimate-tree} follows in this case.

    Now suppose that for each $\fp \in \fT(\fu)$ with $L \cap E(\fp) \ne \emptyset$, we have $\ps(\fp) > s(L)$.
    Since there exists at least one such $\fp$, there exists in particular at least one cube $L'' \in \mathcal{D}$ with $L \subset L''$ and $s(L'') > s(L)$.
    By \eqref{coverdyadic}, there exists $L' \in \mathcal{D}$ with $L \subset L'$ and $s(L') = s(L) + 1$.
    By the definition of $\mathcal{L}(\fT(\fu))$ there exists a tile $\fp'' \in \fT(\fu)$ with $\scI(\fp'') \subset L'$.

    It now suffices to show that there exists a tile $\fp' \in \fP(\fT(\fu))$ with $\scI(\fp') = L'$, $d_{\fp'}(\fcc(\fp'), \fcc(\fu)) < 4$ and $9\fp'' \lesssim 9\fp'$.
    For then, let $\fq \in \fT(\fu)$ with $L \cap E(\fq) \ne \emptyset$. As shown above, this implies $\ps(\fq) \ge s(L')$, so by \eqref{dyadicproperty} $L' \subset \scI(\fq)$.
    If $q \in B_{\fq}(\fcc(\fq), 1)$, then by a similar calculation as above, using the triangle inequality, \Cref{monotone-cube-metrics} and \eqref{forest1}, we obtain
    $$
        d_{\fp'}(\fcc(\fp'), q) \le d_{\fp'}(\fcc(\fp'), \fcc(\fq)) + d_{\fq}(\fcc(\fq), q) \le 9\,.
    $$
    Thus $L \cap G \cap E(\fq) \subset E_2(9, \fp')$. We deduce using the definition \eqref{definedens1} of $\dens_1$
    $$
        \mu(L \cap G \cap \bigcup_{\fq \in \fT(\fu)} E(\fq)) \le \mu(E_2(9, \fp')) \le 9^a \dens_1(\fT(\fu)) \mu(L')\,.
    $$
    Using the doubling property \eqref{doublingx}, \eqref{eq-vol-sp-cube}, and $a \ge 4$ this is estimated by
    $$
        9^a 2^{100a^3 + 5a}\dens_1(\fT(\fu)) \mu(L) \le 2^{101 a^3} \dens_1(\fT(\fu))\mu(L)\,.
    $$
    To show existence of $\fp'$ with the given properties, if $\scI(\fp'') = L'$ we can take $\fp' = \fp''$, which satisfies the distance property by \eqref{forest1} and the other properties trivially.
    Otherwise, let $\fp'$ be the unique tile such that $\scI(\fp') = L'$ and such that $\Omega(\fu) \cap \Omega(\fp') \ne \emptyset$.
    Since $\scI(\fp') \subset \scI(\fp)$ and $\fp \in \fT(\fu)$, we have $\fp' \in \fP(\fT(\fu))$.
    Since by \eqref{forest1} $\ps(\fp') = s(L') \le \ps(\fp) < \ps(\fu)$, we have by \eqref{dyadicproperty} and \eqref{eq-freq-dyadic} that $\Omega(\fu) \subset \Omega(\fp')$, and hence the distance property.
    $9\fp'' \lesssim 9\fp'$ follows by the triangle inequality, \eqref{forest1}, \Cref{monotone-cube-metrics} and \eqref{eq-freq-comp-ball}.
    This completes the proof.
\end{proof}

\begin{proof}[Proof of \Cref{local-dens2-tree-bound}]
    \proves{local-dens2-tree-bound}
    \leanok
    We prove the inequality with the constant $2^{201a^3}$ replaced by $2 ^ {200a^3 + 14a}$; this
    is stronger because $a \geq 4$. It suffices to show the existence of a tile $\fp \in \fT(\fu)$
    and an $r \geq 4 D ^ {\ps(\fp)}$ such that $J \subset B(\pc(\fp), r)$ and
    $\mu(B(\pc(\fp), r)) \le 2 ^ {200a^3 + 14a} \mu(J)$, because then it follows from the definition
    \eqref{definedens2} of $\dens_2$ that
    $$
        \mu(F \cap J) \le \mu(F \cap B(\pc(\fp), r))
    $$
    $$
        \le \dens_2(\fT(\fu)) \mu(B(\pc(\fp), r)) \le 2^{200a^3 + 14a} \dens_2(\fT(\fu))\mu(J)\,.
    $$

    In particular, these criteria are satisfied, with $r = 4 D ^ {\ps(\fp)}$, by any
    $\fp \in \fT(\fu)$ such that $J \subseteq B(\pc(\fp, 4 D ^ {\ps(\fp)}))$ and
    $\mu(\scI(\fp)) \le 2^{100a^3 + 10a} \mu(J)$, because then by the doubling property
    \eqref{doublingx},
    $$
        \mu(B(\pc(\fp), 4 D ^ {\ps(\fp)})) \le 2^{4a} \mu(B(\pc(\fp), D ^ {\ps(\fp)} / 4))
    $$
    $$
        \le 2^{4a} \mu(\scI(\fp)) \le 2^{100a^3 + 14a} \mu(J)\,.
    $$

    Suppose first that $s(J) = S$. Then $J = I_0$, so \eqref{subsetmaxcube} and the fact that
    $J \in \mathcal{J}(\fT(\fu)) \subseteq \mathcal{J}_0(\fT(\fu))$ imply that $s(J) = -S$. Thus
    $S = 0$. It follows that $J$ is the only dyadic cube, so any $\fp \in \fT(\fu)$ has
    $\scI(\fp) = J$, and therefore satisfies $J \subseteq B(\pc(\fp, 4 D ^ {\ps(\fp)}))$ and
    $\mu(\scI(\fp)) \le 2^{100a^3 + 10a} \mu(J)$.

    It remains to consider the case $s(J) < S$. Then, by \eqref{coverdyadic} and
    \eqref{dyadicproperty}, there exists some cube $J' \in \mathcal{D}$ with $s(J') = s(J) + 1$ and
    $J \subset J'$. By definition of $\mathcal{J}(\fT(\fu))$ there exists some $\fp \in \fT(\fu)$
    such that $\scI(\fp) \subset B(c(J'), 100 D^{s(J') + 1})$.

    Since $c(J) \in J \subset J' \subset B(c(J'), 4D^{s(J')})$, the triangle inequality,
    $s(J') = s(J) + 1$ and $D=2^{100a^2}$ imply
    $$
        B(c(J'), 204D^{s(J')+1}) \subset B(c(J), 204D^{s(J') + 1} + 4D^{s(J')})
                                 \subset B(c(J), 2^8 D^{s(J) + 2})\,.
    $$
    From the doubling property \eqref{doublingx}, $D=2^{100a^2}$ and \eqref{eq-vol-sp-cube}, we
    obtain
    \begin{equation}
        \label{measure-comparison}
        \mu(B(c(J'), 204D^{s(J') + 1})) \leq 2 ^ {200a^3 + 10a} \mu(J)\,.
    \end{equation}

    If $J \subset B(\pc(\fp), 4 D^{\ps(\fp)})$, then we need only check that
    $\mu(\scI(\fp)) \le 2^{100a^3 + 10a} \mu(J)$. This follows immediately from
    $\scI(\fp) \subset B(c(J'), 100 D^{s(J') + 1})$ and \eqref{measure-comparison}.

    From now on we assume $J \not \subset B(\pc(\fp), 4 D^{\ps(\fp)})$. Since
    \begin{equation*}
        \pc(\fp) \in \scI(\fp) \subset B(c(J'), 100 D^{s(J') + 1})\, ,
    \end{equation*}
    we have by \eqref{eq-vol-sp-cube} and the triangle inequality
    $$
        J \subset J' \subset B(c(J'), 4D^{s(J')}) \subset B(\pc(\fp), 104 D^{s(J') + 1})\,.
    $$
    In particular this implies $104 D^{s(J') + 1} > 4D^{\ps(\fp)}$. By the triangle inequality
    we also have
    $$
        B(\pc(\fp), 104 D^{s(J') + 1}) \subset B(c(J), 204 D^{s(J') + 1})\,,
    $$
    so from \eqref{measure-comparison},
    $$
        \mu( B(\pc(\fp), 104 D^{s(J') + 1})) \le 2^{200a^3 + 10a} \mu(J)\,,
    $$
    which proves $\fp$ satisfies the needed criteria with $r=104 D^{s(J') + 1}$.

\end{proof}

\subsection{Almost orthogonality of separated trees}

The main result of this subsection is the almost orthogonality estimate for operators associated to distinct trees in a forest in \Cref{correlation-separated-trees} below. We will deduce it from Lemmas \ref{correlation-distant-tree-parts} and \ref{correlation-near-tree-parts}, which are proven in Subsections \ref{subsec-big-tiles} and \ref{subsec-rest-tiles}, respectively. Before stating it, we introduce some relevant notation.

The adjoint of the operator $T_{\fp}$ defined in \eqref{definetp} is given by
\begin{equation}
    \label{definetp*}
    T_{\fp}^* g(x) = \int_{E(\fp)} \overline{K_{\ps(\fp)}(y,x)} e(-\tQ(y)(x)+ \tQ(y)(y)) g(y) \, \mathrm{d}\mu(y)\,.
\end{equation}

\begin{lemma}[adjoint tile support]
    \label{adjoint-tile-support}
    \uses{kernel-summand}
    \leanok
    \lean{TileStructure.Forest.adjoint_tile_support1, TileStructure.Forest.adjoint_tile_support2}
    For each $\fp \in \fP$, we have
    $$
        T_{\fp}^* g = \mathbf{1}_{B(\pc(\fp), 5D^{\ps(\fp)})} T_{\fp}^* \mathbf{1}_{\scI(\fp)} g\,.
    $$
    For each $\fu \in \fU$ and each $\fp \in \fT(\fu)$, we have
    $$
        T_{\fp}^* g = \mathbf{1}_{\scI(\fu)} T_{\fp}^* \mathbf{1}_{\scI(\fu)} g\,.
    $$
\end{lemma}

\begin{proof}
    \leanok
    By \eqref{forest1}, $E(\fp) \subset \scI(\fp) \subset \scI(\fu)$. Thus by \eqref{definetp*}
    $$
         T_{\fp}^* g(x) = T_{\fp}^* (\mathbf{1}_{\scI(\fp)} g)(x)
    $$
    $$
        = \int_{E(\fp)} \overline{K_{\ps(\fp)}(y,x)} e(-\tQ(y)(x) + \tQ(y)(y)) \mathbf{1}_{\scI(\fp)}(y) g(y) \, \mathrm{d}\mu(y)\,.
    $$
    If this integral is not $0$, then there exists $y \in \scI(\fp)$ such that $K_{\ps(\fp)}(y,x) \ne 0$. By \eqref{supp-Ks}, \eqref{eq-vol-sp-cube} and the triangle inequality, it follows that
    \begin{equation*}
        x \in B(\pc(\fp), 5 D^{\ps(\fp)})\, .
    \end{equation*}
    Thus
    $$
        T_{\fp}^* g(x) = \mathbf{1}_{B(\pc(\fp), 5D^{\ps(\fp)})}(x) T_{\fp}^* (\mathbf{1}_{\scI(\fp)} g)(x)\,.
    $$
    The second claimed equation follows now since $\scI(\fp) \subset \scI(\fu)$ and by \eqref{forest6} we have $B(\pc(\fp), 5D^{\ps(\fp)}) \subset \scI(\fu)$.
\end{proof}

\begin{lemma}[adjoint tree estimate]
    \label{adjoint-tree-estimate}
    \leanok
    \lean{TileStructure.Forest.adjoint_tree_estimate}
    \uses{densities-tree-bound}
    For all $g$ with $\text{support}(g) \subseteq G$, we have that
    $$
        \left\| \sum_{\fp \in \fT(\fu)} T_{\fp}^* g\right\|_2 \le 2^{181a^3} \dens_1(\fT(\fu))^{1/2} \|g\|_2\,.
    $$
\end{lemma}

\begin{proof}
    \leanok
    By \Cref{densities-tree-bound}, we have for all bounded $f$ and $g$ with $|g| \le \mathbf{1}_G$ that
    $$
        \left| \int_X \overline{\sum_{\fp\in \fT(\fu)} T_{\fp}^* g} f \,\mathrm{d}\mu \right| = \left| \int_X \overline{g} \sum_{\fp \in \fT(\fu)} T_{\fp} f \,\mathrm{d}\mu \right|
    $$
    \begin{equation}
        \label{eq-adjoint-bound}
        \le 2^{181a^3} \dens_1(\fT(\fu))^{1/2} \|g\|_2 \|f\|_2\,.
    \end{equation}
    Let $f = \sum_{\fp \in \fT(\fu)} T_{\fp}^* g$. Since $|g| \le \mathbf{1}_G$, $f$ is bounded and has bounded support. In particular $\|f\|_2 < \infty$. Dividing \eqref{eq-adjoint-bound} by $\|f\|_2$ completes the proof.
\end{proof}

We define
$$
    S_{2,\fu}g := \left|\sum_{\fp \in \fT(\fu)} T_{\fp}^*g \right| + M_{\mathcal{B},1}g + |g|\,.
$$
\begin{lemma}[adjoint tree control]
    \label{adjoint-tree-control}
    \leanok
    \lean{TileStructure.Forest.adjoint_tree_control}
    \uses{adjoint-tree-estimate}
    We have for all $\fu \in \fU$ and for all $g$ with $\text{support}(g) \subseteq G$
    $$
        \|S_{2, \fu} g\|_2 \le 2^{182a^3} \|g\|_2\,.
    $$
\end{lemma}

\begin{proof}
    \leanok
    This follows immediately from Minkowski's inequality, \Cref{Hardy-Littlewood} and \Cref{adjoint-tree-estimate}, using that $a \ge 4$.
\end{proof}

Now we are ready to state the main result of this subsection.

\begin{lemma}[correlation separated trees]
    \label{correlation-separated-trees}
    \leanok
    \lean{TileStructure.Forest.correlation_separated_trees}
    \uses{correlation-distant-tree-parts,correlation-near-tree-parts}
    For any $\fu_1 \ne \fu_2 \in \fU$ and all bounded $g_1, g_2$ with bounded support, we have
    \begin{equation}
        \label{eq-lhs-sep-tree}
        \left| \int_X \sum_{\fp_1 \in \fT(\fu_1)} \sum_{\fp_2 \in \fT(\fu_2)} T^*_{\fp_1}g_1 \overline{T^*_{\fp_2}g_2 }\,\mathrm{d}\mu \right|
    \end{equation}
    \begin{equation}
        \label{eq-rhs-sep-tree}
        \le 2^{512a^3-4n} \prod_{j =1}^2 \| S_{2, \fu_j} g_j\|_{L^2(\scI(\fu_1) \cap \scI(\fu_2))}\,.
    \end{equation}
\end{lemma}

\begin{proof}[Proof of \Cref{correlation-separated-trees}]
    \leanok
    \proves{correlation-separated-trees}
    By \Cref{adjoint-tile-support} and \eqref{dyadicproperty}, the left hand side \eqref{eq-lhs-sep-tree} is $0$ unless $\scI(\fu_1) \subset \scI(\fu_2)$ or $\scI(\fu_2) \subset \scI(\fu_1)$. Without loss of generality we assume that $\scI(\fu_1) \subset \scI(\fu_2)$.

    Define
    \begin{equation}
        \label{def-Tree-S-set}
         \mathfrak{S} := \{\fp \in \fT(\fu_1) \cup \fT(\fu_2) \ : \ d_{\fp}(\fcc(\fu_1), \fcc(\fu_2)) \ge 2^{Zn/2}\,\}.
    \end{equation}
    \Cref{correlation-separated-trees} follows by combining the definition \eqref{defineZ} of $Z$ with the following two lemmas.
\end{proof}

\begin{lemma}[correlation distant tree parts]
        \label{correlation-distant-tree-parts}
        \leanok
        \lean{TileStructure.Forest.correlation_distant_tree_parts}
        \uses{Holder-van-der-Corput,Lipschitz-partition-unity,Holder-correlation-tree,lower-oscillation-bound}
        We have for all $\fu_1 \ne \fu_2 \in \fU$ with $\scI(\fu_1) \subset \scI(\fu_2)$ and all bounded $g_1, g_2$ with bounded support
      \begin{equation}
            \label{eq-lhs-big-sep-tree}
            \left| \int_X \sum_{\fp_1 \in \fT(\fu_1)} \sum_{\fp_2 \in \fT(\fu_2) \cap \mathfrak{S}} T^*_{\fp_1}g_1 \overline{T^*_{\fp_2}g_2 }\,\mathrm{d}\mu \right|
        \end{equation}
      \begin{equation}
            \label{eq-rhs-big-sep-tree}
            \le 2^{511a^3} 2^{-Zn/(4a^2 + 2a^3)} \prod_{j =1}^2 \| S_{2, \fu_j} g_j\|_{L^2(\scI(\fu_1))}\,.
        \end{equation}
    \end{lemma}
    \begin{lemma}[correlation near tree parts]
        \label{correlation-near-tree-parts}
        \leanok
        \lean{TileStructure.Forest.correlation_near_tree_parts}
        \uses{tree-projection-estimate,dyadic-partition-2,bound-for-tree-projection}
        We have for all $\fu_1 \ne \fu_2 \in \fU$ with $\scI(\fu_1) \subset \scI(\fu_2)$ and all bounded $g_1, g_2$ with bounded support
      \begin{equation}
            \label{eq-lhs-small-sep-tree}
            \left| \int_X \sum_{\fp_1 \in \fT(\fu_1)} \sum_{\fp_2 \in \fT(\fu_2) \setminus \mathfrak{S}} T^*_{\fp_1}g_1 \overline{T^*_{\fp_2}g_2 }\,\mathrm{d}\mu \right|
        \end{equation}
      \begin{equation}
            \label{eq-rhs-small-sep-tree}
            \le 2^{232a^3+21a+5} 2^{-\frac{25}{101a}Zn\kappa} \prod_{j=1}^2 \|S_{2, \fu_j} g_j\|_{L^2(\scI(\fu_1))}\,.
        \end{equation}
    \end{lemma}

In the proofs of both lemmas, we will need the following observation.

\begin{lemma}[overlap implies distance]
    \label{overlap-implies-distance}
    \leanok
    \lean{TileStructure.Forest.𝔗_subset_𝔖₀, TileStructure.Forest.overlap_implies_distance}
    Let $\fu_1 \ne \fu_2 \in \fU$ with $\scI(\fu_1) \subset \scI(\fu_2)$. If $\fp \in \fT(\fu_1) \cup \fT(\fu_2)$ with $\scI(\fp) \cap \scI(\fu_1) \ne \emptyset$, then $\fp \in \mathfrak{S}$. In particular, we have $\fT(\fu_1) \subset \mathfrak{S}$.
\end{lemma}

\begin{proof}
    \leanok
    Suppose first that $\fp \in \fT(\fu_1)$. Then $\scI(\fp) \subset \scI(\fu_1) \subset \scI(\fu_2)$, by \eqref{forest1}. Thus we have by the separation condition \eqref{forest5}, \eqref{eq-freq-comp-ball}, \eqref{forest1} and the triangle inequality
    \begin{align*}
        d_{\fp}(\fcc(\fu_1), \fcc(\fu_2)) &\ge d_{\fp}(\fcc(\fp), \fcc(\fu_2)) - d_{\fp}(\fcc(\fp), \fcc(\fu_1))\\
        &\ge 2^{Z(n+1)} - 4\\
        &\ge 2^{Zn/2}\,,
    \end{align*}
    using that $Z= 2^{12a}\ge 4$. Hence $\fp \in \mathfrak{S}$.

    Suppose now that $\fp \in \fT(\fu_2)$. If $\scI(\fp) \subset \scI(\fu_1)$, then the same argument as above with $\fu_1$ and $\fu_2$ swapped shows $\fp \in \mathfrak{S}$. If $\scI(\fp) \not \subset \scI(\fu_1)$ then, by \eqref{dyadicproperty}, $\scI(\fu_1) \subset \scI(\fp)$. Pick $\fp' \in \fT(\fu_1)$, we have $\scI(\fp') \subset \scI(\fu_1) \subset \scI(\fp)$. Hence, by \Cref{monotone-cube-metrics} and the first paragraph
    $$
        d_{\fp}(\fcc(\fu_1), \fcc(\fu_2)) \ge d_{\fp'}(\fcc(\fu_1), \fcc(\fu_2)) \ge 2^{Zn}\,,
    $$
    so $\fp \in \mathfrak{S}$.
\end{proof}

To simplify the notation, we will write at various places throughout the proof of Lemmas \ref{correlation-distant-tree-parts} and \ref{correlation-near-tree-parts} for a subset $\fC \subset \fP$
$$
    T_{\fC} f := \sum_{\fp \in \fC} T_{\fp} f\,, \quad\quad T_{\fC}^* g := \sum_{\fp\in\fC} T_{\fp}^* g\,.
$$

\subsection{Proof of the Tiles with large separation Lemma}
    \label{subsec-big-tiles}

\Cref{correlation-distant-tree-parts} follows from the van der Corput estimate in \Cref{Holder-van-der-Corput}. We apply this proposition in \Cref{subsubsec-van-der-corput}. To prepare this application, we first, in \Cref{subsubsec-pao}, construct a suitable partition of unity, and show then, in \Cref{subsubsec-holder-estimates} the H\"older estimates needed to apply \Cref{Holder-van-der-Corput}.

\subsubsection{A partition of unity}
\label{subsubsec-pao}
    Define
    $$
        \mathcal{J}' = \{J \in \mathcal{J}(\mathfrak{S}) \ : \ J \subset \scI(\fu_1)\}\,.
    $$

    \begin{lemma}[dyadic partition 1]
        \label{dyadic-partition-1}
        \leanok
        \lean{TileStructure.Forest.union_𝓙₅, TileStructure.Forest.pairwiseDisjoint_𝓙₅}
        \uses{dyadic-partitions}
        We have that
        $$
            \scI(\fu_1) = \dot{\bigcup_{J \in \mathcal{J}'}} J\,.
        $$
    \end{lemma}

    \begin{proof}
        \leanok
        By \Cref{dyadic-partitions}, it remains only to show that each $J \in \mathcal{J}(\mathfrak{S})$ with $J \cap \scI(\fu_1) \ne \emptyset$ is in $\mathcal{J}'$. But if $J \notin \mathcal{J}'$, then by \eqref{dyadicproperty} $\scI(\fu_1) \subset J$. Pick $\fp \in \fT(\fu_1) \subset \mathfrak{S}$. Then $\scI(\fp) \subset J$. This contradicts the definition of $\mathcal{J}(\mathfrak{S})$.
    \end{proof}

    For cubes $J \in \mathcal{D}$, denote
    \begin{equation}
        \label{def-BJ}
        B(J) := B(c(J), 8D^{s(J)}).
    \end{equation}
    The main result of this subsubsection is the following.

    \begin{lemma}[Lipschitz partition unity]
        \label{Lipschitz-partition-unity}
        \leanok
        \lean{TileStructure.Forest.sum_χ, TileStructure.Forest.χ_le_indicator,
          TileStructure.Forest.dist_χ_le}
        \uses{dyadic-partition-1,moderate-scale-change}
        There exists a family of functions $\chi_J$, $J \in \mathcal{J}'$ such that \begin{equation}
            \label{eq-pao-1}
            \mathbf{1}_{\scI(\fu_1)} = \sum_{J \in \mathcal{J}'} \chi_J\,,
        \end{equation}
        and for all $J \in \mathcal{J}'$ and all $y,y' \in \scI(\fu_1)$
      \begin{equation}
            \label{eq-pao-2}
            0 \leq \chi_J(y) \leq \mathbf{1}_{B(J)}(y)\,,
        \end{equation}
      \begin{equation}
            \label{eq-pao-3}
            |\chi_J(y) - \chi_J(y')| \le 2^{227a^3} \frac{\rho(y,y')}{D^{s(J)}}\,.
        \end{equation}
    \end{lemma}

    In the proof, we will use the following auxiliary lemma.

    \begin{lemma}[moderate scale change]
        \label{moderate-scale-change}
        \leanok
        \lean{TileStructure.Forest.moderate_scale_change}

        If $J, J' \in \mathcal{J'}$ with
        $$
            B(J) \cap B(J') \ne \emptyset\,,
        $$
        then $|s(J) - s(J')| \le 1$.
    \end{lemma}

    \begin{proof}[Proof of \Cref{Lipschitz-partition-unity}]
        \leanok
        \proves{Lipschitz-partition-unity}
        For each cube $J \in \mathcal{J}$ let
        $$
            \tilde\chi_J(y) = \mathbf{1}_{\scI(\fu_1)}(y)\max\{0, 8 - D^{-s(J)} \rho(y, c(J))\}\,,
        $$
        and set
        $$
            a(y) = \sum_{J \in \mathcal{J}'} \tilde \chi_J(y)\,.
        $$
        We define
        \[
            \chi_J(y) := \frac{\tilde \chi_J(y)}{a(y)}\,.
        \]
        Then, due to \eqref{forest6} and \eqref{def-BJ}, the properties \eqref{eq-pao-1} and \eqref{eq-pao-2} are clearly true. Estimate \eqref{eq-pao-3} follows from \eqref{eq-pao-2} if $y, y' \notin B(J)$. Thus we can assume that $y \in B(J)$. We have by the triangle inequality
        $$
            |\chi_J(y) - \chi_J(y')| \le \frac{|\tilde \chi_J(y) - \tilde \chi_J(y')|}{a(y)} + \frac{\tilde \chi_J(y')|a(y) - a(y')|}{a(y)a(y')}
        $$
        Since $\tilde \chi_J(z) \ge 4$ for all $z \in B(c(J),4D^{s(J)}) \supset J$ and by \Cref{dyadic-partition-1}, we have that $a(z) \ge 4$ for all $z \in \scI(\fu_1)$. So we can estimate the above further by
        $$
            \le 2^{-2}(|\tilde \chi_J(y) - \tilde \chi_J(y')| + \tilde \chi_J(y')|a(y) - a(y')|)\,.
        $$
        If $y' \notin B(\pc(\fp), 8D^{\ps(\fp)})$ then the second summand vanishes. Else, we can estimate the above, using also that $|\tilde \chi_J(y')| \le 8$, by
        $$
            \le 2^{-2} |\tilde \chi_J(y) - \tilde \chi_J(y')| + 2 \sum_{\substack{J' \in \mathcal{J}'\\ B(J') \cap B(J) \ne \emptyset}}|\tilde \chi_{J'}(y) - \tilde \chi_{J'} (y')|\,.
        $$
        By the triangle inequality, we have for all dyadic cubes $I$
        $$
            |\tilde \chi_I(y) - \tilde \chi_I(y')| \le \rho(y, y') D^{-s(I)}\,.
        $$
        Using this above, we obtain
        $$
            |\chi_J(y) - \chi_J(y')| \le \rho(y,y') \Big( \frac{1}{4} D^{-s(J)} + 2 \sum_{\substack{J' \in \mathcal{J}'\\ B(J') \cap B(J) \ne \emptyset}} D^{-s(J')}\Big)\,.
        $$
        By \Cref{moderate-scale-change}, this is at most
        $$
             \frac{\rho(y,y')}{D^{s(J)}} \left( \frac{1}{4} + 2D |\{J' \in \mathcal{J}' \ : \ B(J') \cap B(J) \ne \emptyset\}|\right)\,.
        $$
        By \eqref{eq-vol-sp-cube} and \Cref{dyadic-partition-1}, the balls $B(c(J'), \frac{1}{4} D^{s(J')})$ are pairwise disjoint. By the triangle inequality and \Cref{moderate-scale-change}, each such ball for $J'$ in the set of the last display is contained in
        $$
            B(c(J), 9 D^{s(J) + 1})\,.
        $$
        By the doubling property \eqref{doublingx}, we further have
        $$
            \mu(B(c(J), 9 D^{s(J) + 1})) \le 2^{200a^3 + 7a} \mu\Big(B(c(J'), \frac{1}{4}D^{s(J')})\Big)
        $$
        for each such ball.
        Thus
        $$
            |\{J' \in \mathcal{J}' \ : \ B(J') \cap B(J) \ne \emptyset\}| \le 2^{200a^3 + 7a}\,.
        $$
        Recalling that $D=2^{100a^2}$, we obtain
        $$\frac{1}{4} + 2D |\{J' \in \mathcal{J}' \ : \ B(J') \cap B(J) \ne \emptyset\}|\leq 2^{200a^3 + 100a^2 + 7a + 2}.$$
        Since $a\ge 4$, \eqref{eq-pao-3} follows.
    \end{proof}

    \begin{proof}[Proof of \Cref{moderate-scale-change}]
        \leanok
        \proves{moderate-scale-change}
        Suppose that $s(J') < s(J) - 1$. Then $s(J) > -S$. Thus, by the definition of $\mathcal{J}'$ there exists no $\fp \in \mathfrak{S}$ with
      \begin{equation}
            \label{eq-tile-incl-1}
            \scI(\fp) \subset B(c(J), 100D^{s(J) + 1})\,.
        \end{equation}
        Since $s(J') < s(J)$, there exists a cube $J'' \in \mathcal{D}$ with $J \subset J''$ and $s(J'') = s(J') + 1$. By the definition of $\mathcal{J}'$, there exists a tile $\fp \in \mathfrak{S}$ with
      \begin{equation}
            \label{tile-incl-2}
            \scI(\fp) \subset B(c(J''), 100 D^{s(J')+2})\,.
        \end{equation}
        But by the triangle inequality and \eqref{defineD}, we have
        $$
            B(c(J''), 100 D^{s(J')+2}) \subset B(c(J), 100D^{s(J) + 1})\,,
        $$
        which contradicts \eqref{eq-tile-incl-1} and \eqref{tile-incl-2}.
    \end{proof}

\subsubsection{H\"older estimates for adjoint tree operators}
\label{subsubsec-holder-estimates}
    Let $g_1, g_2:X \to \mathbb{C}$ be bounded with bounded support.
    Define for $J \in \mathcal{J}'$
    \begin{equation}
        \label{def-hj}
        h_J(y) := \chi_J(y)\cdot(e(\fcc(\fu_1)(y)) T_{\fT(\fu_1)}^* g_1(y)) \cdot \overline{(e(\fcc(\fu_2)(y)) T_{\fT(\fu_2) \cap \mathfrak{S}}^* g_2(y))}\,.
    \end{equation}
    The main result of this subsubsection is the following $\tau$-H\"older estimate for $h_J$, where $\tau = 1/a$.

    \begin{lemma}[Holder correlation tree]
        \label{Holder-correlation-tree}
        \leanok
        \lean{TileStructure.Forest.holder_correlation_tree}
        \uses{global-tree-control-2}
        We have for all $J \in \mathcal{J}'$ that
      \begin{equation}
            \label{hHolder}
            \|h_J\|_{C^{\tau}(B(c(J), 16D^{s(J)}))} \le 2^{485a^3} \prod_{j = 1,2} (\inf_{B(c(J), \frac{1}{8}D^{s(J)})} |T_{\fT(\fu_j)}^* g_j| + \inf_J M_{\mathcal{B}, 1} |g_j|)\,.
        \end{equation}
    \end{lemma}

    We will prove this lemma at the end of this section, after establishing several auxiliary results.

    We begin with the following H\"older continuity estimate for adjoints of operators associated to tiles.
    \begin{lemma}[Holder correlation tile]
        \label{Holder-correlation-tile}
        \leanok
        \lean{TileStructure.Forest.holder_correlation_tile}
        \uses{adjoint-tile-support}
        Let $\fu \in \fU$ and $\fp \in \fT(\fu)$. Then for all $y, y' \in X$ and all bounded $g$ with bounded support, we have
        $$
            |e(\fcc(\fu)(y)) T_{\fp}^* g(y) - e(\fcc(\fu)(y')) T_{\fp}^* g(y')|
        $$
      \begin{equation}
            \label{T*Holder2}
            \le \frac{2^{128a^3}}{\mu(B(\pc(\fp), 4D^{\ps(\fp)}))} \left(\frac{\rho(y, y')}{D^{\ps(\fp)}}\right)^{1/a} \int_{E(\fp)} |g(x)| \, \mathrm{d}\mu(x)\,.
        \end{equation}
    \end{lemma}

    \begin{proof}
        \leanok
        By \eqref{definetp*}, we have
        $$
            |e(\fcc(\fu)(y)) T_{\fp}^* g(y) - e(\fcc(\fu)(y')) T_{\fp}^* g(y')|
        $$
        \begin{multline*}
            =\bigg| \int_{E(\fp)} e(\tQ(x)(x) - \tQ(x)(y) + \fcc(\fu)(y)) \overline{K_{\ps(\fp)}(x, y)} g(x) \\
            - e(\tQ(x)(x) - \tQ(x)(y') + \fcc(\fu)(y')) \overline{K_{\ps(\fp)}(x, y')} g(x) \, \mathrm{d}\mu(x)\bigg|
        \end{multline*}
        \begin{multline*}
            \leq\int_{E(\fp)} |g(x)| |e(\tQ(x)(y) - \tQ(x)(y') - \fcc(\fu)(y) + \fcc(\fu)(y'))\overline{K_{\ps(\fp)}(x, y)}\\
            - \overline{K_{\ps(\fp)}(x, y')}| \, \mathrm{d}\mu(x)
        \end{multline*}
        \begin{multline}
            \leq\int_{E(\fp)} |g(x)| |e(-\tQ(x)(y) + \tQ(x)(y') + \fcc(\fu)(y) - \fcc(\fu)(y')) - 1| \\
            \times |\overline{K_{\ps(\fp)}(x, y)}|\, \mathrm{d}\mu(x) \label{T*Holder1b}
        \end{multline}
      \begin{equation}
            + \int_{E(\fp)} |g(x)| |\overline{K_{\ps(\fp)}(x, y)} - \overline{K_{\ps(\fp)}(x, y')} |\, \mathrm{d}\mu(x)\,.\label{T*Holder1}
        \end{equation}
        By the oscillation estimate \eqref{osccontrol}, we have
        $$
            |-\tQ(x)(y) + \tQ(x)(y') + \fcc(\fu)(y) - \fcc(\fu)(y')|
        $$
      \begin{equation}
            \label{eq-lem-tile-Holder-comp}
            \le d_{B(y, 1.6\rho(y,y'))}(\tQ(x), \fcc(\fu))\,.
        \end{equation}
        Suppose that $y, y' \in B(\pc(\fp), 5D^{\ps(\fp)})$, so that $\rho(y,y') \le 10D^{\ps(\fp)}$.
        Let $k \in \mathbb{Z}$ be such that $2^{ak}\rho(y,y') \le 10D^{\ps(\fp)}$ but $2^{a(k+1)} \rho(y,y') > 10D^{\ps(\fp)}$.
        In particular, $k \ge 0$. Then, using \eqref{seconddb} followed by \eqref{firstdb}, we can bound \eqref{eq-lem-tile-Holder-comp} from above by
        $$
            2^{-k} d_{B(\pc(\fp), 16 D^{\ps(\fp)})}(\tQ(x), \fcc(\fu)) \le 2^{6a - k} d_{\fp}(\tQ(x), \fcc(\fu))\,.
        $$
        Since $x \in E(\fp)$ we have $\tQ(x) \in \Omega(\fp) \subset B_{\fp}(\fcc(\fp), 1)$, and since $\fp \in \fT(\fu)$ we have $\fcc(\fu) \in B_{\fp}(\fcc(\fp), 4)$, so this is estimated by
        $$
            \le 5 \cdot 2^{6a - k}\,.
        $$
        By definition of $k$, we have
        $$
            -k < 1 - \frac{1}{a} \log_2\left(\frac{10 D^{\ps(\fp)}}{\rho(y,y')}\right)\,,
        $$
        which gives
      \begin{equation}
            \label{eq-lem-Tile-holder-im1}
             |-\tQ(x)(y) + \tQ(x)(y') + \fcc(\fu)(y) - \fcc(\fu)(y')| \le 10 \cdot 2^{6a} \left(\frac{\rho(y,y')}{10 D^{\ps(\fp)}}\right)^{1/a}\,.
        \end{equation}
        For all $x \in \scI(\fp)$, we have by \eqref{doublingx} that
        $$
            \mu(B(x, D^{\ps(\fp)})) \ge 2^{-3a} \mu(B(\pc(\fp), 4D^{\ps(\fp)}))\,.
        $$
        Combining the above with \eqref{eq-Ks-size}, \eqref{eq-Ks-smooth} and \eqref{eq-lem-Tile-holder-im1},
        we obtain
        $$
            \eqref{T*Holder1b}+\eqref{T*Holder1} \le \frac{2^{3a}}{\mu(B(\pc(\fp), 4D^{\ps(\fp)}))} \int_{E(\fp)}|g(x)| \, \mathrm{d}\mu(x) \times
        $$
        $$
            (2^{102a^3} \cdot 10 \cdot 2^{6a} \left(\frac{\rho(y,y')}{ D^{\ps(\fp)}}\right)^{1/a} + 2^{127a^3} \left(\frac{\rho(y,y')}{D^{\ps(\fp)}}\right)^{1/a})
        $$
        Since $\rho(y,y') \le 10 D^{\ps(\fp)}$, we conclude
        $$
            \eqref{T*Holder1b}+\eqref{T*Holder1} \le \frac{2^{128a^3}}{\mu(B(\pc(\fp), 4D^{\ps(\fp)}))} \left(\frac{\rho(y,y')}{D^{\ps(\fp)}}\right)^{1/a} \int_{E(\fp)}|g(x)| \, \mathrm{d}\mu(x)\,.
        $$

        Next, if $y,y' \notin B(\pc(\fp), 5D^{\ps(\fp)})$, then $T_{\fp}^*g(y) = T_{\fp}^*g(y') = 0$, by \Cref{adjoint-tile-support}. Then \eqref{T*Holder2} holds.

        Finally, if $y \in B(\pc(\fp), 5D^{\ps(\fp)})$ and $y' \notin B(\pc(\fp), 5D^{\ps(\fp)})$, then
        $$
            |e(\fcc(\fu)(y)) T_{\fp}^* g(y) - e(\fcc(\fu)(y')) T_{\fp}^* g(y')| = |T_{\fp}^* g(y)|
        $$
        $$
            \le \int_{E(\fp)} |K_{\ps(\fp)}(x,y)| |g(x)| \, \mathrm{d}\mu(x)\,.
        $$
        By the same argument used to prove \eqref{eq-Ks-aux}, this is bounded by
      \begin{equation}
            \label{eq-lem-Tile-holder-im2}
            \le 2^{102a^3} \int_{E(\fp)} \frac{1}{\mu(B(x, D^s))} \psi(D^{-s} \rho(x,y)) |g(x)| \, \mathrm{d}\mu(x)\,.
        \end{equation}
        It follows from the definition of $\psi$ that
        $$
            \psi(x) \le \max\{0, (2 - 4x)^{1/a}\}\,.
        $$
        Now for all $x\in E(\fp)$, it follows by the triangle inequality and \eqref{eq-vol-sp-cube} that
        \begin{multline*}
        2 - 4D^{-\ps(\fp)}\rho(x,y)\leq 2 - 4D^{-\ps(\fp)}\rho(y, \pc(\fp)) + 4 D^{-\ps(\fp)}\rho(x, \pc(\fp))\\\leq 18 - 4 D^{-\ps(\fp)} \rho(y, \pc(\fp)) \leq 4 D^{-\ps(\fp)}\rho(y,y') - 2 .
        \end{multline*}
        Combining the above with the previous estimate on $\psi$, we get
        $$
            \psi(D^{-\ps(\fp)}\rho(x,y)) \le 4 (D^{-\ps(\fp)}\rho(y,y'))^{1/a}.
        $$
        Further, we obtain from the doubling property \eqref{doublingx} and \eqref{eq-vol-sp-cube} that
        $$
            \mu(B(x, D^{\ps(\fp)})) \ge 2^{-3a} \mu(B(\pc(\fp), 4D^{\ps(\fp)}))\,.
        $$
        Plugging this into \eqref{eq-lem-Tile-holder-im2} and using $a \ge 4$, we get
        $$
            |T_{\fp}^* g(y)| \le \frac{2^{103a^3}}{\mu(B(\pc(\fp), 4D^{\ps(\fp)}))} \left(\frac{\rho(y,y')}{D^{\ps(\fp)}}\right)^{1/a} \int_{E(\fp)} |g(x)| \, \mathrm{d}\mu(y)\,,
        $$
        which completes the proof of the lemma.
    \end{proof}

    Recall that
    \begin{equation*}
        B(J) := B(c(J), 8D^{s(J)}).
    \end{equation*}
    We also denote
    \begin{equation*}
     B'(J) := B(c(J), 16D^{s(J)}),
    \end{equation*}
    \begin{equation*}
     B^\circ{}(J) := B(c(J), \frac{1}{8}D^{s(J)})\, .
    \end{equation*}

    \begin{lemma}[limited scale impact]
        \label{limited-scale-impact}
        \leanok
        \lean{TileStructure.Forest.limited_scale_impact}
        \uses{overlap-implies-distance}
        Let $\fp \in \fT(\fu_2) \setminus \mathfrak{S}$, $J \in \mathcal{J}'$ and suppose that
        $$
           B(\scI(\fp)) \cap B^\circ(J) \ne \emptyset\,.
        $$
        Then
        $$
            s(J) \le \ps(\fp) \le s(J) +3\,.
        $$
    \end{lemma}

    \begin{proof}
        \leanok
        For the first estimate, assume that $\ps(\fp) < s(J)$, then in particular $\ps(\fp) \le \ps(\fu_1)$. Since $\fp \notin \mathfrak{S}$, we have by \Cref{overlap-implies-distance} that $\scI(\fp) \cap \scI(\fu_1) = \emptyset$.
        Since $B\Big(c(J), \frac{1}{4} D^{s(J)}\Big) \subset \scI(J) \subset \scI(\fu_1)$, this implies
        $$
            \rho(c(J), \pc(\fp)) \ge \frac{1}{4}D^{s(J)}\,.
        $$
        On the other hand
        $$
            \rho(c(J), \pc(\fp)) \le \frac{1}{8} D^{s(J)} + 8 D^{\ps(\fp)}\,,
        $$
        by our assumption. Thus $D^{\ps(\fp)} \ge 64^{-1} D^{s(J)}$, which contradicts \eqref{defineD} and $a \ge 4$.

        For the second estimate, assume that $\ps(\fp) > s(J) + 3$. Since $J \in \mathcal{J}'$, we have $J \subsetneq \scI(\fu_1)$. Thus there exists $J' \in \mathcal{D}$ with $J \subset J'$ and $s(J') = s(J) + 1$, by \eqref{coverdyadic} and \eqref{dyadicproperty}. By definition of $\mathcal{J}'$, there exists some $\fp' \in \mathfrak{S}$ such that $\scI(\fp') \subset B(c(J'), 100 D^{s(J) + 2})$. On the other hand, since $B(\scI(\fp)) \cap B^\circ(J) \ne \emptyset$, by the triangle inequality it holds that
        $$
            B(c(J'), 100 D^{s(J) + 3}) \subset B(\pc(\fp), 10 D^{\ps(\fp)})\,.
        $$
        Using the definition of $\mathfrak{S}$, we have
        $$
            2^{Zn/2} \le d_{\fp'}(\fcc(\fu_1), \fcc(\fu_2)) \le d_{B(c(J'), 100 D^{s(J) + 2})}(\fcc(\fu_1), \fcc(\fu_2))\,.
        $$
        By \eqref{seconddb}, this is
        $$
            \le 2^{-100a} d_{B(c(J'), 100 D^{s(J) + 3})}(\fcc(\fu_1), \fcc(\fu_2))
        $$
        $$
            \le 2^{-100a} d_{B(\pc(\fp), 10 D^{\ps(\fp)})}(\fcc(\fu_1), \fcc(\fu_2))\,,
        $$
        and by \eqref{firstdb} and the definition of $\mathfrak{S}$
        $$
            \le 2^{-94a} d_{\fp}(\fcc(\fu_1), \fcc(\fu_2)) \le 2^{-94a} 2^{Zn/2}\,.
        $$
        This is a contradiction, the second estimate follows.
    \end{proof}

    \begin{lemma}[local tree control]
        \label{local-tree-control}
        \leanok
        \lean{TileStructure.Forest.local_tree_control}
        \uses{limited-scale-impact}
        For all $J \in \mathcal{J}'$ and all bounded $g$ with bounded support
        $$
            \sup_{B^\circ{}(J)} |T_{\mathfrak{T}(\mathfrak{u}_2)\setminus\mathfrak{S}}^* g| \le 2^{104a^3} \inf_J M_{\mathcal{B},1}|g|
        $$
    \end{lemma}

    \begin{proof}
        \leanok
        By the triangle inequality and since $T_{\fp}^* g = \mathbf{1}_{B(\pc(\fp), 5D^{\ps(\fp)})} T_{\fp}^* g$, we have
        $$
            \sup_{B^\circ{}(J)} |T_{\fT(\fu_2) \setminus\mathfrak{S}}^* g|
            \leq \sup_{B^\circ{}(J)} \sum_{\substack{\fp \in \fT(\fu_2) \setminus \mathfrak{S}\\ B(\scI(\fp)) \cap B^\circ(J) \ne \emptyset}} |T_{\fp}^*g|\,.
        $$
        By \Cref{limited-scale-impact}, this is at most
      \begin{equation}
            \label{eq-sep-tree-aux-3}
            \sum_{s = s(J)}^{s(J) + 3} \sum_{\substack{\fp \in \fP, \ps(\fp) = s\\ B(\scI(\fp)) \cap B^\circ(J) \ne \emptyset}} \sup_{B^\circ{}(J)} |T_{\fp}^* g|\,.
        \end{equation}
        If $x \in E(\fp)$ and $B(\scI(\fp)) \cap B^\circ(J) \ne \emptyset$, then
        $$
            B(c(J), 16D^{\ps(\fp)}) \subset B(x, 32 D^{\ps(\fp)})\,,
        $$
        by \eqref{eq-vol-sp-cube} and the triangle inequality. Using the doubling property \eqref{doublingx}, it follows that
        $$
            \mu(B(x, D^{\ps(\fp)})) \ge 2^{-5a} \mu(B(c(J), 16D^{\ps(\fp)}))\,.
        $$
        Using \eqref{definetp*}, \eqref{eq-Ks-size} and that $a \ge 4$, we bound \eqref{eq-sep-tree-aux-3} by
        $$
            2^{103a^3}\sum_{s = s(J)}^{s(J) + 3} \sum_{\substack{\fp \in \fP, \ps(\fp) = s\\B(\scI(\fp)) \cap B^\circ(J) \ne \emptyset}} \frac{1}{\mu(B(c(J), 16 D^s)} \int_{E(\fp)} |g| \, \mathrm{d}\mu\,.
        $$
        For each $I \in \mathcal{D}$, the sets $E(\fp)$ for $\fp \in \fP$ with $\scI(\fp) = I$ are pairwise disjoint by \eqref{defineep} and \eqref{eq-dis-freq-cover}.
        Further, if $B(\scI(\fp)) \cap B^\circ(J) \ne \emptyset$ and $\ps(\fp) \ge s(J)$, then $E(\fp) \subset B(c(J), 16 D^{\ps(\fp)})$. Thus the last display is bounded by
        $$
            2^{103a^3}\sum_{s = s(J)}^{s(J) + 3} \frac{1}{\mu(B(c(J), 16 D^s))} \int_{B(c(J), 16 D^s)} |g| \, \mathrm{d}\mu\,.
        $$
        $$
            \le \inf_{x' \in J} 2^{103a^3 +2} M_{\mathcal{B}, 1} |g|\,.
        $$
        The lemma follows since $a \ge 4$.
    \end{proof}

    \begin{lemma}[scales impacting interval]
        \label{scales-impacting-interval}
        \leanok
        \lean{TileStructure.Forest.scales_impacting_interval}
        \uses{overlap-implies-distance}
        Let $\fC = \fT(\fu_1)$ or $\fC = \fT(\fu_2) \cap \mathfrak{S}$. Then for each $J \in \mathcal{J}'$ and $\fp \in \fC$ with $B(\scI(\fp)) \cap B'(J) \neq \emptyset$, we have $\ps(\fp) \ge s(J)$.
    \end{lemma}

    \begin{proof}
        \leanok
        By \Cref{overlap-implies-distance}, we have that in both cases, $\fC \subset \mathfrak{S}$. If $\fp \in \fC$ with $B(\scI(\fp)) \cap B'(J) \neq \emptyset$ and $\ps(\fp) < s(J)$, then $\scI(\fp) \subset B(c(J), 100 D^{s(J) + 1})$. Since $\fp \in \mathfrak{S}$, it follows from the definition of $\mathcal{J}'$ that $s(J) = -S$, which contradicts $\ps(\fp) < s(J)$.
    \end{proof}

    \begin{lemma}[global tree control 1]
        \label{global-tree-control-1}
        \leanok
        \lean{TileStructure.Forest.global_tree_control1_edist_left, TileStructure.Forest.global_tree_control1_edist_right, TileStructure.Forest.global_tree_control1_supbound}
        \uses{Holder-correlation-tile,scales-impacting-interval}
        Let $\fC_1 = \fT(\fu_1)$ and $\fC_2 = \fT(\fu_2) \cap \mathfrak{S}$. Then for $i = 1,2$ and each $J \in \mathcal{J}'$ and all bounded $g$ with bounded support, we have
        \begin{align}
            \label{TreeUB}
            \sup_{B'(J)} |T_{\fC_i}^*g| \leq \inf_{B^\circ{}(J)} |T^*_{\fC_i} g| + 2^{128a^3+4a+3} \inf_{J} M_{\mathcal{B}, 1} |g|
        \end{align}
        and for all $y,y' \in B'(J)$
        $$
            |e(\fcc(\fu_i)(y)) T_{\fC_i}^* g(y) - e(\fcc(\fu_i)(y')) T_{\fC_i}^* g(y')|
        $$
        \begin{equation}
            \label{TreeHolder}
             \le 2^{128a^3+4a+1} \left(\frac{\rho(y,y')}{D^{s(J)}}\right)^{1/a} \inf_J M_{\mathcal{B},1} |g|\,.
        \end{equation}
    \end{lemma}

    \begin{proof}
        \leanok
        Note that \eqref{TreeUB} follows from \eqref{TreeHolder}, since for $y'\in B^\circ{}(J)$, by the triangle inequality,
        $$\left(\frac{\rho(y,y')}{D^{s(J)}}\right)^{1/a}\le \Big(16 + \frac{1}8\Big)^{1/a}\le 2^2.$$

        By the triangle inequality, \Cref{adjoint-tile-support} and \Cref{Holder-correlation-tile}, we have for all $y, y' \in B'(J)$
        \begin{equation}
            \label{eq-C-Lip}
            |e(\fcc(\fu_i)(y)) T_{\fC_i}^* g(y) - e(\fcc(\fu_i)(y')) T_{\fC_i}^* g(y')|
        \end{equation}
        $$
            \le \sum_{\substack{\fp \in \fC_i\\ B(\scI(\fp)) \cap B'(J) \neq \emptyset}} |e(\fcc(\fu_i)(y)) T_{\fp}^* g(y) - e(\fcc(\fu_i)(y')) T_{\fp}^* g(y')|
        $$
        $$
            \le 2^{128a^3}\rho(y,y')^{1/a} \sum_{\substack{\fp \in \fC_i\\ B(\scI(\fp)) \cap B'(J) \neq \emptyset}} \frac{D^{- \ps(\fp)/a}}{\mu(B(\pc(\fp), 4D^{\ps(\fp)}))} \int_{E(\fp)} |g| \, \mathrm{d}\mu\,.
        $$
        By \Cref{scales-impacting-interval}, we have $\ps(\fp) \ge s(J)$ for all $\fp$ occurring in the sum. Further, for each $s \ge s(J)$, the sets $E(\fp)$ for $\fp \in \fP$ with $\ps(\fp) = s$ are pairwise disjoint by \eqref{defineep} and \eqref{eq-dis-freq-cover}, and contained in $B(c(J), 32D^{s})$ by \eqref{eq-vol-sp-cube} and the triangle inequality. Using also the doubling estimate \eqref{doublingx}, we obtain that the expression in the last display can be estimated by
        $$
            2^{128a^3}\rho(y,y')^{1/a} \sum_{S \ge s \ge s(J)} D^{-s/a} \frac{2^{4a}}{\mu(B(c(J), 32D^{s}))} \int_{B(c(J), 32D^{s})} |g| \, \mathrm{d}\mu
        $$
        $$
            \le 2^{128a^3+4a} \left(\frac{\rho(y,y')}{D^{s(J)}}\right)^{1/a} \sum_{S \ge s \ge s(J)} D^{(s(J) - s)/a} \inf_J M_{\mathcal{B},1} |g|\,.
        $$
        Since $D^{-1/a}\le\frac12$, we have
        $$
            \sum_{S \ge s \ge s(J)} D^{(s(J) - s)/a} \le 2.
        $$
        Estimate \eqref{TreeHolder}, and therefore the lemma, follow.
    \end{proof}

    \begin{lemma}[global tree control 2]
        \label{global-tree-control-2}
        \leanok
        \lean{TileStructure.Forest.global_tree_control2}
        \uses{global-tree-control-1, local-tree-control}
        We have for all $J \in \mathcal{J}'$ and all bounded $g$ with bounded support
        $$
            \sup_{B'(J)} |T^*_{\fT(\fu_2) \cap \mathfrak{S}} g| \le \inf_{B^\circ{}(J)} |T^*_{\fT(\fu_2)} g| + 2^{129a^3} \inf_{J} M_{\mathcal{B},1}|g|\,.
        $$
    \end{lemma}

    \begin{proof}
        \leanok
        By \Cref{global-tree-control-1}
        $$
            \sup_{B'(J)} |T^*_{\fT(\fu_2) \cap \mathfrak{S}} g| \le \inf_{B^\circ{}(J)} |T_{\fT(\fu_2) \cap \mathfrak{S}}^* g| + 2^{128a^3+4a+3} \inf_{J} M_{\mathcal{B}, 1} |g|
        $$
        $$
            \le \inf_{B^\circ{}(J)} |T_{\fT(\fu_2)}^* g| + \sup_{B^\circ{}(J)} |T_{\fT(\fu_2) \setminus \mathfrak{S}}^* g| + 2^{128a^3+4a+3} \inf_{J} M_{\mathcal{B}, 1} |g|\,,
        $$
        and by \Cref{local-tree-control}
        $$
            \le \inf_{B^\circ{}(J)} |T_{\fT(\fu_2)}^* g| + (2^{104a^3} + 2^{128a^3+4a+3}) \inf_{J} M_{\mathcal{B}, 1} |g|\,.
        $$
        This completes the proof.
    \end{proof}

    \begin{proof}[Proof of \Cref{Holder-correlation-tree}]
        \leanok
        \proves{Holder-correlation-tree}
        Let $P$ be the product on the right hand side of \eqref{hHolder}, and $h_J$ be as defined in \eqref{def-hj}.
        By \eqref{eq-pao-2} and \Cref{adjoint-tile-support}, the function $h_J$ is supported in $B'(J) \cap \scI(\fu_1)$.
        By \eqref{eq-pao-2}, \Cref{global-tree-control-1} and \Cref{global-tree-control-2}, we have for all $y \in B'(J)$:
        $$
            |h_J(y)| \le 2^{257a^3+4a+3} P\,.
        $$
        We have by the triangle inequality
        \begin{align}
            &|h_J(y) - h_J(y')|\nonumber\\
            \label{eq-h-Lip-1}
            &\le |\chi_J(y) - \chi_J(y')| |T_{\fT(\fu_1)}^* g_1(y)| |T_{\fT(\fu_2) \cap \mathfrak{S}}^* g_2(y)|\\
            \label{eq-h-Lip-2}
            & + |\chi_J(y')| |e(\fcc(\fu_1)(y)) T_{\fT(\fu_1)}^* g_1(y) - e(\fcc(\fu_1)(y')) T_{\fT(\fu_1)}^* g_1(y')| |T_{\fT(\fu_2) \cap \mathfrak{S}}^* g_2(y)|\\
            \label{eq-h-Lip-3}
            & + |\chi_J(y')| |T_{\fT(\fu_1)}^* g_1(y')| |e(\fcc(\fu_2)(y)) T_{\fT(\fu_2) \cap \mathfrak{S}}^* g_2(y) - e(\fcc(\fu_2)(y')) T_{\fT(\fu_2) \cap \mathfrak{S}}^* g_2(y')|\,.
        \end{align}

        As $h_J$ is supported in $\scI(\fu_1)$, we can assume without loss of generality that $y' \in \scI(\fu_1)$.
        If $y \notin \scI(\fu_1)$, then \eqref{eq-h-Lip-1} vanishes. If $y \in \scI(\fu_1)$ then we have by \eqref{eq-pao-3}, \Cref{global-tree-control-1} and \Cref{global-tree-control-2}
        $$
            \eqref{eq-h-Lip-1} \le 2^{484a^3+4a+3} \frac{\rho(y,y')}{D^{s(J)}} P\,,
        $$
        where $P$ denotes the product on the right hand side of \eqref{hHolder}.

        By \eqref{eq-pao-2}, \Cref{global-tree-control-1} and \Cref{global-tree-control-2}, we have
        $$
            \eqref{eq-h-Lip-2} \le 2^{257a^3+4a+1} \left(\frac{\rho(y,y')}{D^{s(J)}}\right)^{1/a} P\,.
        $$

        By \eqref{eq-pao-2}, and twice \Cref{global-tree-control-1}, we have
        $$
            \eqref{eq-h-Lip-3} \le 2^{256a^3+8a+5} \left(\frac{\rho(y,y')}{D^{s(J)}}\right)^{1/a} P\,.
        $$
        Using that $\rho(y,y') \le 32D^{s(J)}$ and $a \ge 4$, the lemma follows.
    \end{proof}

\subsubsection{The van der Corput estimate}
\label{subsubsec-van-der-corput}
    \begin{lemma}[lower oscillation bound]
        \label{lower-oscillation-bound}
        \leanok
        \lean{TileStructure.Forest.lower_oscillation_bound}
        \uses{overlap-implies-distance}
        For all $J \in \mathcal{J}'$, we have that
        $$
            d_{B(J)}(\fcc(\fu_1), \fcc(\fu_2)) \ge 2^{-201a^3} 2^{Zn/2}\,.
        $$
    \end{lemma}

    \begin{proof}
    \leanok
    Since $\emptyset \ne \fT(\fu_1) \subset \mathfrak{S}$ by \Cref{overlap-implies-distance}, there exists at least one tile $\fp \in \mathcal{S}$ with $\scI(\fp) \subsetneq \scI(\fu_1)$. Thus $\scI(\fu_1) \notin \mathcal{J}'$, so $J \subsetneq \scI(\fu_1)$. Thus there exists a cube $J' \in \mathcal{D}$ with $J \subset J'$ and $s(J') = s(J) + 1$, by \eqref{coverdyadic} and \eqref{dyadicproperty}. By definition of $\mathcal{J'}$ and the triangle inequality, there exists $\fp \in \mathfrak{S}$ such that
    $$
        \scI(\fp) \subset B(c(J'), 100 D^{s(J') + 1}) \subset B(c(J), 128 D^{s(J)+2})\,.
    $$
    Thus, by definition of $\mathfrak{S}$:
    \begin{align*}
        2^{Zn/2} \le d_{\fp}(\fcc(\fu_1), \fcc(\fu_2)) \le d_{B(c(J), 128 D^{s(J)+2})}(\fcc(\fu_1), \fcc(\fu_2))\,.
    \end{align*}
    By the doubling property \eqref{firstdb}, this is
    $$
        \le 2^{200a^3 + 4a} d_{B(J)}(\fcc(\fu_1), \fcc(\fu_2))\,,
    $$
    which gives the lemma using $a \ge 4$.
    \end{proof}

    Now we are ready to prove \Cref{correlation-distant-tree-parts}.
    \begin{proof}[Proof of \Cref{correlation-distant-tree-parts}]
    \leanok
    \proves{correlation-distant-tree-parts}
    We have
    $$
        \eqref{eq-lhs-big-sep-tree} = \left| \int_{X} T_{\fT(\fu_1)}^* g_1 \overline{T_{\fT(\fu_2) \cap \mathfrak{S}}^* g_2 }\right|\,.
    $$
    By \Cref{adjoint-tile-support}, the right hand side is supported in $\scI(\fu_1)$. Using \eqref{eq-pao-1} of \Cref{Lipschitz-partition-unity} and the definition \eqref{def-hj} of $h_J$, we thus have
    $$
        \le \sum_{J \in \mathcal{J}'} \left|\int_{X} e(\fcc(\fu_2)(y) - \fcc(\fu_1)(y)) h_J(y) \, \mathrm{d}\mu(y) \right|\,.
    $$
    Using \Cref{Holder-van-der-Corput} with the ball $B(J)$, we bound this by
    $$
        \le 2^{7a} \sum_{J \in \mathcal{J}'} \mu(B(J)) \|h_J\|_{C^{\tau}(B'(J))} (1 + d_{B(J)}(\fcc(\fu_2), \fcc(\fu_1)))^{-1/(2a^2+a^3)}\,.
    $$
    Using \Cref{Holder-correlation-tree}, \Cref{lower-oscillation-bound} and $a \ge 4$, we bound the above by
    \begin{multline}
        \label{eq-big-sep-1}
        \le 2^{485a^3+7a+3a^3+3a} 2^{-Zn/(4a^2 + 2a^3)} \sum_{J \in \mathcal{J}'} \mu(B(J)) \\
        \times \prod_{j=1}^2 (\inf_{B^\circ{}(J)} |T_{\fT(\fu_j)}^* g_j| + \inf_J M_{\mathcal{B},1} g_j)\,.
    \end{multline}
    By the doubling property \eqref{doublingx}
    $$
        \mu(B(J)) \le 2^{6a} \mu(B^\circ{}(J))\,,
    $$
    thus
    $$
        \mu(B(J)) \prod_{j=1}^2 (\inf_{B^\circ{}(J)} |T_{\fT(\fu_j)}^* g_j| + \inf_J M_{\mathcal{B},1} g_j)
    $$
    $$
        \le 2^{6a} \int_{B^\circ{}(J)} \prod_{j=1}^2 ( |T_{\fT(\fu_j)}^* g_j|(x) + M_{\mathcal{B},1} g_j(x)) \, \mathrm{d}\mu(x)
    $$
    $$
        \le 2^{6a} \int_J \prod_{j=1}^2 ( |T_{\fT(\fu_j)}^* g_j|(x) + M_{\mathcal{B},1} g_j(x)) \, \mathrm{d}\mu(x)\,.
    $$
    Summing over $J \in \mathcal{J}'$, we obtain
    $$
        \eqref{eq-big-sep-1} \le 2^{499a^3} 2^{-Zn/(4a^2 + 2a^3)} \int_X \prod_{j=1}^2 ( |T_{\fT(\fu_j)}^* g_j|(x) + M_{\mathcal{B},1} g_j(x)) \, \mathrm{d}\mu(x)\,.
    $$
    Applying the Cauchy-Schwarz inequality, \Cref{correlation-distant-tree-parts} follows.
    \end{proof}

\subsection{Proof of The Remaining Tiles Lemma}
    \label{subsec-rest-tiles}
    We define
    $$
        \mathcal{J}' = \{J \in \mathcal{J}(\fT(\fu_1)) \, : \, J \subset \scI(\fu_1)\}\,,
    $$
    note that this is different from the $\mathcal{J}'$ defined in the previous subsection.
    \begin{lemma}[dyadic partition 2]
        \label{dyadic-partition-2}
        \leanok
        \lean{TileStructure.Forest.union_𝓙₆, TileStructure.Forest.pairwiseDisjoint_𝓙₆}
        \uses{dyadic-partitions}
        We have
        $$
            \scI(\fu_1) = \dot{\bigcup_{J \in \mathcal{J}'}} J\,.
        $$
    \end{lemma}

    \begin{proof}\leanok
        By \Cref{dyadic-partitions}, it remains only to show that each $J \in \mathcal{J}(\fT(\fu_1))$ with $J \cap \scI(\fu_1) \ne \emptyset$ is in $\mathcal{J}'$. But if $J \notin \mathcal{J}'$, then by \eqref{dyadicproperty} $\scI(\fu_1) \subsetneq J$. Pick $\fp \in \fT(\fu_1)$. Then $\scI(\fp) \subsetneq J$. This contradicts the definition of $\mathcal{J}(\fT(\fu_1))$.
    \end{proof}

    \Cref{correlation-near-tree-parts} follows from the following key estimate.

    \begin{lemma}[bound for tree projection]
        \label{bound-for-tree-projection}
        \leanok
        \lean{TileStructure.Forest.bound_for_tree_projection}
        \uses{adjoint-tile-support,overlap-implies-distance,dyadic-partition-2,thin-scale-impact,square-function-count}
        We have for all bounded $f$ with bounded support
        $$
            \|P_{\mathcal{J}'}|T_{\fT(\fu_2) \setminus \mathfrak{S}}^* g_2|\|_2
            \le 2^{102a^3+21a+5} 2^{-\frac{25}{101a}Zn\kappa} \|\mathbf{1}_{\scI(\fu_1)} M_{\mathcal{B},1} |g_2|\|_2
        $$
    \end{lemma}

    We prove this lemma below. First, we deduce \Cref{correlation-near-tree-parts}.

    \begin{proof}[Proof of \Cref{correlation-near-tree-parts}]
        \leanok
        \proves{correlation-near-tree-parts}
        By \Cref{tree-projection-estimate} and \Cref{adjoint-tile-support}, we have
        \begin{align*}
            \eqref{eq-lhs-small-sep-tree} \le 2^{130a^3} \|P_{\mathcal{L}(\fT(\fu_1))} |\mathbf{1}_{\scI(\fu_1)}g_1| \|_2 \|P_{\mathcal{J}(\fT(\fu_1) )}|\mathbf{1}_{\scI(\fu_1)} T_{\fT(\fu_2) \setminus \mathfrak{S}}^* g_2|\|_2\,.
        \end{align*}
        It follows from the definition of the projection operator $P$ and Jensen's inequality that
        $$
            \|P_{\mathcal{L}(\fT(\fu_1))} |g_1\mathbf{1}_{\scI(\fu_1)}| \|_2 \le \|g_1 \mathbf{1}_{\scI(\fu_1)}\|_2\,.
        $$
        By \Cref{dyadic-partition-2}, a cube $J \in \mathcal{J}(\fT(\fu_1))$ intersects $\scI(\fu_1)$ if and only if $J \in \mathcal{J}'$. Thus
        $$
            P_{\mathcal{J}(\fT(\fu_1) )}|\mathbf{1}_{\scI(\fu_1)} T_{\fT(\fu_2) \setminus \mathfrak{S}}^* g_2| = P_{\mathcal{J}'}|T_{\fT(\fu_2) \setminus \mathfrak{S}}^* g_2|\,.
        $$
        Combining this with \Cref{bound-for-tree-projection}, the definition \eqref{definekappa} and $a \ge 4$ proves the lemma.
    \end{proof}

    We need two more auxiliary lemmas before we prove \Cref{bound-for-tree-projection}.

    \begin{lemma}[thin scale impact]
        \label{thin-scale-impact}
        \leanok
        \lean{TileStructure.Forest.thin_scale_impact}
        If $\fp \in \fT(\fu_2) \setminus \mathfrak{S}$ and $J \in \mathcal{J'}$ with $B(\scI(\fp)) \cap B(J) \ne \emptyset$, then
        $$
            \ps(\fp) \le s(J) + 2 - \frac{Zn}{202a^3}\,.
        $$
    \end{lemma}

    \begin{proof}
        \leanok
        Suppose that $\ps(\fp) > s(J) + 2 -\frac{Zn}{202a^3} =: s(J) - s_1$. Then, we have $s_1 + 2 \ge 0$ so
        $$
            \rho(\pc(\fp), c(J)) \le 8D^{s(J)}+8D^{\ps(\fp)} \le 16 D^{\ps(\fp) + s_1 + 2}\,.
        $$
        There exists a tile $\fq \in \fT(\fu_1)$. By \eqref{forest1}, it satisfies $\scI(\fq) \subsetneq \scI(\fu_1)$. Thus $\scI(\fu_1) \notin \mathcal{J}'$. It follows that $J \subsetneq \scI(\fu_1)$. By \eqref{coverdyadic} and \eqref{dyadicproperty}, there exists a cube $J' \in \mathcal{D}$ with $J \subset J'$ and $s(J') = s(J) + 1$. By definition of $\mathcal{J}'$, there exists a tile $\fp' \in \fT(\fu_1)$ with
        $$
            \scI(\fp') \subset B(c(J'), 100 D^{s(J') + 1})\,.
        $$
        By the triangle inequality, the definition \eqref{defineD} and $a \ge 4$, we have
        $$
            B(c(J'), 100 D^{s(J')+1}) \subset B(\pc(\fp), 128 D^{\ps(\fp) + s_1 + 2})\,.
        $$
        Since $\fp' \in \fT(\fu_1)$ and $\scI(\fu_1) \subset \scI(\fu_2)$, we have by \eqref{forest5}
        $$
            d_{\fp'}(\fcc(\fp'), \fcc(\fu_2)) > 2^{Z(n+1)}\,.
        $$
        Hence, by \eqref{forest1}, the triangle inequality and using that by \eqref{defineZ} $Z(n+1) = 2^{12a}(n+1) \ge 3$
        $$
            d_{\fp'}(\fcc(\fu_1), \fcc(\fu_2)) > 2^{Z(n+1)} - 4 \ge 2^{Z(n+1) - 1}\,.
        $$
        It follows that
        $$
            2^{Z(n+1)-1} \le d_{\fp'}(\fcc(\fu_1), \fcc(\fu_2)) \le d_{B(\pc(\fp), 128 D^{\ps(\fp) + s_1+ 2})}(\fcc(\fu_1), \fcc(\fu_2))\,.
        $$
        Using \eqref{firstdb}, we obtain
        $$
            \le 2^{9a + 100a^3 (s_1+3)} d_{\fp}(\fcc(\fu_1), \fcc(\fu_2))\,.
        $$
        Since $\fp' \notin \mathfrak{S}$ this is bounded by
        $$
            \le 2^{9a + 100a^3 (s_1+3)} 2^{Zn/2}\,.
        $$
        Thus
        $$
            Z n/2 + Z - 1 \le 9a + 100a^3(s_1 + 3)\,,
        $$
        contradicting the definition of $s_1$.
    \end{proof}

    \begin{lemma}[square function count]
        \label{square-function-count}
        \leanok
        \lean{TileStructure.Forest.square_function_count}
        For each $J \in \mathcal{J}'$ and all $s$, we have
        $$
            \frac{1}{\mu(J)} \int_J \Bigg(\sum_{\substack{I \in \mathcal{D}, s(I) = s(J) - s\\ I \cap \scI(\fu_1) = \emptyset\\
        J \cap B(I) \ne \emptyset}} \mathbf{1}_{B(I)}\bigg)^2 \, \mathrm{d}\mu \le 2^{14a+1} (8 D^{-s})^\kappa\,.
        $$
    \end{lemma}

    \begin{proof}
        \leanok
        Since $J \in \mathcal{J}'$ we have $J \subset \scI(\fu_1)$. Thus, if $B(I) \cap J \ne \emptyset$ then
    \begin{equation}
        \label{eq-sep-small-incl}
        B(I) \cap J \subset \{x \in J \ : \ \rho(x, X \setminus J) \le 8D^{s(I)}\}\,.
    \end{equation}
    Furthermore, for each $s$ the balls $B(I)$ with $s(I) = s$ have bounded overlap: Consider the collection $\mathcal{D}_{s,x}$ of all $I \in \mathcal{D}$ with $x \in B(I)$ and $s(I) = s$. By \eqref{eq-vol-sp-cube} and \eqref{dyadicproperty}, the balls $B(c(I), \frac{1}{4} D^{s(I)})$, $I \in \mathcal{D}_{s,x}$ are disjoint, and by the triangle inequality, they are contained in $B(x, 9 D^{s})$. By the doubling property \eqref{doublingx}, we have
    $$
        \mu(B(x, 9D^{s})) \le 2^{7a} \mu(B(c(I), \frac{1}{4} D^{s(I)}))
    $$
    for each $I \in \mathcal{D}_{s,x}$.
    Thus
    $$
        \mu(B(x, 9D^{s})) \ge \sum_{I \in \mathcal{D}_{s,x}} \mu(B(c(I), \frac{1}{4} D^{s(I)})) \ge 2^{-7a} |\mathcal{D}_{s,x}| \mu(B(x, 9D^{s}))\,.
    $$
    Dividing by the positive $\mu(B(x, 9D^{s}))$, we obtain that for each $x$
    \begin{equation}
        \label{eq-sep-small-bound}
        \Bigg(\sum_{\substack{I \in \mathcal{D}, s(I) = s(J) - s\\ I \cap \scI(\fu_1) = \emptyset\\
        J \cap B(I) \ne \emptyset}} \mathbf{1}_{B(I)}(x) \bigg)^2 \le |\mathcal{D}_{s(J) - s,x}|^2 \le 2^{14a} \,.
    \end{equation}
    Combining \eqref{eq-sep-small-incl}, \eqref{eq-sep-small-bound} and the small boundary property \eqref{eq-small-boundary}, noting that $8D^{s(I)}=8D^{-s}D^{s(J)}$, the lemma follows.
\end{proof}

\begin{proof}[Proof of \Cref{bound-for-tree-projection}]
    \leanok
    \proves{bound-for-tree-projection}
    Expanding the definition of $P_{\mathcal{J}'}$, we have
    $$
        \|P_{\mathcal{J}'}|T_{\fT(\fu_2) \setminus \mathfrak{S}}^* g_2|\|_2
    $$
    $$
        = \left(\sum_{J \in \mathcal{J}'} \frac{1}{\mu(J)} \left(\int_J \left| \sum_{\fp \in \fT(\fu_2) \setminus \mathfrak{S}} T_{\fp}^* g_2(y) \right| \, \mathrm{d}\mu(y) \right)^2 \right)^{1/2}\,.
    $$
    By \Cref{adjoint-tile-support}, the innermost sum in the last display is $0$ if $J \cap B(\scI(\fp)) = \emptyset$.
    Then we split that sum according to the scale of $\fp$ relative to the scale of $J$.
    By \Cref{thin-scale-impact}, $s_1 \le s(J) - \ps(\fp) \le 2S$ with $s_1 := \lfloor\frac{Zn}{202a^3} - 2\rfloor$:
    $$
        = \left(\sum_{J \in \mathcal{J}'} \frac{1}{\mu(J)} \left(\int_J \left| \sum_{s=s_1}^{2S} \sum_{\substack{\fp \in \fT(\fu_2) \setminus \mathfrak{S}\\ \ps(\fp) = s(J) - s\\
        J \cap B(\scI(\fp)) \ne \emptyset}} T_{\fp}^* g_2(y) \right| \, \mathrm{d}\mu(y) \right)^2 \right)^{1/2}\,.
    $$
    Then we apply the triangle inequality and Minkowski's inequality to get
    \begin{equation}
    \label{eq-sep-tree-small-1}
        \le \sum_{s=s_1}^{2S} \left( \sum_{J \in \mathcal{J}'} \frac{1}{\mu(J)} \left(\int_J \sum_{\substack{\fp \in \fT(\fu_2) \setminus \mathfrak{S}\\ \ps(\fp) = s(J) - s\\
        J \cap B(\scI(\fp)) \ne \emptyset}} |T_{\fp}^* g_2(y)| \, \mathrm{d}\mu(y) \right)^2\right)^{1/2}\,.
    \end{equation}
    We have by \Cref{adjoint-tile-support} and \eqref{eq-Ks-size}
    $$
        |T_{\fp}^* g_2(y)| \le 2^{102a^3} \mathbf{1}_{B(\pc(\fp), 8D^{\ps(\fp)})}(y) \int_{E(\fp)} \frac{1}{\mu(B(x, D^{\ps(\fp)}))} |g_2(x)| \, \mathrm{d}\mu(x)\,.
    $$
    Using the doubling property \eqref{doublingx}, it follows that
    $$
        \mu(B(\pc(\fp), 8D^{\ps(\fp)})) \le 2^{4a} \mu(B(x, D^{\ps(\fp)}))\,.
    $$
    Thus, using also $a \ge 4$
    $$
        |T_{\fp}^* g_2(y)| \le 2^{102a^3+4a} \mathbf{1}_{B(\pc(\fp), 8D^{\ps(\fp)})}(y) \frac{1}{\mu(B(\pc(\fp), 8D^{\ps(\fp)}))} \int_{E(\fp)} |g_2(x)| \, \mathrm{d}\mu(x)\,.
    $$
    Since for each $I \in \mathcal{D}$ the sets $E(\fp), \fp \in \fP(I)$ are disjoint, it follows that
    $$
        \int_J \sum_{\substack{\fp \in \fT(\fu_2) \setminus \mathfrak{S}\\ \scI(\fp) = I\\
        J \cap B(\scI(\fp)) \ne \emptyset}} |T_{\fp}^* g_2(y)| \, \mathrm{d}\mu
    $$
    $$
        \le 2^{102a^3+4a} \int_J \mathbf{1}_{B(I)} \frac{1}{\mu(B(\pc(\fp), 8D^{\ps(\fp)}))} \int_{B(\pc(\fp), 8D^{\ps(\fp)})} |g_2(x)| \, \mathrm{d}\mu(x)
    $$
    $$
        \le 2^{102a^3+4a} \int_J M_{\mathcal{B},1} |g_2|(y) \mathbf{1}_{B(I)}(y) \, \mathrm{d}\mu(y)\,.
    $$
    By \Cref{overlap-implies-distance}, we have $\scI(\fp) \cap \scI(\fu_1) = \emptyset$ for all $\fp \in \fT(\fu_2) \setminus \mathfrak{S}$.
    Thus we can estimate \eqref{eq-sep-tree-small-1} by
    $$
        2^{102a^3+4a} \sum_{s=s_1}^{2S} \left( \sum_{J \in \mathcal{J}'} \frac{1}{\mu(J)} \left(\int_J \sum_{\substack{I \in \mathcal{D}, s(I) = s(J) - s\\ I \cap \scI(\fu_1) = \emptyset\\
        J \cap B(I) \ne \emptyset}} M_{\mathcal{B},1} |g_2| \mathbf{1}_{B(I)} \, \mathrm{d}\mu \right)^2\right)^{1/2}\,,
    $$
    which is by Cauchy-Schwarz at most
    \begin{equation}
    \label{eq-sep-tree-small-2}
        2^{102a^3+4a} \sum_{s=s_1}^{2S} \left( \sum_{J \in \mathcal{J}'} \int_J ( M_{\mathcal{B},1} |g_2|)^2 \, \mathrm{d}\mu \times
        \frac{1}{\mu(J)} \int_J \left(\sum_{\substack{I \in \mathcal{D}, s(I) = s(J) - s\\ I \cap \scI(\fu_1) = \emptyset\\
        J \cap B(I) \ne \emptyset}} \mathbf{1}_{B(I)}\right)^2 \, \mathrm{d}\mu \right)^{1/2}\,.
    \end{equation}
    Using \Cref{square-function-count}, we bound \eqref{eq-sep-tree-small-2} by
    $$
        2^{102a^3+4a} \sum_{s=s_1}^{2S} \left(2^{14a+1} (8 D^{-s})^\kappa \sum_{J \in \mathcal{J}'} \int_J (M_{\mathcal{B},1} |g_2|)^2\right)^{1/2}\,,
    $$
    and since dyadic cubes in $\mathcal{J}'$ form a partition of $\scI(\fu_1)$ by \Cref{dyadic-partition-2}, $\kappa \le 1$ by \eqref{definekappa}, and $a \ge 4$
    $$
        \le 2^{102a^3+11a+2} \sum_{s=s_1}^{2S} D^{-s\kappa/2} \|\mathbf{1}_{\scI(\fu_1)} M_{\mathcal{B},1} |g_2|\|_2
    $$
    $$
        \le 2^{102a^3+11a+2} D^{-s_1 \kappa /2} \frac{1}{1 - D^{-\kappa/2}} \|\mathbf{1}_{\scI(\fu_1)} M_{\mathcal{B},1} |g_2|\|_2\,.
    $$
    By convexity of $t \mapsto D^{-t}$ and since $D \ge 2$, we have for all $0 \le t \le 1$
    $$
        D^{-t} \le 1 - t(1 - 1/D) \le 1 - \frac{1}{2}t\,.
    $$
    Using this for $t = \kappa/2$ and using that $s_1 = \frac{Zn}{202a^3} - 2$ and the definitions \eqref{defineD} and \eqref{definekappa} of $\kappa$ and $D$
    $$
        \le 2^{102a^3+11a+2} 2^{-100a^2(\frac{Zn}{202a^3} - 3) \frac{\kappa}{2}} \frac{2}{\kappa} \|\mathbf{1}_{\scI(\fu_1)} M_{\mathcal{B},1} |g_2|\|_2
    $$
    $$
        = 2^{102a^3+21a+4} 2^{150a^2\kappa} 2^{-\frac{100}{404a}Zn\kappa} \|\mathbf{1}_{\scI(\fu_1)} M_{\mathcal{B},1} |g_2|\|_2\,.
    $$
    Using the definition \eqref{definekappa} of $\kappa$ and $a \ge 4$, the lemma follows.
\end{proof}

\subsection{Forests}
In this subsection, we complete the proof of \Cref{forest-operator} from the results of the previous subsections.

Define an $n$-row to be an $n$-forest $(\fU, \fT)$, i.e. satisfying conditions \eqref{forest1} - \eqref{forest6}, such that in addition the sets $\scI(\fu), \fu \in \fU$ are pairwise disjoint.

\begin{lemma}[forest row decomposition]
    \label{forest-row-decomposition}
    \leanok
    \lean{TileStructure.Forest.rowDecomp, TileStructure.Forest.biUnion_rowDecomp,
    TileStructure.Forest.pairwiseDisjoint_rowDecomp}
    Let $(\fU, \fT)$ be an $n$-forest. Then there exists a decomposition
    $$
        \fU = \dot{\bigcup_{1 \le j \le 2^n}} \fU_j
    $$
    such that for all $j = 1, \dotsc, 2^n$ the pair $(\fU_j, \fT|_{\fU_j})$ is an $n$-row.
\end{lemma}

\begin{proof}
    \leanok
    Define recursively $\fU_j$ to be a maximal disjoint set of tiles $\fu$ in
    $$
        \fU \setminus \bigcup_{j' < j} \fU_{j'}
    $$
    with inclusion maximal $\scI(\fu)$. Properties \eqref{forest1}, -\eqref{forest6} for $(\fU_j, \fT|_{\fU_k})$ follow immediately from the corresponding properties for $(\fU, \fT)$, and the cubes $\scI(\fu), \fu \in \fU_j$ are disjoint by definition. The collections $\fU_j$ are also disjoint by definition.

    Now we show by induction on $j$ that each point is contained in at most $2^n - j$ cubes $\scI(\fu)$ with $\fu \in \fU \setminus \bigcup_{j' \le j} \fU_{j'}$. This implies that $\bigcup_{j = 1}^{2^n} \fU_j = \fU$, which completes the proof of the Lemma. For $j = 0$ each point is contained in at most $2^n$ cubes by \eqref{forest3}. For larger $j$, if $x$ is contained in any cube $\scI(\fu)$ with $\fu \in \fU \setminus \bigcup_{j' < j} \fU_{j'}$, then it is contained in a maximal such cube. Thus it is contained in a cube in $\scI(\fu)$ with $\fu \in \fU_j$. Thus the number $\fu \in \fU \setminus \bigcup_{j' \le j} \fU_{j'}$ with $x\in \scI(\fu)$ is zero, or is less than the number of $\fu \in \fU \setminus \bigcup_{j' \le j-1} \fU_{j'}$ with $x \in \scI(\fu)$ by at least one.
\end{proof}

We pick a decomposition of the forest $(\fU, \fT)$ into $2^n$ $n$-rows
\begin{equation*}
(\fU_j, \fT_j) := (\fU_j, \fT|_{\fU_j})
\end{equation*}
as in \Cref{forest-row-decomposition}. To save some space in the proofs of the remaining lemmas in this section we will write
    $$
        T_{\fC} = \sum_{\fp \in \fC} T_{\fp},\qquad T_{\fC}^* = \sum_{\fp \in \fC} T_{\fp}^*,
    $$
    $$
        T_{\mathfrak{R}_j} = \sum_{\fu \in \fU_j} T_{\fT(\fu)},\qquad T_{\mathfrak{R}_j}^* = \sum_{\fu \in \fU_j} T_{\fT(\fu)}^*.
    $$
\begin{lemma}[row bound]
    \label{row-bound}
    \leanok
    \lean{TileStructure.Forest.row_bound, TileStructure.Forest.indicator_row_bound}
    \uses{densities-tree-bound,adjoint-tile-support}
    For each $1 \le j \le 2^n$ and each bounded $g$ with bounded support with $\text{support}(g) \subseteq G$, we have
    \begin{equation}
        \label{eq-row-bound-1}
        \left\| T_{\mathfrak{R}_j}^*g \right\|_2 \le 2^{182a^3} 2^{-n/2} \|g\|_2
    \end{equation}
    and
    \begin{equation}
        \label{eq-row-bound-2}
        \left\| \mathbf{1}_F T_{\mathfrak{R}_j}^*g \right\|_2 \le 2^{283a^3} 2^{-n/2} \dens_2(\bigcup_{\fu\in \fU}\fT(\fu))^{1/2} \|g\|_2\,.
    \end{equation}
\end{lemma}

\begin{proof}
    \leanok
    By \Cref{densities-tree-bound} and the density assumption \eqref{forest4}, we have for each $\fu \in \fU$ and all bounded $f$ of bounded support that
    \begin{equation}
        \label{eq-explicit-tree-bound-1}
        \left\|\mathbf{1}_G\sum_{\fp \in \fT(\fu)} T_{\fp} f \right\|_{2} \le 2^{181a^3} 2^{(4a+1-n)/2} \|f\|_2\,
    \end{equation}
    and
    \begin{equation}
        \label{eq-explicit-tree-bound-2}
        \left\|\mathbf{1}_G \sum_{\fp \in \fT(\fu)} T_{\fp} \mathbf{1}_F f \right\|_{2} \le 2^{282a^3} 2^{(4a + 1-n)/2} \dens_2(\fT(\fu))^{1/2} \|f\|_2\,.
    \end{equation}
    Since for each $j$ the top cubes $\scI(\fu)$, $\fu \in \fU_j$ are disjoint, we further have for all bounded $g$ of bounded support and $|g| \le \mathbf{1}_G$ by \Cref{adjoint-tile-support}
    $$
        \left\|\mathbf{1}_F \sum_{\fu \in \fU_j} \sum_{\fp \in \fT(\fu)} T_{\fp}^* g\right\|_2^2 = \left\|\mathbf{1}_F \sum_{\fu \in \fU_j} \sum_{\fp \in \fT(\fu)} \mathbf{1}_{\scI(\fu)} T_{\fp}^* \mathbf{1}_{\scI(\fu)} g\right\|_2^2
    $$
    $$
        = \sum_{\fu \in \fU_j} \int_{\scI(\fu)} \left| \mathbf{1}_F \sum_{\fp \in \fT(\fu)} T_{\fp}^* \mathbf{1}_{\scI(\fu)} g\right|^2 \, \mathrm{d}\mu
        \le \sum_{\fu \in \fU_j} \left\|\sum_{\fp \in \fT(\fu)} \mathbf{1}_F T_{\fp}^* \mathbf{1}_G \mathbf{1}_{\scI(\fu)}  g\right\|_2^2\,.
    $$
    Applying the estimate for the adjoint operator following from equation \eqref{eq-explicit-tree-bound-2}, we obtain
    $$
        \le 2^{282a^3} 2^{(4a+1-n)/2} \max_{\fu \in \fU_j}\dens_2(\fT(\fu))^{1/2} \sum_{\fu \in \fU_j} \left\| \mathbf{1}_{\scI(\fu)} g\right\|_2^2\,.
    $$
    Again by disjointedness of the cubes $\scI(\fu)$, this is estimated by
    $$
        2^{282a^3} 2^{(4a+1-n)/2} \max_{\fu \in \fU_j}\dens_2(\fT(\fu))^{1/2} \|g\|_2^2\,.
    $$
    Thus, \eqref{eq-row-bound-2} follows, since $a \ge 4$.
    The proof of \eqref{eq-row-bound-1} from \eqref{eq-explicit-tree-bound-1} is the same up to replacing $F$ by $X$.
\end{proof}

Because all involved operators are linear, the inequality
\eqref{eq-row-bound-2} is also true for uniformly bounded functions $g$
supported on $G$. This generalization will be needed for
\Cref{forest-operator}.

\begin{lemma}[row correlation]
    \label{row-correlation}
    \leanok
    \lean{TileStructure.Forest.row_correlation}
    \uses{adjoint-tree-control,correlation-separated-trees}
    For all $1 \le j,j' \le 2^n$ with $j\ne j'$ and for all $g_1, g_2$ with $\text{support}(g_i) \subseteq G$, it holds that
    $$
        \left| \int T_{\mathfrak{R}_j}^*g_1 \overline{T_{\mathfrak{R}_{j'}}^*g_2} \, \mathrm{d}\mu \right| \le
        2^{876a^3-4n}\|g_1\|_2 \|g_2\|_2\,.
    $$
\end{lemma}

\begin{proof}
    \leanok
    We have by \Cref{adjoint-tile-support} and the triangle inequality that
    $$
        \left| \int T_{\mathfrak{R}_j}^*g_1 \overline{T_{\mathfrak{R}_{j'}}^*g_2} \, \mathrm{d}\mu \right|
    $$
    $$
        \le \sum_{\fu \in \fU_j} \sum_{\fu' \in \fU_{j'}} \left| \int T^*_{\fT_j(\fu)} (\mathbf{1}_{\scI(\fu)} g_1) \overline{T^*_{\fT_{j'}(\fu')} (\mathbf{1}_{\scI(\fu')} g_2)} \, \mathrm{d}\mu \right|\,.
    $$
    By \Cref{correlation-separated-trees}, this is bounded by
    \begin{equation}
        \label{eq-S2uu'}
         2^{512a^3-4n} \sum_{\fu \in \fU_j} \sum_{\fu' \in \fU_{j'}} \|S_{2,\fu} (\mathbf{1}_{\scI(\fu)}g_1)\|_{L^2(\scI(\fu')\cap \scI(\fu)} \|S_{2, \fu'} (\mathbf{1}_{\scI(\fu')}g_2)\|_{L^2(\scI(\fu')\cap\scI(\fu))}\,.
    \end{equation}
    We apply the Cauchy-Schwarz inequality in the form
    \begin{equation*}
        \sum_{i \in M} a_i b_i \le (\sum_{i \in M} a_i^2 )^{1/2}(\sum_{i \in M} b_i^2 )^{1/2}
    \end{equation*} to the outer two sums:
    $$
        \le 2^{512a^3-4n} \left(\sum_{\fu \in \fU_j} \sum_{\fu' \in \fU_{j'}} \|S_{2,\fu} (\mathbf{1}_{\scI(\fu)}g_1)\|_{L^2(\scI(\fu')\cap \scI(\fu))}^2 \right)^{1/2}
    $$
    $$
        \left(\sum_{\fu \in \fU_j} \sum_{\fu' \in \fU_{j'}} \|S_{2,\fu'} (\mathbf{1}_{\scI(\fu')}g_2)\|_{L^2(\scI(\fu')\cap\scI(\fu))}^2 \right)^{1/2}\,.
    $$
    We can now estimate the factor involving $g_1$ as follows:
    \begin{multline*}
        \sum_{\fu \in \fU_j}\sum_{\fu' \in \fU_{j'}} \|S_{2,\fu} (\mathbf{1}_{\scI(\fu)}g_1)\|_{L^2(\scI(\fu')\cap \scI(\fu))}^2
        \\ = \sum_{\fu \in \fU_j}\sum_{\fu' \in \fU_{j'}} \int_{\scI(\fu) \cap \scI(\fu')} |S_{2,\fu} (\mathbf{1}_{\scI(\fu)}(y)g_1(y))|^2 \, \mathrm{d}\mu(y)
    \end{multline*}
    By pairwise disjointedness of the sets $\scI(\fu')$ for $\fu' \in \fU_{j'}$, we have
    $$
        \le \sum_{\fu \in \fU_j}\int_{\scI(\fu)} |S_{2,\fu} (\mathbf{1}_{\scI(\fu)}(y)g_1(y))|^2 \, \mathrm{d}\mu(y)
        \le \sum_{\fu \in \fU_j}\int_{X} |S_{2,\fu} (\mathbf{1}_{\scI(\fu)}(y)g_1(y))|^2 \, \mathrm{d}\mu(y)
    $$
    $$
       = \sum_{\fu \in \fU_j}\|S_{2,\fu} (\mathbf{1}_{\scI(\fu)}g_1)\|_2^2
    $$
    By \Cref{adjoint-tree-control} we now estimate:
    $$
        \le \sum_{\fu \in \fU_j}(2^{182a^3})^2 \|\mathbf{1}_{\scI(\fu)}g_1\|_2^2
    $$
    By pairwise disjointedness of the sets $\scI(\fu)$ for $\fu \in \fU_j$ (and writing out the definition of $L^2$-norms), we have
    $$
        \le (2^{182a^3})^2 \|g_1\|_2^2
    $$
    Arguing similarly for $g_2$, we obtain the desired inequality.
\end{proof}

Define for $1 \le j \le 2^n$
$$
    E_j := \bigcup_{\fu \in \fU_j} \bigcup_{\fp \in \fT(\fu)} E(\fp)\,.
$$

\begin{lemma}[disjoint row support]
    \label{disjoint-row-support}
    \leanok
    \lean{TileStructure.Forest.pairwiseDisjoint_rowSupport}
    The sets $E_j$, $1 \le j \le 2^n$ are pairwise disjoint.
\end{lemma}

\begin{proof}
    \leanok
    Suppose that $\fp \in \fT(\fu)$ and $\fp' \in \fT(\fu')$ with $\fu \ne \fu'$ and $x \in E(\fp) \cap E(\fp')$. Suppose without loss of generality that $\ps(\fp) \le \ps(\fp')$. Then $x \in \scI(\fp) \cap \scI(\fp') \subset \scI(\fu')$. By \eqref{dyadicproperty} it follows that $\scI(\fp) \subset \scI(\fu')$. By \eqref{forest5}, it follows that
    $$
        d_{\fp}(\fcc(\fp), \fcc(\fu')) > 2^{Z(n+1)}\,.
    $$
    By the triangle inequality. \Cref{monotone-cube-metrics} and \eqref{forest1} it follows that
    \begin{align*}
        d_{\fp}(\fcc(\fp), \fcc(\fp')) &\ge d_{\fp}(\fcc(\fp), \fcc(\fu')) - d_{\fp}(\fcc(\fp'), \fcc(\fu'))\\
        &> 2^{Z(n+1)} - d_{\fp'}(\fcc(\fp'), \fcc(\fu'))\\
        &\ge 2^{Z(n+1)} - 4\,.
    \end{align*}
    Since $Z \ge 3$ by \eqref{defineZ}, it follows that $\fcc(\fp') \notin B_{\fp}(\fcc(\fp), 1)$, so $\Omega(\fp') \not\subset \Omega(\fp)$ by \eqref{eq-freq-comp-ball}. Hence, by \eqref{eq-freq-dyadic}, $\Omega(\fp) \cap \Omega(\fp') = \emptyset$. But if $x \in E(\fp) \cap E(\fp')$ then $Q(x) \in \Omega(\fp) \cap \Omega(\fp')$. This is a contradiction, and the lemma follows.
\end{proof}

Now we prove \Cref{forest-operator}.

\begin{proof}[Proof of \Cref{forest-operator}]
    \leanok
    \proves{forest-operator}
    By \eqref{definetp*}, we have for each $j$
    $$
        T_{\mathfrak{R}_j}^*g = \sum_{\fu \in \fU_j} \sum_{\fp \in \fT(\fu)} T_{\fp}^* g = \sum_{\fu \in \fU_j} \sum_{\fp \in \fT(\fu)} T_{\fp}^* \mathbf{1}_{E_j} g = T_{\mathfrak{R}_j}^* \mathbf{1}_{E_j} g\,.
    $$
    Hence, by \Cref{forest-row-decomposition} and the triangle inequality,
    $$
        \left\|\sum_{\fu \in \fU} \sum_{\fp \in \fT(\fu)} T^*_{\fp} g\right\|_2^2 = \left\|\sum_{j = 1}^{2^n} T^*_{\mathfrak{R}_{j}} g\right\|_2^2 = \left\|\sum_{j=1}^{2^n} T^*_{\mathfrak{R}_{j}} \mathbf{1}_{E_j} g\right\|_2^2
    $$
    $$
        = \int_X \left|\sum_{j=1}^{2^n} T^*_{\mathfrak{R}_{j}} \mathbf{1}_{E_j} g\right|^2 \, \mathrm{d}\mu
    $$
    $$
        \le \sum_{j=1}^{2^n} \|T_{\mathfrak{R}_j}^* \mathbf{1}_{E_j} g\|_2^2 + \sum_{j =1}^{2^n} \sum_{\substack{j' = 1\\j' \ne j}}^{2^n} \left| \int_X \overline{ T_{\mathfrak{R}_j}^* \mathbf{1}_{E_j} g} T_{\mathfrak{R}_{j'}}^* \mathbf{1}_{E_{j'}} g \right| \, \mathrm{d}\mu\,.
    $$
    We use \Cref{row-bound} to estimate each term in the first sum, and \Cref{row-correlation} to bound each term in the second sum:
    $$
        \le 2^{566a^3-n} \sum_{j = 1}^{2^n} \|\mathbf{1}_{E_j} g\|_2^2 + 2^{876a^3-4n}\sum_{j=1}^{2^n}\sum_{j' = 1}^{2^n} \|\mathbf{1}_{E_j} g\|_2 \|\mathbf{1}_{E_{j'}}g\|_2\,.
    $$
    By Cauchy-Schwarz in the second two sums, this is at most
    $$
        2^{876a^3} (2^{-n} + 2^{n}2^{-4n}) \sum_{j = 1}^n \|\mathbf{1}_{E_j} g\|_2^2\,,
    $$
    and by disjointedness of the sets $E_j$, this is at most
    $$
        2^{877a^3 - n} \|g\|_2^2\,.
    $$
    Taking square roots, it follows that for all $g$
    \begin{equation}
        \label{eq-forest-bound-1}
        \left\|\sum_{\fu \in \fU} \sum_{\fp \in \fT(\fu)} T_{\fp}^* g\right\|_2 \le 2^{439a^3-\frac{n}{2}} \|g\|_2\,.
    \end{equation}
    On the other hand, we have by disjointedness of the sets $E_j$ from \Cref{disjoint-row-support} and the triangle inequality
    $$
        \left\|\mathbf{1}_G \sum_{\fu \in \fU} \sum_{\fp \in \fT(\fu)} T_{\fp} f\right\|_2^2 = \left\|\sum_{j=1}^{2^n} \mathbf{1}_{E_j} \mathbf{1}_G T_{\mathfrak{R}_{j}} f\right\|_2^2 \le \sum_{j = 1}^{2^n} \|\mathbf{1}_{E_j} \mathbf{1}_G T_{\mathfrak{R}_{j}} f\|_2^2 \le \sum_{j = 1}^{2^n} \|\mathbf{1}_G T_{\mathfrak{R}_{j}} f\|_2^2\,.
    $$
    Now with $|f| \le \mathbf{1}_F$ and \Cref{row-bound} we obtain
    $$
        \| \mathbf{1}_G T_{\mathfrak{R}_j} f \|_2^2 = \left| \int_X \overline{\mathbf{1}_G T_{\mathfrak{R}_j} f} T_{\mathfrak{R}_j} f \right| = \left| \int_X \overline{T_{\mathfrak{R}_j}^* \mathbf{1}_G T_{\mathfrak{R}_j} f} \mathbf{1}_F f \right|
    $$
    $$
        \le \|f\|_2 \| \mathbf{1}_F T_{\mathfrak{R}_j}^* \mathbf{1}_G T_{\mathfrak{R}_j} f \|_2 \le 2^{283a^3-n/2} \dens_2(\bigcup_{\fu\in \fU}\fT(\fu))^{1/2} \| \mathbf{1}_G T_{\mathfrak{R}_j} f \|_2 \|f\|_2.
    $$
    Dividing this last inequality by the finite $\| \mathbf{1}_G T_{\mathfrak{R}_j} f \|_2$, substituting back and taking square roots we get
    $$
        \left\|\mathbf{1}_G \sum_{\fu \in \fU} \sum_{\fp \in \fT(\fu)} T_{\fp} f\right\|_2 \le 2^{283a^3} \dens_2(\bigcup_{\fu\in \fU}\fT(\fu))^{\frac{1}{2}} 2^{-\frac{n}{2}} (\sum_{j = 1}^{2^n} \|f\|_2^2)^{\frac{1}{2}}
    $$
    \begin{equation}
        \label{eq-forest-bound-2}
        = 2^{283a^3} \dens_2(\bigcup_{\fu\in \fU}\fT(\fu))^{\frac{1}{2}} \|f\|_2\,.
    \end{equation}
    \Cref{forest-operator} follows by taking the product of the $(2 - \frac{2}{q})$-th power of \eqref{eq-forest-bound-1} and the $(\frac{2}{q} - 1)$-st power of \eqref{eq-forest-bound-2}.
\end{proof}

\section{Proof of the H\"older cancellative condition}
\label{liphoel}

We need the following auxiliary lemma.
Recall that $\tau = 1/a$.

\begin{lemma}[Lipschitz Holder approximation]
    \label{Lipschitz-Holder-approximation}
    \leanok
    \lean{support_holderApprox_subset, enorm_holderApprox_sub_le, iLipENorm_holderApprox_le}
    Let $z\in X$ and $R>0$. Let $\varphi: X \to \mathbb{C}$ be a function supported in the ball
    $B:=B(z,R)$ with finite norm $\|\varphi\|_{C^\tau(B(z, 2R))}$.
    Let $0<t \leq 1$. There exists a function $\tilde \varphi : X \to \mathbb{C}$, supported in $B(z,2R)$, such that for every $x\in X$
    \begin{equation}\label{eq-firstt}
        |\varphi(x) - \tilde \varphi(x)| \leq (t/2)^{\tau} \|\varphi\|_{C^\tau(B(z,2R))}
    \end{equation}and
   \begin{equation}\label{eq-secondt}
       \|\tilde \varphi\|_{\Lip(B(z,2R))} \leq 2^{4a}t^{-1-a} \|\varphi\|_{C^{\tau}(B(z,2R))}\, .
   \end{equation}
\end{lemma}

\begin{proof}\leanok
    Define for $x,y\in X$ the Lipschitz and thus measurable function
    \begin{equation}
        L(x,y) := \max\{0, 1 - \frac{\rho(x,y)}{tR}\}\, .
    \end{equation}
We have that $L(x,y)\neq 0$ implies
\begin{equation}\label{eql01}
    y\in B(x, tR)\, .
\end{equation}
We have for $y\in B(x, 2^{-1}tR)$ that
\begin{equation}\label{eql30}
    |L(x,y)|\ge 2^{-1} \ .
\end{equation}
Hence
\begin{equation}
    \int L(x,y) \, \mathrm{d}\mu(y)\ge 2^{-1}\mu(B(x, 2^{-1}tR))\, .
\end{equation}
 Let $n$ be the smallest integer so that
 \begin{equation}\label{2nt1}
    2^n t\ge 1\, .
 \end{equation}
 Iterating $n+2$ times the doubling condition \eqref{doublingx}, we obtain
 \begin{equation}\label{eql32}
    \int L(x,y) \, \mathrm{d}\mu(y)\ge 2^{-1-a(n+2)}\mu(B(x, 2R))\, .
 \end{equation}

Now define
    $$
        \tilde \varphi(x) := \left(\int L(x,y) \, \mathrm{d}\mu(y)\right)^{-1}\int L(x,y) \varphi(y) \, \mathrm{d}\mu(y)\, .
    $$
Using that $\varphi$ is supported in $B(z,R)$ and
\eqref{eql01}, we have that $\tilde{\varphi}$ is supported in $B(z,2R)$.

We prove \eqref{eq-firstt}.
 For any $x\in X$, using that $L$ is nonnegative,
   \begin{equation}\label{eql1}
    \left(\int L(x,y) \, \mathrm{d}\mu(y)\right)
        |\varphi(x) - \tilde \varphi(x)|
    \end{equation}
 \begin{equation}\label{eql2}
 = \left| \int L(x,y)(\varphi(x) - \varphi(y)) \, \mathrm{d}\mu(y)\right|\, .
    \end{equation}
Using \eqref{eql01}, we estimate the last display by
 \begin{equation}\label{eql3}
         \le \int_{B(x, tR)} L(x,y)|\varphi(x) - \varphi(y)| \, \mathrm{d}\mu(y)\, .\end{equation}
We claim that in this integral, $|\varphi(x) - \varphi(y)|\le \rho(x,y)^\tau \|\varphi\|_{C^\tau(B(z, 2R))}(2R)^{-\tau}$.
Indeed, if $x$ or $y$ does not belong to $B(z, 2R)$, then the other point is not in $B(z,R)$ as $\rho(x,y)\le tR \le R$. Therefore,
$\varphi(x)=\varphi(y)=0$ since $\varphi$ is supported in $B (z, R)$. Otherwise, both points are in $B(z, 2R)$, and the inequality
follows from the definition of $\|\varphi\|_{C^\tau(B(z, 2R))}$.

Therefore, we can estimate the last display further by
       \begin{equation}\label{eql4}
         \le \left(\int_{B(x, tR)} L(x,y)
          \rho(x,y)^\tau \, \mathrm{d}\mu(y) \right)\|\varphi\|_{C^\tau(B(z, 2R))}(2R)^{-\tau}\, .
    \end{equation}

  Using the condition on the domain of integration to estimate $\rho(x,y)$ by $tR$ and then expanding the domain by positivity of the integrand, we estimate this further by

   \begin{equation}\label{eql5}
         \le \left(\int L(x,y) \, \mathrm{d}\mu(y)\right)
         \|\varphi\|_{C^\tau(B(z, 2R))} (t/2)^{\tau} \, .
    \end{equation}
 Dividing the string of inequalities from \eqref{eql1} to
\eqref{eql5} by the positive integral of $L$ proves \eqref{eq-firstt}.

We turn to \eqref{eq-secondt}. For every $x\in X$, we have
\begin{equation}
    \left|\int L(x,y) \, \mathrm{d}\mu(y)\right||\tilde{\varphi}(x)|
    =\left|\int L(x,y) {\varphi}(y)\, \mathrm{d}\mu(y)\right|
\end{equation}
 \begin{equation}
    \le \left|\int L(x,y) \, \mathrm{d}\mu(y)\right| \sup_{x'\in X}
    |{\varphi}(x')|\ .
\end{equation}
As $\varphi$ is supported on $B$, dividing by the integral of $L$, we obtain
\begin{equation}\label{eql42}
 |\tilde{\varphi}(x)|\le \sup_{x'\in B}
    |{\varphi}(x')|\le \|\varphi\|_{C^\tau(B(z, 2R))}\ .
\end{equation}
If $\rho(x,x')\ge R$, then we have by the triangle inequality
  \begin{equation}\label{eql52}
 R\frac{|\tilde{\varphi}(x') - \tilde \varphi(x)|}{\rho(x,x')} \le
 2\sup_{x''\in X} |\tilde{\varphi}(x'')|\le 2\|\varphi\|_{C^\tau(B(z, 2R))}\, .
\end{equation}
Now assume $\rho(x,x')< R$. For $y\in X$ we have by the triangle inequality and a two fold case distinction
for the maximum in the definition of $L$,
\begin{equation}\label{eql10}
   |L(x,y) - L(x',y)| \le \frac{\rho(x,x')}{tR}.
\end{equation}
We compute with \eqref{eql10}, first adding and subtracting a term in the integral,
\begin{equation}
    \left(\int L(x,y) \, \mathrm{d}\mu(y)\right)
    |\tilde{\varphi}(x') - \tilde \varphi(x)|=
\end{equation}
\begin{equation}
    \left|\int L(x,y) \tilde{\varphi}(x')
    -L(x,y) \tilde{\varphi}(x)
    +L(x',y) \tilde{\varphi}(x')-
     L(x',y) \tilde{\varphi}(x')
    \, \mathrm{d}\mu(y)\right|\,.
\end{equation}
Grouping the second and third and the first and fourth term, we obtain using the definition of $\tilde \varphi$ and Fubini,
\begin{equation}\label{eql21}
    \le \left| \int (L(x',y)-L(x,y)) \varphi(y) \, \mathrm{d}\mu(y)\right|
\end{equation}
\begin{equation}\label{eql22}
    + \left| \int L(x,y) \, \mathrm{d}\mu(y)-\int L(x',y) \, \mathrm{d}\mu(y)\right||\tilde{\varphi}(x')|
\end{equation}
\begin{equation}\label{eql23}
    \le 2 \int |L(x,y) -L(x',y)| \, \mathrm{d}\mu(y)
    \|\varphi\|_{C^\tau(B(z, 2R))}\, ,
\end{equation}
where in the last inequality we have used \eqref{eql42}.
Using further \eqref{eql10} and the support of $L$, we estimate the last display by
\begin{equation}\label{eql224}\le 2 \frac{\rho(x,x')} {tR}\mu(B(x,tR)\cup B(x',tR))
\|\varphi\|_{C^\tau(B(z, 2R))}\, .
    \end{equation}
  Using $\rho(x,x')<R$ and the triangle inequality, we estimate the last display by
\begin{equation}\label{eql225}\le 2
\frac{\rho(x,x')} {tR}
\mu(B(x,2R))
\|\varphi\|_{C^\tau(B(z, 2R))}\, .
    \end{equation}
Dividing by the integral over $L$ and using \eqref{eql32} and \eqref{2nt1}, we obtain
\begin{equation}\label{eql226}
 \frac {R |\tilde{\varphi}(x') - \tilde \varphi(x)|}{\rho(x,x')}
 \le 2^{2+a(n+2)}t^{-1}\|\varphi\|_{C^\tau(B(z, 2R))} \le
 2^{2+3a} t^{-1-a} \|\varphi\|_{C^\tau(B(z, 2R))}\, .
\end{equation}
Combining \eqref{eql52} and \eqref{eql226} using $a\ge 4$ and $t\le 1$ and
adding \eqref{eql42} proves \eqref{eq-secondt} and completes the proof
of \Cref{Lipschitz-Holder-approximation}.
\end{proof}

We turn to the proof of \Cref{Holder-van-der-Corput}.
\begin{proof}[Proof of \Cref{Holder-van-der-Corput}]
    \proves{Holder-van-der-Corput}\leanok
Let $z\in X$ and $R>0$ and set $B=B(z,R)$. Let $\varphi$
be given as in \Cref{Holder-van-der-Corput}.
Set
\begin{equation}\label{eql69}
    t:=(1+d_B(\mfa,\mfb))^{-\frac{\tau}{2+a}}
\end{equation}
and define $\tilde{\varphi}$ as in \Cref{Lipschitz-Holder-approximation}. Let $\mfa$ and $\mfb$ be in $\Mf$.
Then
   \begin{equation}\label{eql60}
       \left|\int e(\mfa(x)-{\mfb(x)}) \varphi (x)\, \mathrm{d}\mu(x)\right|
   \end{equation}
    \begin{equation}\label{eql61}
   \le \left|\int e(\mfa(x)-{\mfb(x)}) \tilde{\varphi} (x)\, \mathrm{d}\mu(x)\right|
   \end{equation}
         \begin{equation}\label{eql62}
     + \left|\int e(\mfa(x)-{\mfb(x)}) (\varphi (x)-\tilde{\varphi}(x))\, \mathrm{d}\mu(x)\right|
   \end{equation}
Using the cancellative condition \eqref{eq-vdc-cond} of $\Mf$ on the ball $B(z,2R)$, the term \eqref{eql61} is bounded above by
 \begin{equation}\label{eql63}
       2^a \mu(B(z,2R)) \|\tilde{\varphi}\|_{\Lip(B(z,2R))} (1 + d_{B(z,2R)}(\mfa,\mfb))^{-\tau} \, .
 \end{equation}

Using the doubling condition \eqref{doublingx},
the inequality \eqref{eq-secondt}, and the estimate
$d_B\le d_{B(z,2R)}$ from the definition,
we estimate \eqref{eql63} from above by
\begin{equation}\label{eql64}
       2^{6a}t^{-1-a} \mu(B) \|{\varphi}\|_{C^\tau(B)}
       (1 + d_{B}(\mfa,\mfb))^{-\tau} \, .
 \end{equation}

The term \eqref{eql62} we estimate using
\eqref{eq-firstt} and that
$\mfa$ and $\mfb$ are real and thus $e(\mfa)$ and
$e(\mfb)$ bounded in absolute value by $1$.
We obtain for \eqref{eql62} with \eqref{doublingx}
the upper bound
  \begin{equation}\label{eql65}
      \mu(B(z,2R)) (t/2)^{\tau} \|\varphi\|_{C^\tau(B)}
      \le 2^a \mu(B) t^{\tau} \|\varphi\|_{C^\tau(B)}
      \,.
 \end{equation}
Using the definition \eqref{eql69} of $t$ and adding
\eqref{eql64} and \eqref{eql65} estimates
\eqref{eql60} from above by
\begin{equation}
       2^{6a} \mu(B) \|{\varphi}\|_{C^\tau(B)}
       (1 + d_{B}(\mfa,\mfb))^{-\frac{\tau}{2+a}}
       \end{equation}
\begin{equation} +
        2^a \mu(B) \|{\varphi}\|_{C^\tau(B)}
       (1 + d_{B}(\mfa,\mfb))^{-\frac{\tau^2}{2+a}}\, .
 \end{equation}
\begin{equation}\label{eql66}
      \le 2^{1+6a} \mu(B) \|{\varphi}\|_{C^\tau(B)}
       (1 + d_{B}(\mfa,\mfb))^{-\frac{\tau^2}{2+a}} \, ,
 \end{equation}
where we used $\tau\le 1$.
This completes the proof of \Cref{Holder-van-der-Corput}.
\end{proof}

\section{Proof of Vitali covering and Hardy--Littlewood}
\label{sec-hlm}

We begin with a classical representation of the Lebesgue norm.
\begin{lemma}[layer cake representation]\label{layer-cake-representation}
\leanok
\lean{MeasureTheory.eLpNorm_pow_eq_distribution}
Let $1\le p< \infty$. Then for any measurable function $u:X\to [0,\infty)$ on the measure space $X$
relative to the measure $\mu$
we have
\begin{equation}\label{eq-layercake}
    \|u\|_p^p=p\int_0^\infty \lambda^{p-1}\mu(\{x: u(x)\ge \lambda\})\, d\lambda\, .
\end{equation}
\end{lemma}
\begin{proof}
    \leanok
    The left-hand side of \eqref{eq-layercake} is by definition
\begin{equation}
    \int_X u(x)^p \, d\mu(x)\, .\end{equation}
    Writing $u(x)$ as an elementary integral in $\lambda$ and then using Fubini, we write for the last display
    \begin{equation}
    =\int_X \int _0^{u(x)}
    p \lambda^{p-1} d\lambda\, d\mu(x)
\end{equation}
\begin{equation}
 =p\int _0^{\infty}
    \lambda^{p-1} \mu(\{x: u(x)\ge \lambda\}) d\lambda\, .
\end{equation}
This proves the lemma.
\end{proof}

The following lemma will be used to define $M$ in the proof of \Cref{Hardy-Littlewood}.
\begin{lemma}[covering separable space]
    \label{covering-separable-space}
    \leanok
    \lean{covering_separable_space}
    For each $r > 0$, there exists a countable collection $C(r) \subset X$ of points such that
    $$
        X \subset \bigcup_{c \in C(r)} B(c, r)\,.
    $$
\end{lemma}

\begin{proof}
    \leanok
    It clearly suffices to construct finite collections $C(r,k)$ such that
    $$
        B(o, r2^k) \subset \bigcup_{c \in C(r,k)} B(c,r)\,,
    $$
    since then the collection $C(r) = \bigcup_{k \in \mathbb{N}} C(r,k)$ has the desired property.

    Suppose that $Y \subset B(o, r2^k)$ is a collection of points such that for all $y, y' \in Y$ with $y \ne y'$, we have $\rho(y,y') \ge r$. Then the balls $B(y, r/2)$ are pairwise disjoint and contained in $B(o, r2^{k+1})$. If $y \in B(o, r)$, then $B(o, r2^{k+1}) \subset B(y, r2^{k+2})$. Thus, by the doubling property \eqref{doublingx},
    $$
        \mu(B(y, \frac{r}{2})) \ge 2^{-(k+2)a} \mu(B(o, r2^{k+1}))\,.
    $$
    Thus, we have
    $$
        \mu(B(o, r2^{k+1})) \ge \sum_{y \in Y} \mu(B(y, \frac{r}{2})) \ge |Y| 2^{-(k+2)a} \mu(B(o, r2^{k+1}))\,.
    $$
    We conclude that $|Y| \le 2^{(k+2)a}$. In particular, there exists a set $Y$ of maximal cardinality. Define $C(r,k)$ to be such a set.

    If $x \in B(o, r2^k)$ and $x \notin C(r,k)$, then there must exist $y \in C(r,k)$ with $\rho(x,y) < r$. Thus $C(r,k)$ has the desired property.
\end{proof}

We turn to the proof of \Cref{Hardy-Littlewood}.
\begin{proof}[Proof of \Cref{Hardy-Littlewood}]
\leanok
\proves{Hardy-Littlewood}
Let the collection $\mathcal{B}$ be given.
We first show \eqref{eq-besico}.

We recursively choose a finite sequence $B_i\in \mathcal{B}$
for $i\ge 0$ as follows. Assume $B_{i'}$
is already chosen for $0\le i'<i$.
If there exists a ball $B_{i}\in \mathcal{B}$ so that $B_{i}$
is disjoint from all $B_{i'}$
with $0\le i'<i$, then choose
such a ball $B_i=B(x_i,r_i)$ with maximal $r_i$.

If there is no such ball, stop the selection and set
$i'':=i$.

By disjointedness of the chosen balls and since $0 \le u$, we have
\begin{equation}
\sum_{0\le i<i''}\int_{B_i} u(x)\, d\mu(x) \le \int_X u(x)\, d\mu(x)\, .
\end{equation}
By \eqref{eq-ball-assumption}, we conclude
\begin{equation}\label{eqbes1}
\lambda \sum_{0\le i<i''}\mu(B_i)
\le \int_X u(x)\, d\mu(x)\, .
\end{equation}
Let $x\in \bigcup \mathcal{B}$.
Choose a ball $B'=B(x',r')\in \mathcal{B}$
such that $x\in B'$.
If $B'$ is one of the selected balls, then
\begin{equation}\label{3rone}
    x\in \bigcup _{0\le i< i''}B_i\subset \bigcup _{0\le i< i''}B(x_i,3r_i)\, .
\end{equation}
If $B'$ is not one of the selected balls, then as it is not selected at time $i''$, there is a selected ball $B_i$ with
$B'\cap B_i\neq \emptyset$.
Choose such $B_i$ with minimal index $i$. As $B'$ is therefore disjoint from all
balls $B_{i'}$ with $i'<i$ and
as it was not selected in place of $B_i$, we have $r_i\ge r'$.

Using a point $y$ in the intersection of $B_i$ and $B'$,
we conclude by the triangle inequality
\begin{equation}
   \rho(x_i,x')\le \rho(x_i,y)+\rho(x',y)\le r_i+r'\le 2r_i \, .
\end{equation}
By the triangle inequality again, we further conclude
\begin{equation}
   \rho(x_i,x)\le \rho(x_i,x')+\rho(x',x)\le 2r_i+r'\le 3r_i \, .
\end{equation}
It follows that
\begin{equation}\label{3rtwo}
    x\in \bigcup _{0\le i< i''}B(x_i,3r_i)\, .
\end{equation}
With \eqref{3rone} and \eqref{3rtwo}, we conclude
\begin{equation}
\bigcup \mathcal{B}\subset
\bigcup _{0\le i< i''}B(x_i,3r_i)\, .
\end{equation}
With the doubling property
\eqref{doublingx} applied twice, we conclude
\begin{equation}\label{eqbes2}
    \mu(\bigcup{\mathcal{B}})
    \le \sum _{0\le i< i''}\mu (B(x_i,3r_i))
    \le 2^{2a}\sum _{0\le i< i''}\mu (B_i)\, .
\end{equation}
With \eqref{eqbes1} and \eqref{eqbes2} we conclude
\eqref{eq-besico}.

We turn to the proof of \eqref{eq-hlm}. We first consider the case $p_1=1$ and recall $M_{\mathcal{B}}=M_{\mathcal{B},1}$.
We write for the $p_2$-th power of left-hand side of \eqref{eq-hlm}
with \Cref{layer-cake-representation}
and a change of variables
\begin{equation}
    \|M_{\mathcal{B}}u(x)\|_{p_2}^{p_2}
   =p_2\int _0^{\infty}
    \lambda^{p_2-1} \mu(\{x: M_{\mathcal{B}}u(x)\ge \lambda\}) d\lambda\,
\end{equation}
\begin{equation} \label{eqbesi11}
   =2^{p_2} p_2\int _0^{\infty}
    \lambda^{p_2-1} \mu(\{x: M_{\mathcal{B}}u(x)\ge 2\lambda\}) d\lambda\, .
\end{equation}
Fix $\lambda\ge 0$ and let $x\in X$ satisfy $M_{\mathcal{B}}u(x)\ge 2\lambda$. By definition of $M_{\mathcal{B}}$, there is a ball
$B'\in \mathcal{B}$ such that
$x\in B'$ and
\begin{equation}\label{eqbesi10}
\int_{B'} u(y)\, d\mu(y)\ge 2\lambda \mu({B'}) \, .
\end{equation}
Define
$u_\lambda(y):=0$ if $|u(y)|<\lambda$ and $u_\lambda(y):=u(y)$ if $|u(y)|\ge \lambda$.
Then with \eqref{eqbesi10}
\begin{equation}
\int_{B'} u_\lambda (y)\, d\mu(y)
=\int_{B'} u (y)\, d\mu(y)-
\int_{B'} (u-u_\lambda) (y) d\mu(y)\,
\end{equation}
\begin{equation}
\ge 2\lambda \mu({B'})-
\int_{B'} (u-u_\lambda) (y) d\mu(y)\, .
\end{equation}
As $(u-u_\lambda)(y)\le \lambda$
by definition, we can estimate the last display by
\begin{equation}
\ge 2\lambda \mu({B'})-
\int_{B'} \lambda \, d\mu(y)
=\lambda \mu({B'})\, .
\end{equation}
Hence $x$ is contained in
$\bigcup(\mathcal{B}_\lambda)$,
where $\mathcal{B}_\lambda$
is the collection of balls $B''$ in $\mathcal{B}$ such that
\begin{equation}
    \int_{B''} u_\lambda (y)\, d\mu(y)\ge \lambda \mu(B'')\, .
\end{equation}
We have thus seen
\begin{equation}
    \{x: M_{\mathcal{B}}u(x)\ge 2\lambda\}\subset
    \bigcup \mathcal{B}_\lambda
\, .
\end{equation}
Applying \eqref{eq-besico} to the collection $\mathcal{B}_\lambda$
gives
\begin{equation}
    \lambda \mu(\{x: M_{\mathcal{B}}u(x)\ge 2\lambda\})\le
   2^{2a}
    \int u_\lambda (x)\, dx\, .
\end{equation}
With \Cref{layer-cake-representation},
\begin{equation}\label{eqbesi12}
    \lambda \mu(\{x: M_{\mathcal{B}}u(x)\ge 2\lambda\})\le
   2^{2a}
    \int_0^\infty \mu (\{x: |u_\lambda (x)|\ge \lambda'\})\, d\lambda'\, .
\end{equation}
By definition of $u_\lambda$, making a case distinction between $\lambda\ge \lambda'$ and $\lambda <\lambda'$, we see that
\begin{equation}\label{eqbesi13}
   \mu (\{x: |u_\lambda (x)|\ge \lambda'\})
   \le
   \mu (\{x: |u (x)|\ge \max(\lambda,\lambda')\})\, .
\end{equation}
We obtain with \eqref{eqbesi11},
\eqref{eqbesi12}, and \eqref{eqbesi13}
\begin{equation}
    \|M_{\mathcal{B}}u(x)\|_{p_2}^{p_2}
 \end{equation}
 \begin{equation}
   \le 2^{p_2+2a} p_2
   \int_0^\infty \lambda^{p_2-2}
   \int_0^\infty
   \mu (\{x: |u (x)|\ge \max(\lambda,\lambda')\})
   \, d\lambda'd\lambda\, .
\end{equation}
We split the integral into $\lambda\ge \lambda'$ and $\lambda<\lambda'$ and resolve the
maximum correspondingly.
We have for $\lambda\ge \lambda'$
with \Cref{layer-cake-representation}
\begin{equation}
    \int_0^\infty \lambda^{p_2-2}
   \int_0^\lambda
   \mu (\{x: |u (x)|\ge \lambda\})
   \, d\lambda'd\lambda
\end{equation}
\begin{equation}
   =\int_0^\infty \lambda^{p_2-1}
     \mu (\{x: |u (x)|\ge \lambda\})
d\lambda.
\end{equation}
\begin{equation}\label{eqbesi14}
   =p_2^{-1} \|u\|_{p_2}^{p_2}\, .
\end{equation}
We have for $\lambda< \lambda'$
with Fubini and \Cref{layer-cake-representation}
\begin{equation}
    \int_0^\infty \lambda^{p_2-2}
   \int_\lambda^\infty
   \mu (\{x: |u(x)|\ge \lambda'\})
   \, d\lambda'd\lambda.
\end{equation}
\begin{equation}
   =\int_0^\infty \int_0^{\lambda'}\lambda^{p_2-2}
     \mu (\{x: |u (x)|\ge \lambda'\})
d\lambda d\lambda'.
\end{equation}
\begin{equation}
   =(p_2-1)^{-1}\int_0^\infty (\lambda')^{p_2-1}
     \mu (\{x: |u(x)|\ge \lambda'\})
d\lambda'.
\end{equation}
\begin{equation}\label{eqbesi15}
   =(p_2-1)^{-1} p_2^{-1}\|u\|_{p_2}^{p_2}\, .
\end{equation}
Adding the two estimates
\eqref{eqbesi14} and \eqref{eqbesi15} gives
\begin{equation}
    \|M_{\mathcal{B}}u(x)\|_{p_2}^{p_2}
   \le 2^{p_2+2a} (1+(p_2-1)^{-1})\|u\|_{p_2}^{p_2}
   = 2^{p_2+2a} p_2(p_2-1)^{-1}\|u\|_{p_2}^{p_2}
   \, .
   \end{equation}
With $a\ge 1$ and $p_2>1$, taking the $p_2$-th root, we obtain \eqref{eq-hlm}.
We turn to the case of general
$1\le p_1<p_2$.
We have
\begin{equation}
    M_{\mathcal{B},p_1}u=(M_{\mathcal{B}} (|u|^{p_1}))^{\frac 1{p_1}}\, .
\end{equation}
Applying the special case of \eqref{eq-hlm} for $M_{\mathcal{B}}$ gives
\begin{equation}
    \|M_{\mathcal{B},p_1}u\|_{p_2}=
    \|M_{\mathcal{B}} (|u|^{p_1})\|_{p_2/p_1}^{\frac 1{p_1}}
\end{equation}
\begin{equation}
    \le 2^{2a} (p_2/p_1) (p_2/p_1-1)^{-1}
    \|(|u|^{p_1})\|_{p_2/p_1}^{\frac 1{p_1}}
    =2^{2a} p_2(p_2-p_1)^{-1}\|u\|_{p_2}\, .
\end{equation}
This proves \eqref{eq-hlm} in general.

Now we construct the operator $M$ satisfying \eqref{eq-ball-av} and \eqref{eq-hlm-2}.
For each $k \in \mathbb{Z}$ we choose a countable set $C(2^k)$ as in \Cref{covering-separable-space}.
Define
$$
    \mathcal{B}_\infty = \{B(c, 2^k) \ : \ c \in C(2^k), k \in \mathbb{Z}\}\,.
$$
By \Cref{covering-separable-space}, this is a countable collection of balls. We choose an enumeration $\mathcal{B}_\infty = \{B_1, \dotsc\}$ and define
$$
    \mathcal{B}_n = \{B_1, \dotsc, B_n\}\,.
$$
We define
$$
    Mw := 2^{2a}\sup_{n \in \mathbb{N}} M_{\mathcal{B}_n}w\,.
$$
This function is measurable for each measurable $w$, since it is a countable supremum of measurable functions. Estimate \eqref{eq-hlm-2} follows immediately from \eqref{eq-hlm} and the monotone convergence theorem.

It remains to show \eqref{eq-ball-av}. Let $B = B(x, r) \subset X$. Let $k$ be the smallest integer such that $2^k \ge r$, in particular we have $2^k < 2r$. By definition of $C(2^k)$, there exists $c \in C(2^k)$ with $x \in B(c, 2^k)$. By the triangle inequality, we have $B(c, 2^k) \subset B(x, 4r)$, and hence by the doubling property \eqref{doublingx}
$$
    \mu(B(c, 2^k)) \le 2^{2a} \mu(B(x,r))\,.
$$
It follows that for each $z \in B(x,r)$
\begin{align*}
    \frac{1}{\mu(B(x,r))}\int_{B(x,r)} |w(y)| \, \mathrm{d}\mu(y) &\le \frac{2^{2a}}{\mu(B(c,2^k))}\int_{B(c,2^k)} |w(y)| \, \mathrm{d}\mu(y) \\
    &\le Mw(z)\,.
\end{align*}
This completes the proof.
\end{proof}

\section{Two-sided Metric Space Carleson}

We prove a variant of \Cref{metric-space-Carleson} for a two-sided Calder\'on--Zygmund kernel on the doubling metric measure space $(X,\rho,\mu,a)$, i.e. a one-sided Calder\'on--Zygmund kernel $K$ which additionally satisfies for all $x,x',y\in X$ with $x\neq y$ and $2\rho(x,x') \leq \rho(x,y)$,
\begin{equation}
    \label{eqkernel-x-smooth}
    |K(x,y) - K(x',y)| \leq \left(\frac{\rho(x,x')}{\rho(x,y)}\right)^{\frac{1}{a}}\frac{2^{a^3}}{V(x,y)}.
\end{equation}
By the additional regularity, we can weaken the assumption \eqref{nontanbound} to a family of operators that is easier to work with in applications.
Namely, for $r > 0$, $x\in X$, and a bounded, measurable function $f:X\to\C$ supported on a set of finite measure, we define
\begin{equation}
\label{def-T-r}
T_r f(x):= \int_{r\le\rho(x,y)} K(x,y) f(y) \, d\mu(y) = \int_{X\setminus B(x,r)} K(x,y) f(y) \, d\mu(y).
\end{equation}
\begin{theorem}[two-sided metric space Carleson]
    \label{two-sided-metric-space-Carleson}
    \uses{metric-space-Carleson, nontangential-from-simple}
    \leanok
    \lean{two_sided_metric_carleson}
        For all  integers $a \ge  4$ and real numbers $1<q\le 2$
        the following holds.
        Let $(X,\rho,\mu,a)$ be a doubling metric measure space. Let  $\Mf$ be a
        cancellative compatible  collection of functions and let $K$ be a two-sided Calder\'on--Zygmund kernel on $(X,\rho,\mu,a)$. Assume  that for every bounded measurable function $g$ on $X$ supported on a set of finite measure and all $r>0$ we have
      \begin{equation}\label{two-sided-Hr-bound-assumption}
            \|T_r g\|_{2} \leq 2^{a^3} \|g\|_2\,.
        \end{equation}
        Then for all Borel sets $F$ and $G$ in $X$ and
        all Borel functions $f:X\to \C$ with
        $|f|\le \mathbf{1}_F$, we have, with $T$ defined in \eqref{def-main-op},
      \begin{equation}
        \label{two-sided-resweak}
            \left|\int_{G} T f \, \mathrm{d}\mu\right| \leq \frac{2^{474a^3}}{(q-1)^6} \mu(G)^{1-\frac{1}{q}} \mu(F)^{\frac{1}{q}}\, .
        \end{equation}
\end{theorem}

For the remainder of this chapter, fix an integer $a\ge 4$, a doubling metric measure space $(X,\rho,\mu,a)$ and a two-sided Calder\'on--Zygmund kernel $K$ as in \Cref{two-sided-metric-space-Carleson}.

The following lemma is proved in \Cref{subsec-cotlar}.
\begin{lemma}[nontangential-from-simple]\label{nontangential-from-simple}
    \leanok
    \lean{nontangential_from_simple}
    \uses{simple-nontangential-operator, nontangential-operator-boundary}
    Assume \eqref{two-sided-Hr-bound-assumption} holds.
    Then, for every bounded measurable function $g : X \to \C$ supported on a set of finite measure we have
    \begin{equation}\label{concretetstarbound}
        \|T_*g\|_2\le 2^{3a^3}\|g\|_2.
    \end{equation}
\end{lemma}

\begin{proof}[Proof of \Cref{two-sided-metric-space-Carleson}]
    \proves{two-sided-metric-space-Carleson}
    \leanok
    Let $1<q\le 2$ be a real number. Let $\Theta$ be a cancellative compatible collection of functions.
    By the assumption \eqref{two-sided-Hr-bound-assumption}, we can apply \Cref{nontangential-from-simple} to obtain for every bounded measurable $g:X\to\C$ supported on a set of finite measure,
    \begin{equation}\label{original-operator-assumption}
        \|T_*g\|_2\le 2^{3a^3}\|g\|_2.
    \end{equation}
    Define
    \begin{equation*}
        K'(x,y):= 2^{-2a^3} K(x,y) \,.
    \end{equation*}
    Then $K'$ is a two-sided Calder\'on--Zygmund kernel on $(X,\rho,\mu,a)$. Denote the corresponding maximally truncated non-tangential singular operator by $T_*'$ and the corresponding generalized Carleson operator by $T'$.
    With \eqref{original-operator-assumption}, we obtain for $g$ as above,
    \begin{equation}\label{modified-operator-assumption}
        \|T_*'g\|_2\le 2^{a^3}\|g\|_2.
    \end{equation}
    Applying \Cref{metric-space-Carleson} for $K'$ yields that for all Borel sets $F$ and $G$ in $X$ and
    all Borel functions $f:X\to \C$ with
    $|f|\le \mathbf{1}_F$, we have
    \begin{equation*}
        \left|\int_{G} T' f \, \mathrm{d}\mu\right| \leq \frac{2^{450a^3}}{(q-1)^6} \mu(G)^{1-\frac{1}{q}} \mu(F)^{\frac{1}{q}}\, .
    \end{equation*}
    This finishes the proof since for all $x\in X$,
    \begin{equation*}
        T'f(x) = 2^{-2a^3} Tf(x) \,.
    \end{equation*}
\end{proof}

The proof of \Cref{nontangential-from-simple} relies on the following auxiliary lemma which is proved in \Cref{subsec-CZD}.
\begin{lemma}[Calderon-Zygmund Weak (1, 1)]
    \label{calderon-zygmund-weak-1-1}
    \leanok
    \lean{czOperator_weak_1_1}
    \uses{Calderon-Zygmund-decomposition,estimate-bad,estimate-good}
    Let $f:X\to\C$ be a bounded measurable function supported on a set of finite measure and assume for some $r>0$ that for every bounded measurable function $g:X\to\C$ supported on a set of finite measure,
    \begin{equation}
        \label{eq-strong-2-2-assumption}
        \|T_rg\|_{2}\le 2^{a^3} \|g\|_2.
    \end{equation}
    Then for all $\alpha>0$, we have
    \begin{equation}
        \label{eq-weak-1-1}
        \mu\left(\{x\in X: |T_r f(x)|>\alpha\}\right)\le \frac{2^{a^3 + 21a}}{\alpha} \int |f(y)|\, d\mu(y).
    \end{equation}
\end{lemma}
Throughout \Cref{subsec-CZD} and \Cref{subsec-cotlar}, for any measurable bounded function $w: X \to \C$, let $Mw: X \to [0, \infty)$ denote the corresponding Hardy--Littlewood maximal function defined in \Cref{Hardy-Littlewood}.
Note: In some cases, a stronger formalization of \Cref{Hardy-Littlewood} is used, where $w$ is not necessarily bounded. In particular it is applied to $T_rg$.
Apart from \Cref{Hardy-Littlewood}, \Cref{subsec-CZD} and \Cref{subsec-cotlar} have no dependencies in the previous chapters.

\subsection{Proof of Cotlar's Inequality}
\label{subsec-cotlar}

\begin{lemma}[geometric series estimate]
    \label{geometric-series-estimate}
    \leanok
    \lean{geometric_series_estimate}
    For all real numbers $x\ge 4$,
    \begin{equation*}
        \sum_{n=0}^\infty 2^{-\frac{n}{x}} \le 2^x.
    \end{equation*}
\end{lemma}
\begin{proof}
    \leanok
    By convexity, for all $0\le\lambda\le1$
    \begin{equation*}
        2^{\lambda(-\frac{1}{4})} \le \lambda 2^{-\frac{1}{4}} + (1-\lambda)2^0.
    \end{equation*}
    For $\lambda:=\frac{4}{x}$, we obtain
    \begin{equation*}
        2^{-\frac 1 x} \le 1 - (1-2^{-\frac 1 4}) \frac{4}{x}.
    \end{equation*}
    We conclude
    \begin{equation*}
        \sum_{n=0}^\infty 2^{-\frac{n}{x}} = \frac{1}{1-2^{-\frac 1 x}} \le \frac{1}{4(1-2^{-\frac 1 4})} x \le 2^x.
    \end{equation*}
\end{proof}

\begin{lemma}[estimate x shift]
\label{estimate-x-shift}
\leanok
\lean{estimate_x_shift}
\uses{geometric-series-estimate}
Let $0<r$ and  $x\in X$. Let $g:X\to\C$ be a bounded measurable function supported on a set of finite measure.
Then for all $x'$ with $\rho(x,x')\le r$.
\begin{equation*}
\left| T_r g(x) - T_r g(x') \right|
\le 2^{a^3 + 2a + 2} Mg(x)\, .
\end{equation*}
\end{lemma}
\begin{proof}
\proves{estimate-x-shift}
\leanok
By definition,
\begin{equation}
\label{xx'difference}
\left| T_r g(x) - T_r g(x') \right|
=\left|\int_{r\le\rho(x,y)} K(x,y) g(y) \,d\mu(y)
-\int_{r\le\rho(x',y)}
K(x',y) g(y)\, d\mu(y)
\right| \, .
\end{equation}
We split the first integral in \eqref{xx'difference}
into the domains $r\le\rho(x,y)<2r$
and $2r\le\rho(x,y)$. The integral over the first domain we estimate
by $\eqref{firstxx'}$ below.
For the second domain, we
observe with $\rho(x,x')\le r$ and the triangle inequality that $r\le\rho(x',y)$. We therefore combine on this domain with the
corresponding part of the second integral in \eqref{xx'difference} and estimate that by $\eqref{secondxx'}$
below. The remaining part of the second integral in
\eqref{xx'difference} we estimate by $\eqref{thirdxx'}$.
Overall, we have estimated \eqref{xx'difference} by
\begin{equation}
\label{firstxx'}
\int_{r\le\rho(x,y)< 2r} |K(x,y)| |g(y)| \,d\mu(y)
\end{equation}
\begin{equation}
\label{secondxx'}
  +  \left|\int_{2r\le\rho(x,y)} (K(x,y)-K(x',y)) g(y) \,d\mu(y)
\right|
\end{equation}
\begin{equation}
\label{thirdxx'}
  +\int_{r\le\rho(x',y), r\le\rho(x,y)<2r}
|K(x',y)| |g(y)|\, d\mu(y)\,.
\end{equation}
Using the bound on $K$ in \eqref{eqkernel-size} and the doubling condition \eqref{doublingx}, we estimate \eqref{firstxx'} by
\begin{align}
 \int_{r\le\rho(x,y)<2r} \frac{2^{a^3}}{V(x,y)} |g(y)| \,d\mu(y)\,
&\le
\frac{2^{a^3}}{\mu(B(x,r))} \int_{r\le\rho(x,y)<2r}
|g(y)| \, d\mu(y)\, \\
\label{firstxx'b}
&\le
\frac{2^{a^3} \cdot 2^a}{\mu(B(x,2r))} \int_{\rho(x,y)<2r}
|g(y)| \, d\mu(y)\,
.
\end{align}
Using the definition of $Mg$, we estimate
\eqref{firstxx'b} by
\begin{equation}
\label{firstxx'c}
 \le 2^{a^3 + a} {Mg(x)}\, .
\end{equation}
Similarly, in the domain of \eqref{thirdxx'}
we note by the triangle inequality
and assumption on $x'$ that $\rho(x',y)<3r$ and thus we estimate
\eqref{thirdxx'}
by
\begin{equation}
\label{thirdxx'b}
\frac{2^{a^3}}{\mu(B(x',r))} \int_{\rho(x',y)<4r}
|g(y)| \, d\mu(y)\le 2^{a^3 + 2a} Mg(x)
\end{equation}

We turn to the remaining term. Using \eqref{eqkernel-x-smooth}, we estimate \eqref{secondxx'} by
\begin{equation}\label{secondxx'b}
 \int_{2r\le\rho(x,y)} \left(\frac{\rho(x,x')}{\rho(x,y)}\right)^{\frac{1}{a}}\frac{2^{a^3}}{V(x,y)} |g(y)| \,d\mu(y)
\end{equation}
We decompose and estimate
\eqref{secondxx'}
with the triangle inequality by
\begin{align}
 &\sum_{j=1}^\infty  \int_{2^jr\le\rho(x,y)< 2^{j+1}r}
 \left(\frac{\rho(x,x')}{\rho(x,y)}\right)^{\frac{1}{a}}\frac{2^{a^3}}{V(x,y)} |g(y)| \,d\mu(y)\, \\
 \le&
 \sum_{j=1}^\infty \left( 2^{-j} \right)^{\frac{1}{a}} \int_{2^jr\le\rho(x,y)< 2^{j+1}r} \frac{2^{a^3}}{\mu(B(x,2^j r))} |g(y)| \,d\mu(y) \\
 \le&
\sum_{j=1}^\infty 2^{-\frac{j}{a}} \frac{2^{a^3 + a}}{\mu(B(x,2^{j+1} r))} \int_{\rho(x,y)<2^{j+1}r}
|g(y)| \, d\mu(y) \\
\label{secondxx'c}
\le& 2^{a^3 + a} \sum_{j=1}^\infty 2^{-\frac{j}{a}} Mg(x)\,.
\end{align}
Using \Cref{geometric-series-estimate}, we estimate
\eqref{secondxx'c} by
\begin{equation}\label{secondxx'd}
\le 2^{a^3 + 2a} Mg(x) \,.
\end{equation}
Summing the estimates
for \eqref{firstxx'},
\eqref{secondxx'}, and
\eqref{thirdxx'}
proves the lemma.
\end{proof}

\begin{lemma}[Cotlar control]
\label{Cotlar-control}
\leanok
\lean{cotlar_control}
\uses{estimate-x-shift}
Let $0<r\le R$ and  $x\in X$. Let $g:X\to\C$ be a bounded measurable function supported on a set of finite measure.
Then for all $x'\in X$ with $\rho(x,x')\le\frac {R}{4}$ we have
\begin{equation}\label{eq-cotlar-control}
\left|T_R g(x)
\right|\le
|T_r(g-g\mathbf{1}_{B(x,\frac {R} 2)})(x')| +
2^{a^3 + 4a + 1} Mg(x)\, .
\end{equation}
\end{lemma}

\begin{proof}
\leanok
Let $x$ and $x'$ be given with $\rho(x,x')\le\frac {R}{4}$.
By \Cref{estimate-x-shift}, we estimate
the left-hand-side of
\eqref{eq-cotlar-control} by
\begin{equation}\label{eqcotlar0}
|T_R(g)(x')|+
2^{a^3 + 2a + 2} Mg(x)\,.
\end{equation}
We  have
\begin{equation}
\label{eqcotlar-1}
T_R(g)(x')=
\int_{R\le\rho(x',y)} K(x',y) g(y) \, d\mu(y)\, .
\end{equation}
On the domain $R\le\rho(x',y)$, we have $\frac {R}2\le\rho(x,y)$. Hence we may write
for \eqref{eqcotlar-1}
\begin{equation*}
T_R(g)(x')=\int_{R\le\rho(x',y)} K(x',y) (g-g\mathbf{1}_{B(x,\frac {R} 2)})(y) \,d\mu(y)
\end{equation*}
\begin{equation}\label{eqcotlar1}
=T_R(g-g\mathbf{1}_{B(x,\frac {R} 2)})(x')\, .
\end{equation}
Combining the estimate \eqref{eqcotlar0} with the identification \eqref{eqcotlar1}, we obtain
\begin{equation}\label{eqcotlar5}
\left|T_R g(x)
\right|\le
|T_R(g-g\mathbf{1}_{B(x,\frac {R} 2)})(x')|+
2^{a^3 + 2a + 2} Mg(x)\, .
\end{equation}
We  have
\begin{equation*}
(T_r-T_R)(g-g\mathbf{1}_{B(x,\frac {R} 2)})(x')
\end{equation*}
\begin{equation*}
= \int_{B(x',R)\setminus B(x',r)} K(x',y) (g-g\mathbf{1}_{B(x,\frac {R} 2)})(y) \,d\mu(y)
\end{equation*}
\begin{equation}\label{eqcotlar2}
= \int_{B(x',R)\setminus (B(x',r) \cup B(x,\frac{R}{2}))} K(x',y) g(y) \,d\mu(y)
\end{equation}
As $\frac{R}{2}\le\rho(x,y)$ together with $\rho(x,x')\le\frac {R}{4}$ implies
$\frac {R}{4}\le\rho(x',y)$, we can estimate the absolute value of \eqref{eqcotlar2} with \eqref{eqkernel-size} by
\begin{align*}
\le &\frac{2^{a^3}}{\mu(B(x',\frac{R}{4}))} \int_{B(x,2R)\setminus B(x',\frac{R}{4})}
|g(y)|\, d\mu(y) \\
\le &\frac{2^{a^3+a}}{\mu(B(x',\frac{R}{2}))} \int_{B(x,2R)}
|g(y)|\, d\mu(y) \\
\le &2^{a^3+a} \frac{\mu(B(x,2R))}{\mu(B(x,\frac{R}{4}))} Mg(x) \le 2^{a^3 + 4a} Mg(x)\, .
\end{align*}

By the triangle inequality, \eqref{eq-cotlar-control} follows now from \eqref{eqcotlar5} and the estimate for \eqref{eqcotlar2}.
\end{proof}

\begin{lemma}[Cotlar sets]
\label{Cotlar-sets}
\leanok
\lean{cotlar_set_F₁, cotlar_set_F₂}
\uses{calderon-zygmund-weak-1-1}
Assume that \eqref{two-sided-Hr-bound-assumption} holds. Let $0<r\le R$ and  $x\in X$. Let $g:X\to\C$ be a bounded measurable function supported on a set of finite measure.
Then the measure $|F_1|$ of the set $F_1$ of all $x'\in B(x,\frac {R} 4)$ such that
\begin{equation}
\label{first-cotlar-exception}
    |T_rg(x')|> 4 M(T_rg)(x)
\end{equation}
is less than or equal to $\mu(B(x,\frac{R}{4}))/4$.
Moreover,  the measure $|F_2|$ of the set $F_2$ of all $x'\in
B(x,\frac {R} 4)$ such that
\begin{equation}
\label{second-cotlar-exception}
    |T_r(g\mathbf{1}_{B(x,\frac {R} 2)})(x')| > 2^{a^3 + 22a + 2} Mg(x)
\end{equation}
is less than or equal to  $\mu(B(x,\frac{R}{4}))/4$.
\end{lemma}

\begin{proof}
\leanok
Let $r$, $R$, $x$ and $g$ be given.
If $M(T_rg)(x)=0$, then $T_rg$ is zero almost everywhere and the estimate on $|F_1|$ is trivial.
Assume $M(T_rg)(x)>0$.
We have with \eqref{first-cotlar-exception}
\begin{equation}
    M(T_rg)(x)\ge
    \frac 1{\mu(B(x,\frac{R}{4}))}\int_{B(x,\frac{R}{4})}|T_rg(x')|\, dx'
\end{equation}
\begin{equation}
    \ge
    \frac 1{\mu(B(x,\frac{R}{4}))}\int_{F_1} 4 M(T_rg)(x)\, dx' \,.
\end{equation}
Dividing by $M(T_rg)(x)$ gives
\begin{equation}
    1\ge \frac{4}{\mu(B(x,\frac{R}{4}))} |F_1|\, .
\end{equation}
This gives the desired bound for the measure of $F_1$.
We turn to the set $F_2$. Similarly as above we may assume $Mg(x)>0$.
The set $F_2$ is then estimated with \Cref{calderon-zygmund-weak-1-1} by
\begin{equation}
   \frac {2^{a^3+21a}}{2^{a^3 + 22a + 2}Mg(x)}\int |g\mathbf{1}_{B(x,\frac {R}2)}|(y)\, d\mu(y)
\end{equation}
\begin{equation}
   \le \frac {1}{2^{a+2} Mg(x)}\mu(B(x,\frac {R}{2})) Mg(x) \le \frac{\mu(B(x,\frac {R}{4}))}{4} \,.
\end{equation}
This gives the desired bound for the measure of $F_2$.
\end{proof}

\begin{lemma}[Cotlar estimate]
\label{Cotlar-estimate}
\leanok
\lean{cotlar_estimate}
\uses{Cotlar-control, Cotlar-sets}
Assume that \eqref{two-sided-Hr-bound-assumption} holds.
Let $0<r\le R$ and  $x\in X$. Let $g:X\to\C$ be a bounded measurable function supported on a set of finite measure. Then
\begin{equation}\label{eq-cotlar-estimate}
|T_Rg(x)| \le 2^{2}M(T_rg)(x)+ 2^{a^3+22a+3} Mg(x)
\, .
\end{equation}
\end{lemma}

\begin{proof}
\leanok
By \Cref{Cotlar-sets}, the set of all $x'\in B(x,\frac {R} 4)$
such that at least one of the conditions
\eqref{first-cotlar-exception} and
\eqref{second-cotlar-exception} is satisfied has measure less than or equal to $\mu(B(x,\frac{R}{4}))/2$ and hence is not all of $B(x,\frac {R} 4)$.
Pick an $x'\in B(x,\frac {R} 4)$ such that both conditions are not satisfied.
Applying \Cref{Cotlar-control} for this $x'$ and using the triangle inequality
estimates the left-hand side of \eqref{eq-cotlar-estimate}
by
\begin{equation}
    4 M(T_rg)(x) + 2^{a^3 + 22a + 2} Mg(x) + 2^{a^3 + 4a + 1} Mg(x)\, .
\end{equation}
This proves the lemma.
\end{proof}

For the next Lemma, we define the following operation.
\begin{equation}\label{eq-simple--nontangential}
    T_{*}^r g(x):=\sup_{r<R}\sup_{x'\in B(x,R)} |T_R(g)(x')| \, .
\end{equation}

\begin{lemma}[simple nontangential operator]\label{simple-nontangential-operator}
\leanok
\lean{simple_nontangential_operator, lowerSemicontinuous_simpleNontangentialOperator}
\uses{Hardy-Littlewood,Cotlar-estimate}
Assume that \eqref{two-sided-Hr-bound-assumption} holds.
For every $r>0$ and every bounded measurable function $g$ supported on a set of finite measure
the function $T_{*}^r g$ is lower-semicontinuous and we have
\begin{equation}\label{trzerobound}
    \|T_{*}^rg\|_2\le 2^{a^3+26a+6}\|g\|_2,
\end{equation}
\end{lemma}
\begin{proof}
\leanok
For a fixed $\lambda$ we note that (rewriting $x'\in B(x,R)$ as $x\in B(x',R)$):
\begin{equation}
    \{x : T_{*}^r g(x) > \lambda\} = \bigcup_{r<R}\bigcup_{\{x':|T_R(g)(x')|>\lambda\}} B(x',R)
\end{equation}
The right-hand side is open, hence $T_{*}^r g$ is lower-semicontinuous by definition.

For the second part, with \Cref{estimate-x-shift} and the triangle inequality, we estimate for every $x\in X$
\begin{equation}
     T_{*}^r g(x)
     \le 2^{a^3 + 2a + 2} Mg(x)+\sup_{r<R} |T_R(g)(x)|\, .
\end{equation}
Using further \Cref{Cotlar-estimate}, we estimate
\begin{equation}
      T_{*}^r g(x)
     \le 2^{a^3+2a+2}Mg(x) + 2^{a^3+22a+3}Mg(x) + 2^{2}M(T_rg)(x)\, .
\end{equation}
Taking the $L^2$ norm and using \Cref{Hardy-Littlewood} with $a=4$  and $p_2=2$ and $p_1=1$ , we obtain
\begin{equation}
      \|T_{*}^r g\|_2
     \le 2^{a^3+22a+4} \|Mg\|_2 + 2^{2}\|M(T_rg)\|_2
\end{equation}
\begin{equation}
     \le 2^{a^3+26a+5} \|g\|_2 + 2^{4a+3}\|T_r g\|_2\, .
\end{equation}
Applying \eqref{two-sided-Hr-bound-assumption} gives
\begin{equation}
      \|T_{*}^r g\|_2\le 2^{a^3+26a+5}\|g\|_2 + 2^{a^3+4a+3}\|g\|_2\, .
\end{equation}
This shows \eqref{trzerobound} and completes the proof of the inequality.
\end{proof}
In order to pass from the one-sided truncation in $T_r$ and $T_{*}^r$ to the two-sided truncation in $T_*$, we show in the following two lemmas that the integral in \eqref{def-tang-unm-op} can be exchanged for an integral over the difference of two balls.
\begin{lemma}[small annulus] \label{small-annulus}
    \leanok
    \lean{small_annulus_right, small_annulus_left}
    Let $f:X\to\C$ be a bounded measurable function supported on a set of finite measure.
    Let $x\in X$ and $R>0$.
    Then, for all $\epsilon>0$, there exists some $\delta>0$ such that
    \begin{equation}
        \left| \int_{R<\rho(x,y)<R+\delta} K(x,y) f(y) \, d\mu(y) \right| \le \epsilon
    \end{equation}
    and
    \begin{equation}
        \left| \int_{R-\delta<\rho(x,y)<R} K(x,y) f(y) \, d\mu(y) \right| \le \epsilon \,.
    \end{equation}
\end{lemma}
\begin{proof}
\leanok
We only prove the second inequality, the first one is analogous.
Note that the integrand is bounded in $X\setminus B(x,\frac{R}{2})$. So for $0<\delta\le\frac{R}{2}$,
\begin{align*}
    &\left| \int_{R-\delta<\rho(x,y)<R} K(x,y) f(y) \, d\mu(y) \right| \\
    \le &\frac{2^{a^3}}{\mu(B(x,\frac{R}{2}))} \, \sup_{y\in X}|f(x)| \cdot \mu(\{y\in X: R-\delta<\rho(x,y)<R\}).
\end{align*}
By continuity from above of $\mu$, the right factor becomes arbitrarily small as $\delta\rightarrow 0$. Thus, for small enough $\delta$, the whole expression is $\le\epsilon$.
\end{proof}

\begin{lemma}[nontangential operator boundary]
    \label{nontangential-operator-boundary}
    \leanok
    \lean{nontangential_operator_boundary}
    \uses{small-annulus}
    Let $f:X\to\C$ be a bounded measurable function supported on a set of finite measure.
    For all $x\in X$,
    \begin{equation}
    \label{tang-unm-op-eq}
        T_*f(x) = \sup_{R_1 < R_2} \sup_{x'\in B(x,R_1)} \left|\int_{B(x',R_2)\setminus B(x',R_1)} K(x',y) f(y) \, \mathrm{d}\mu(y) \right|
    \end{equation}
\end{lemma}
\begin{proof}
\leanok
We show two inequalities. Let $\epsilon>0$.
Let $R_1<R_2$ and $x'\in B(x,R_1)$. Then for small enough $\delta>0$,
\begin{alignat}{3}
    \label{eq-without-suprema-1}
    &&&\left|\int_{R_1<\rho(x',y)<R_2} K(x',y) f(y) \, \mathrm{d}\mu(y) \right| \\
    \label{eq-diff-small-1}
    \le &&&\left|\int_{R_1<\rho(x',y)<R_1+\delta} K(x',y) f(y) \, \mathrm{d}\mu(y) \right| \\
    \label{eq-other-1}
    &+&&\left|\int_{R_1+\delta\le\rho(x',y)<R_2} K(x',y) f(y) \, \mathrm{d}\mu(y) \right| .
\end{alignat}
By \Cref{small-annulus}, we can choose $\delta$ such that \eqref{eq-diff-small-1} is bounded by $\epsilon$. Without loss of generality, we can assume $R_1+\delta<R_2$. Then \eqref{eq-other-1} is bounded by the right hand side of \eqref{tang-unm-op-eq} and we obtain
\begin{equation*}
    \le \epsilon + \sup_{R_1 < R_2} \sup_{x'\in B(x,R_1)} \left|\int_{B(x',R_2)\setminus B(x',R_1)} K(x',y) f(y) \, \mathrm{d}\mu(y) \right| .
\end{equation*}
The inequality still holds when taking the suprema over $R_1<R_2$ and $\rho(x,x')<R_1$ in \eqref{eq-without-suprema-1}. Since $\epsilon>0$ was arbitrary, this proves the first inequality.

The other direction is similar. Let $\epsilon>0$.
Let $R_1<R_2$ and $x'\in B(x,R_1)$. Then for $\delta>0$,
\begin{alignat}{3}
    \label{eq-without-suprema-2}
    &&&\left|\int_{B(x',R_2)\setminus B(x',R_1)} K(x',y) f(y) \, \mathrm{d}\mu(y) \right| \\
    \label{eq-diff-small-2}
    \le &&&\left|\int_{R_1-\delta<\rho(x',y)< R_1} K(x',y) f(y) \, \mathrm{d}\mu(y) \right| \\
    \label{eq-other-2}
    &+&&\left|\int_{R_1-\delta<\rho(x',y)<R_2} K(x',y) f(y) \, \mathrm{d}\mu(y) \right| .
\end{alignat}

By \Cref{small-annulus}, we can choose $\delta$ such that \eqref{eq-diff-small-2} is bounded by $\epsilon$. Without loss of generality, we can assume $\rho(x,x')<R_1-\delta$. Then \eqref{eq-other-2} is bounded by the left hand side of \eqref{tang-unm-op-eq} and we obtain
\begin{equation*}
    \le \epsilon + \sup_{R_1 < R_2} \sup_{x'\in B(x,R_1)} \left|\int_{R_1<\rho(x',y)<R_2} K(x',y) f(y) \, \mathrm{d}\mu(y) \right| .
\end{equation*}
The inequality still holds when taking the suprema over $R_1<R_2$ and $\rho(x,x')<R_1$ in \eqref{eq-without-suprema-1}. Since $\epsilon>0$ was arbitrary, this proves the second inequality.
\end{proof}

\begin{proof}[Proof of \Cref{nontangential-from-simple}]
    \proves{nontangential-from-simple}
\leanok
Fix $g$ as in the Lemma.
Applying \Cref{simple-nontangential-operator} with a
sequence of $r$ tending to $0$ and using Lebesgue monotone convergence shows
\begin{equation}\label{tzerobound}
    \|T_*^{0}g\|_2\le 2^{a^3+26a+6}\|g\|_2,
\end{equation}
where
\begin{equation}\label{eq-simpler--nontangential}
    T_*^{0} g(x):=\sup_{0<R}\sup_{x'\in B(x,R)} \left|\int_{X\setminus B(x',R)}
K(x',y) g(y)\, d\mu(y)\right|\, .
\end{equation}
We now write using \Cref{nontangential-operator-boundary} and the triangle inequality,
\begin{alignat*}{3}\label{concretetstartriangle}
    T_* g(x)\le&&&\sup_{0<R_1<R_2}\sup_{x'\in B(x,R_1)} \left|\int_{X\setminus B(x',R_1)}K(x',y) g(y)\, d\mu(y)\right| \\
&+&&\sup_{0<R_1<R_2}\sup_{x'\in B(x,R_1)} \left|\int_{X\setminus B(x',R_2)} K(x',y) g(y)\, d\mu(y)\right|\, .
\end{alignat*}
Noting that the first integral does not depend on $R_2$ and
estimating the second summand by the larger supremum over all
$x'\in B(x,R_2)$, at which time the integral does not depend on $R_1$, we estimate further
\begin{alignat*}{3}
    \le&&&\sup_{0<R_1}\sup_{x'\in B(x,R_1)} \left|\int_{X\setminus B(x',R_1)}K(x',y) g(y)\, d\mu(y)\right| \\
    &+  &&\sup_{0<R_2}\sup_{x'\in B(x,R_2)} \left|\int_{X\setminus B(x',R_2)} K(x',y) g(y)\, d\mu(y)\right|\, .
\end{alignat*}
    Applying the triangle inequality on the left-hand side
    of \eqref{concretetstarbound} and applying
     \eqref{tzerobound} twice
    proves \eqref{concretetstarbound}.
    This completes the proof of \Cref{nontangential-from-simple}.
\end{proof}

\subsection{Calder\'on-Zygmund Decomposition}
\label{subsec-CZD}
Calder\'on-Zygmund decomposition is a tool to extend $L^2$ bounds to $L^p$ bounds with $p<2$ or to the so-called weak $(1, 1)$ type endpoint bound.
It is classical and can be found in \cite{stein-book}.

The following lemma is Theorem 3.1(b) in \cite{stein-book}. The proof uses \Cref{Hardy-Littlewood}.
\begin{lemma}[Maximal theorem]
    \label{maximal-theorem}
    \leanok
    \lean{maximal_theorem}
    \uses{Hardy-Littlewood}
    Let $f: X \to \C$ be bounded, measurable, supported on a set of finite measure, and let $\alpha > 0$. Then
    \begin{equation}
        \label{maximal-theorem-equation}
        \mu(\{x\in X : Mf(x) > \alpha\}) \le \frac{2^{2a}}{\alpha} \int |f(y)|\, d\mu(y).
    \end{equation}
\end{lemma}
\begin{proof}
    \leanok
    By definition, for each $x\in X$ with $Mf(x) > \alpha$, there exists a ball $B_x$ such that $x\in B_x$ and
    \begin{equation}
        \label{maximal-theorem-a}
        \alpha \mu(B_x) < \int_{B_x} |f(y)|\, d\mu(y).
    \end{equation}
    Since $\{x\in X : Mf(x) > \alpha\}$ is open and $\mu$ is inner regular on open sets, it suffices to show that
    \begin{equation*}
        \mu(E) \le \frac{2^{2a}}{\alpha} \int |f(y)|\, d\mu(y)
    \end{equation*}
    for every compact $E\subset \{x\in X : Mf(x) > \alpha\}$.
    For such an $E$, by compactness, we can select a finite subcollection $\mathcal{B} \subset \{B_x: x\in E\}$ that covers $E$.
    By \eqref{eq-besico} applied to \eqref{maximal-theorem-a},
    \begin{equation}
        \alpha \mu(\bigcup \mathcal{B}) \le 2^{2a} \int |f(y)|\,d\mu(y)
    \end{equation}
    and hence
    \begin{equation*}
        \mu(E) \le \mu(\bigcup \mathcal{B}) \le \frac{2^{2a}}{\alpha} \int |f(y)|\,d\mu(y).
    \end{equation*}
\end{proof}

\begin{lemma}[Lebesgue differentiation]
    \label{Lebesgue-differentiation}
    \leanok
    \lean{lebesgue_differentiation}
    Let $f$ be a bounded measurable function supported on a set of finite measure. Then for $\mu$ almost every $x$, we have
    $$\lim_{n\to \infty} \frac{1}{\mu(B_n)}\int_{B_n} f(y)\, dy= f(x),$$
    where $\{B_n\}_{n\geq 1}$ is a sequence of balls with radii $r_n>0$ such that $x\in B_n$ for each $n\geq 1$ and
    \begin{equation*}
        \lim_{n\to \infty} r_n=0 \,.
    \end{equation*}
\end{lemma}
\begin{proof}
    \leanok
    This follows from the Lebesgue differentiation theorem, which is already formalized in Lean.
\end{proof}

\begin{lemma}[Disjoint family countable]
    \leanok
    \lean{Pairwise.countable_of_isOpen_disjoint}
    \label{disjoint-family-countable}
    In a doubling metric measure space $(X,\rho,\mu, a)$, every disjoint family of balls $B_j = B(x_j, r_j)$, $j\in J$, is countable.
\end{lemma}
\begin{proof}
    \leanok
    Choose an arbitrary $x\in X$ as reference point. For $q, Q\in\Q_+$, let $J_{q,Q}$ denote the set of all $j\in J$ such that $B_j\subset B(x, Q)$ and $r_j \ge q$. It suffices to show that all the $J_{q,Q}$ are finite.
    Indeed, for all $j\in J_{q,Q}$,
    \begin{equation*}
        \mu(B(x, Q)) \le \mu(B(x_j, 2Q)) = \mu(B(x_j, \frac{2 Q}{r_j} r_j))
        \le 2^{a\log_2{\lceil \frac{2 Q}{r_j}\rceil}} \mu(B_j).
    \end{equation*}
    Since the $B_j$ are disjoint,
    \begin{equation}
        |J_{q,Q}| \mu(B(x, Q)) \le 2^{a\log_2{\lceil\frac{2Q}{q}\rceil}} \sum_{j\in J_{q,Q}} \mu(B_j) \le 2^{a\log_2{\lceil\frac{2Q}{q}\rceil}} \mu(B(x,Q))
    \end{equation}
    and hence $|J_{q,Q}| \le 2^{a\log_2{\lceil\frac{2Q}{q}\rceil}}$.
\end{proof}

The following lemma corresponds to Lemma 3.2 in \cite{stein-book} with additional proof of the bounded intersection property taken from the proof of Proposition 7.1.
It uses the following notation, which will be used throughout the rest of this section: whenever $B_j$ denotes a ball $B_j = B(x_j, r_j)$, define
\begin{equation*}
    B_{n,j} = B(x_j, nr_j).
\end{equation*}
\begin{lemma}[Ball covering]
    \label{ball-covering}
    \leanok
    \lean{ball_covering}
    \uses{disjoint-family-countable}
    Given an open set $O\ne X$, there exists a countable family of balls $B_j = B(x_j, r_j)$ such that
    \begin{equation}
        \label{balls-disjoint}
        B_j \cap B_{j'} = \emptyset \quad \text{ for } j \ne j',
    \end{equation}
    and
    \begin{equation}
        \label{balls-covering}
        \bigcup_j B_{3,j} = O,
    \end{equation}
    and
    \begin{equation}
        \label{enlarged-balls-intersect-complement}
        B_{7,j} \cap (X \setminus O) \ne \emptyset \quad \text{ for all } j
    \end{equation}
    and we have the bounded intersection property that each $x\in O$ is contained in at most $2^{6a}$ of the $B_{3,j}$.
\end{lemma}
\begin{proof}
\leanok
Define for $x\in O$,
\begin{equation}
    \delta(x):= \sup \{\delta\in\R : B(x,\delta)\subset O\}.
\end{equation}
Since $O$ is open, and $O\ne X$, we have
\begin{equation}
    0 < \delta(x) < \infty \,.
\end{equation}
Using Zorn's Lemma, we select a maximal disjoint subfamily of $\{B(x,\frac{\delta(x)}{6}) : x \in O\}$.
We obtain a (by \Cref{disjoint-family-countable} countable) family of balls $B_j = B(x_j, \frac{\delta(x_j)}{6}), j \in J$ such that \eqref{balls-disjoint}, \eqref{enlarged-balls-intersect-complement}, and $\bigcup_j B_{3,j} \subset O$ are also immediate. For the other inclusion, first observe that for $x,y\in X$,  if $B(x,\frac{\delta(x)}{6}) \cap B(y,\frac{\delta(y)}{6}) \ne \emptyset$, then
\begin{equation*}
    \delta(x) \le \rho(x,y) + \delta(y) \le (\frac{\delta(x)}{6} + \frac{\delta(y)}{6}) + \delta(y) = \frac{\delta(x)}{6} + \frac{7\delta(y)}{6},
\end{equation*}
so
\begin{equation}
    \label{control-distance-growth}
    \delta(x) \le 2 \delta(y).
\end{equation}
Now let $z\in O$. By maximality, there exists some $j\in J$ with $B(z,\frac{\delta(z)}{6}) \cap B_j \ne \emptyset$. By \eqref{control-distance-growth},
\begin{equation*}
    \rho(z,x_j)< \frac{\delta(z)}{6} + \frac{\delta(x_j)}{6} \le \frac{3\delta(x_j)}{6} = 3r_j
\end{equation*}
and thus $z\in B_{3,j}$.

We now turn to the bounded intersection property. Assume that for some $j_1,\dots,j_N$,
\begin{equation}
    z\in \bigcap_{k=1}^N B_{3,j_k}.
\end{equation}
Similarly as above, observe for $1\le k \le N$,
\begin{equation}
    \label{control-distance-growth-b}
    \delta(z) \le \rho(z,x_{j_k}) + \delta(x_{j_k}) \le \frac{\delta(x_{j_k})}{2} + \delta(x_{j_k}) = \frac{3\delta(x_{j_k})}{2}
\end{equation}
and
\begin{equation*}
    \delta(x_{j_k}) \le \rho(x_{j_k},z) + \delta(z) \le \frac{\delta(x_{j_k})}{2} + \delta(z),
\end{equation*}
so
\begin{equation}
    \label{control-distance-growth-c}
    \delta(x_{j_k}) \le 2 \delta(z).
\end{equation}

By \eqref{control-distance-growth-b} and \eqref{control-distance-growth-c}, for all $1\le k \le N$, $B(z,\frac{\delta(z)}{6}) \subset B(x_{j_k}, 5r_{j_k})$ and $B_{j_k} \subset B(z,\frac{8\delta(z)}{6})$. Using this and \eqref{balls-disjoint}, we obtain
\begin{align}
    N \mu(B(z,\frac{\delta(z)}{6})) &\le \sum_{k=1}^N \mu(B(x_{j_k}, 5r_j)) \le 2^{3a} \sum_{k=1}^N \mu(B_{j_k}) \\
    &= 2^{3a} \mu(\bigcup_{k=1}^N B_{j_k}) \le 2^{3a} \mu(B(z,\frac{8\delta(z)}{6})) \le 2^{6a} \mu(B(z,\frac{\delta(z)}{6}))
\end{align}
and conclude $N\le 2^{6a}$.
\end{proof}

Most of the next lemma and its proof is taken from Theorem 4.2 in \cite{stein-book}.
\begin{lemma}[Calderon Zygmund decomposition]
    \label{Calderon-Zygmund-decomposition}
    \leanok
    \lean{encard_czBall3_le, tsum_czRemainder', aemeasurable_czApproximation, czApproximation_add_czRemainder, enorm_czApproximation_le, enorm_czApproximation_le_infinite, eLpNorm_czApproximation_le, support_czRemainder'_subset, integral_czRemainder', integral_czRemainder, eLpNorm_czRemainder'_le, eLpNorm_czRemainder_le, tsum_volume_czBall3_le, volume_univ_le, tsum_eLpNorm_czRemainder'_le, tsum_eLpNorm_czRemainder_le}
    \uses{ball-covering,Lebesgue-differentiation,maximal-theorem}
    Let $f$ be a bounded, a.e. measurable function supported on a set of finite measure and let $\alpha>\frac{1}{\mu(X)}\int |f|\,d\mu$.
    Then there exists a bounded a.e. measurable function $g$ supported on a set of finite measure, a countable family of balls $B_{3,j}$ (where we allow $B_{3,1} = X$ in the special case that $\mu(X)<\infty$)
    such that each $x\in X$ is contained in at most $2^{6a}$ of the $B_{3,j}$, and a countable family of a.e. measurable functions $\{b_j\}_{j\in J}$ such that for all $x \in X$
    \begin{equation}
       \label{eq-gb-dec}
       f(x)= g(x) + \sum_{j} b_j(x)
    \end{equation}
    and such that the following holds. For almost every $x\in X$,
    \begin{equation}
        \label{eq-g-max}
       |g(x)|\leq 2^{3a} \alpha\,.
    \end{equation}
    We have
    \begin{equation}
        \label{eq-g-L1-norm}
        \int |g(y)|\, d\mu(y)\leq \int |f(y)|\, d\mu(y).
    \end{equation}
    For every $j$
    \begin{equation}
        \label{eq-supp-bj}
        \operatorname{supp} b_j \subset B_{3,j}\,.
    \end{equation}
    For every $j$
    \begin{equation}
        \label{eq-bad-mean-zero}
        \int_{B_{3,j}} b_j(x)\, d\mu(x)=0,
    \end{equation}
    and
     \begin{equation}
        \label{eq-bj-L1}
        \int_{B_{3,j}} |b_j(x)|\, d\mu(x) \leq 2^{2a+1} \alpha \mu(B_{3,j}).
    \end{equation}
    We have
    \begin{equation}
        \label{eq-bset-length-sum}
        \sum_j \mu(B_{3,j})\leq \frac{2^{6a}}{\alpha}\int |f(y)|\, d\mu(y)
    \end{equation}
    and
    \begin{equation}
    \label{eq-b-L1}
    \sum_{j}\int_{B_{3,j}} |b_j(y)|\, d\mu(y)\leq 2 \int |f(y)|\, d\mu(y)\,.
    \end{equation}
\end{lemma}

\begin{proof}
\leanok
Let $E_\alpha:=\{x\in X: Mf(x)>\alpha\}$.
Then $E_\alpha$ is open. Assume first that $E_\alpha \ne X$. We apply \Cref{ball-covering} with $O=E_\alpha$ to obtain the family $B_j, j\in J,$. Without loss of generality, we can assume $J=\N$. Define inductively
\begin{equation}
    Q_j := B_{3,j} \setminus \left(\bigcup_{i<j} Q_i \cup \bigcup_{i>j} B_i \right).
\end{equation}
Then $B_j\subset Q_j\subset B_{3,j}$, the $Q_j$ are pairwise disjoint and $\bigcup_j Q_j = E_\alpha$.
Define
\begin{equation}
    \label{eq-g-def}
    g(x):=\begin{cases}
     f(x), & x\in X\setminus E_\alpha,\\
     \frac{1}{\mu(Q_j)}\int_{Q_j} f(y)\, d\mu(y), &x\in Q_j,
    \end{cases}
\end{equation}
and, for each $j$,
\begin{equation}
    b_j(x):= \mathbf{1}_{Q_j}(x) \left(f(x)-\frac{1}{\mu(Q_j)}\int_{Q_j} f(y)\, d\mu(y) \right).
\end{equation}
Then \eqref{eq-gb-dec}, \eqref{eq-supp-bj} and \eqref{eq-bad-mean-zero} are true by construction.
Boundedness of $g$ is immediate from the definition, as is the fact that $\supp g$ is contained within
$\supp f \cup \bigcup_j Q_j$. Now the fact that $\supp g$ has finite measure follows from
$\bigcup_j Q_j = E_\alpha$ and \eqref{eq-hlm-2}.

For \eqref{eq-g-max}, we first do the case $x\in X\setminus E_\alpha$. By definition of $Mf$,
\begin{equation}
    \frac{1}{\mu(B)}\int_B |f(y)|\,d\mu(y) \le \alpha
\end{equation}
for every ball $B\subset X$ with $x\in B$. It follows by \Cref{Lebesgue-differentiation} that for almost every $x\in X\setminus E_\alpha$, $|f(x)|\le \alpha$.
In the case $x\in E_\alpha$, there exists some $j\in J$ with $x\in Q_j$ and we have that
\begin{equation}
    \label{large-ball-estimate}
    \frac{1}{\mu(B_{7,j})} \int_{B_{7,j}} |f(y)| \,d\mu(y) \le \alpha
\end{equation}
because $B_{7,j}\cap (X\setminus E_\alpha) \ne \emptyset$. We get
\begin{equation}
    |g(x)| \le \frac{1}{\mu(Q_j)}\int_{Q_j} |f(y)|\, d\mu(y) \le \frac{1}{\mu(B_j)}\int_{B_{7,j}} |f(y)| \, d\mu(y) \le 2^{3a}\alpha .
\end{equation}

To prove \eqref{eq-g-L1-norm}, we estimate
\begin{align*}
    \int |g(z)|\, d\mu(z) &\le \int_{X\setminus E_\alpha} |f(z)|\, d\mu(z) + \sum_{j} \int_{Q_j}\frac{1}{\mu(Q_j)}\int_{Q_j}|f(y)|\, d\mu(y)\,d\mu(z) \\
    &= \int |f(z)|\,d\mu(z).
\end{align*}
Using the triangle inequality, we have that
\begin{align}
    \label{eq-bj-int}
    \int_{B_{3,j}} |b_j(y)|\, dy &\le \int_{Q_j} |f(y)|\, d\mu(y) + \int_{Q_j} \frac{1}{\mu(Q_j)}\int_{Q_j} |f(x)|\, d\mu(x)\, d\mu(y) \\
    &= 2 \int_{Q_j} |f(y)|\, dy.
\end{align}
With \eqref{large-ball-estimate}, we estimate further
\begin{equation}
    \le 2 \int_{B_{7,j}} |f(y)|\, dy \le 2\mu(B_{7,j})\alpha \le 2^{2a+1} \alpha \mu(B_{3,j})
\end{equation}
to obtain \eqref{eq-bj-L1}.
Further, summing up \eqref{eq-bj-int} in $j$ yields \eqref{eq-b-L1}.
At last, we estimate with \Cref{maximal-theorem}
\begin{equation}
    \sum_j \mu(B_{3,j}) \le 2^{2a} \sum_j \mu(B_j) \le 2^{2a} \mu(E_\alpha) \le \frac{2^{6a}}{\alpha}\int |f(y)|\, d\mu(y),
\end{equation}
proving \eqref{eq-bset-length-sum}.

Assume now that $E_\alpha = X$. It follows from \Cref{maximal-theorem} that then $\mu(X)<\infty$.
Define
\begin{equation*}
    g := \frac{1}{\mu(X)} \int |f(y)|\,d\mu(y)
\end{equation*}
and
\begin{equation*}
    b_1 := f - g.
\end{equation*}
Then $f = g + b_1$ and $\supp b_1 \subset B_{3,1}:=X$ and \eqref{eq-gb-dec}, \eqref{eq-g-L1-norm}, \eqref{eq-supp-bj}, \eqref{eq-bad-mean-zero} all hold immediately. By assumption, $\alpha>\frac{1}{\mu(X)}\int |f|\,d\mu = g$, so \eqref{eq-g-max} holds.
We also have, using the definitions and the same assumption,
\begin{equation}
    \int |b_1(y)|\, d\mu(y) \le 2 \int |f(y)|\,d\mu(y) \le 2\alpha\mu(X),
\end{equation}
which verifies both \eqref{eq-b-L1} and \eqref{eq-bj-L1}.
Finally, by \Cref{maximal-theorem},
\begin{equation*}
    \mu(X) \le \frac{2^{2a}}{\alpha} \int |f(y)|\,d\mu(y),
\end{equation*}
which shows \eqref{eq-bset-length-sum}.
\end{proof}

We use \Cref{Calderon-Zygmund-decomposition} to prove \Cref{calderon-zygmund-weak-1-1}. For the remainder of this section, let $f:X\to\C$, $r>0$ and $\alpha>0$ as in \Cref{calderon-zygmund-weak-1-1}.
We define the constant
\begin{equation} \label{weak-1-1-proof-cz-const}
    c:= 2^{-a^3-12a-4}
\end{equation}
and $\alpha' := c\alpha$. If $\alpha'\le\frac{1}{\mu(X)}\int |f|\,d\mu$, then we directly have
\begin{equation*}
    \mu\left(\{x\in X: |T_r f(x)|>\alpha\}\right)\le \mu(X) \le \frac{1}{\alpha'} \int |f(y)|\, d\mu(y)
    \le \frac{2^{a^3 + 19a}}{\alpha} \int |f(y)|\, d\mu(y),
\end{equation*}
which proves \eqref{eq-weak-1-1}.
So assume from now on that $\alpha'>\frac{1}{\mu(X)}\int |f|\,d\mu$.
Using \Cref{Calderon-Zygmund-decomposition} for $f$ and $\alpha'$, we obtain the decomposition
\begin{equation*}
    f=g+b=g+\sum_j b_j
\end{equation*}
such that the properties \eqref{eq-gb-dec}-\eqref{eq-b-L1} are satisfied (with $\alpha'$ replacing $\alpha$). Let
\begin{equation}
    \label{eq-Ij-cj}
    B_{3,j}=B(x_j, 3r_j)
\end{equation}
as in the lemma. Then
\begin{equation}
    \label{eq-Ij*}
    B_{6,j}=B(x_j, 6r_j)
\end{equation}
is a ball with the same center as $B_{3,j}$ but with
\begin{equation}
    \label{eq-Ij*-dim}
    \mu(B_{6,j})\le 2^{a} \mu(B_{3,j}).
\end{equation}
Let
\begin{equation}
    \label{eq-omega}
    \Omega:=\bigcup_j B_{6,j}.
\end{equation}
We deal with $T_rg$ and $T_rb$ separately in the following lemmas.

\begin{lemma}[Estimate good]
    \label{estimate-good}
    \leanok
    \lean{estimate_good}
    \uses{Calderon-Zygmund-decomposition}
    \begin{equation*}
        \mu\left(\{x\in X: |T_r g(x)|>{\alpha}/2\}\right)
        \le \frac{2^{2a^3+3a+2}c}{\alpha} \int |f(y)|\, d\mu(y).
    \end{equation*}

\end{lemma}

\begin{proof}
    \leanok
    We estimate using monotonicity of the integral
\begin{equation*}
     \mu\left(\{x\in X: |T_r g(x)|>{\alpha}/2\}\right)\leq \frac{4}{\alpha^2} \int |T_r g(y)|^2\, d\mu(y).
\end{equation*}
Using \eqref{eq-strong-2-2-assumption} followed by \eqref{eq-g-max} and \eqref{eq-g-L1-norm}, we estimate the right hand side above by
\begin{equation}
    \label{eq-Hr-g}
    \leq \frac{4\cdot 2^{2a^3}}{\alpha^2} \int |g(y)|^2\, d\mu(y)\leq \frac{2^{2a^3+3a+2}c}{\alpha} \int |g(y)|\, dy \le \frac{2^{2a^3+3a+2}c}{\alpha} \int |f(y)|\, d\mu(y).
\end{equation}
\end{proof}

\begin{lemma}[Estimate bad partial]
    \label{estimate-bad-partial}
    \leanok
    \lean{estimate_bad_partial}
    \uses{Calderon-Zygmund-decomposition}
    Let $x\in X\setminus\Omega$. Then
    \begin{equation*}
        |T_rb(x)| \le 3F(x)+\alpha/8,
    \end{equation*}
    where
    \begin{equation*}
        F(x) := 2^{a^3+2a+1} c\alpha \sum_{j\in J} \left(\frac{3r_j}{\rho(x,x_j)}\right)^{\frac{1}{a}}\frac{\mu(B_{3,j})}{V(x,x_j)}.
    \end{equation*}
\end{lemma}

\begin{proof}
\leanok
We decompose the index set $J$ into the following disjoint sets:
\begin{align*}
    \mathcal{J}_1(x)&:=\{j\,: r+3r_j \le \rho(x,x_j) \},\\
    \mathcal{J}_2(x)&:=\{j\,: r-3r_j \le \rho(x,x_j) < r+3r_j\},\\
    \mathcal{J}_3(x)&:=\{j\,: \rho(x,x_j) < r-3r_j\}.
\end{align*}
Then
\begin{alignat}{3}
    \label{eq-b-dec-1}
    |T_r b(x)|\le&&&\sum_{j\in \mathcal{J}_1(x)} |T_rb_j(x)| \\
    \label{eq-b-dec-2}
                &+&&\sum_{j\in \mathcal{J}_2(x)} |T_rb_j(x)| \\
    \label{eq-b-dec-3}
                &+&&\sum_{j\in \mathcal{J}_3(x)} |T_rb_j(x)|.
\end{alignat}
For all $j\in \mathcal{J}_3(x)$, $\supp b_j\subset B_{3,j}\subset B(x,r)$, and thus $T_rb_j(x)=0$, so $\eqref{eq-b-dec-3} = 0$.

Next, for $j\in \mathcal{J}_1(x)$, $\supp b_j\subset B_{3,j}\subset X \setminus B(x,r)$, and we have
\begin{equation*}
    T_rb_j(x)=\int_{X\setminus B(x,r)} K(x,y) b_j(y)\,d\mu(y)=\int_{B_{3,j}} K(x,y) b_j(y)\,d\mu(y)\,.
\end{equation*}

Using \eqref{eq-bad-mean-zero}, the above is equal to
\begin{equation*}
    \int_{{3,j}} (K(x,y)-K(x,x_j)) b_j(y)\,d\mu(y)\,.
\end{equation*}
Since $x\in X\setminus\Omega$, we have for each $y\in B_{3,j}$ that
\begin{equation}
    \label{eq-Om-cj}
    \rho(x,x_j)\ge 6r_j > 2\rho(x_j,y),
\end{equation}
so we can apply \eqref{eqkernel-y-smooth} to estimate
\begin{equation*}
    \eqref{eq-b-dec-1} \le \sum_{j\in \mathcal{J}_1(x)}\int_{B_{3,j}}\left(\frac{\rho(x_j,y)}{\rho(x,x_j)}\right)^{\frac{1}{a}}\frac{2^{a^3}}{V(x,x_j)} |b_j(y)|\, d\mu(y)
\end{equation*}
\begin{equation}
    \le 2^{a^3} \sum_{j} \left(\frac{3r_j}{\rho(x,x_j)}\right)^{\frac{1}{a}}\frac{1}{V(x,x_j)}\int_{B_{3,j}} |b_j(y)|\, dy
\end{equation}
and by \eqref{eq-bj-L1},
\begin{equation}
    \label{eq-J1-diff-est}
    \le 2^{a^3+2a+1} c\alpha \sum_{j} \left(\frac{3r_j}{\rho(x,x_j)}\right)^{\frac{1}{a}}\frac{\mu(B_{3,j})}{V(x,x_j)} = F(x).
\end{equation}
Next, we estimate \eqref{eq-b-dec-2}. For each $j\in \mathcal{J}_2(x)$, set
\begin{equation*}
    d_j:=\frac{1}{\mu(B_{3,j})}\int_{B_{3,j}} \mathbf{1}_{X\setminus B(x,r)}(y) b_j(y)\, dy.
\end{equation*}
Then by \eqref{eq-bj-L1}
\begin{equation}
    \label{eq-dj-est}
    |d_j|\le 2^{2a+1} c\alpha.
\end{equation}
For each $j\in \mathcal{J}_2(x)$, we have
\begin{align*}
    T_r b_j(x)&=\int_{B_{3,j}} K(x,y) (\mathbf{1}_{X\setminus B(x,r)}(y)b_j(y)-d_j)\, dy + \int_{B_{3,j}} d_j K(x,y) \, dy \\
    &= \int_{B_{3,j}} (K(x,y)-K(x,x_j)) (\mathbf{1}_{X\setminus B(x,r)}(y)b_j(y)-d_j)\, dy + \int_{B_{3,j}} d_j K(x,y) \, dy.
\end{align*}
Thus, using the triangle inequality, the equation above and \eqref{eq-dj-est}, we obtain
\begin{equation*}
    |T_r b_j(x)|\le
\end{equation*}
\begin{equation}
    \label{eq-J2-diff-est}
    \int_{B_{3,j}} |K(x,y)-K(x,x_j)| \left(|b_j(y)|+2^{2a+1} c\alpha\right)\, dy +2^{2a+1} c\alpha \int_{B_{3,j}}  |K(x,y)| \, dy.
\end{equation}
By \eqref{eq-Om-cj}, we can apply \eqref{eqkernel-y-smooth} and arguing as in \eqref{eq-J1-diff-est}, we get that
\begin{equation}
    \label{eq-J2-diff-est-2}
    \eqref{eq-b-dec-2} \le 2F(x) + 2^{2a+1}c\alpha \sum_{j\in\mathcal{J}_2(x)} \int_{B_{3,j}}  |K(x,y)|\,d\mu(y),
\end{equation}
with $F$ as in \eqref{eq-J1-diff-est}.
Define
\begin{equation*}
    A := \bigcup_{j \in \mathcal{J}_2(x)} B_{3,j}.
\end{equation*}
We claim that
\begin{equation}
    \label{eq-J2-union-subset}
    A\subset B(x,3r) \setminus B(x,\frac{r}{3}).
\end{equation}
Indeed, for each $j\in \mathcal{J}_2(x)$ and $y\in B_{3,j}$, using again \eqref{eq-Om-cj},
\begin{equation*}
    \rho(x,x_j) < r+r_{3,j} \le r + \frac{1}{2} \rho(x,x_j) \implies \rho(x,x_j) < 2r
\end{equation*}
and hence
\begin{equation*}
    \rho(x,y) \le \rho(x,x_j) + \rho(x_j,y) < 2r + 3r_j \le 2r + \frac{1}{2}\rho(x,x_j) < 3r.
\end{equation*}
For the lower bound, we observe
\begin{equation*}
    \rho(x,x_j) \ge r-3r_j \ge r - \frac{1}{2} \rho(x,x_j) \implies \rho(x,x_j) \ge \frac{2}{3}r,
\end{equation*}
and conclude
\begin{equation*}
    \rho(x,y) \ge \rho(x,x_j) - \rho(y,x_j) \ge \rho(x,x_j) - 3r_j \ge \rho(x,x_j) - \frac{1}{2}\rho(x,x_j) \ge \frac{1}{3} r.
\end{equation*}

Using the bounded intersection property of the $B_{3,j}$, \eqref{eq-J2-union-subset} and \eqref{eqkernel-size}, we get
\begin{align}
    \sum_{j\in\mathcal{J}_2(x)} \int_{B_{3,j}}  |K(x,y)|\,d\mu(y) &\le 2^{6a} \int_{A}  |K(x,y)|\,d\mu(y) \\
    &\le 2^{6a} \int_{B(x,3r) \setminus B(x,\frac{r}{3})}  |K(x,y)|\,d\mu(y) \\
    &\le 2^{6a} \int_{B(x,3r) \setminus B(x,\frac{r}{3})} \frac{2^{a^3}}{V(x,y)} \,d\mu(y) \\
    &\le 2^{a^3+6a} \int_{B(x,3r) \setminus B(x,\frac{r}{3})} \frac{1}{\mu(B(x,\frac{r}{3}))} \,d\mu(y) \\
    &\le 2^{a^3+6a} \frac{\mu(B(x,3r))}{\mu(B(x,\frac{r}{3}))} \\
    \label{eq-J2-diff-est-4}
    &\le 2^{a^3+10a}.
\end{align}

If we combine the estimates \eqref{eq-J1-diff-est} for \eqref{eq-b-dec-1}, \eqref{eq-J2-diff-est-2} for \eqref{eq-b-dec-2}, and \eqref{eq-J2-diff-est-4}, we get
\begin{equation*}
    |T_rb(x)|\leq 3F(x)+2^{a^3+12a+1}c\alpha.
\end{equation*}
By the definition \eqref{weak-1-1-proof-cz-const} of $c$, this equals
\begin{equation*}
    3F(x)+\alpha/8.
\end{equation*}
\end{proof}

\begin{lemma}[Estimate F set]
    \label{estimate-F-set}
    \leanok
    \lean{distribution_czOperatorBound}
    \uses{Calderon-Zygmund-decomposition,geometric-series-estimate}
    For $F$ as defined in \Cref{estimate-bad-partial}, we have
    \begin{equation}
        \label{eq-F-X-minus-Omega}
        \mu(\{x\in X\setminus\Omega: F(x)>\alpha/8\}) \le \frac{2^{a^3+11a+4}}{\alpha} \int |f(y)|\,d\mu(y)\,.
    \end{equation}
\end{lemma}

\begin{proof}
    \leanok
    We estimate
    \begin{align}
        \mu(\{x\in X\setminus\Omega&: F(x)> \alpha/8\})
        \le \frac{8}{\alpha} \int_{X\setminus \Omega} F(x)\,d\mu(x) \\
        &\le \frac{8}{\alpha} \int_{X\setminus \Omega} 2^{a^3+2a+1} c\alpha \sum_{j} \left(\frac{3r_j}{\rho(x,x_j)}\right)^{\frac{1}{a}}\frac{\mu(B_{3,j})}{V(x,x_j)}\,d\mu(x) \\
        \label{eq-F-est-1}
        &\le 2^{a^3+2a+4} c \sum_{j} \mu(B_{3,j}) \int_{X\setminus B_{6,j}} \left(\frac{3r_j}{\rho(x,x_j)}\right)^{\frac{1}{a}}\frac{1}{V(x,x_j)}\,d\mu(x)
    \end{align}
    Using
    \begin{equation*}
        V(x,x_j) = \mu(B(x,\rho(x,x_j))) \ge 2^{-a}\mu(B(x,2\rho(x,x_j))) \ge 2^{-a} \mu(B(x_j,\rho(x_j,x))),
    \end{equation*}
    we have for all $j\in J$,
    \begin{align*}
        &\int_{X\setminus B_{6,j}} \left(\frac{3r_j}{\rho(x,x_j)}\right)^{\frac{1}{a}}\frac{1}{V(x,x_j)}\,d\mu(x) \\
        \le& 2^a \int_{X\setminus B_{6,j}} \left(\frac{3r_j}{\rho(x,x_j)}\right)^{\frac{1}{a}}\frac{1}{\mu(B(x_j,\rho(x_j,x)))}\,d\mu(x) \\
        \le&2^a \sum_{n=1}^\infty \int_{B(x_j,2^{n+1}\cdot3r_j)\setminus B(x_j,2^n\cdot3r_j)} \left(\frac{3r_j}{2^n \cdot3r_j}\right)^{\frac{1}{a}}\frac{1}{\mu(B(x_j,2^n \cdot3r_j))}\,d\mu(x) \\
        \le&2^a \sum_{n=1}^\infty 2^{-\frac{n}{a}} \frac{\mu(B(x_j,2^{n+1}\cdot3r_j))}{\mu(B(x_j,2^n \cdot3r_j))} \\
        \le&2^{3a},
    \end{align*}
    where we used \Cref{geometric-series-estimate} in the last step.
    Plugging this into \eqref{eq-F-est-1} and using \eqref{eq-bset-length-sum}, we conclude that
    \begin{equation*}
        \mu(\{x\in X\setminus\Omega: F(x)> \alpha/8\}) \le \frac{2^{a^3+11a+4}}{\alpha} \int |f(y)|\,d\mu(y)\,.
    \end{equation*}
\end{proof}

\begin{lemma}[Estimate bad]
    \label{estimate-bad}
    \leanok
    \lean{estimate_bad}
    \uses{Calderon-Zygmund-decomposition,estimate-bad-partial,estimate-F-set}
    We have
    \begin{equation*}
        \mu\left({\{x\in X: |T_r b(x)|>\alpha/2\}}\right) \le  \frac{\frac{2^{7a}}{c} + 2^{a^3+11a+4}}{\alpha} \int |f(y)|\,d\mu(y) \,.
    \end{equation*}
\end{lemma}

\begin{proof}
    \leanok
    We estimate
    \begin{equation*}
        \mu\left(\{x\in X: |T_r b(x)|>\alpha/2\}\right)
    \end{equation*}
    \begin{equation}
        \label{eq-set-dec-2}
        \le \mu (\Omega) +  \mu\left(\{x\in X\setminus\Omega: |T_r b(x)|>{\alpha}/2\}\right)\,.
    \end{equation}
    Using \eqref{eq-Ij*-dim} and \eqref{eq-bset-length-sum}, we conclude that
    \begin{equation}
        \label{eq-omega-bd}
        \mu(\Omega) \le \sum_{j} \mu (B_{6,j})
        \le 2^a \sum_j \mu(B_{3,j}) \le \frac{2^{7a}}{c\alpha} \int |f(y)|\, d\mu(y)\,.
    \end{equation}
    It follows from \Cref{estimate-bad-partial} and the triangle inequality that
    \begin{equation}
        \label{eq-set-dec-3}
        \mu({\{x\in X\setminus\Omega: |T_r b(x)|>\alpha/2\}}) \le \mu(\{x\in X\setminus\Omega: F(x)> \alpha/8\})\,.
    \end{equation}
    The claim now follows from \Cref{estimate-F-set}.
\end{proof}

\begin{proof}[Proof of \Cref{calderon-zygmund-weak-1-1}]
\leanok
\proves{calderon-zygmund-weak-1-1}
It follows by the triangle inequality and subadditivity of $\mu$ that
\begin{equation*}
    \mu\left(\{x\in X: |T_r f(x)|>\alpha\}\right)
\end{equation*}
\begin{equation*}
\label{eq-set-dec-1}
 \le \mu\left(\{x\in X: |T_r g(x)|>{\alpha}/2\}\right) +  \mu\left(\{x\in X: |T_r b(x)|>{\alpha}/2\}\right).
\end{equation*}
Using \Cref{estimate-good}, \Cref{estimate-bad} and the definition \eqref{weak-1-1-proof-cz-const} of $c$, we get
\begin{align*}
    \le& \frac{2^{2a^3+3a+2}c + \frac{2^{7a}}{c} + 2^{a^3+11a+4}}{\alpha} \int |f(y)|\, d\mu(y) \\
    =& \frac{2^{a^3-9a-2} + 2^{a^3+19a+4} + 2^{a^3+11a+4}}{\alpha} \int |f(y)|\, d\mu(y) \\
    \le& \frac{2^{a^3 + 21a}}{\alpha} \int |f(y)|\, d\mu(y).
\end{align*}
\end{proof}

\section{Proof of The Classical Carleson Theorem}

The convergence of partial Fourier sums is proved in
\Cref{10classical} in two steps. In the first step,
we establish convergence on a suitable dense subclass of functions. We choose smooth functions as subclass, the convergence is stated in \Cref{convergence-for-smooth} and proved in \Cref{10smooth}.
In the second step, one controls the relevant error of approximating a general function by a function in
the subclass. This is stated in \Cref{control-approximation-effect} and proved
in \Cref{10difference}.
The proof relies on a bound on the real Carleson maximal operator stated in \Cref{real-Carleson} and proved in \Cref{10carleson}, which involves showing that the real line fits into the setting of \Cref{overviewsection}.
This latter proof refers to the two-sided variant of the Carleson \Cref{two-sided-metric-space-Carleson}. Two assumptions in \Cref{metric-space-Carleson} require more work. The boundedness of the operator $T_r$ defined in \eqref{def-T-r} is established in \ref{Hilbert-strong-2-2}. This lemma is proved in \Cref{10hilbert}.
The cancellative property is verified by \Cref{van-der-Corput}, which is proved in \Cref{10vandercorput}.
Several further auxiliary lemmas are stated and proved in \Cref{10classical}, the proof of one of these auxiliary lemmas, \Cref{spectral-projection-bound}, is done in \Cref{10projection}.

All subsections past \Cref{10classical} are mutually independent.

\subsection{The classical Carleson theorem}
\label{10classical}

Let a uniformly continuous $2\pi$-periodic function $f:\R\to \mathbb{C}$ and $\epsilon>0$ be given.
Let
\begin{equation}
    C_{a,q} := \frac{2^{474a^3}}{(q-1)^6}
\end{equation}
denote the constant from \Cref{two-sided-metric-space-Carleson}.
Define
\begin{equation}
    \epsilon' := \frac {\epsilon} {4 C_\epsilon} \,,
\end{equation}
where
\begin{equation*}
    C_\epsilon = \left(\frac{8}{\pi\epsilon}\right)^\frac{1}{2} C_{4,2} + \pi \,.
\end{equation*}
Since $f$ is continuous and periodic, $f$ is uniformly continuous.
Thus, there is a $0<\delta<\pi$
such that for all $x,x' \in \R$ with $|x-x'|\le \delta$
we have
\begin{equation}\label{uniconbound}
|f(x)-f(x')|\le \epsilon' \, .
\end{equation}
Define
\begin{equation}\label{def-fzero}
f_0:=f \ast \phi_\delta,
\end{equation}
where $\phi_\delta$ is a nonnegative smooth bump function with $\supp (\phi_\delta) \subset (-\delta, \delta)$ and $\int _\R \phi_\delta (x) \, dx = 1$.

\begin{lemma}[smooth approximation]
\label{smooth-approximation}
\leanok
\lean{close_smooth_approx_periodic}
    The function $f_0$ is $2\pi$-periodic.
    The function $f_0$ is smooth (and therefore measurable).
    The function $f_0$ satisfies for all $x\in \R$:
    \begin{equation}\label{eq-ffzero}
    |f(x)-f_0(x)|\le \epsilon' \, ,
    \end{equation}
\end{lemma}

\begin{proof}
    \leanok
    Periodicity follows directly from the definitions. The other properties are part of the Lean library.
\end{proof}

We prove in \Cref{10smooth}:
\begin{lemma}[convergence for smooth]
\label{convergence-for-smooth}
\uses{convergence-for-twice-contdiff}
\leanok
\lean{fourierConv_ofTwiceDifferentiable}
    There exists some $N_0 \in \N$ such that for all $N>N_0$ and $x\in [0,2\pi]$ we have
    \begin{equation}
        |S_N f_0 (x)- f_0(x)|\le \frac \epsilon 4\, .
    \end{equation}
\end{lemma}

We prove in \Cref{10difference}:
\begin{lemma}[control approximation effect]
\label{control-approximation-effect}
\uses{partial-Fourier-sums-of-small}
\leanok
\lean{control_approximation_effect}
    There is a set $E \subset \R$ with Lebesgue measure
    $|E|\le \epsilon$ such that for all
    \begin{equation}
        x\in [0,2\pi)\setminus E
    \end{equation}
   we have
   \begin{equation}
    \label{eq-max-partial-sum-diff}
    \sup_{N\ge 0} |S_Nf(x)-S_Nf_0(x)| \le \frac \epsilon 4\,.
    \end{equation}
\end{lemma}

We are now ready to prove \Cref{classical-carleson}. We first prove a version with explicit exceptional sets.

\begin{theorem}[classical Carleson with exceptional sets]
    \label{exceptional-set-carleson}
    \leanok
    \lean{exceptional_set_carleson}
    \uses{smooth-approximation, convergence-for-smooth,
    control-approximation-effect}
    Let $f$ be a $2\pi$-periodic complex-valued continuous function on $\mathbb{R}$. For all $\epsilon>0$, there exists a Borel set $E\subset [0,2\pi]$ with Lebesgue measure $|E|\le \epsilon$ and a positive integer $N_0$ such that for all $x\in [0,2\pi]\setminus E$ and all integers $N>N_0$, we have
    \begin{equation}\label{aeconv}
    |f(x)-S_N f(x)|\le \epsilon.
    \end{equation}
\end{theorem}
\begin{proof}
\leanok
Let $N_0$ be as in \Cref{convergence-for-smooth}.
For every
\begin{equation}
x\in [0, 2\pi) \setminus E\, ,
\end{equation}
and every $N>N_0$ we have by the triangle inequality
\begin{equation*}
    |f(x)-S_Nf(x)|
    \end{equation*}
    \begin{equation}\label{epsilonthird}
    \le |f(x)-f_0(x)|+ |f_0(x)-S_Nf_0(x)|+|S_Nf_0(x)-S_N f(x)|\, .
\end{equation}
Using \Cref{smooth-approximation,convergence-for-smooth,control-approximation-effect}, we estimate \eqref{epsilonthird} by
\begin{equation}
    \le \epsilon' +\frac \epsilon 4 +\frac \epsilon 4\le \epsilon\, .
\end{equation}
This shows \eqref{aeconv} for the given $E$ and $N_0$.
\end{proof}

Now we turn to the proof of the statement of Carleson theorem given in the introduction.
\begin{proof}[Proof of \Cref{classical-carleson}]
\proves{classical-carleson}
\leanok
By applying \Cref{exceptional-set-carleson} with a sequence of $\epsilon_n:= 2^{-n}\delta$ for $n\ge 1$ and taking the union
of corresponding exceptional sets $E_n$, we see that outside a set of measure $\delta$, the partial Fourier
sums converge pointwise for $N\to \infty$. Applying this with a sequence of $\delta$ shrinking to zero and
taking the intersection of the corresponding exceptional sets, which has measure zero, we see that the Fourier
series converges outside a set of measure zero.
\end{proof}

Let $\kappa:\R\to \C$ be the function defined by
$\kappa(0)=0$ and for $0<|x|<1$
\begin{equation}\label{eq-hilker}
\kappa(x)=\frac { 1-|x|}{1-e^{ix}}\,
\end{equation}
and for $|x|\ge 1$,
\begin{equation}\label{eq-hilker1}
\kappa(x)=0\, .
\end{equation}
Note that this function is continuous at every point $x$ with $|x|>0$.

The proof of \Cref{control-approximation-effect} will
use the following \Cref{real-Carleson}, which itself is proven
 in \Cref{10carleson} as an application of
 \Cref{metric-space-Carleson}.
\begin{lemma}[real Carleson]\label{real-Carleson}
\uses{two-sided-metric-space-Carleson,real-line-doubling,real-line-metric,real-line-measure,oscillation-control,frequency-monotone,frequency-ball-doubling,frequency-ball-growth,integer-ball-cover,real-van-der-Corput,Hilbert-strong-2-2,Hilbert-kernel-regularity,Hilbert-kernel-bound}
\leanok
\lean{rcarleson}
    Let $F,G$ be Borel subsets of $\R$ with finite measure. Let $f$ be a bounded measurable function on $\R$ with $|f|\le \mathbf{1}_F$. Then
\begin{equation}
    \left|\int _G Tf(x) \, dx\right| \le C_{4,2} |F|^{\frac 12} |G|^{\frac 12} \, ,
\end{equation}
where
\begin{equation}
    \label{define-T-carleson}
    T f(x)=\sup_{n\in \mathbb{Z}}
    \sup_{r>0}\left|\int_{r<|x-y|<1} f(y)\kappa(x-y) e^{iny}\, dy\right|\, .
\end{equation}
\end{lemma}

One of the main assumptions of \Cref{two-sided-metric-space-Carleson}, concerning the operator $T_r$ defined in \eqref{def-T-r}, is verified by the following lemma, which is proved in \Cref{10hilbert}.
\begin{lemma}[Hilbert strong 2 2]
    \label{Hilbert-strong-2-2}
    \uses{modulated-averaged-projection,integrable-bump-convolution,Dirichlet-approximation}
    \leanok
    \lean{Hilbert_strong_2_2}
    Let $0<r$. Let $f$ be a bounded, measurable function on $\mathbb{R}$. Then
    \begin{equation}
        \label{eq-Hr-L2-bound}
        \|H_rf\|_{2}\leq 2^{9} \|f\|_2,
    \end{equation}
    where
    \begin{equation}
        \label{def-H-r}
        H_r f(x) := T_r f(x) = \int_{r\le |x-y|} \kappa(x-y) f(y) \, dy
    \end{equation}
\end{lemma}

The next lemma will be used to verify that the collection $\Mf$ of modulation functions in our application of \Cref{metric-space-Carleson} satisfies the condition \eqref{eq-vdc-cond}.
It is proved in \Cref{10vandercorput}.

\begin{lemma}[van der Corput]
\label{van-der-Corput}
\leanok
\lean{van_der_Corput}
    Let $\alpha\le\beta$ be real numbers. Let $g:\R\to \C$ be a measurable function and assume
    \begin{equation}
        \|g\|_{Lip(\alpha,\beta)}:=\sup_{\alpha\le x\le \beta}|g(x)|+\frac{|\beta-\alpha|}{2}
        \sup_{\alpha\le x<y\le \beta} \frac {|g(y)-g(x)|}{|y-x|}<\infty\, .
    \end{equation}
    Then for any $\alpha \le \beta$ and $n\in\Z$ we have
    \begin{equation}
        \int _{\alpha}^{\beta} g(x) e^{inx}\, dx\le 2\pi |\beta-\alpha|\|g\|_{Lip(\alpha,\beta)}(1+|n||\beta-\alpha|)^{-1}\, .
    \end{equation}

\end{lemma}

We close this section with five lemmas that are used
across the following subsections.

\begin{lemma}[Dirichlet kernel]
\label{Dirichlet-kernel}
\leanok
\lean{dirichletKernel_eq, partialFourierSum_eq_conv_dirichletKernel}
We have for every $2\pi$-periodic bounded measurable $f$ and every $N\ge 0$
\begin{equation}
    S_Nf(x)=\frac 1{2\pi}\int_{0}^{2\pi}f(y) K_N(x-y)\, dy
\end{equation}
where $K_N$ is the $2\pi$-periodic continuous function of
$\R$ given by
\begin{equation}\label{eqksumexp}
\sum_{n=-N}^N e^{in x'}\, .
\end{equation}
We have for $e^{ix'}\neq 1$ that
\begin{equation}\label{eqksumhil}
    K_N(x')=\frac{e^{iNx'}}{1-e^{-ix'}}
      +\frac {e^{-iNx'}}{1-e^{ix'}} \, .
\end{equation}

\end{lemma}

\begin{proof}
\leanok
We have by definitions and interchanging sum and integral
   \begin{equation*}
        S_Nf(x)=\sum_{n=-N}^N \widehat{f}_n e^{inx}
    \end{equation*}
       \begin{equation*}
    =\sum_{n=-N}^N \frac 1{2\pi}\int_{0}^{2\pi}
    f(x) e^{in(x-y)}\, dy
    \end{equation*}
 \begin{equation}\label{eq-expsum}
     =\frac 1{2\pi}\int_{0}^{2\pi}
    f(y) \sum_{n=-N}^N e^{in(x-y)}\, dy\, .
 \end{equation}
This proves the first statement of the lemma.
By a telescoping sum, we have for every $x'\in \R$
\begin{equation}
    \left( e^{\frac 12 ix'}-e^{-\frac 12 ix'}\right) \sum_{n=-N}^N e^{inx'}= e^{(N+\frac 12) ix'}-e^{-(N+\frac 12) ix'}\, .
\end{equation}
If $e^{ix'}\neq 1$, the first factor on the left-hand side is not $0$ and we may divide by this factor to obtain
\begin{equation}
      \sum_{n=-N}^N e^{inx'}= \frac{e^{i(N+\frac 1 2)x'}}{e^{\frac 12 ix'}-e^{-\frac 12ix'}}
      -\frac{e^{-i(N+\frac 1 2)x'}}{e^{\frac 12 ix'}-e^{-\frac 12ix'}}
       =\frac{e^{iNx'}}{1-e^{-ix'}}
      +\frac {e^{-iNx'}}{1-e^{ix'}}\, .
\end{equation}
This proves the second part of the lemma.
\end{proof}

\begin{lemma}[lower secant bound]
    \label{lower-secant-bound}
    \leanok
    \lean{lower_secant_bound'}
    Let $\eta>0$ and $-2\pi +\eta \le x\le 2\pi-\eta$ with $|x|\ge \eta$. Then
    \begin{equation}
        |1-e^{ix}|\ge \frac{2}{\pi} \eta
    \end{equation}
\end{lemma}
\begin{proof}
    \leanok
    We have
    $$
        |1 - e^{ix}| = \sqrt{(1 - \cos(x))^2 + \sin^2(x)} \ge |\sin(x)|\,.
    $$
    If $0 \le x \le \frac{\pi}{2}$, then we have from concavity of $\sin$ on $[0, \pi]$ and $\sin(0) = 0$ and $\sin(\frac{\pi}{2}) = 1$
    $$
        |\sin(x)| \ge \frac{2}{\pi} x \ge \frac{2}{\pi} \eta\,.
    $$
    When $x\in \frac{m\pi}{2} + [0, \frac{\pi}{2}]$ for $m \in \{-4, -3, -2, -1, 1, 2, 3\}$ one can argue similarly.
\end{proof}

The following lemma will be proved in \Cref{10projection}.

\begin{lemma}[spectral projection bound]
\label{spectral-projection-bound}
\leanok
\lean{spectral_projection_bound}
    Let $f$ be a bounded $2\pi$-periodic measurable function. Then, for all $N\ge 0$
   \begin{equation}\label{snbound}
   \|S_Nf\|_{L^2[0, 2\pi]} \le \|f\|_{L^2[0, 2\pi]}.
   \end{equation}
\end{lemma}

\begin{lemma}[Hilbert kernel bound]
\label{Hilbert-kernel-bound}
\leanok
\lean{Hilbert_kernel_bound}
\uses{lower-secant-bound}
    For $x,y\in \R$ with $x\neq y$ we have
    \begin{equation}\label{eqcarl30}
        |\kappa(x-y)|\le 2^2(2|x-y|)^{-1}\, .
    \end{equation}
\end{lemma}
\begin{proof}
\leanok
    Fix $x\neq y$. If $\kappa(x-y)$ is zero, then \eqref{eqcarl30} is evident. Assume $\kappa(x-y)$ is not zero, then $0<|x-y|<1$.
    We have
\begin{equation}\label{eqcarl31}
|\kappa(x-y)|=\left|\frac {1-|x-y|}{1-e^{i(x-y)}}\right|\, .
\end{equation}
We estimate
with \Cref{lower-secant-bound}
\begin{equation}\label{eqcarl311}
|\kappa(x-y)|\le \frac {1}{|1-e^{i(x-y)}|}\le \frac 2{|x-y|}\, .
\end{equation}
This proves \eqref{eqcarl30} in the given case and completes the proof of the lemma.
\end{proof}

\begin{lemma}[Hilbert kernel regularity]
\label{Hilbert-kernel-regularity}
\leanok
\lean{Hilbert_kernel_regularity}
\uses{lower-secant-bound}
    For $x,y,y'\in \R$ with $x\neq y,y'$ and
    \begin{equation}
        \label{eq-close-hoelder}
        2|y-y'|\le |x-y|\, ,
    \end{equation}
    we have
    \begin{equation}\label{eqcarl301}
        |\kappa(x-y) - \kappa(x-y')|\le 2^{8}\frac{1}{|x-y|} \frac{|y-y'|}{|x-y|}\, .
    \end{equation}
\end{lemma}

\begin{proof}
\leanok
    Upon replacing $y$ by $y-x$ and $y'$ by $y'-x$ on the left-hand side of \eqref{eq-close-hoelder}, we can assume that $x = 0$. Then the assumption \eqref{eq-close-hoelder} implies that $y$ and $y'$ have the same sign. Since $\kappa(y) = \bar \kappa(-y)$ we can assume that they are both positive. Then it follows from \eqref{eq-close-hoelder} that
    $$
        \frac{y}{2} \le y' \,.
    $$
    We distinguish four cases. If $y, y' \le 1$, then we have
    $$
        |\kappa(-y) - \kappa(-y')| = \left| \frac{1 - y}{1- e^{-iy}} - \frac{1 - y'}{1- e^{-iy'}}\right|
    $$
    and by the fundamental theorem of calculus
    $$
        = \left| \int_{y'}^{y} \frac{-1 + e^{-it} + i(1-t)e^{it}}{(1 - e^{-it})^2} \,dt \right|\,.
    $$
    Using $y' \ge \frac{y}{2}$ and \Cref{lower-secant-bound}, we bound this by
    $$
        \le |y - y'| \sup_{\frac{y}{2} \le t \le 1} \frac{3}{|1 - e^{-it}|^2} \le 3 |y-y'| (2 \frac{2}{y})^2 \le 2^{6} \frac{|y-y'|}{|y|^2}\,.
    $$
    If $y \le 1$ and $y' > 1$, then $\kappa(-y') = 0$ and we have from the first case
    $$
        |\kappa(-y) - \kappa(-y')| = |\kappa(-y) - \kappa(-1)| \le 2^{6} \frac{|y-1|}{|y|^2} \le 2^{6} \frac{|y-y'|}{|y|^2}\,.
    $$
    Similarly, if $y > 1$ and $y' \le 1$, then $\kappa(-y) = 0$ and we have from the first case
    $$
        |\kappa(-y) - \kappa(-y')| = |\kappa(-y') - \kappa(-1)| \le 2^{6} \frac{|y'-1|}{|y'|^2} \le 2^{6} \frac{|y-y'|}{|y'|^2}\,.
    $$
    Using again $y' \ge \frac{y}{2}$, we bound this by
    $$
        \le 2^{6} \frac{|y-y'|}{|y / 2|^2} = 2^{8} \frac{|y-y'|}{|y|^2}
    $$
    Finally, if $y, y' > 1$ then
    $$
        |\kappa(-y) - \kappa(-y')| = 0 \le 2^{8} \frac{|y-y'|}{|y|^2}\,.
    $$
\end{proof}

\subsection{Smooth functions.}
\label{10smooth}
\begin{lemma}
\label{fourier-coeff-derivative}
\leanok
\lean{fourierCoeffOn_of_hasDerivAt}
Let $f:\R \to \C$ be $2\pi$-periodic and continuously differentiable, and let $n \in \Z \setminus \{0\}$. Then
\begin{equation}
    \widehat{f}_n = \frac{1}{i n} \widehat{f'}_n.
\end{equation}
\end{lemma}
\begin{proof}
\leanok
This is part of the Lean library.
\end{proof}

\begin{lemma}
\label{convergence-of-coeffs-summable}
\leanok
\lean{hasSum_fourier_series_of_summable}
Let $f:\R \to \C$ such that
\begin{equation}
    \sum_{n\in \Z} |\widehat{f}_n| < \infty.
\end{equation}
Then
\begin{equation}
    \sup_{x\in [0,2\pi]} |f(x) - S_Nf(x)| \rightarrow 0
\end{equation}
as $N \rightarrow \infty$.
\end{lemma}

\begin{proof}
\leanok
    This is part of the Lean library.
\end{proof}

\begin{lemma}
\label{convergence-for-twice-contdiff}
\uses{fourier-coeff-derivative,convergence-of-coeffs-summable}
\leanok
\lean{fourierConv_ofTwiceDifferentiable}
    Let $f:\R \to \C$ be $2\pi$-periodic and twice continuously differentiable. Then
    \begin{equation}
        \sup_{x\in [0,2\pi]} |f(x) - S_Nf(x)| \rightarrow 0
    \end{equation}
    as $N \rightarrow \infty$.
\end{lemma}
\begin{proof}
\leanok
By \Cref{convergence-of-coeffs-summable}, it suffices to show that the Fourier coefficients $\widehat{f}_n$ are summable.
Applying \Cref{fourier-coeff-derivative} twice and using the fact that $f''$ is continuous and thus bounded on $[0,2\pi]$ , we compute
\begin{equation*}
    \sum_{n\in \Z} |\widehat{f}_n| = |\widehat{f}_0| + \sum_{n\in \Z \setminus \{0\}} \frac {1}{n^2} |\widehat{f''}_n|
    \le |\widehat{f}_0| + \left(\sup_{x\in [0,2\pi]} |f(x)| \right) \cdot \sum_{n\in \Z \setminus \{0\}} \frac {1}{n^2}
    < \infty.
\end{equation*}
\end{proof}

\begin{proof}[Proof of \Cref{convergence-for-smooth}]
\leanok
\proves{convergence-for-smooth}
    \Cref{convergence-for-smooth} now follows directly from \Cref{convergence-for-twice-contdiff}.
\end{proof}

\subsection{The truncated Hilbert transform}
\label{10hilbert}

Let $M_n$ be the modulation operator
acting on measurable $2\pi$-periodic functions
defined by
\begin{equation}
    M_ng(x)=g(x) e^{inx}\, .
\end{equation}
Define the approximate Hilbert transform by
\begin{equation}
    L_N g=\frac 1N\sum_{n=0}^{N-1}
       M_{n+N} S_{N+n}M_{-N-n}g\, .
\end{equation}

\begin{lemma}[modulated averaged projection]
\label{modulated-averaged-projection}
\leanok
\lean{modulated_averaged_projection}
\uses{spectral-projection-bound}
We have for every bounded measurable $2\pi$-periodic function $g$
\begin{equation}\label{lnbound}
    \|L_Ng\|_{L^2[0, 2\pi]}\le \|g\|_{L^2[0, 2\pi]}\,.
\end{equation}
\end{lemma}
\begin{proof}
    \leanok
    We have
    \begin{equation}\label{mnbound}
        \|M_ng\|_{L^2[0, 2\pi]}^2=\int_{0}^{2\pi} |e^{inx}g(x)|^2\, dx
        =\int_{0}^{2\pi} |g(x)|^2\, dx=\|g\|_{L^2[0, 2\pi]}^2\, .
    \end{equation}
     We have by the triangle inequality, the square root of the identity in \eqref{mnbound}, and \Cref{spectral-projection-bound}
    \begin{equation*}
        \|L_ng\|_{L^2[0, 2\pi]}=\|\frac 1N\sum_{n=0}^{N-1}
       M_{n+N} S_{N+n}M_{-N-n}g\|_{L^2[0, 2\pi]}
    \end{equation*}
    \begin{equation*}
        \le \frac 1N\sum_{n=0}^{N-1} \|
       M_{n+N} S_{N+n}M_{-N-n}g\|_{L^2[0, 2\pi]}
         = \frac 1N\sum_{n=0}^{N-1} \|
    S_{N+n}M_{-N-n}g\|_{L^2[0, 2\pi]}
    \end{equation*}
     \begin{equation}
     \le \frac 1N\sum_{n=0}^{N-1} \|
 M_{-N-n}g\|_{L^2[0, 2\pi]} = \frac 1N\sum_{n=0}^{N-1} \|
g\|_{L^2[0, 2\pi]} =\|g\|_{L^2[0, 2\pi]}\, .
    \end{equation}
This proves \eqref{lnbound} and completes the proof of the lemma.
\end{proof}

\begin{lemma}[periodic domain shift]\label{periodic-domain-shift}
\leanok
\lean{Function.Periodic.intervalIntegral_add_eq, intervalIntegral.integral_comp_sub_right}
Let $f$ be a bounded $2\pi$-periodic function. We have for any
$0 \le x\le 2\pi$ that
\begin{equation}
 \int_0^{2\pi} f(y)\, dy= \int_{-x}^{2\pi -x} f(y)\, dy
 =\int_{-\pi}^{\pi} f(y-x)\, dy\,.
\end{equation}
\end{lemma}
\begin{proof}
\leanok
    We have by periodicity and change of variables
    \begin{equation}\label{eqhil9}
 \int_{-x}^{0} f(y)\, dy=\int_{-x}^{0} f(y+2\pi)\, dy= \int_{2\pi -x}^{2\pi} f(y)\, dy\, .
\end{equation}
We then have by breaking up the domain of integration
and using \eqref{eqhil9}
\begin{equation*}
 \int_0^{2\pi} f(y)\, dy= \int_0^{2\pi -x} f(y)\, dy+
 \int_{2\pi -x}^{2\pi} f(y)\, dy
 \end{equation*}
\begin{equation}
= \int_0^{2\pi -x} f(y)\, dy+
 \int_{ -x}^{0} f(y)\, dy
 = \int_{-x}^{2\pi-x} f(y)\, dy\, .
 \end{equation}
This proves the first identity of the lemma. The second identity follows by substitution of $y$ by $y-x$.
\end{proof}

\begin{lemma}[Young convolution]\label{Young-convolution}
\leanok
\lean{young_convolution}
\uses{periodic-domain-shift}
    Let $f$ and $g$ be two bounded non-negative measurable $2\pi$-periodic functions on $\R$. Then
    \begin{equation}\label{eqyoung}
        \left(\int_0^{2\pi} \left(\int_0^{2\pi}
        f(y)g(x-y)\, dy\right)^2\, dx\right)^{\frac 12}\le \|f\|_{L^2[0, 2\pi]} \|g\|_{L^1[0, 2\pi]}\, .
    \end{equation}
    \end{lemma}
\begin{proof}
\leanok
    Using Fubini and \Cref{periodic-domain-shift}, we observe
\begin{equation*}
  \int_0^{2\pi}\int_0^{2\pi}f(y)^2g(x-y)\, dy
    \, dx=\int_0^{2\pi}f(y)^2\int_0^{2\pi}g(x-y)\, dx
    \, dy
\end{equation*}
\begin{equation}\label{eqhil4}
=\int_0^{2\pi}f(y)^2\int_0^{2\pi}g(x) \, dx
     dy
=\|f\|_{L^2[0, 2\pi]}^2\|g\|_{L^1[0, 2\pi]}\, .
\end{equation}
   Let $h$ be the nonnegative square root of $g$, then
   $h$ is bounded and $2\pi$-periodic with $h^2=g$.
   We estimate the square of the left-hand side of
   \eqref{eqyoung} with Cauchy-Schwarz and then with
   \eqref{eqhil4} by
       \begin{equation*}
         \int_0^{2\pi} (\int_0^{2\pi}f(y)h(x-y)h(x-y)\, dy)^2\, dx
   \end{equation*}
\begin{equation*}
    \le \int_0^{2\pi}\left(\int_0^{2\pi}f(y)^2g(x-y)\, dy\right)
    \left(\int_0^{2\pi}g(x-y)\, dy\right)\, dx
\end{equation*}
\begin{equation*}
    = \|f\|_{L^2[0, 2\pi]}^2\|g\|_{L^1[0, 2\pi]}^2\, .
\end{equation*}
Taking square roots, this proves the lemma.
\end{proof}

For $0<r<1$, Define the kernel $k_r$ to be the $2\pi$-periodic function
\begin{equation}
    k_r(x):=\min \left(r^{-1}, 1+\frac r{|1-e^{ix}|^2}\right)\, ,
\end{equation}
where the minimum is understood to be $r^{-1}$ in case $1=e^{ix}$.
\begin{lemma}[integrable bump convolution]
\label{integrable-bump-convolution}
\leanok
\lean{integrable_bump_convolution}
\uses{Young-convolution}
Let $g,f$ be bounded measurable $2\pi$-periodic functions. Let $0<r<\pi$.
Assume we have for all $x$
\begin{equation}\label{ebump1}
    |g(x)|\le k_r(x)\, .
\end{equation}
Let
\begin{equation}
   h(x)= \int_0^{2\pi} f(y)g(x-y)\, dy \, .
\end{equation}
Then
\begin{equation}
   \|h\|_{L^2[0, 2\pi]}\le 17\|f\|_{L^2[-\pi, \pi]} \, .
\end{equation}

\end{lemma}

\begin{proof}
\leanok
From monotonicity of the integral and \eqref{ebump1},
\begin{equation}
    \|g\|_{L^1[0, 2\pi]} \le \int_0^{2\pi}k_r(x)\, dx\,.
\end{equation}
Using the symmetry
$k_r(x)=k_r(-x)$, the assumption, and \Cref{lower-secant-bound}, the last display
is equal to
\begin{equation*}
    = 2 \int_0^\pi \min\left(\frac 1r, 1+\frac r{|1-e^{ix}|^2}\right)\, dx
\end{equation*}
\begin{equation*}
    \le 2\int_0^{r} \frac 1r \, dx+2\int_r^{\pi}1+\frac {64r}{x^2}\, dx
\end{equation*}
\begin{equation}
    \le 2+2\pi + 2\left(\frac {4r}r-\frac {4r}{\pi}\right)
    \le 17\, .
\end{equation}
    Together with \Cref{Young-convolution}, this proves the lemma.
\end{proof}

\begin{lemma}[Dirichlet approximation]\label{Dirichlet-approximation}
\uses{Dirichlet-kernel,lower-secant-bound}
\leanok
\lean{continuous_dirichletApprox, periodic_dirichletApprox, approxHilbertTransform_eq_dirichletApprox, dist_dirichletApprox_le}
Let $0<r<1$. Let $N$ be the smallest
integer larger than $\frac 1r$.
There is a $2\pi$-periodic continuous function
 ${L'}$ on $\R$ that satisfies for all $0\le x\le 2\pi$
and all $2\pi$-periodic bounded measurable functions $f$ on $\R$
\begin{equation}\label{lthroughlprime}
    L_Nf(x)=\frac 1{2\pi}\int_{0}^{2\pi}f(y) {L'}(x-y)\, dy\,.
\end{equation}
Moreover, for all $-\pi \le x \le \pi$,
\begin{equation}\label{eqdifflhil}
    \left|L'(x)-\mathbf{1}_{\{y:\, r<|y|<1\}} \kappa(x)\right|\le 12 k_r(x)\, .
\end{equation}
\end{lemma}

\begin{proof}
\leanok
We have by definition and \Cref{Dirichlet-kernel}
\begin{equation}
    L_Ng(x)=
    \frac 1N\sum_{n=0}^{N-1}
       \int_0^{2\pi} e^{-i(N+n)x} K_{N+n}(x-y) e^{i(N+n)y}g(y)
\, dy \, .\end{equation}
This is of the form \eqref{lthroughlprime} with
the continuous function
\begin{equation}
    {L'}(x)= \frac 1N\sum_{n=0}^{N-1}
      K_{N+n}(x) e^{-i(N+n)x}\, .
\end{equation}
With \eqref{eqksumexp} of \Cref{Dirichlet-kernel}
we have $|K_N(x)|\le 2N+1$ for every $x$ and thus
\begin{equation}\label{eqhil13}
    |{L'}(x)|\le \frac 1N\sum_{n=0}^{N-1}
      (2N+2n+1) \le 4N\le 2^3 r^{-1}\, .
\end{equation}
Therefore, for $|x|\in [0, r)\cup (1, \pi]$, we have
\begin{equation}
    \label{eqdiffzero}
    \left|L'(x)-\mathbf{1}_{\{y:\, r<|y|<1\}}(x)\kappa(x)\right|=|L'(x)|\leq 2^{3} r^{-1}.
\end{equation}
This proves \eqref{eqdifflhil} for $|x|\in [0, r)$ since $k_r(x)=r^{-1}$ in this case.

For $e^{ix}\neq 1$
one may use the expression
\eqref{eqksumhil} for $K_N$
in \Cref{Dirichlet-kernel} to obtain
\begin{equation*}
    {L'}(x)= \frac 1N\sum_{n=0}^{N-1}
     \left(\frac{e^{i(N+n)x}}{1-e^{-ix}}
      +\frac {e^{-i(N+n)x}}{1-e^{ix}}\right) e^{i(N+n)x}
\end{equation*}
\begin{equation*}
    = \frac 1N\sum_{n=0}^{N-1}
    \left(\frac{e^{i2(N+n)x}}{1-e^{-ix}}
      +\frac {1}{1-e^{ix}}\right)
\end{equation*}
\begin{equation}\label{eqhil3}
    = \frac{1}{1-e^{ix}} +
     \frac 1N \frac {e^{i2Nx}}{1-e^{-ix}}
     \sum_{n=0}^{N-1}
    {e^{i2nx}}
\end{equation}
and thus
\begin{equation}
\label{eq-L'L''}
  {L'}(x) -\mathbf{1}_{\{y:\, r<|y|<1\}}\kappa(x)=L''(x)+ \frac{1-\mathbf{1}_{{\{y:\, r<|y|<1\}}}(x)(1-|x|)}{1-e^{ix}},
\end{equation}
where
$$L''(x):=\frac 1N \frac {e^{i2Nx}}{1-e^{-ix}}
     \sum_{n=0}^{N-1}
    {e^{i2nx}}.$$
For $x\in [-\pi, r]\cup [r, \pi]$, we have using \Cref{lower-secant-bound} that
\begin{equation*}
   \left|\frac{1-\mathbf{1}_{{\{y:\, r<|y|<1\}}}(x)(1-|x|)}{1-e^{ix}} \right|=\left|\frac{\min(|x|, 1)}{1-e^{ix}} \right|\leq \frac{2\min(|x|, 1)}{|x|}
\end{equation*}
\begin{equation}
 \label{eq-diffzero2}
    \leq 2\cdot 1\leq 2 k_r(x).
\end{equation}
Next, we need to estimate $L''(x)$. If the real part of
$e^{ix}$ is negative, we have
\begin{equation}
  1\le |1-e^{-ix}|\le 2\, .
\end{equation}
and hence
\begin{equation}\label{eqhil12}
    |L''(x)|\le
     \frac 1N
     \sum_{n=0}^{N-1}
    1=1\le 1+\frac r{|1-e^{ix}|^2}\, .
\end{equation}
If the real part of $e^{ix}$ is positive and in particular while still $e^{ix}\neq \pm 1$, then we have by telescoping
\begin{equation}
 (1-e^{2ix})
     \sum_{n=0}^{N-1}
    {e^{i2nx}}=1-e^{i2Nx}\, .
\end{equation}
As $e^{2ix}\neq 1$, we may divide by $1-e^{2ix}$ and insert this into
\eqref{eqhil3} to obtain
\begin{equation}
 L''(x)=
           \frac 1N \frac {e^{i2Nx}}{1-e^{-ix}}
     \frac{1-e^{i2Nx}}{1-e^{2ix}}\, .
\end{equation}
Hence, using the nonnegativity of the real part of $e^{ix}$
\begin{equation*}
    |L''(x)|
 \le \frac 2 N \frac {1}{|1-e^{ix}|}
     \frac{1}{|1-e^{2ix}|}
 \end{equation*}
\begin{equation}\label{eqhil11}
    = \frac 2 N \frac {1}{|1-e^{ix}|^2}
     \frac{1}{|1+e^{ix}|}\le
 \frac {2r}{|1-e^{ix}|^2}\le 2 \left(1+\frac {r}{|1-e^{ix}|^2}\right)
\end{equation}
Moreover, for $|x| \in [r, \pi]$, one checks easily that
\begin{equation*}
  1+\frac {r}{|1-e^{ix}|^2} \leq 5 k_r(x).
\end{equation*}
Therefore,~\eqref{eqhil12} and~\eqref{eqhil11} show that, in this range of $x$,
\begin{equation*}
  |L''(x)| \le 10 k_r(x).
\end{equation*}
Together with~\eqref{eq-diffzero2}, this prove \eqref{eqdifflhil} for $|x| \in [r, \pi]$.
As the range $|x| \in [0,r)$ is covered by~\eqref{eqdiffzero}, this completes the proof of the lemma.
\end{proof}

We now prove \Cref{Hilbert-strong-2-2}.

\begin{proof}[Proof of \Cref{Hilbert-strong-2-2}]
    \proves{Hilbert-strong-2-2}
    \leanok
    We first show that if $f$ is supported in $[1, 4]$, then
    \begin{equation}
        \label{eq-Hr-short-support}
        \|H_r f\|_{L^2[2, 3]} \le 2^{8} \|f\|_{L^2(\R)}\,.
    \end{equation}
    Let $\tilde{f}$ be the $2\pi$-periodic extension of $f$ to $\mathbb{R}$. Let $N$ be the smallest
    integer larger than $\frac 1r$. Then, \Cref{eqdifflhil} shows that the kernels of $H_r$ and $2\pi L_N$ differ by at
    most $12k_r$ on $[-\pi, \pi]$. Consider $x \in [2,3]$. When computing $H_r f(x)$ and $2\pi L_N f(x)$, the kernels are
    computed at points of the form $x-y$ with $f(y) \ne 0$, i.e., $y \in [1,4]$. As $x\in [2,3]$, the difference
    $x-y$ is bounded in absolute value by $2$, and therefore belongs to $[-\pi, \pi]$, where the above bound holds.

    It follows that, for $x \in [2,3]$,
    \begin{equation*}
        |H_r \tilde{f}(x)|\leq 2\pi |L_N \tilde{f}(x)|+12\left|\int_{0}^{2\pi}k_r(x-y)\tilde{f}(y)\, dy\right|.
    \end{equation*}
    Taking $L^2$ norm and using its sub-additivity, we get
    \begin{equation*}
       \|H_r \tilde{f}\|_{L^2[2, 3]}\leq 2\pi \|L_N \tilde{f}\|_{L^2[0, 2\pi]}\, + 12\left(\int_{0}^{2\pi} \left|\int_{0}^{2\pi}k_r(x-y)\tilde{f}(y)\, dy\right|^2\, dx\right)^{\frac{1}{2}}.
    \end{equation*}
    Since $\kappa$ is supported in $[-1,1]$, we have that $H_r\tilde{f}$ agrees with $H_r f$ on $[2,3]$.
    Using \Cref{modulated-averaged-projection} and \Cref{integrable-bump-convolution}, we conclude
    \begin{equation}
        \|H_r f\|_{L^2[2, 3]} \le 2\pi \|f\|_{L^2[0, 2\pi]} + 12 \cdot 17 \|f\|_{L^2[0, 2\pi]}\,,
    \end{equation}
    which gives \eqref{eq-Hr-short-support}.

    Suppose now that $f$ is supported in $[c, c+3]$ for some $c \in \R$. Then the function $g(x) = f(x-c+1)$ is supported in $[1, 4]$.
    By a change of variables in \eqref{def-H-r}, we have $H_r g(x ) = H_r f(x - c + 1)$. Thus, by \eqref{eq-Hr-short-support}
    \begin{equation}
        \label{eq-Hr-short-support-2}
        \|H_rf\|_{L^2[c+1, c+2]} = \|H_r g\|_{L^2[2,3]} \le 2^{8} \|g\|_2 = 2^8 \|f\|_2\,.
    \end{equation}

    Let now $f$ be arbitrary.
    Since $\kappa(x) = 0$ for $|x| > 1$, we have for all $x \in [c+1, c+2]$
    $$
        H_rf(x) = H_r(f \mathbf{1}_{[c, c+3]})(x)\,.
    $$
    Thus
    $$
        \int_{c+1}^{c+2} |H_r f(x)|^2 \, \mathrm{d}x = \int_{c+1}^{c+2}|H_r(f \mathbf{1}_{[c, c+3]})(x)|^2 \, \mathrm{d}x\,.
    $$
    Applying the bound \eqref{eq-Hr-short-support-2}, this is
    $$
        \le 2^{16} \int_{c}^{c+3} |f(x)|^2 \, \mathrm{d}x\,.
    $$
    Summing over all $c \in \mathbb{Z}$, we obtain
    $$
        \int_{\R} |H_rf(x)|^2 \, \mathrm{d}x \le 3 \cdot 2^{16} \int_{\R} |f(x)|^2 \, \mathrm{d}x\,.
    $$
    This completes the proof.
\end{proof}

\subsection{The proof of the van der Corput Lemma}
\label{10vandercorput}

\begin{proof}[Proof of \Cref{van-der-Corput}]
\proves{van-der-Corput}
\leanok
Let $g$ be a Lipschitz continuous function as in the lemma.
Assume first that $n=0$. Then
\begin{equation*}
    \int_\alpha^\beta g(x) \, \mathrm{d}x \le |\beta - \alpha| \sup_{\alpha\le x\le \beta}|g(x)|
    \le |\beta-\alpha|\|g\|_{Lip(\alpha,\beta)}(1+|n||\beta-\alpha|)^{-1}
\end{equation*}
Assume now $n\ne 0$. Without loss of generality, we may assume $n>0$.
We distinguish two cases. If $\beta-\alpha < \frac{\pi}{n}$, we have by the triangle inequality
\begin{equation*}
    \left|\int_\alpha^\beta g(x) e^{inx} \, \mathrm{d}x\right|
    \le |\beta -\alpha| \sup_{x \in [\alpha,\beta]} |g(x)|
    \le 2\pi |\beta-\alpha|\|g\|_{Lip(\alpha,\beta)}(1+|n||\beta-\alpha|)^{-1} \,.
\end{equation*}
We turn to the case $\frac{\pi}{n} \le \beta-\alpha$.
We have
$$
    e^{in(x + \pi/n)} = -e^{inx}\,.
$$
Using this, we write
$$
    \int_\alpha^\beta g(x) e^{inx} \, \mathrm{d}x
    = \frac{1}{2} \int_\alpha^\beta g(x) e^{inx} \, \mathrm{d}x - \frac{1}{2} \int_\alpha^\beta g(x) e^{in(x + \pi/n)}) \, \mathrm{d}x\,.
$$
We split the first integral at $\alpha + \frac{\pi}{n}$ and the second one at $\beta - \frac{\pi}{n}$, and make a change of variables in the second part of the first integral to obtain
$$
    = \frac{1}{2} \int_{\alpha}^{\alpha + \frac{\pi}{n}} g(x) e^{inx} \, \mathrm{d}x - \frac{1}{2} \int_{\beta - \frac{\pi}{n}}^{\beta} g(x) e^{in(x + \pi/n)} \, \mathrm{d}x
$$
$$
    + \frac{1}{2} \int_{\alpha + \frac{\pi}{n}}^{\beta} (g(x) - g(x - \frac{\pi}{n})) e^{inx} \, \mathrm{d}x\,.
$$
The sum of the first two terms is by the triangle inequality bounded by
$$
    \frac{\pi}{n} \sup_{x \in [\alpha,\beta]} |g(x)|\,.
$$
The third term is by the triangle inequality at most
$$
    \frac{1}{2} \int_{\alpha + \frac{\pi}{n}}^\beta |g(x) - g(x - \frac{\pi}{n})| \, \mathrm{d}x
$$
$$
    \le \frac{|\beta-\alpha|}{2} \frac{\pi}{n} \sup_{\alpha \le x < y \le \beta} \frac{|g(x) - g(y)|}{|x-y|}\,.
$$
Adding the two terms, we obtain
$$
    \left|\int_\alpha^\beta g(x) e^{-inx} \, \mathrm{d}x\right| \le \frac{\pi}{n} \|g\|_{\mathrm{Lip}(\alpha,\beta)}\,.
$$
This completes the proof of the lemma, using that with $\frac{\pi}{n} \le \beta-\alpha$,
$$
    \frac{\pi}{n} = \frac{2 \pi |\beta-\alpha|}{2n|\beta-\alpha|} \le 2 \pi |\beta-\alpha|(1 + n|\beta-\alpha|)^{-1}\,.
$$
\end{proof}

\subsection{Partial sums as orthogonal projections}
\label{10projection}

This subsection contains the proof of \Cref{spectral-projection-bound}

\begin{proof}[Proof of \Cref{spectral-projection-bound}]
\proves{spectral-projection-bound}
\leanok

The functions $e_n:x\mapsto e^{inx}$ form an orthonormal basis in $L^2[-\pi, \pi]$
(this is already in Mathlib).
Therefore we have
\begin{align*}
    \|S_Nf\|^2_{L^2[-\pi, \pi]}
    &= \|\sum_{n=-N}^N \langle\widehat{f}_n, e_n\rangle_{L^2[-\pi, \pi]}\|^2_{L^2[-\pi, \pi]} \\
    &= \sum_{n=-N}^N |\widehat{f}_n| \\
    &\le \sum_{n\in \mathbb{Z}} |\widehat{f}_n| \\
    &= \|f\|_{L^2[-\pi, \pi]}.
\end{align*}

This completes the proof of the lemma.
\end{proof}

\subsection{The error bound}
\label{10difference}

\begin{lemma}[Dirichlet kernel - Hilbert kernel relation]
    \label{Dirichlet-Hilbert}
    \leanok
    \lean{Dirichlet_Hilbert_diff}
    \uses{Dirichlet-kernel,lower-secant-bound}
    For all $N\in\Z$ and $x\in [-\pi,\pi] \setminus \{0\}$,
    \begin{equation*}
        \left|K_N(x) - (e^{-iNx}\kappa(x) + \overline{e^{-iNx}\kappa(x)})\right| \le \pi \,.
    \end{equation*}
\end{lemma}

\begin{proof}
    \leanok
    Let $N\in\Z$ and $x\in [-\pi,\pi] \setminus \{0\}$. With \Cref{Dirichlet-kernel}, we obtain
    \begin{equation*}
        K_N(x) - (e^{-iNx}\kappa(x) + \overline{e^{-iNx}\kappa(x)})
        = e^{-iNx} \frac{\min(|x|, 1) }{1 - e^{ix}} + e^{iNx} \frac{\min(|x|, 1) }{1 - e^{-ix}} \,.
    \end{equation*}
    Using \Cref{lower-secant-bound} with $\eta = \min(|x|, 1)$, we bound
    \begin{equation*}
        \left|K_N(x) - (e^{-iNx}\kappa(x) + \overline{e^{-iNx}\kappa(x)})\right|
        \le \frac{\min(|x|, 1) }{|1 - e^{ix}|} + \frac{\min(|x|, 1)}{|1 - e^{-ix}|}
        \le \frac{\pi}{2} + \frac{\pi}{2} = \pi \,.
    \end{equation*}
\end{proof}

\begin{lemma}[partial Fourier sum bound]
    \label{partial-Fourier-sum-bound}
    \uses{Dirichlet-kernel, Dirichlet-Hilbert}
    \leanok
    \lean{partialFourierSum_bound}
    Let $g:\R\to\C$ be a measurable $2\pi$-periodic function such that for some $\delta>0$ and every $x\in\R$,
    \begin{equation}
        |g(x)|\le \delta \,.
    \end{equation}
    Then for every $x\in [0,2\pi]$ and $N>0$,
    \begin{equation*}
        |S_N g(x)| \le \frac{1}{2\pi} (Tg(x) + T\bar{g}(x)) + \pi\delta.
    \end{equation*}
\end{lemma}

\begin{proof}
    \leanok
    Let $x\in [0,2\pi]$ and $N>0$. We have with \Cref{Dirichlet-kernel}
    \begin{equation*}
        |S_N g(x)| = \frac{1}{2\pi} \left| \int_0^{2\pi} g(y) K_N(x-y) \, dy\right|\,.
    \end{equation*}
    We use $2\pi$-periodicity of $g$ and $K_N$ to shift the domain of integration to obtain
    \begin{equation*}
        = \frac{1}{2\pi} \left|\int_{x-\pi}^{x+\pi} g(y) K_N(x-y) \, dy\right|\,.
    \end{equation*}
    Using the triangle inequality, we split this as
    \begin{equation}
        \label{eq-diff-integrable}
        \le \frac{1}{2\pi} \left|\int_{x-\pi}^{x+\pi} g(y) \left(K_N(x-y) - \max(|x-y|,0) K_N(x-y)\right) \, dy \right|
    \end{equation}
    \begin{equation}
        \label{eq-diff-singular}
        + \frac{1}{2\pi} \left|\int_{x-\pi}^{x+\pi} g(y) \max(|x-y|,0) K_N(x-y) \, dy\right|\,.
    \end{equation}
    Note that all integrals are well defined, since $K_N$ is by \eqref{eqksumexp} bounded by $2N+1$.
    Using that
    \begin{equation}
        \label{eq-Dirichlet-Hilbert}
        \max(|x-y|,0) K_N(x-y) = e^{-iN(x-y)}\kappa(x-y) + \overline{e^{-iN(x-y)}\kappa(x-y)} \,,
    \end{equation}
    \Cref{Dirichlet-Hilbert} and \eqref{g-small}, we bound \eqref{eq-diff-integrable} by
    \begin{equation*}
        \frac{1}{2\pi} \int_{x-\pi}^{x+\pi} |g(y)| \left|K_N(x-y) - e^{-iN(x-y)}\kappa(x-y) + \overline{e^{-iN(x-y)}\kappa(x-y)}\right|\, dy
        \le \pi\delta \,.
    \end{equation*}
    By dominated convergence and since $\kappa(x-y) = 0$ for $|x-y| > 1$, \eqref{eq-diff-singular} equals
    \begin{equation*}
        \frac{1}{2\pi} \lim_{r \to 0^+} \left| \int_{r < |x-y| < 1} g(y) \max(|x-y|,0) K_N(x-y) \, dy\right|\,.
    \end{equation*}
    We bound the limit by a supremum and rewrite using \eqref{eq-Dirichlet-Hilbert},
    \begin{equation*}
        \le \frac{1}{2\pi} \sup_{r > 0} \left| \int_{r < |x-y| < 1} g(y) \left(e^{-iN(x-y)}\kappa(x-y) + \overline{e^{-iN(x-y)}\kappa(x-y)}\right) \, dy\right|
    \end{equation*}
    Using the triangle inequality, we further bound this by
    \begin{alignat*}{3}
        \le&&&\frac{1}{2\pi} \sup_{r > 0} \left| \int_{r < |x-y| < 1} g(y) e^{-iNy} \kappa(x-y) \, dy\right| \\
        &+ &&\frac{1}{2\pi} \sup_{r > 0} \left| \int_{r < |x-y| < 1} \overline{g}(y) e^{-iNy} \kappa(x-y) \, dy\right|\,.
    \end{alignat*}
    By the definition \eqref{define-T-carleson} of $T$, this is
    \begin{equation*}
        \le \frac{1}{2\pi} (Tg(x) + T\bar{g}(x))\,.
    \end{equation*}
\end{proof}

\begin{lemma}[real Carleson operator measurable]
    \label{real-Carleson-operator-measurable}
    \uses{Hilbert-kernel-bound}
    \leanok
    \lean{carlesonOperatorReal_measurable}
    Let $f$ be a bounded measurable function on $\R$. Then $Tf$ as defined in \eqref{define-T-carleson} is measurable.
\end{lemma}
\begin{proof}
    \leanok
    Since a countable supremum of measurable functions is measurable, it suffices to show that for every $n\in\Z$,
    \begin{equation*}
        x \mapsto \sup_{r>0}\left|\int_{r<|x-y|<1} f(y)\kappa(x-y) e^{iny}\, dy\right|
    \end{equation*}
    is measurable. So let $n\in\Z$.
    Note that for each $x\in\R$, the function
    \begin{equation*}
        r \mapsto \left|\int_{r<|x-y|<1} f(y)\kappa(x-y) e^{iny}\, dy\right|
    \end{equation*}
    is continuous on $(0,\infty)$ since the integrand is locally bounded on the domain $0<|x-y|<1$ by the assumptions on $f$ and \Cref{Hilbert-kernel-bound}.
    Thus, for each $x\in\R$,
    \begin{equation*}
        \sup_{r>0}\left|\int_{r<|x-y|<1} f(y)\kappa(x-y) e^{iny}\, dy\right|
        =\sup_{r\in\Q_{>0}}\left|\int_{r<|x-y|<1} f(y)\kappa(x-y) e^{iny}\, dy\right|
    \end{equation*}
    The right hand side is again a countable supremum so it remains to show that for every $r\in\Q_{>0}$,
    \begin{equation*}
        x \mapsto \left|\int_{r<|x-y|<1} f(y)\kappa(x-y) e^{iny}\, dy\right| = \left|\int \mathbf{1}_{\{r<|x-\cdot|<1\}}(y) f(y)\kappa(x-y) e^{iny}\, dy\right|
    \end{equation*}
    is measurable, which follows from the fact that the integrand is measurable in $(x,y)$.
\end{proof}

\begin{lemma}[partial Fourier sums of small]
    \label{partial-Fourier-sums-of-small}
    \uses{real-Carleson-operator-measurable,real-Carleson,partial-Fourier-sum-bound}
    \leanok
    \lean{control_approximation_effect} 
    Let $g:\R\to\C$ be a measurable $2\pi$-periodic function such that for some $\delta>0$ and every $x\in\R$,
    \begin{equation} \label{g-small}
        |g(x)|\le \delta \,.
    \end{equation}
    Then for every $\epsilon>0$, there exists a measurable set $E\subset [0,2\pi]$ with $|E|<\epsilon$ such that for every $x\in [0,2\pi]\setminus E$ and $N>0$,
    \begin{equation}  \label{S-Ng-small}
        |S_N g(x)|\le C_\epsilon \delta,
    \end{equation}
    where
    \begin{equation} \label{C-epsilon-def}
        C_\epsilon = \left(\frac{8}{\pi\epsilon}\right)^\frac{1}{2} C_{4,2} + \pi \,.
    \end{equation}
\end{lemma}

\begin{proof}
\leanok
Define
\begin{equation*}
    E := \{x \in [0, 2\pi] \ : \ \sup_{N > 0} |S_N g (x)| > C_\epsilon \delta \} \,.
\end{equation*}
Then \eqref{S-Ng-small} clearly holds, and it remains to show that $|E| \le \epsilon$.
With \Cref{partial-Fourier-sum-bound}, we obtain
\begin{equation*}
    E \subset \{x\in [0,2\pi] \ : \ C_\epsilon \delta < \frac{1}{2\pi} (Tg(x) + T\bar{g}(x)) + \pi\delta\} \subset E_1 \cup E_2,
\end{equation*}
where
\begin{align*}
    E_1 :=& \{x\in [0,2\pi] \ : \ \pi(C_\epsilon - \pi) \delta < Tg(x)\} \\
    E_2 :=& \{x\in [0,2\pi] \ : \ \pi(C_\epsilon - \pi) \delta < T\bar{g}(x)\}.
\end{align*}
By \Cref{real-Carleson-operator-measurable}, $E_1$ and $E_2$ are measurable. Thus,
\begin{equation*}
    \pi(C_\epsilon - \pi) \delta |E_1| \le \int_{E_1} Tg(x) \, dx = \delta \int_{E_1} T(\delta^{-1} g\mathbf{1}_{[-\pi,3\pi]})(x) \, dx \,.
\end{equation*}
Applying \Cref{real-Carleson} with $F = [-\pi, 3\pi]$ and $G = E'$, it follows that this is
\begin{equation*}
    \le \delta \cdot C_{4,2} |F|^{\frac{1}{2}} |E_1|^{\frac{1}{2}} \le (4\pi)^\frac{1}{2} C_{4,2} \delta \cdot |E'|^{\frac{1}{2}}\,.
\end{equation*}
Rearranging, we obtain
\begin{equation*}
    |E_1| \le \left(\frac{(4\pi)^\frac{1}{2} C_{4,2}}{\pi(C_\epsilon - \pi)}\right)^2 = \frac{\epsilon}{2}\,.
\end{equation*}
Analogously, we get the same estimate for $|E_2|$. This completes the proof using $|E| \le |E_1| + |E_2|$.
\end{proof}

\begin{proof}[Proof of \Cref{control-approximation-effect}]
\proves{control-approximation-effect}
\leanok
    \Cref{control-approximation-effect} follows directly from \Cref{partial-Fourier-sums-of-small} with $\delta:=\epsilon'$.
\end{proof}

\subsection{Carleson on the real line}
\label{10carleson}

We prove \Cref{real-Carleson}.

Consider the standard distance function
\begin{equation}
    \rho(x,y)=|x-y|
\end{equation}
on the real line $\R$.
\begin{lemma}[real line metric]
\label{real-line-metric}
\leanok
\lean{instProperSpaceReal, locally_compact_of_proper, Real.instCompleteSpace}
The space $(\R,\rho)$ is a complete locally compact metric space.
\end{lemma}
\begin{proof}
    \leanok
    This is part of the Lean library.
\end{proof}
\begin{lemma}[real line ball]
\label{real-line-ball}
    \leanok
    \lean{Real.ball_eq_Ioo}
    For $x\in R$ and $R>0$, the ball $B(x,R)$ is the interval $(x-R,x+R)$
\end{lemma}
\begin{proof}
\leanok
Let $x'\in B(x,R)$. By definition of the ball,
$|x'-x|<R$. It follows that $x'-x<R$ and $x-x'<R$.
It follows $x'<x+R$ and $x'>x-R$. This implies
$x'\in (x-R,x+R)$.
Conversely, let $x'\in (x-R,x+R)$. Then
$x'<x+R$ and $x'>x-R$. It follows that
$x'-x<R$ and $x-x'<R$. It follows that $|x'-x|<R$,
hence $x'\in B(x,R)$.
This proves the lemma.
\end{proof}
We consider the Lebesgue measure $\mu$ on $\R$.
\begin{lemma}[real line measure]
\label{real-line-measure}
\leanok
\lean{instIsAddHaarMeasureVolume}
    The measure $\mu$ is a sigma-finite non-zero
    Radon-Borel measure on $\R$.
\end{lemma}
\begin{proof}
    \leanok
    This is part of the Lean library.
\end{proof}
\begin{lemma}[real line ball measure]
\label{real-line-ball-measure}
\uses{real-line-ball}
    \leanok
    \lean{Real.volume_ball}
    We have for every $x\in \R$ and $R>0$
    \begin{equation}
        \mu(B(x,R))=2R\, .
    \end{equation}
\end{lemma}
\begin{proof}
\leanok
We have with \Cref{real-line-ball}
\begin{equation}
    \mu(B(x,R))=\mu((x-R,x+R))=2R\, .
\end{equation}
\end{proof}

\begin{lemma}[real line doubling]
\label{real-line-doubling}
\leanok
\lean{MeasureTheory.InnerProductSpace.IsDoubling}
\uses{real-line-ball-measure}
    We have for every $x\in \R$ and $R>0$
    \begin{equation}
        \mu(B(x,2R))=2\mu(B(x,R))\, .
    \end{equation}
\end{lemma}
\begin{proof}
    \leanok
    We have with \Cref{real-line-ball-measure}
\begin{equation}
    \mu(B(x,2R)=4R=2\mu(B(x,R))\, .
\end{equation}
This proves the lemma.
\end{proof}

The preceding four lemmas show that $(\R, \rho, \mu, 4)$ is a doubling metric measure space. Indeed, we even show
that $(\R, \rho, \mu, 1)$ is a doubling metric measure space, but we may relax the estimate in \Cref{real-line-doubling} to conclude that $(\R, \rho, \mu, 4)$
is a doubling metric measure space.

For each $n\in \mathbb{Z}$ define
$\mfa_n:\R\to \R$ by
\begin{equation}
    \mfa_n(x)=nx\, .
\end{equation}

Let $\Mf$ be the collection $\{\mfa_n, n\in \mathbb{Z}\}$.
Note that for every $n\in \mathbb{Z}$ we have $\mfa_n(0)=0$.
Define
\begin{equation}
    \label{eqcarl4}
    d_{B(x,R)}(\mfa_n, \mfa_m) := 2R|n-m|\,.
\end{equation}

\begin{lemma}[frequency metric]
    \label{frequency-metric}
    \leanok
    \lean{instFunctionDistancesReal}
    For every $R > 0$ and $x \in X$, the function $d_{B(x,R)}$ is a metric on $\Mf$.
\end{lemma}

\begin{proof}
    \leanok
    This follows immediately from the fact that the standard metric on $\mathbb{Z}$ is a metric.
\end{proof}

\begin{lemma}[oscillation control]
\label{oscillation-control}
\leanok
\lean{oscillation_control}
    For every $R > 0$ and $x \in X$, and for all $n, m \in \mathbb{Z}$, we have
    \begin{equation}\label{eqcarl2}
       \sup_{y,y'\in B(x,R)}|ny-ny'-my+my'|\le 2|n-m|R\, .
    \end{equation}
\end{lemma}

\begin{proof}
\leanok
    The right hand side of \eqref{eqcarl2} equals
    $$
        \sup_{y,y'\in B(x,R)}|(n-m)(y-x)-(n-m)(y'-x)|\,.
    $$
    The lemma then follows from the triangle inequality.
\end{proof}

\begin{lemma}[frequency monotone]
\label{frequency-monotone}
\leanok
\lean{frequency_monotone}

    For any $x, x' \in X$ and $R, R' > 0$ with $B(x,R) \subset B(x, R')$, and for any $n, m \in \mathbb{Z}$
    $$
        d_{B(x,R)}(\mfa_n, \mfa_m) \le d_{B(x',R')}(\mfa_n, \mfa_m)\,.
    $$
\end{lemma}

\begin{proof}
\leanok
    This follows immediately from the definition \eqref{eqcarl4} and $R \le R'$.
\end{proof}

\begin{lemma}[frequency ball doubling]
\label{frequency-ball-doubling}
\leanok
\lean{frequency_ball_doubling}
\uses{frequency-metric}
  For any $x,x'\in \R$ and $R>0$ with
   $x\in B(x',2R)$ and any $n,m\in \mathbb{Z}$, we have
\begin{equation}\label{firstdb1}
    d_{B(x',2R)}(\mfa_n,\mfa_m)\le 2 d_{B(x,R)}(\mfa_n,\mfa_m) \, .
\end{equation}
\end{lemma}

\begin{proof}
\leanok
With \eqref{eqcarl4}, both sides of \eqref{firstdb1} are equal to $4R|n-m|$. This proves the lemma.
\end{proof}

\begin{lemma}[frequency ball growth]
\label{frequency-ball-growth}
\leanok
\lean{frequency_ball_growth}
    For any $x,x'\in \R$ and $R>0$ with
   $B(x,R)\subset B(x',2R)$ and any $n,m\in \mathbb{Z}$, we have
\begin{equation}\label{seconddb1}
    2d_{B(x,R)}(\mfa_n,\mfa_m)\le d_{B(x',2R)}(\mfa_n,\mfa_m) \, .
\end{equation}
\end{lemma}

\begin{proof}
\leanok
    With \eqref{eqcarl4}, both sides of \eqref{firstdb1} are equal to $4R|n-m|$. This proves the lemma.
\end{proof}

\begin{lemma}[integer ball cover]
    \label{integer-ball-cover}
    \leanok
    \lean{integer_ball_cover}
    For every $x\in \R$ and $R>0$ and every
    $n\in \mathbb{Z}$ and $R'>0$,
    there exist $m_1, m_2, m_3\in \mathbb{Z}$
    such that
    \begin{equation}\label{eqcarl5}
        B'\subset B_1\cup B_2\cup B_3\, ,
    \end{equation}
where
\begin{equation}
B'= \{ \mfa \in \Mf: d_{B(x,R)}(\mfa, \mfa_n)<2R'\}
\end{equation}
and for $j=1,2,3$
\begin{equation}
  B_j=
     \{ \mfa \in \Mf: d_{B(x,R)}(\mfa, \mfa_{m_j})<R'\}
     \, .
\end{equation}
\end{lemma}

\begin{proof}
\leanok
Let $m_1$ be the largest integer smaller than
or equal to
$n- \frac {R'} {2R}$.
Let $m_2=n$.
Let $m_3$ be the smallest integer larger than
or equal to $n+ \frac {R'} {2R}$.

Let $\mfa_{n'}\in B'$, then with \eqref{eqcarl4},
we have
\begin{equation}\label{eqcarl6}
    2R|n-n'|< 2R'\, .
\end{equation}

Assume first $n'\le m_1$. With \eqref{eqcarl6}
we have
\begin{equation*}
    R|m_1-n'|=R(m_1-n')=R(m_1-n)+R(n-n')
\end{equation*}
\begin{equation}
    < -\frac{R'}2+R'=-\frac{R'}2\, .
\end{equation}
We conclude $\mfa_{n'}\in B_1$.

Assume next $m_1<n'<m_3$. Then
$\mfa_{n'}\in B_2$.

Assume finally that $m_3\le n'$. With \eqref{eqcarl6}
we have
\begin{equation*}
    R|m_3-n'|=R(n'-m_3)=R(n'-n)+R(n-m_3)
\end{equation*}
\begin{equation}
    < R' -\frac{R'}2=-\frac{R'}2\, .
\end{equation}
We conclude $\mfa_{n'}\in B_1$.
This completes the proof of the lemma.
\end{proof}

\begin{lemma}[real van der Corput]
\label{real-van-der-Corput}
\leanok
\lean{real_van_der_Corput}
\uses{van-der-Corput}
    For any $x\in \R$ and $R>0$ and any
    function $\varphi: X\to \C$ supported on $B'=B(x,R)$
    such that
\begin{equation}
    \|\varphi\|_{\Lip(B')} = \sup_{x \in B'} |\varphi(x)| + R \sup_{x,y \in B', x \neq y} \frac{|\varphi(x) - \varphi(y)|}{\rho(x,y)}
\end{equation}
is finite and for any $n,m\in \mathbb{Z}$, we have
\begin{equation}
    \label{eq-vdc-cond1}
    \left|\int_{B'} e(\mfa_n(x)-{\mfa_m(x)}) \varphi(x) d\mu(x)\right|\le 2\pi \mu(B')\frac{\|\varphi\|_{\Lip(B')}}{1+d_{B'}(\mfa_n,\mfa_m)}
\, .
\end{equation}
\end{lemma}
\begin{proof}
    \leanok
Set $n'=n-m$. Then we have to prove
\begin{equation}
    \label{eq-vdc-cond2}
    \left|\int_{x-R}^{x+R} e^{in'y}\varphi(y) dy\right|\le 4\pi R\|\varphi\|_{\Lip(B')}
(1+2R|n'|)^{-1}\, .
\end{equation}
This follows from \Cref{van-der-Corput} with $\alpha = x - R$ and $\beta = x + R$.
\end{proof}

\begin{proof}[Proof of \Cref{real-Carleson}]
\proves{real-Carleson}
\leanok

The preceding chain of lemmas establishes that $\Mf$ is a cancellative, compatible collection of functions on $(\R, \rho, \mu, 4)$. Again, some of the statements in these lemmas are stronger than what is needed for $a=4$, but can be relaxed to give the desired conclusion for $a=4$.

With $\kappa$ as near \eqref{eq-hilker}, define
the function $K:\R\times \R\to \mathbb{C}$ as in Theorem
\ref{metric-space-Carleson} by
\begin{equation}
    K(x,y):=\kappa(x-y)\, .
\end{equation}
The function $K$ is continuous outside the diagonal
$x=y$ and vanishes on the diagonal. Hence it is measurable.

By \Cref{Hilbert-kernel-bound,Hilbert-kernel-regularity}, it follows that $K$ is a two-sided Calder\'on--Zygmund kernel on $(\R,\rho,\mu,4)$.
\Cref{Hilbert-strong-2-2} verifies \eqref{two-sided-Hr-bound-assumption}.
Thus the assumptions of \Cref{two-sided-metric-space-Carleson} are all satisfied. Applying the Theorem, \Cref{real-Carleson} follows.

\end{proof}

\printbibliography
\end{document}